\pdfoutput=1
\RequirePackage{ifpdf}
\ifpdf 
\documentclass[pdftex]{sigma}
\else
\documentclass{sigma}
\fi

\usepackage{psfrag}
\usepackage[active]{srcltx}
\usepackage[all]{xy}

\numberwithin{equation}{section}
\numberwithin{figure}{section}

\newtheorem{teo}{Theorem}[section]
\newtheorem{lem}[teo]{Lemma}
\newtheorem{cor}[teo]{Corollary}

\newtheorem{prop}[teo]{Proposition}

 { \theoremstyle{definition}
\newtheorem{defi}[teo]{Definition}

\newtheorem{remk}[teo]{Remark}
}

\newcommand{\fonc}[5]{ 
 \begin{array}{lcll}#1 :& #2 & \longrightarrow & #3 \\ %
 &#4 &\longmapsto & #5 %
 \end{array}}

\newcommand{\Ue}{U_\epsilon}

\newcommand{\Oq}{\Oo_{q}}
\newcommand{\OqD}{\Oo_{q}\big(q^{1/D}\big)}

\newcommand{\mUq}{\mathbb{U}_q}

\newcommand{\e}{\epsilon}

\newcommand{\mr}{\mathbb{R}}
\newcommand{\mc}{\mathbb{C}}
\newcommand{\mz}{\mathbb{Z}}

\newcommand{\mn}{\mathbb{N}}

\newcommand{\Gg}{{\mathcal G}}

\newcommand{\Ll}{{\mathcal L}}

\newcommand{\Oo}{{\mathcal O}}

\newcommand{\lra}{\longrightarrow}

\newcommand{\ra}{\rightarrow}

\begin{document}
\allowdisplaybreaks

\newcommand{\arXivNumber}{1912.02440}

\renewcommand{\PaperNumber}{025}

\FirstPageHeading

\ShortArticleName{Unrestricted Quantum Moduli Algebras. I.~The Case of Punctured Spheres}

\ArticleName{Unrestricted Quantum Moduli Algebras.\\ I.~The Case of Punctured Spheres}

\Author{St\'ephane BASEILHAC and Philippe ROCHE}

\AuthorNameForHeading{S.~Baseilhac and P.~Roche}

\Address{IMAG, Univ Montpellier, CNRS, Montpellier, France}
\Email{\href{mailto:stephane.baseilhac@umontpellier.fr}{stephane.baseilhac@umontpellier.fr}, \href{mailto:philippe.roche@umontpellier.fr}{philippe.roche@umontpellier.fr}}

\ArticleDates{Received April 09, 2021, in final form March 07, 2022; Published online March 29, 2022}

\Abstract{Let $\Sigma$ be a finite type surface, and $G$ a complex algebraic simple Lie group with Lie algebra $\mathfrak{g}$. The quantum moduli algebra of $(\Sigma,G)$ is a quantization of the ring of functions of $X_G(\Sigma)$, the variety of $G$-characters of $\pi_1(\Sigma)$, introduced by Alekseev--Grosse--Schomerus and Buffenoir--Roche in the mid '90s. It can be realized as the invariant subalgebra of so-called graph algebras, which are $U_q(\mathfrak{g})$-module-algebras associated to graphs on~$\Sigma$, where~$U_q(\mathfrak{g})$ is the quantum group corresponding to $G$. We study the structure of the quantum moduli algebra in the case where $\Sigma$ is a sphere with $n+1$ open disks removed, $n\geq 1$, using the graph algebra of the ``daisy'' graph on $\Sigma$ to make computations easier. We provide new results that hold for arbitrary $G$ and generic $q$, and develop the theory in the case where $q=\e$, a primitive root of unity of odd order, and $G={\rm SL}(2,\mc)$. In such a situation we introduce a Frobenius morphism that provides a natural identification of the center of the daisy graph algebra with a finite extension of the coordinate ring~$\mathcal{O}(G^n)$. We extend the quantum coadjoint action of De-Concini--Kac--Procesi to the daisy graph algebra, and show that the associated Poisson structure on the center corresponds by the Frobenius morphism to the Fock--Rosly Poisson structure on $\mathcal{O}(G^n)$. We show that the set of fixed elements of the center under the quantum coadjoint action is a finite extension of $\mc[X_G(\Sigma)]$ endowed with the Atiyah--Bott--Goldman Poisson structure. Finally, by using Wilson loop operators we identify the Kauffman bracket skein algebra $K_{\zeta}(\Sigma)$ at $\zeta:={\rm i}\epsilon^{1/2}$ with this quantum moduli algebra specialized at $q=\e$. This allows us to recast in the quantum moduli setup some recent results of Bonahon--Wong and Frohman--Kania-Bartoszy\'nska--L\^e on $K_{\zeta}(\Sigma)$.\looseness=-1}


\Classification{16R30; 17B37; 20G42; 57M27; 57R56; 81R50}

\section{Introduction}

Let $\Sigma$ be an oriented surface of finite type, and $G$ a complex algebraic simple Lie group with Lie algebra $\mathfrak{g}$. In this paper we begin our investigation of the quantum moduli algebra defined by quantum lattice gauge field theory (qLGFT) on $\Sigma$ with gauge algebra $U_\e(\mathfrak{g})$, where $U_\e(\mathfrak{g})$ is the adjoint unrestricted quantum group $U_\e(\mathfrak{g})$ at a primitive root of unity $\e$. For technical simplicity we focus in this paper on the case where $\Sigma$ has genus $0$ and $n\geq 2$ boundary components, we assume that $\e$ has odd order, and we prove our main results in the case of $G = {\rm SL}(2,\mc)$. Their formulation for arbitrary $G$ has qualitatively the same form.

\looseness=-1 Our main motivation comes from quantum topology. We aim at showing that the quantum moduli algebras make a very efficient and unifying setting by which quantum invariant theory for manifolds equipped with $G$-characters can be studied. As an example, in this paper we will verify this postulate on the Kauffman bracket skein algebra $K_\zeta(\Sigma)$, where $\zeta = {\rm i}\e^{1/2}$ is a~primitive root of unity of order~$4l$, with $l\geq 3$ odd. Namely, we will recast some recent results of Bonahon--Wong \cite{BW1,BW2} and Frohman--Kania-Bartoszy\'nska--L\^e \cite{FKL, FKL2} on $K_\zeta(\Sigma)$ in the setup of quantum moduli algebras, where they follow from our general results applied to the case of $\mathfrak{g}={\mathfrak{sl}}(2)$.

Our approach consists genuinely of doing geometric invariant theory for quantum groups. In order to present our results, let us recall a few facts about qLGFTs.

\looseness=-1 The qLGFTs were introduced in the mid '$90$s by Alekseev--Grosse--Schomerus \cite{A,AGS1,AGS2,AS} and Buffenoir--Roche \cite{BR1,BR2}, who used as gauge algebras the quantum groups $H=U_q(\mathfrak{g})$ with $q$ generic, or semi-simplifications thereof when $q$ is a root of unity. Assuming that $\Sigma$ has non-empty boundary, which simplifies this presentation and is the case studied in this paper, the qLGFT on~$\Sigma$ with gauge algebra $H$ associates a $H$-module-algebra~$\Ll_\Gamma(H)$, called graph algebra, to any ribbon graph $\Gamma$ embedded in~$\Sigma$ and onto which~$\Sigma$ deformation retracts. The $H$-invariant subalgebra $\Ll_\Gamma(H)^H$ is independent up to isomorphism of the choice of~$\Gamma$, and so is canonically associated to~$\Sigma$. Abusing of notations we call $\Ll_\Gamma(H)^H$ ``the'' quantum moduli algebra of the qLGFT.

In these papers $\Ll_\Gamma(H)$ was defined by the method of ``combinatorial quantization", which yields presentations by generators and relations given in matrix form. Assuming $\Gamma$ has one vertex and no edge contractible in $\Sigma$, these presentations make $\Ll_\Gamma(H)$ a natural deformation quantization of the ring of regular functions $\mathcal{O}\big(G^E\big)$ endowed with the Fock--Rosly Poisson structure \cite{AM0,AM,FR}, where $E$ is the number of edges (loops) of $\Gamma$. Hence $\Ll_\Gamma(H)^H$ is a deformation quantization of $\mathcal{O}\big(G^E\big)^G \subset \mathcal{O}\big(G^E\big)$, the ring of regular functions invariant under the coadjoint action of $G$, i.e., the ring of regular functions on the variety $X_G(\Sigma)$ of characters of representations $\pi_1(\Sigma) \ra G$, endowed with the Atiyah--Bott--Goldman Poisson structure.

Bullock--Frohman--Kania-Bartoszy\'nska provided in \cite{BFK} a coordinate free construction of the qLGFTs that works for any ribbon Hopf gauge algebra~$H$, based on the Reshetikhin--Turaev functor. They also related the Wilson loop elements of $\Ll_\Gamma(H)^H$, introduced in \cite{BR2} and associated to the isotopy classes of framed oriented links in $\Sigma \times [0,1]$, to the Kauffman bracket skein algebra of~$\Sigma$. In particular, by adapting their approach to $H=U_q({\mathfrak{sl}}(2))$ (which is strictly speaking ribbon in a certain completion), they showed in \cite{BFK2} that the construction of Wilson loop elements yields an isomorphism from the skein algebra $K_q(\Sigma)$ defined over $\mc(q)$ to the quantum moduli algebra for $H=U_q({\mathfrak{sl}}(2))$. This isomorphism, that we call Wilson loop map, explained in a~very natural way the emergence in qLGFTs of the Witten--Reshetikhin--Turaev mapping class group representations and of the Jones polynomial of links, already discovered in~\cite{AS,BR2}. We refer to~\cite{MW} for a comprehensive account of the axiomatic and algebraic structures of qLGFTs, and to~\cite{B-BZ-J} for their re-appearance in the context of factorization homology.

With the exception of the work of Frolov \cite{Fro}, until recently the qLGFTs for (non semi-simplified) quantum groups at roots of unity were not studied, certainly because of their apparent complicated definitions. Although the combinatorial quantization approach is the less intrinsic, it has the advantage of providing computationally transparent connections with representation theory. The more topological approach of~\cite{BFK} makes some invariance statements obvious, but encapsulates part of representation theory in a diagrammatic calculus, leaving many aspects rather implicit.

A major progress has been achieved recently by Faitg in \cite{Faitg3,Faitg,Faitg2,Faitg4}. He defined the qLGFTs for arbitrary finite-dimensional ribbon factorizable Hopf gauge algebras $H$, not necessarily semisimple, using combinatorial quantization. He showed that the mapping class group representations associated to such qLGFTs coincide with those of Lyubashenko--Majid \cite{Lyu, LM}, which were originally defined by categorical means, and provided explicit and ready-to-use formulas for Dehn twists. In particular this includes (with little adaptation) the restricted quantum group for ${\mathfrak{sl}}(2)$ at a primitive root of unity of even order. In this case, he generalized the Wilson loop map, obtained new non semisimple representations of $K_q(\Sigma)$, and established new relationships with the stated skein algebras. He showed also that the corresponding qLGFT mapping class group representations coincide in genus $g=1$ with those derived from logarithmic conformal field theory in~\cite{FGST}.

\looseness=-1 In the present paper we make a further step in the root of unity case; we note that another approach is being developed in \cite{G-J-S}. We consider the qLGFTs on $\Sigma$ with gauge algebra an unrestricted (adjoint) quantum group $U_q(\mathfrak{g})$; when $q=\e$ is a root of unity we mainly focus on the case $\mathfrak{g}={\mathfrak{sl}}(2)$. Also for technical simplicity we assume $\e$ is a primitive root of unity of odd order $l$; the case of even order can be treated similarly. As mentioned above, for simplicity also we focus on the case where $\Sigma$ is a sphere with $n+1$ open disks removed, $n\geq 1$ (sometimes we say for short that $\Sigma$ is punctured), but qualitatively similar results hold for surfaces of non zero genus.

First we construct the graph algebra $\Ll_{0,n} := \Ll_{0,n}(\mathfrak{g})$ associated to the ``daisy graph'' in $\Sigma$, made of one vertex, one loop encircling each deleted disk, and one ``cilium" at the vertex, which provides an ordering of the loops by using the orientation of $\Sigma$. The embedding in $\Sigma$ gives the daisy graph a structure of ribbon graph, as shown in the picture below:

\begin{figure}[h!]\centering
\vspace*{-30mm}
\def\svgwidth{0.35\textwidth}
\begingroup%
 \makeatletter%
 \providecommand\color[2][]{%
 \errmessage{(Inkscape) Color is used for the text in Inkscape, but the package 'color.sty' is not loaded}%
 \renewcommand\color[2][]{}%
 }%
 \providecommand\transparent[1]{%
 \errmessage{(Inkscape) Transparency is used (non-zero) for the text in Inkscape, but the package 'transparent.sty' is not loaded}%
 \renewcommand\transparent[1]{}%
 }%
 \providecommand\rotatebox[2]{#2}%
 \ifx\svgwidth\undefined%
 \setlength{\unitlength}{595.27559055bp}%
 \ifx\svgscale\undefined%
 \relax%
 \else%
 \setlength{\unitlength}{\unitlength * \real{\svgscale}}%
 \fi%
 \else%
 \setlength{\unitlength}{\svgwidth}%
 \fi%
 \global\let\svgwidth\undefined%
 \global\let\svgscale\undefined%
 \makeatother%
 \begin{picture}(1,1.41428571)%
 \put(0,0){\includegraphics[width=\unitlength,page=1]{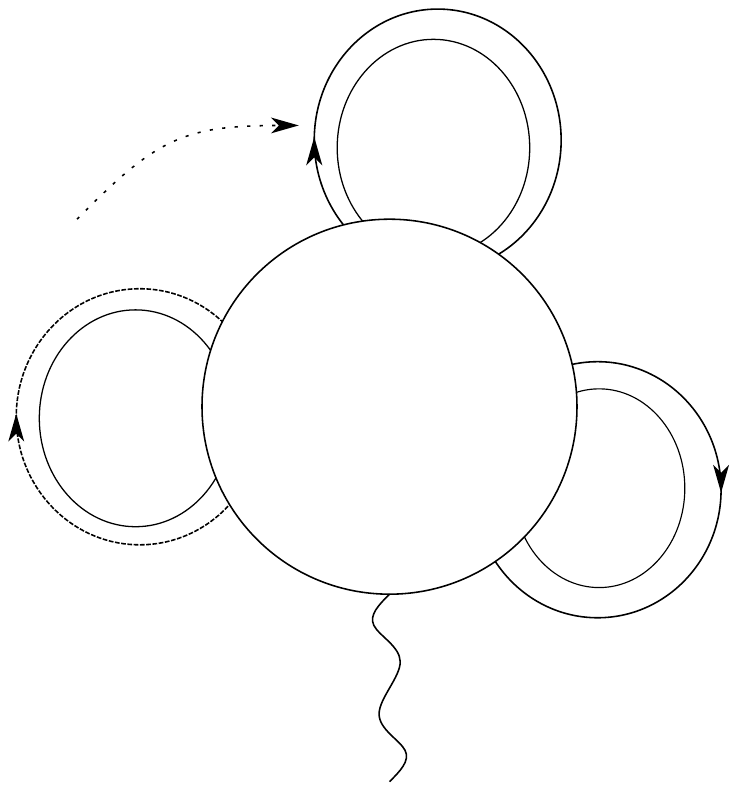}}%
 \put(0.65515871,0.30238095){\color[rgb]{0,0,0}\makebox(0,0)[lb]{\smash{$n$}}}%
 \put(0.15119048,0.37797619){\color[rgb]{0,0,0}\makebox(0,0)[lb]{\smash{$1$}}}%
 \put(2.66666667,-0.97866516){\color[rgb]{0,0,0}\makebox(0,0)[lt]{\begin{minipage}{0.66666667\unitlength}\raggedright \end{minipage}}}%
 \put(0.6047619,0.83154762){\color[rgb]{0,0,0}\makebox(0,0)[lb]{\smash{$a$}}}%
 \put(0,0){\includegraphics[width=\unitlength,page=2]{dessin5.pdf}}%
 \end{picture}%
\endgroup%

\vspace*{-15mm}

\caption{The daisy graph.}
\end{figure}

The graph algebra $\Ll_{0,n}$ is a module-algebra over the quantum group $U_q:=U_q(\mathfrak{g})$ with ground ring $\mc(q)$. We define $\Ll_{0,n}$ by means of combinatorial quantization based on $U_q$, that we reformulate also in terms of twists of module-algebras and braided tensor product. More precisely, $U_q$~is not a ribbon Hopf algebra, but a suitable extension of the category of finite-dimensional $U_q$-modules is ribbon. So, to make sense of the construction of $\Ll_{0,n}$ we replace $U_q$ by a categorical completion~$\mathbb{U}_q$.

In particular, $\Ll_{0,1} = \mathcal{O}_q$ as a $U_q$-module, where $\mathcal{O}_q$ is the restricted dual of $U_q$ endowed with the right coadjoint action of $U_q$, and the algebra structure of $\Ll_{0,1}$ is compatible with that action. Eventually, we find that a map due to Alekseev \cite{A} yields an equivariant embedding of $U_q$-module-algebras
\[\Phi_n\colon \ \Ll_{0,n}\ra \tilde U_q^{\otimes n}, \]
where $\tilde U_q$ is the simply-connected quantum group associated to $\mathfrak{g}$, and $\Ll_{0,n}$ and $\tilde U_q^{\otimes n}$ are endowed with a coadjoint and an adjoint action of $U_q$ respectively. In the case $n=1$, $\Phi_1$ coincides with a celebrated isomorphism of $U_q$-modules $\mathcal{O}_q \ra \tilde U_q^{\rm lf}$, where $\tilde U_q^{\rm lf}$ is the subalgebra of locally finite elements of~$\tilde U_q$, which was first introduced by Drinfeld and Reshetikhin--Semenov-Tian-Shansky~\cite{RSTS}, and further studied by Caldero, Joseph--Letzter and others (see, e.g., \cite{Bau1, Ca,JL2}).

For the purpose of defining specializations at $q=\e$, we then introduce an integral form~$\Ll_{0,n}^A$ of~$\Ll_{0,n}$, that is, an $A$-algebra satisfying $\Ll_{0,n} = \Ll_{0,n}^A\otimes_A \mc(q)$, where $A=\mc\big[q,q^{-1}\big]$. It is a~module-algebra over the {\it unrestricted} integral form~$U_A$ of~$U_q$, as defined by De Concini--Kac--Procesi \cite{DC-K,DC-K-P2}. The construction of~$\Ll_{0,n}^A$ is based on Lusztig's~\cite{Lusztig} {\it restricted} integral form~$U_A^{\rm res}$ of $U_q$ and some integrality properties of the $R$-matrix on $U_A^{\rm res}$-modules. The Alekseev map yields an equivariant embedding of $U_A$-module-algebras $\Phi_n\colon \Ll_{0,n}^A\ra \tilde U_A^{\otimes n}$, and the iterated coproduct~$\Delta^{(n-1)}$ of a (natural) integral form $\mathcal{O}_A$ of $\mathcal{O}_q$ defines a morphism of algebras $\Delta^{(n-1)}\colon \Ll_{0,1}^A\rightarrow \Ll_{0,n}^A$. We show (see Proposition~\ref{nozeroq} and Theorem~\ref{CenterLinv} for details):
\begin{teo}\label{teointrogen} The invariant subalgebra $\Ll_{0,n}^{U_q}$ does not have non trivial zero divisors, and its center is a polynomial algebra, generated by $\Delta^{(n-1)}(\mathcal{Z}(\Ll_{0,1}))$ and $\mathcal{Z}(\Ll_{0,n})$.
\end{teo}

In the case $\mathfrak{g} = {\mathfrak{sl}}(2)$ this result also follows from the isomorphism $\mathcal{W}$ of Theorem~\ref{teointroOBS} below, and the fact, proved in \cite{BW0,PS}, that the skein algebras satisfy the corresponding properties. In~\cite{PS} it is also shown that the skein algebras are finitely generated and Noetherian.

The constructions and results above are developed from Section~\ref{CATCOMP} to Section~\ref{Lgnalg}, for any of the quantum groups $U_q = U_q(\mathfrak{g})$. As they form the basis of all results that follow we give full details, though part of this material has already been considered in some ways in the litterature. Along the text and especially in Section~\ref{sl2} we consider in detail the case of $\mathfrak{g}={\mathfrak{sl}}(2)$.

Section \ref{specialization} is devoted to the center of the specializations of $\Ll_{0,n}^A$ at roots of unity. For every $\e\in \mc^\times$, set $U_\e := U_A\otimes_A {\mathbb C}_{\e}$ and
\[\Ll_{0,n}^\e = \Ll_{0,n}^A\otimes_A {\mathbb C}_{\e},\]
where ${\mathbb C}_{\e}={\mathbb C}$ as a vector space, and as an $A$-module, $q$ is evaluated as $\e$ on ${\mathbb C}_{\e}$. The Alekseev map affords an embedding of $U_\e$-modules $\Phi_n\colon \Ll_{0,n}^\e\ra \tilde{U}_\e^{\otimes n}$. We study this map when $\e$ is a~primitive root of unity of odd order~$l$.

De Concini--Kac--Procesi \cite{DC-K,DC-K-P1,DC-K-P2} showed that the center $\mathcal{Z}(\tilde U_\e)$ can be identified with the coordinate ring of a Poisson--Lie group $G^*$ (dual to $G$ endowed with the standard Poisson--Lie structure), and that certain Hamiltonian vector fields on $\operatorname{Spec}\big(\mathcal{Z}\big(\tilde U_\e\big)\big)\cong G^*$ can be integrated to define an {\it infinite dimensional} group $\Gg_{\rm DCK}$ acting by {\it analytic} automorphisms on $\operatorname{Spec}\big(\mathcal{Z}\big(\tilde U_\e\big)\big)$, and by automorphisms on a suitable completion of~$U_\e$. The orbits of this action, called {\it quantum coadjoint action}, lift the orbits of the conjugation action via the natural covering map $G^*\ra G^0$, where $G^0$ is the big cell of $G$. They used this action to obtain a series of fundamental results on the simple $U_\e$-modules.

It is not hard to make the quantum coadjoint action explicit for $U_\e({\mathfrak{sl}}(2))$. Therefore, starting from Section~\ref{specialization} we restrict to the case of $\mathfrak{g}={\mathfrak{sl}}(2)$, thus omitting~$\mathfrak{g}$ from the notations (denoting, e.g., $\Ll_{0,n}^\e({\mathfrak{sl}}(2))$ by $\Ll_{0,n}^\e$), and we put $G={\rm SL}(2,\mc)$. In the case of $\mathfrak{g}={\mathfrak{sl}}(2)$ the Alekseev map affords an isomorphism
\[\Phi_n\colon \ {}_{\rm loc}\Ll_{0,n}^\e\ra U_\e^{\otimes n}, \]
where ${}_{\rm loc}\Ll_{0,n}^\e$ is the specialization at $q=\e$ of a localization of $\Ll_{0,n}^A({\mathfrak{sl}}(2))$ introduced in Section~\ref{Lgnalg}.

After proving preliminary results on the center $\mathcal{Z}(\Ll_{0,n}^\e)$ in Section~\ref{centersl2}, we extend the quantum coadjoint action by means of the Alekseev map, to get in Section~\ref{QCAL0nsec} an action of a {\it finite dimensional} Lie group $\Gg$ on the fraction ring of $\mathcal{Z}(\Ll_{0,n}^\e)$ (hence a {\it partial action} on $\operatorname{Spec}(\mathcal{Z}(\Ll_{0,n}^\e))$ by rational transformations) and on a suitable completion of $\Ll_{0,n}^\e$. Then we study the invariant subalgebra $\mathcal{Z}(\Ll_{0,n}^\epsilon)^{\mathcal{G}}$. The groups $\Gg$ and $\Gg_{\rm DCK}$ are different. In fact, it is necessary to adapt the constructions of De Concini--Kac--Procesi in the case $n>1$ (see Remark \ref{explanation}). We note that the problem of extending the quantum coadjoint action to graph algebras has already been considered by Frolov in \cite{Fro}.

We can summarize the main results of Section \ref{specialization} as follows (see Corollary~\ref{equivcor} for a precise and more complete statement). Denote by $X_G(\Sigma)$ the variety of characters of representations $\pi_1(\Sigma)\ra G$. We define an $l$-fold branched covering space $\tilde{G}$ of $G$, and an $l^n$-fold branched covering space $\tilde X_G'(\Sigma)$ of $X_G(\Sigma)$, related by a natural branched covering identification map $\tilde{\mathfrak{c}}\colon \tilde X_G'(\Sigma)\ra \tilde{G}^{n}$. Then we prove:
\begin{teo} \label{teointrocentre} The center $\mathcal{Z}(\Ll_{0,n}^\epsilon)$ is naturally endowed with a Poisson bracket $\{\ ,\}_{\rm QCA}$ inhe\-ri\-ted from the algebra structure of $\Ll_{0,n}^\epsilon$, so that there is an isomorphism of Poisson algebras
\[\widetilde{\rm Fr}\colon \ \big(\mathcal{O}\big(\tilde{G}^{n}\big),\{\, ,\, \}_{\rm FR}\big) \ra \big(\mathcal{Z}(\Ll_{0,n}^\epsilon),\{\ ,\}_{\rm QCA}\big),\]
where $\{\, ,\, \}_{\rm FR}$ is the trivial extension to $\mathcal{O}\big(\tilde{G}^{n}\big)$ of the Fock--Rosly Poisson bracket on $\mathcal{O}(G^n)$. Moreover, $\widetilde{\rm Fr}$ yields an isomorphism of Poisson algebras
\[\widetilde{\rm Fr} \circ \tilde{\mathfrak{c}}^{-1*}\colon \ \big(\mathcal{O}\big(\tilde X_G'(\Sigma)\big),\{\, ,\, \}_{\rm Gold}\big) \ra \big(\mathcal{Z}(\Ll_{0,n}^\epsilon)^{\mathcal{G}},\{\, ,\, \}_{\rm QCA}\big),\]
where $\{\, ,\, \}_{\rm Gold}$ is the trivial extension to $\mathcal{O}\big(\tilde X_G'(\Sigma)\big)$ of the Atiyah--Bott--Goldman Poisson bracket on $\mathcal{O}(X_G(\Sigma))$.
\end{teo}

The isomorphism $\widetilde{\rm Fr}$ provides a precise formulation of what means combinatorial quantization of the $\operatorname{Ad}(G)$-module-algebra $\mathcal{O}(G^n)$ at roots of unity. It maps the elements generating $\mathcal{O}\big(\tilde{G}^{n}\big)$ as an extension of $\mathcal{O}(G^n)$ to analogs of Casimir elements in $\mathcal{Z}(\Ll_{0,n}^\epsilon)$. On the subalgebra $\mathcal{O}(G^n)$, it is given by a {\it Frobenius map} ${\rm Fr}$, analogous in the case $n=1$ to the one defined for $\mathcal{O}_\e = {\rm SL}_\e(2)$ by Parshall--Wang in~\cite{PW}, though more complicated, see Definitions~\ref{defFr1} and~\ref{defFrn}. It satisfies the remarkable identity (and a similar one for $n>1$, see Proposition~\ref{Thcentral})
\[
{\rm Fr}\big(\!\operatorname{Tr}\big(\underline{\stackrel{V_2}{M}}\big)\big) = T_l\big(\!\operatorname{qTr}\big(\!\stackrel{V_2}{M}\!\big)\big)= \operatorname{qTr}\big({T_l(V_2)}{M}\big),
\]
where $\underline{\stackrel{V_2}{M}}$ is the matrix of coordinate functions of $G$ in its fundamental representation, $\stackrel{V_2}{M}$ is the matrix of generators of $\Ll_{0,1}^\e$ in the fundamental representation $V_2$ of $U_\e$, $T_l$ is the $l$-th Chebyshev polynomial of the first type (suitably normalized), $T_l(V_2)$ the corresponding virtual representation in the Grothendieck ring of $U_q$-modules, and $\operatorname{qTr}$ and $\operatorname{Tr}$ are the quantum trace and classical trace of $2\times 2$ matrices respectively. This identity shows how $Fr$ relates invariant functions on~$G$ to $\Gg$-invariant central elements. The appearance of the $l$-th Chebyshev polynomial~$T_l$ in this context relies on the fact that it generates the defining relation of $\mathcal{Z}(U_\e)$, between the Casimir element $\Omega$ and the generators $E^l$, $F^l$, $K^{\pm l}$ of the ``small'' center $\mathcal{Z}_0(U_\e) \subset \mathcal{Z}(U_\e)$.

In the context of the quantum function algebra ${\rm SL}_\e(2)$, identities similar to the above one and its extension to $n>1$ have been obtained by Bonahon in~\cite{Bon17}.

In Section \ref{SKEIN} we develop a topological (i.e., skein theoretic) formulation of some of the previous results.
In Section \ref{Wfunctor} we give two definitions of a {\it Wilson loop functor}~$\mathbb{W}$, defined on a category of ribbon oriented graphs in $\Sigma \times [0,1]$ colored by $U_A^{\rm res}$-modules, extending the Wilson loop map of~\cite{BR2} and defined for any~$\mathfrak{g}$. One of these formulations uses the Reshetikhin--Turaev functor, and is close to the one of~\cite{BFK}.

In Section \ref{Wiso} we consider the restriction $W$ of $\mathbb{W}$ to {\it closed} colored ribbon oriented graphs. The image of $W$ is $\big(\Ll_{0,n}^A\big)^{U_A}\otimes_A \mc\big[q^{1/D},q^{-1/D}\big]$, where we recall that $\big(\Ll_{0,n}^A\big)^{U_A}$ is the $A$-algebra of $U_A$-invariant elements of $\Ll_{0,n}^A$, and $A=\mc\big[q,q^{-1}\big]$. Moreover, we prove the following result, which is an integral version (i.e., over the ring $A$) of a combination of \cite[Theorem~10]{BFK} and \cite[Theorem~1]{BFK2}. Denote by $K_\zeta(\Sigma)$ the Kauffman bracket skein algebra of $\Sigma$, defined over the ring $\mc\big[\zeta,\zeta^{-1}\big]$. We have \big(see Theorem~\ref{teoOBS}, and Remark~\ref{teoOBSrem}(3) for a statement over $\mz\big[\zeta,\zeta^{-1}\big]$\big):
\begin{teo}\label{teointroOBS} When $\mathfrak{g}={\mathfrak{sl}}(2)$, the linear map defined by $\mathcal{W}(L) = {\rm i}^{{\rm lk}(L)}W(L)$ on ribbon oriented links~$L$ colored by the fundamental representation, where ${\rm lk}(L)$ is the linking number of~$L$, descends to an isomorphism of algebras $\big($where $\zeta := {\rm i}q^{1/2}\big)$:
\[
\mathcal{W}\colon \ K_\zeta(\Sigma)\ra \big(\Ll_{0,n}^A\big)^{U_A}\otimes_A \mc\big[\zeta,\zeta^{-1}\big].
\]
\end{teo}
By using the image by $\mathcal{W}$ of the multicurve basis of $K_\zeta(\Sigma)$ we prove in Theorem~\ref{centerprop1} the following result, about the specialization $\big(\Ll_{0,n}^A\big)^{U_A}_\e = \big(\Ll_{0,n}^A\big)^{U_A} \otimes_A \mc_\e$.
\begin{teo} \label{teointroder} The algebra $\mathcal{Z}(\Ll_{0,n}^\epsilon)^{\mathcal{G}}$ is contained in $\big(\Ll_{0,n}^A\big)^{U_A}_\e$.
Moreover, the bracket $\{\, ,\, \}_{\rm QCA}$ extends to an action by derivations of $\mathcal{Z}(\Ll_{0,n}^\epsilon)^{\mathcal{G}}$ on $\big(\Ll_{0,n}^A\big)^{U_A}_\e$.
\end{teo}

Theorems \ref{teointrogen}, \ref{teointrocentre}, \ref{teointroOBS} and \ref{teointroder} imply:
\begin{cor}\quad
\begin{enumerate}\itemsep=0pt
\item[$(1)$] The skein algebra $K_\zeta(\Sigma)$ does not have non trivial zero divisors, and its center is the polynomial ring generated by the skein classes of the boundary components of~$\Sigma$.

\item[$(2)$] When $\zeta$ is specialized to a root of unity $\e'$ of order $4l$, with $l\geq 3$ odd, the center of $K_{\e'}(\Sigma)$ contains a subalgebra isomorphic to $\mathcal{O}(X_G(\Sigma))$, endowed with the image of the Poisson bracket $\{\, ,\, \}_{\rm Gold}$, which extends to an action by derivations on $K_{\e'}(\Sigma)$.
\end{enumerate}
\end{cor}

Details are given in Section~\ref{appskein}. As already mentioned above, (1) has been proved in~\cite{BW0} and~\cite{PS}. The claim (2) belongs to a corpus of results proved in \cite{BW1,BW2} and~\cite{FKL,FKL2}, for any finite type surface and root of unity $\e'$. One interest of our method is to be intrinsically algebro-geometric, and valid for any complex simple Lie algebra~$\mathfrak{g}$. We note that the Frobenius map ${\rm Fr}$ discussed after Theorem~\ref{teointrocentre} provides an explicit, geometric realization of the threading map ${\rm Ch}\colon K_{\epsilon'{}^{l^2}}(\Sigma) \ra \mathcal{Z}_{\epsilon'}(\Sigma)$ of Bonahon--Wong, see~\cite{BW1} and also~\cite{FKL2} (note that $\epsilon'{}^{l^2}\in\{\pm 1, \pm {\rm i}\}$).

In \cite{BaRo2} and works in preparation we study the structure of the algebra $\Ll_{0,n}^\e$ and its subalgebras $(\Ll_{0,n}^\e)^{U_\e}$ and $\big(\Ll_{0,n}^A\big)^{U_A}_\e$, their representations at roots of unity, and we extend the results of this paper to arbitrary finite type surfaces.

We note that the quantum coadjoint action implies remarkable properties of the intertwiners of quantum moduli algebras. In the ${\mathfrak{sl}}(2)$ case these properties should eventually recast the quantum hyperbolic field theories (see \cite{BB1,BB2} and the references therein), and therefore quantum Teichm\"uller theory (by the results of \cite{BB3}) within the theory of quantum moduli algebras. In~another direction, integrating the action by derivations of $\mathcal{Z}(\Ll_{0,n}^\epsilon)^{\mathcal{G}}$ should provide inte\-res\-ting information on the space of finite-dimensional $\big(\Ll_{0,n}^A\big)^{U_A}_\e$-modules. There should be no major difficulties in generalizing these results to the higher rank case.

\section{Categorical completions}\label{CATCOMP}
We recall here the notion of categorical completion of a Hopf algebra, which is connected to the theory of multiplier Hopf algebras, see \cite{VD}. We need it because our main object of interest, the algebra $\Ll_{0,n}(\mathfrak{g})$, is built from the Hopf algebra $U_q(\mathfrak{g})$ and its braided structure which exists only in some completion. The categorical one is suited to algebras defined over $\mc(q)$.

Let $k$ be a field (in the sequel it will always be ${\mathbb C}$ or ${\mathbb C}(q)$), and $U$ a $k$-associative algebra (not necessarily with unit). We denote by $\mu$ the multiplication map of $U$, and by $\operatorname{Mod}_{U}$ the category of left $U$-modules.

Let $F_U\colon \operatorname{Mod}_{U}\rightarrow \operatorname{Vect}$ be the forgetful functor from the category $\operatorname{Mod}_{U}$ to the category of $k$-vector spaces. Denote by ${\mathcal U}$ the set of natural transformations from~$F_U$ to~$F_U.$ An element of~${\mathcal U}$ is a collection $(a_X)_{X\in \operatorname{Mod}_{U}}$, where $a_X\in \operatorname{End}_k(X)$ satisfies $F_U(f)\circ a_X=a_Y\circ F_U(f)$ for any objects $X$, $Y$ in $\operatorname{Mod}_{U}$ and any arrow $f\in \operatorname{Hom}_U(X,Y)$. The direct product $\prod_{X\in \operatorname{Mod}_{U}} \operatorname{End}_k(X)$ is canonically endowed with a structure of unital $k$-algebra, and ${\mathcal U}$ is a unital subalgebra of $\prod_{X\in \operatorname{Mod}_{U}} \operatorname{End}_k(X)$. The multiplication map is given in each factor by the composition map $\mu_{X}\colon \operatorname{End}_k(X)\otimes \operatorname{End}_k(X)\to \operatorname{End}_k(X), u\otimes v\mapsto u\circ v$. The map
\[\iota_U\colon \quad U\to {\mathcal U}, \qquad a\mapsto (a_X)_{X\in \operatorname{Mod}_{U}},\] where $a_X$ is the endomorphism defined by the action of $a\in U$ on $X$, is a morphism of algebra. When $U$ has a unit $1$, $\iota_U$ is an isomorphism with inverse
\[
\iota_U^{-1}\colon \quad {\mathcal U}\to U,\qquad (a_X)_{X\in \operatorname{Mod}_{U}}\mapsto a_{U}(1),
\]
 where $U$ is endowed with its structure of left-regular representation.

We will ``enlarge'' ${\mathcal U}$ by considering only $k$-finite dimensional $U$-modules; we stress that, in this situation, the map corresponding to $\iota_U$ is not necessarily surjective nor injective (see below for the case of $U=U_q({\mathfrak g})$). Let $\operatorname{FinVect}$ and $\operatorname{FinMod}_{U}$ be respectively the full subcategories of $\operatorname{Vect}$ and $\operatorname{Mod}_{U}$ whose objects consist of the finite dimensional $k$-vector spaces and $U$-modules. Let $\operatorname{FinF}_{U}\colon \operatorname{FinMod}_{U}\to \operatorname{FinVect}$ be the forgetful functor, and $\hat{ U}$ the set of natural transformations from $\operatorname{FinF}_{U}$ to $\operatorname{FinF}_{U}$. As $\mathcal{U}$ above, $\hat{ U}$ is a unital subalgebra of $\textstyle \prod_{X\in \operatorname{FinMod}_{U}} \operatorname{End}_k(X)$, and the map (keeping the same notation)
\[\iota_U\colon \quad U\to \hat{ U}, \qquad a\mapsto (a_X)_{X\in \operatorname{FinMod}_{U}}\]
is a morphism of algebra. We will call~$\hat{U}$ the $\operatorname{FinMod}_{U} $-{\it categorical completion} of~$U$.

We will often use special elements in $\hat{U}$ defined by series of elements of $U$. Consider a sequence of elements $x_j$ of~$U$, and assume that for any object $X$ of $\operatorname{FinMod}_{U}$, the set of indices $j$ such that the endomorphism $(x_j)_X$ is not zero is finite. Then we can define an element $\sum_{j} \iota_U(x_j)$ of~$\hat{U}$ by
\begin{equation}\label{finsum}
\bigg(\sum_{j} \iota_U(x_j)\bigg)_X=\sum_{j} (x_j)_X.
\end{equation}
When $\iota_U$ is injective this element is also denoted by $\sum_{j} x_j$.

Let $U$, $V$ be $k-$algebras, and $\hat{U}$, $\hat{V}$ their $\operatorname{FinMod}$-categorical completions. We define the categorical completed tensor product of $\hat{U}$ and $\hat{ V}$, denoted by $\hat{U}\hat{\otimes} \hat{V}$, as the space of natural transformations from $F_{U,V}$ to $F_{U,V}$, where $F_{U,V}\colon \operatorname{FinMod}_{U}\times \operatorname{FinMod}_{V}\to \operatorname{FinVect}, (X,Y)\mapsto X\otimes_k Y$. An element of $\hat{U}\hat{\otimes}\hat{V}$ is a collection of linear maps $a_{X,Y}\in \operatorname{End}_k(X\otimes Y)$, for $X\in \operatorname{FinMod}_{U}$ and $ Y\in \operatorname{FinMod}_{V}$, such that for any arrows
$f\in \operatorname{Hom}_U(X,X')$ and $g\in \operatorname{Hom}_V(Y,Y')$ one has $(f\otimes g)\circ a_{X,Y}=a_{X',Y'}\circ (f\otimes g)$ (identifying $F_{U,V}(f,g)$ with the linear map $f\otimes g$). Again the componentwise composition map endows $\hat{U}\hat{\otimes}\hat{ V}$ with a structure of associative algebra, and the map $\iota_{U,V}\colon U\otimes V\to \hat{U}\hat{\otimes}\hat{ V}, u\otimes v\mapsto (u_X\otimes v_Y)_{X,Y}$ is a morphism of algebra. If $\mu\colon U\otimes U\to U$ is the product of $U$, we define $\hat{\mu}\colon \hat{U}\hat{\otimes }\hat{U}\to\hat{U}$ by
\[
(\hat{\mu}(a))_X=\mu_X(a_{X,X})
\]
for all $a\in \hat{U}\hat{\otimes }\hat{U}$ and objects $X$ of $\operatorname{FinMod}_{U}$. This construction is generalized straightforwardly to $n$-tuples of algebras $U_1,\dots, U_n$. In particular, if $U_1=\cdots =U_n=U$ we put
\[
\hat{U}^{\hat{\otimes}n}:=\hat{U}_1\widehat{\otimes}\cdots \widehat{\otimes}\hat{U}_n.
\]
Adapting Sweedler's coproduct notation $\Delta(x)=\sum_{(x)}x_{(1)}\otimes x_{(2)}$, we find convenient to write a~sum $ T=\sum_{j}u_j\hat{\otimes} v_j$ in $\hat{U}\hat{\otimes}\hat{V}$ as \[
T=\sum_{(T)}T_{(1)}\hat \otimes T_{(2)}.
\]

Assume now that $U$ is a $k$-bialgebra, with coproduct $\Delta\colon U\to U\otimes U$ and counit $\epsilon\colon U\to k$. Denote by $X_0$ the $U$-module structure on $k$ defined by the counit $\epsilon$. Set
\[
\hat{\Delta}\colon \quad \hat{U}\to \hat{U}\hat{\otimes} \hat{U},\qquad (a_{Z})_{Z\in \operatorname{FinMod}_{U}}\mapsto (a_{X,Y})_{X,Y\in \operatorname{FinMod}_{U}},
\]
where $a_{X,Y}:=a_{X\otimes Y}$, and
\[
\hat{\epsilon}\colon \ \hat{U}\rightarrow k, \qquad a\mapsto a_{X_0}.
\]
These maps are morphisms of algebras, and satisfy
\[
\hat{\Delta}\circ\iota_{U}=\iota_{U,U}\circ \Delta ,\qquad (\hat{\epsilon}\hat{\otimes} {\rm id})\big(\hat{\Delta}(a)\big)=({\rm id}\hat{\otimes} \hat{\epsilon})\big(\hat{\Delta}(a)\big)={\rm id}, \qquad \big(\hat{\Delta}\hat{\otimes} {\rm id}\big)\hat{\Delta}=\big({\rm id}\hat{\otimes} \hat{\Delta}\big)\hat{\Delta}.
\]
We will still call $\hat{U}$ a $k$-bialgebra although the tensor product is $\hat{\otimes} $ and not the algebraic one. In particular, given $a=(a_{Z})_{Z\in \operatorname{FinMod}_{U}}\in \hat{U}$ we have
\[
a_{(X,Y),Z}:=\big(\big(\hat{\Delta}\hat{\otimes} {\rm id}\big)\hat{\Delta}(a)\big)_{X,Y,Z} = \big(\big({\rm id}\hat{\otimes} \hat{\Delta}\big)\hat{\Delta}(a)\big)_{X,Y,Z}=: a_{X,(Y,Z)}.
\]

We say that $U$ is {\it quasicocommutative in its categorical completion} if there exists an invertible element $R\in \hat{U}\hat{\otimes} \hat{U}$ such that
\[
\hat{\Delta}^{\rm cop}(a) = R\hat{\Delta}(a)R^{-1}
\]
for all $a\in \hat{U}$, where $\hat{\Delta}^{\rm cop} := \sigma\circ \hat{\Delta}$, and $\sigma(a\hat{\otimes} b) = b\hat{\otimes} a$ for every $a,b\in \hat{U}$. In particular, this implies
\[
\iota_{U\otimes U}\big( \Delta^{\rm cop}(a)\big)=\iota_{U\otimes U}(\sigma(\Delta(a)))=R\iota_{U\otimes U}(\Delta(a))R^{-1}
\]
for all $a\in U$.

Assume that $U$ is quasicocommutative in its categorical completion. We say it is quasitriangular, or braided, if moreover we have
\begin{equation}\label{Rdef0}
\big(\hat{\Delta} \hat{\otimes} {\rm id}\big)(R)=R_{13}R_{23} ,\qquad \big({\rm id}\hat{\otimes} \hat{\Delta}\big)(R)=R_{13}R_{12}.
\end{equation}
Then, we call $R$ a universal $R$-matrix. We use the following notations:
\begin{equation*}
R^+:=R, \qquad R^-:= \sigma(R)^{-1},\qquad R':=\sigma(R).
\end{equation*}
Finally assume that $U$ is a Hopf algebra. For $X$ a $U$-module, denote by $X^*$ the dual module (with the standard action $a\cdot \lambda = \lambda \circ S(a)$, for all $a\in U$ and $\lambda\in X^*$). If $X$ is finite dimensional, let $\psi_X\colon X^{**}\to X$ be the canonical isomorphism of vector spaces. We define an antimorphism of algebras $\hat{S}\colon \hat{U}\to \hat{U}$ by
\[
\hat{S}(a)_{X}=\psi_X \circ {}^t(a_{X^*})\circ \psi_X^{-1}
\]
for all $X\in \operatorname{FinMod}_U$. It satisfies $\hat{\mu}\big(\hat{S}\hat{\otimes} {\rm id} \big)\hat{\Delta}(a)=\hat{\mu}\big({\rm id}\hat{\otimes}\hat{S} \big)\hat{\Delta}(a)=\hat{\eta}\hat{\epsilon}$, where $\eta\colon k\to U$, $\lambda\mapsto \lambda 1_U$, and $\hat{\eta}\colon k\to \hat{U}$, $\lambda\mapsto \lambda 1_{\hat{U}}$, are the unit maps of~$U$ and~$\hat U$ respectively. Hence~$\hat{U}$ is a Hopf algebra in the generalized sense where the tensor product $\hat{\otimes}$ is used in place of the algebraic one.

When $U$ is a Hopf algebra braided in its categorical completion, we say it is {\it ribbon in its categorical completion} if there exists $\theta\in \hat{U}$ such that $\theta$ is central,
$\hat{\Delta}(\theta)=(R'R)^{-1}(\theta\otimes \theta)$, $\epsilon(\theta)=1$ and $\hat{S}(\theta)=\theta$.

For the sake of notational simplicity, from now on we will omit the ``$\hat{\ \ }$'' symbols from the structure morphisms of the categorical completions $\hat U$ under consideration, thus denoting $\hat \Delta$ by~$\Delta$, $\hat \otimes$ by~$\otimes$, and so on, like for~$U$.

\section[The case of U\_q(g)]{The case of $\boldsymbol{U_q({\mathfrak g})}$}\label{Uq}

\subsection{Notations}
Let ${\mathfrak g}$ be a finite dimensional simple complex Lie algebra of rank~$m$. Denote by $(a_{ij})$ the Cartan matrix of ${\mathfrak g}$, and by $d_i$ the unique coprime positive integers such that the matrix $(d_ia_{ij})$ is symmetric. Fix a Cartan subalgebra $\mathfrak{h}$ and a basis of simple roots $\alpha_i \in \mathfrak{h}_{\mr}^{*}$. Denote by $(\,,\,)$ the unique inner product on $\mathfrak{h}_{\mr}^{*}$ such that $d_ia_{ij} = (\alpha_i,\alpha_j)$. The root lattice $Q$ is the $\mz$-lattice in $\mathfrak{h}_{\mr}^{*}$ defined by $\textstyle Q = \sum_{i=1}^m \mathbb{Z} \alpha_i$. The weight lattice $P$ is the $\mz$-lattice formed by all $\lambda\in \mathfrak{h}_{\mr}^{*}$ such that
\begin{gather*}
\langle \lambda,\alpha_i\rangle := 2\frac{(\lambda,\alpha_i)}{(\alpha_i,\alpha_i)}\in \mz
\end{gather*}
for every $i=1,\dots,m$. So $P= \sum_{i=1}^m \mathbb{Z} \omega_i$, where the $\omega_i$ are the fundamental weights corresponding to the simple roots $\alpha_i$, satisfying $\langle \omega_i,\alpha_j\rangle = \delta_{i,j}$. Put $Q_+ := \sum_{i=1}^m \mathbb{Z}_{\geq 0} \alpha_i$ and $P_+:= \sum_{i=1}^m \mathbb{Z}_{\geq 0} \omega_i$, the cone of dominant integral weights. Denote by $\rho$ half the sum of the positive roots, by~$N$ the number of positive roots, and by $D$ the smallest positive integer such that $D(\lambda,\mu)\in {\mathbb Z}$ for every $\lambda,\mu\in P$. Note that $(\lambda,\alpha_i)\in \mz$ for every $\lambda\in P$, $\alpha\in Q$, and $D$ is also the smallest positive integer such that $DP\subset Q$.

Let $q^{1/D}$ be a new variable, and set $q=\big(q^{1/D}\big)^D$, $q_i=q^{d_i}$. The quantum group $U_q = U_q(\mathfrak{g})$ is the algebra over $k=\mathbb {C}(q)$ with generators $E_i$, $F_i$, $K_i$, $K_i^{-1}$, $1\leq i \leq m$, and defining relations (see, e.g., \cite[Chapter~9.1]{CP}):
\begin{gather}
K_iK_j=K_jK_i ,\qquad K_iK_i^{-1}=K_i^{-1}K_i=1 ,\nonumber\\
K_iE_jK_i^{-1}=q_i^{a_{ij}}E_j ,\qquad K_iF_jK_i^{-1}=q_i^{-a_{ij}}F_j,\nonumber\\
E_iF_j-F_jE_i=\delta_{i,j}\frac{K_i-K_i^{-1}}{q_i-q_i^{-1}},\label{EFK}\\
\sum_{r=0}^{1-a_{ij}} (-1)^r \begin{bmatrix} 1-a_{ij} \\ r \end{bmatrix}_{q_i} E_i^{1-a_{ij}-r}E_jE_i^{r} = 0 \qquad {\rm if}\quad i\ne j,\label{Serre1}\\
\sum_{r=0}^{1-a_{ij}} (-1)^r \begin{bmatrix} 1-a_{ij} \\ r \end{bmatrix}_{q_i} F_i^{1-a_{ij}-r}F_jF_i^{r} = 0 \qquad {\rm if}\quad i\ne j.\label{Serre2}
\end{gather}
Here we put for $p$, $k$ non-negative integers with $0\leq k\leq p$,
\[
[0]_q! =1 ,\qquad [p]_q! =[1]_q[2]_q\cdots[p]_q ,\qquad [p]_q = \frac{q^p-q^{-p}}{q- q^{-1}},\qquad \begin{bmatrix} p \\ k \end{bmatrix}_{q} =\frac{[p]_q! }{[p-k]_q![k]_q!}.
\]
The algebra $U_q$ is a Hopf algebra with the coproduct $\Delta$, antipode $S$, and counit $\varepsilon$ defined by
\begin{gather*}
\Delta\big(K_i^{\pm 1}\big)=K_i^{\pm 1}\otimes K_i^{\pm 1} ,\qquad \Delta(E_i)=E_i\otimes K_i+1\otimes E_i ,\qquad \Delta(F_i)=K_i^{-1}\otimes F_i + F_i\otimes 1,\\ S(E_i) = -E_iK_i^{-1} ,\qquad S(F_i) = -K_iF_i ,\qquad S\big(K_i^{\pm 1}\big) = K_i^{\mp 1},\\ \varepsilon(E_i) = \varepsilon(F_i)=0,\qquad \varepsilon(K_i)=1.
\end{gather*}
In the case of $\mathfrak{g}={\mathfrak{sl}}(2)$ we simply write $E=E_1$, $F=F_1$ and $K^{\pm 1} = K_1^{\pm 1}$ the generators of~$U_q({\mathfrak{sl}}(2))$. For every $\alpha\in Q$, $\alpha=\sum_{i=1}^m m_i\alpha_i$, we put
\[K_\alpha=\prod_{i=1}^m K_i^{m_i}.\]
The element $\ell =K_{2\rho }$ is group like and satisfies $S^2(x)=\ell x \ell^{-1}$, $x\in U_q$. Hence~$U_q$ is a {\it pivotal} Hopf algebra, with the pivotal element $\ell$.

We will also consider the {\it simply-connected version} $\tilde U_q$ of $U_q$. It is obtained by adjoining invertible elements $\ell_i$, $1\leq i \leq m$, such that
\begin{gather*}
K_i = \prod_{j=1}^m \ell_j^{a_{ji}} ,\qquad \ell_iE_j\ell_i^{-1}=q_i^{\delta_{i,j}}E_j ,\qquad \ell_iF_j\ell_i^{-1}=q_i^{-\delta_{i,j}}F_j ,\qquad \Delta\big(\ell_i^{\pm 1}\big)=\ell_i^{\pm 1}\otimes \ell_i^{\pm 1}.
\end{gather*}
The expression of the pivotal element is
\[ \ell=\prod_{j=1}^m \ell_j^2.\]
For instance, $\tilde U_q({\mathfrak{sl}}(2))$ is obtained from $U_q({\mathfrak{sl}}(2))$ by adjoining a square root of $K$. We denote $K^{\frac{1}{2}}=\ell_1$.

We fix a reduced expression $s_{i_1}\cdots s_{i_N}$ of the longest element of the Weyl group of $\mathfrak{g}$. It~indu\-ces a total ordering of the positive roots,
\[\beta_1 = \alpha_{i_1}, \qquad \beta_2 = s_{i_1}(\alpha_{i_2}), \qquad \dots, \qquad \beta_N = s_{i_1}\cdots s_{i_{N-1}}(\alpha_{i_N}).\]
The root vectors of $\tilde U_q$ with respect to such an ordering are defined by
\[E_{\beta_k} = T_{i_1}\cdots T_{i_{k-1}}(E_{i_k}) ,\qquad F_{\beta_k} = T_{i_1}\cdots T_{i_{k-1}}(F_{i_k}), \]
where $T_i$ is Lusztig's algebra automorphism of~$\tilde U_q$ associated to the simple root $\alpha_i$ \cite{Lusztig2, Lusztig}, see also \cite[Chapter~8]{CP}). Let us just recall here that the monomials $F_{\beta_1}^{r_1}\cdots F_{\beta_N}^{r_N}K_\lambda E_{\beta_N}^{t_N}\cdots E_{\beta_1}^{t_1}$ ($r_i,t_i\in \mn$, $\lambda\in P$) form a basis of~$\tilde U_q$.

For every positive root $\beta$, we will denote $q_\beta = q^{(\beta,\beta)/2}$ and
\[\bar{E}_\beta = \big(q_\beta-q_\beta^{-1}\big)E_\beta ,\qquad \bar{F}_\beta = \big(q_\beta-q_\beta^{-1}\big)F_\beta .\]

\subsection{Categorical completions} \label{catcomb} By a theorem of Harish-Chandra, the map $\iota_U\colon U\to{\hat U}$ of Section~\ref{CATCOMP} is injective when $U$ is $U({\mathfrak g})$, the universal enveloping algebra of~${\mathfrak g}$ \cite[Theorem 2.5.7]{Dix}. It is also injective when $k={\mathbb C}(q)$ and $U=U_q$ \cite[Lemma~7.1.9]{Jos}.

In the case of $k={\mathbb C}(q)$, only specific categorical completions of $U_q$ are ribbon. We will use the one, that we denote by $\mUq$, defined as follows.

Recall that every finite dimensional $U_q$-module $X$ is completely reducible (see, e.g., \cite[Theorem~10.1.7]{CP}). If~$X$ is irreducible, it is a highest weight module parametrized by a dominant integral weight, say $\lambda\in P_+$, and signs $\sigma_i = \pm 1$, $i\in \{1,\dots,m\}$. Then $X = \oplus_{\e'} X_{\e'}$, where the sum ranges over the tuples $\e' = (\e_1',\dots,\e_m')$ such that $\e'_i{}^{-1}\e_i = q^{(\alpha,\alpha_i)}$ for some $\alpha\in Q_+$, where $\e = (\e_1,\dots,\e_m)$ and $\e_i = \sigma_iq^{(\lambda,\alpha_i)}$, and $X_{\e'}=\{x\in X\, \vert\, K_i\cdot x= \e_i' x,\, i=1,\dots,m\}$ is the weight space of $X$ of weight $\e'$. We say that a $U_q$-module $X$ has {\it type $1$} if it is finite dimensional and the signs $\sigma_i$ of all its irreducible components are equal to~$1$. Equivalently, the generators~$K_i$ are diagonalizable on~$X$, with eigenvalues in~$q_i^{\mz}$. The category~${\mathcal C}$ with objects the $U_q$-modules of type~$1$ is a~semisimple tensor category (where by semisimple we mean that any object is completely reducible). We will systematically denote by $V_\lambda$ the type~$1$ simple $U_q$-module parametrized by~$\lambda\in P_+$.

The algebra $\tilde U_q$ is defined over $\mc(q)$, and therefore also over~$\mc\big(q^{1/D}\big)$. The finite dimensional $\tilde U_q$-modules over $\mc\big(q^{1/D}\big)$ are completely reducible. Similarly as for~$U_q$, one defines the $\tilde U_q$-modules of type~$1$ as the finite dimensional $\tilde U_q$-modules over $\mc\big(q^{1/D}\big)$ on which the $\ell_i$ are diagonalizable with eigenvalues in $q_i^{(1/D)\mz}$. We denote by $\tilde {\mathcal C}$ the category of $\tilde U_q$-modules of type~$1$. Also, we denote by $\mathcal C \otimes \mc\big(q^{1/D}\big)$ the category obtained from~$\mathcal C$ by extending coefficients of objects and morphisms to $\mc\big(q^{1/D}\big)$. The restriction functor $\tilde{\mathcal C} \ra \mathcal C \otimes \mc\big(q^{1/D}\big)$ is an equivalence of categories.

This can be made explicit in the following way. As in Section~\ref{CATCOMP}, we can define the ${\mathcal C}$-categorical completion $\mUq'=\mUq(\mathfrak{g})'$ of $U_q$ as the Hopf algebra of natural transformations from $F_{\mathcal C}$ to $F_{\mathcal C}$, where $F_{\mathcal C}\colon {\mathcal C}\to \operatorname{FinVect}$ is the forgetful functor. Set
\[\mUq = \mUq' \otimes_k {\mathbb C}\big(q^{1/D}\big).\]
Extending the coefficient ring of $\mathcal{C}$ from $\mc(q)$ to $\mc\big(q^{1/D}\big)$ allows one to embed $\tilde U_q$ in $\mUq$. Let us explain this. The type $1$ simple $U_q$-module $V_\lambda$ can be considered as a $q$-deformation of the finite dimensional simple $\mathfrak{g}$-module of highest weight $\lambda$. Therefore, the generators $H_i\in \mathfrak{g}$ such that $\alpha_i(H_j)=a_{ji}$ define elements of $\mUq'$: $H_i$ acts on a weight space $X_{\epsilon'}$ of weight $\epsilon'=(\epsilon'_1,\dots,\epsilon'_m)$, where $\epsilon'_i = q^{(\mu,\alpha_i)}$ and $\mu\in P$, by multiplication by $(\mu,\alpha_i)$. Passing to the coefficient ring $\mc\big(q^{1/D}\big)$, one can define an element $q^{H_i/D}\in \mUq$, acting on $X_{\epsilon'}\otimes \mc\big(q^{1/D}\big)$ by multiplication by $q^{(\mu,\alpha_i)/D}$. Similarly, recalling that $(\mu,\omega_i)\in (1/D)\mz$, we can define the action of the generator $\ell_i \in \tilde U_q$ on $X_{\epsilon'}\otimes \mc\big(q^{1/D}\big)$ as the multiplication by $q^{(\mu,\omega_i)}$. This provides the claimed embedding $\tilde U_q\subset\mUq$.

Extending the coefficient ring of $\mathcal{C}$ from $\mc(q)$ to $\mc\big(q^{1/D}\big)$ allows also to make sure that $\mathcal{C}\otimes \mc\big(q^{1/D}\big)$ is a braided and ribbon category, and that $\mUq$ is a braided and ribbon Hopf algebra. To see this, denote by $B\in M_m(\mathbb{Q})$ the matrix with entries $B_{ij}:=d_j^{-1}a_{ij}$. By the same arguments as above one can define
\begin{equation}\label{exptermR}
\Theta := q^{\sum_{i,j=1}^m (B^{-1})_{ij} H_i\otimes H_j} \in \mUq^{\otimes 2}
\end{equation}
as the operator acting on $X_\epsilon \otimes X_{\epsilon'}\otimes \mc\big(q^{1/D}\big)$ as the multiplication by $q^{(\mu,\nu)} \in q^{(1/D)\mz}$, where $\epsilon=(\epsilon_1,\dots,\epsilon_m)$, $\epsilon'=(\epsilon_1',\dots,\epsilon_m')$, $\epsilon_i = q^{(\mu,\alpha_i)}$, $\epsilon_i' = q^{(\nu,\alpha_i)}$ and $\mu,\nu \in P$. Recall that the Drinfeld universal $R$-matrix of the quantized universal enveloping algebra $U_h(\mathfrak{g})$, defined, e.g., in \cite[Chapter~8.3.C]{CP} acts on tensor products $X\otimes Y$ of $U_q$-modules $X$ and $Y$ of type $1$, where $q=e^h$, and that it can be written in the form
\[R := \Theta \check{R},\]
where $\Theta$ is defined in \eqref{exptermR} and $\check{R}$, called the {\it quasi $R$-matrix}, is
\begin{equation}\label{Rcheck}
\check{R} = \sum_{t_1,\dots,t_N=0}^\infty \prod_{r=1}^N q_{\beta_r}^{\frac{1}{2}t_r(t_r+1)} \frac{\big(1-q_{\beta_r}^{-2}\big)^{t_r}}{[t_r]_{q_{\beta_r}} !} (E_{\beta_r})^{t_r} \otimes (F_{\beta_r})^{t_r}.
\end{equation}
Note that all terms in the above sum are zero but a finite number of them, since $E_{\beta_r}$, $F_{\beta_r}$ act nilpotently on objects of $\mathcal{C}$, so from \eqref{finsum} one derives immediately that
\[R = (R_{X,Y})\in \mUq^{\otimes 2}.\]
Observe that $R_{X,\cdot} \in \operatorname{End}(X) \otimes \tilde U_q$. Indeed, $\check{R}_{X,\cdot} \in \operatorname{End}(X) \otimes U_q$. Moreover, take weight spaces~$X_\epsilon$,~$X_{\epsilon'}$ as in \eqref{exptermR}, with $X_\epsilon$ a weight subspace of~$X$, and with basis vectors $\{e_i\}$ of $X_\epsilon$, and dual basis $\{e^i\}$. Then, setting $\lambda = \sum_{t=1}^m k_t\omega_t$, we have
\begin{equation}\label{Thetaeval}
\big(\langle e^j \mid \cdot \mid e_i\rangle \otimes {\rm id}\big)(\Theta) = \delta_{i,j}\prod_{t=1}^m \ell_t^{k_t} \in \tilde U_q.
\end{equation}
For instance, in the case of $\mathfrak{g}={\mathfrak{sl}}(2)$, we have $D=2$ and $\Theta = q^{H\otimes H/2}$. Hence, identifying elements of $U_q$ with their images by $\iota_{U_q}\colon U_q \ra \mUq$ we can write
\begin{equation}\label{Rsl2}
R=q^{H\otimes H/2} \sum_{n=0}^{+\infty}\frac{\big(q-q^{-1}\big)^n}{[n]_q!} q^{n(n-1)/2}E^n\otimes F^n.
\end{equation}
One checks similarly that the ribbon element of~$U_h(\mathfrak{g})$ induces an element of~$\mUq$. In conclusion, $U_q$ is naturally a braided and ribbon Hopf algebra in the categorical completion~$\mUq$.

\subsection{Integral forms and specializations}

Let ${A}={\mathbb C}\big[q,q^{-1}\big]$. The (De Concini--Kac) {\it unrestricted integral form} $U_A$ is the $A$-subalgebra of~$U_q$ generated by the elements (see, e.g., \cite[Chapter~9.2]{CP})
\begin{equation*}
E_i, \ F_i, \ K_i^{\pm 1}, \ \frac{K_i - K_i^{-1}}{q_i - q_i^{-1}}\qquad {\rm for\quad } 1\leq i\leq m.
\end{equation*}
We will also consider the smallest subalgebra $U_A’\subset U_A$ invariant
under the Lusztig automorphisms $T_i$ and containing the elements $\big(q_i-q_i^{-1}\big)E_i$, $\big(q_i-q_i^{-1}\big)F_i$ and $K_i^{\pm 1}$, ${\rm for\ } 1\leq i\leq m$ (which is the unrestricted integral form considered by De Concini--Procesi~\cite{DC-K-P2}).

The (Lusztig) {\it restricted integral form} $U_A^{\rm res}$ is the $A$-subalgebra of $U_q$ generated by the elements (see, e.g., \cite[Chapter~9.3]{CP})
\begin{equation*}
\frac{E_i^r}{[r]_{q_i} !} , \ \frac{F_i^r}{[r]_{q_i} !} , \ K_i^{\pm 1} \qquad \text{for any} \quad i, \, r \quad \text{with} \quad 1\leq i\leq m, r\geq 1.
\end{equation*}
Note that $U_A\subset U_A^{\rm res}$. Both $U_A$ and $U_A^{\rm res}$ have structures of Hopf $A$-algebras inherited from $U_q$, and they satisfy
\[U_A\otimes _{A}{\mathbb C}(q)=U_q ,\qquad U_A^{\rm res}\otimes _{A}{\mathbb C}(q)=U_q.\]
Moreover they are free $A$-modules, with certain basis of PBW type, described, e.g., in \cite[Chapter~9]{CP}. One defines $\tilde U_A$ in a similar way, replacing the generators $K_i^{\pm 1}$ by $\ell_i^{\pm 1}$.

We say a $U_A^{\rm res}$-module has type $1$ if it is a free $A$-module of finite rank with a basis made of eigenvectors for the elements $K_i$, with eigenvalues in $q^{\mz}$.
Denote by ${\mathcal C}_A$ the category with objects the $U_A^{\rm res}$-modules of type $1$. It is a tensor category, which is not semisimple.

\begin{remk}{\rm The fact that ${\mathcal C}_A$ is not semi simple can be seen on the following elementary example. Let $V_2$ be the two-dimensional type $1$ simple $U_q({\mathfrak{sl}}(2))$-module, with basis vectors~$v_+$,~$v_-$ such that $K.v_+ = qv_+$, $E.v_+=0$, $v_- = F.v_+$. Define $v' := v_-\otimes v_+ - qv_+\otimes v_-$ and $v'' := v_-\otimes v_+ + q^{-1}v_+\otimes v_-$. Then $v_+\otimes v_+$ generates a $U_A^{\rm res}({\mathfrak{sl}}(2))$-module $M$ of type $1$ and rank~$3$, namely $M=A v_+\otimes v_+ \oplus A v' \oplus Av_-\otimes v_-$, and $v''$ generates a $U_A^{\rm res}({\mathfrak{sl}}(2))$-module $N$ of type $1$ and rank $1$. The $A$-module $M\oplus N$ is stricly contained in the $U_A^{\rm res}({\mathfrak{sl}}(2))$-module $V_2 \otimes V_2$, which is clearly of type $1$, and it is not a direct summand thereof, for $(M\oplus N)\otimes \mc(q) = V_2 \otimes V_2$. So $V_2 \otimes V_2$ provides an example of $U_A^{\rm res}({\mathfrak{sl}}(2))$-module of type $1$ which is not completely reducible.}
\end{remk}

We denote by ${\mathcal C}_A \otimes \mc\big[q^{1/D},q^{-1/D}\big]$ the category obtained by extending coefficients of objects and morphisms to $\mc\big[q^{1/D},q^{-1/D}\big]$. We have seen above that ${\mathcal C} \otimes \mc\big(q^{1/D}\big)$ is a ribbon category; let us now explain why this ribbon structure descends to ${\mathcal C}_A\otimes {\mathbb C}\big[q^{1/D},q^{-1/D}\big]$. Define the categorical completion $\mathbb{U}_A^{\rm res}$ of $U_A^{\rm res}$ similarly as $\mUq$ above, that is
\[\mathbb{U}_A^{\rm res} = \mathbb{U}_A^{\rm res}{}' \otimes_A {\mathbb C}\big[q^{1/D},q^{-1/D}\big],\]
where $\mathbb{U}_A^{\rm res}{}'$ is the Hopf algebra of natural transformations $F_{\mathcal C_A}\to F_{\mathcal C_A}$, and $F_{\mathcal C_A}\colon {\mathcal C_A}\to \operatorname{FinMod}_A$ is the forgetful functor ($\operatorname{FinMod}_A$ being the category of finite dimensional $A$-modules).

We will use the following property of the $R$-matrix. Recall that $V_\lambda$ is the type $1$ simple $U_q$-module parametrized by~$\lambda\in P_+$. By a result of Lusztig \cite[Proposition~4.2]{Lusztig1}, $V_\lambda$ contains an $A$-submodule ${}_AV_\lambda$, which is a $U_A^{\rm res}$-module of type $1$ such that ${}_AV_\lambda\otimes \mc(q) = V_\lambda$ (hence ${}_AV_\lambda$ is a~full $A$-sublattice). Moreover ${}_AV_\lambda$ is endowed with a {\it canonical} basis $\mathcal{B}_\lambda$ (the {\it Kashiwara--Lusztig} basis, see, e.g., \cite[Chapter~14]{CP}).

\begin{teo}[integrality property] \label{integteo}
\quad
\begin{enumerate}\itemsep=0pt
\item[$1.$] For any $U_A^{\rm res}$-module $X$, $Y$ of type $1$, the quasi-$R$-matrix $\check{R}_{X,Y}$ is an automorphism of the $A$-module $X\otimes Y$, and $\Theta_{X,Y}$, whence $R_{X,Y}$, is an automorphism of $X\otimes Y \otimes {\mathbb C}\big[q^{1/D},q^{-1/D}\big]$.

\item[$2.$] Denote by $\mathcal{B}_{\lambda\mu}$ the basis of $V_\lambda\otimes V_\mu$ formed by the vectors $x\otimes y$, with $x\in \mathcal{B}_\lambda$, $y\in \mathcal{B}_\mu$. For every $\lambda,\mu\in P_+$, the matrix entries in the basis $\mathcal{B}_{\lambda\mu}$ of the endomorphisms $R_{V_\lambda,V_\mu}^{\pm 1} \in \operatorname{End}_{\mathbb{C}(q^{1/D})}(V_\lambda \otimes V_\mu)$ belong to $q^{\pm (\lambda,\mu)}\mc\big[q,q^{-1}\big]$.

\item[$3.$] For every $U_A^{\rm res}$-module $X$ of type $1$ we have $R_{X,\cdot} \in \operatorname{End}_A(X) \otimes \tilde{U}_A$.
\end{enumerate}
\end{teo}

\begin{proof} The claims (1) and (2) are classical results of Lusztig \cite[Proposition~24.1.4 and Corollary~24.1.6]{Lusztig}; see also \cite{Le}. They follow also from the formulas~\eqref{exptermR} and~\eqref{Rcheck}. The claim~(3) is a consequence of the fact that $(E_{\beta_r})^{t_r}/[t_r]_{q_{\beta_r}} !$ acts on the $A$-lattice $X$, and finally that $\Theta_{X,\cdot} \in \operatorname{End}_A(X) \otimes \tilde{U}_A$, proved via~\eqref{Thetaeval} above. This concludes the proof.
\end{proof}

By Theorem \ref{integteo}(1) the braiding $R$ of $\mathcal{C}\otimes \mc\big(q^{1/D}\big)$ yields a braiding on ${\mathcal C}_A \otimes {\mathbb C}\big[q^{1/D},q^{-1/D}\big]$, and from this one derives easily that the same is true of the ribbon transformation. As a~consequence, $\mathbb{U}_A^{\rm res}$ is a braided and ribbon Hopf algebra.

If $\epsilon\in {\mathbb C}^\times$, the {\it unrestricted specialisation} of $U_q$ at $q=\epsilon$ is the $\mc$-Hopf algebra
 \begin{equation}\label{restrictU}
 U_\epsilon =U_A\otimes_A {\mathbb C}_{\epsilon},
 \end{equation}
 where ${\mathbb C}_{\epsilon}={\mathbb C}$ as a~vector space, and as an $A$-module, $q$ acts on ${\mathbb C}_{\epsilon}$ by $\epsilon$. One defines $\tilde U_\epsilon$ in a~similar way.
 The map $\iota_U\colon U\to{\hat U}$ of Section~\ref{CATCOMP} is injective when $U=U_\epsilon$ or $\tilde U_\e$
 and $\epsilon$ is not a~root of unity \cite[Proposition~5.11]{Jant}.

When $\epsilon$ is not a root of unity, the finite dimensional $U_\epsilon$-modules are completely reducible (see, e.g., \cite[Theorem~10.1.14]{CP}). As above one can define the category ${\mathcal C}_\epsilon$ of $U_\epsilon$-modules of type $1$. The categorical ${\mathcal C}_\epsilon$-completion of $U_\epsilon$ is a braided and ribbon Hopf algebra. Note that, in the case of ${\mathfrak g }={\mathfrak{sl}}(2)$, again when $\epsilon$ is not a root of unity, a similar construction of categorical completion of $U_\epsilon({\mathfrak{sl}}(2))$ has been done in~\cite{FK}.

When $\epsilon$ is a root of unity, the category of finite dimensional $U_\epsilon$-modules is not semisimple and not braided. When ${\mathfrak g}={\mathfrak{sl}}(2)$ a classification of the simple modules is known \cite{DC-K, RA}. We~will use it extensively in a sequel to this paper.

\subsection{Quantum coordinate algebra}\label{OQ} Let $U$ be a Hopf algebra over a field $k$, $X$ an object of $\operatorname{FinMod}_U$, and $\pi_X$ the representation of~$U$ associated to~$X$. For every $v\in X$ and $w\in X^*$ we denote by ${}_X\phi{}^w_v\colon U \ra k$ the linear form defined by
\begin{equation*}
{}_X\phi{}^w_v(a) = w(\pi_X(a)v), \qquad a\in U.
\end{equation*}
The linear form ${}_X\phi{}^w_v$ is called a {\it matrix coefficient} of $\pi_X$. Denote by $U^\circ$ the {\it restricted dual} of~$U$, that is, the subspace of $U^*$ generated by the matrix coefficients of the finite dimensional representations of $U$. It is naturally endowed with a structure of Hopf algebra, such that the bilinear form $\langle\cdot,\cdot\rangle \colon U^\circ \times U \ra k$ is a Hopf pairing (see, e.g., \cite[Chapter~4]{CP}).

When $U=U_q$ and $k=\mc(q)$, we denote by $\Oq$ the Hopf subalgebra of $U_q^\circ$ generated by the matrix coefficients of the representations associated to the objects of ${\mathcal C}$, i.e., the $U_q$-modules of type~$1$. Equivalently one can define $\Oq$ as the set of ${\mathbb C}(q)$-linear maps
\[f\colon \ U_q\ra {\mathbb C}(q)\]
such that $\operatorname{Ker}(f)$ contains a two sided ideal $I\subset U_q$ of finite codimension, and there is an $r\in \mathbb{N}$ such that $\prod_{s=-r}^r (K_i - q_i^s) \in I$ for every $i$. We denote by $\star$ the product on~$\Oq$, induced from~$U_q^\circ$. Since the $U_q$-modules of type $1$ are completely reducible, and the simple ones are the highest weight modules $V_\lambda$, $\lambda\in P_+$, the set of matrix coefficients ${}_{V_{\lambda}}\phi{}^{e^j}_{e_i}$, where $\{e_i\}$ is any basis of $V_\lambda$ and $\{e^j\}$ the dual basis, provide a basis of~$\Oq$ over~$\mc(q)$. The $\mc(q)$-algebra $\Oq$ is also finitely generated by a finite number of elements, e.g., the matrix coefficients of the modules~$V_{\omega_i}$ associated to the fundamental weights $\omega_i$, $i=1,\dots,m$ (see, e.g., \cite[Proposition~10.1.16]{CP}); when $\mathfrak{g}$ is of type~$A$ or~$C$, the set of matrix coefficients of the natural representation $V_{\omega_1}$ already generates~$\Oq$ (see~\cite{Vara} in the classical setup).

Because the morphism $\iota_{U_q}\colon U_q\to \mathbb{U}_q$ is injective (see Section~\ref{catcomb}), the Hopf pairing $\langle\cdot ,\cdot \rangle \colon \Oq \times U_q \ra \mc(q)$ is non degenerate. By extending the coefficient ring from $\mc(q)$ to $\mc\big(q^{1/D}\big)$, we can uniquely extend it to a bilinear pairing
\[\langle\cdot ,\cdot \rangle\colon \ \big(\Oq \otimes_{{\mathbb C}(q)} {\mathbb C}\big(q^{1/D}\big)\big)\times \mUq \ra \mc\big(q^{1/D}\big)\]
such that the following diagram is commutative:
\[
\xymatrix{
 \Oq \otimes U_q \ar[r]^{\langle\cdot,\cdot \rangle}\ar[d]_{{\rm id} \otimes \iota_{ U_q}} & {\mathbb C}(q) \ar[d]\\
 \big(\Oq \otimes_{{\mathbb C}(q)} {\mathbb C}\big(q^{1/D}\big)\big)\otimes \mUq \ar[r]^{\;\;\;\;\;\;\;\ \;\ \;\;\;\ \langle\cdot,\cdot\rangle}& {\mathbb C}\big(q^{1/D}\big).
 }
\]
This pairing is defined by
\[\langle {}_Y\phi{}^w_v , (a_{X})\rangle=w(a_{Y}v)\]
for every $(a_{X})\in \mUq$ and ${}_Y\phi{}^w_v\in \Oq$. It is non degenerate. From now on, we will denote $\OqD := \Oq \otimes_{{\mathbb C}(q)} {\mathbb C}\big(q^{1/D}\big)$, and therefore the pairing above as $\langle\cdot ,\cdot \rangle\colon \OqD \times \mUq \ra \mc\big(q^{1/D}\big)$.

\subsection{Integral quantum coordinate algebra} Finally, we will use the integral form ${\mathcal O}_A$ of~$\Oq$ introduced in~\cite{Lusztig2}, and further studied in \cite{Lusztig} and \cite{DC-L} (see \cite[Remark 4.1]{DC-L} for the equivalence of the two definitions). One defines~${\mathcal O}_A$ as the $A$-module generated by the set of matrix coefficients of the $U_A^{\rm res}$-modules of type $1$. Clearly, it is contained in~$\Oq$; in fact, ${\mathcal O}_A$ is the subset of~$\Oq$ formed by the linear maps $f\colon U_q\ra {\mathbb C}(q)$ such that $f(U_A^{\rm res}) \subset \mc\big[q,q^{-1}\big]$ (see~\cite[Section~7]{Lusztig2}). Equivalently one can define~${\mathcal O}_A$ as the set of $A$-linear maps
\[f\colon \ U_A^{\rm res}\ra \mc\big[q,q^{-1}\big]\]
such that $\operatorname{Ker}(f)$ contains a cofinite two sided ideal $I\subset U_A^{\rm res}$ and there is an $r\in \mathbb{N}$ such that $ \prod_{s=-r}^r \big(K_i - q_i^s\big) \in I$ for every $i$ ($I$ cofinite meaning that there exists a free $A$-module $M$ of finite rank such that $I\oplus M = U_A^{\rm res}$).

The Hopf algebra structure of $\Oq$ descends to a structure of Hopf $A$-algebra on ${\mathcal O}_A$. It is a~full $A$-lattice in $\Oq$, that is, we have ${\mathcal O}_A\otimes_A {\mathbb C}(q) =\Oq$. The Hopf pairing $\langle\cdot,\cdot \rangle \colon \Oq \times U_q \ra \mc(q)$ restricts to a pairing of Hopf $A$-algebras $\Oo_A \otimes_A U_A^{\rm res} \ra A$. Because $U_A^{\rm res} \otimes_A \mc(q) = U_q$ and ${\mathcal O}_A\otimes_A {\mathbb C}(q) =\Oq$, this pairing is also non degenerate. We have described generators of $\Oo_q$ in Section \ref{OQ}. The $A$-algebra $\Oo_A$ is generated by the matrix coefficients of the lattices ${}_AV_\lambda$, $\lambda\in {P_+}$, introduced before Theorem \ref{integteo} (see \cite[Section~7.1]{Lusztig2}).

In \cite{Lusztig}, Lusztig proved that $\Oo_A$ is a free module over $A$ (in fact over $\mz[q,q^{-1}]$). He provided a basis of $\Oo_A$ by considering a multiplier Hopf $\mc(q)$-algebra $\stackrel{.}{\mathbf{U}}$ , containing a multiplier Hopf $A$-algebra $\stackrel{.}{\mathbf{U}}_A$ as a full lattice, such that the so-called {\it unital} $\stackrel{.}{\mathbf{U}}$-modules of finite dimension are exactly the $U_q$-modules of type $1$, and the unital $\stackrel{.}{\mathbf{U}}_A$-modules of finite rank are exactly the $U_A^{\rm res}$-modules of type $1$. The $A$-algebra $\stackrel{.}{\mathbf{U}}_A$ is a free $A$-module, it has a canonical basis $\stackrel{.}{\mathbf{B}}$, and the dual basis is a basis of $\Oo_A$.

When $\mathfrak{g}={\mathfrak{sl}}(2)$, $\Oo_A$ is generated by the matrix coefficients $a$, $b$, $c$, $d$ of the $2$-dimensional type~$1$ simple $U_q$-module $V_2$ in the basis $\{v_-,v_+\}$ introduced before Theorem~\ref{integteo}, and an $A$-basis is formed by the monomials $a^{\star r}\star b^{\star s} \star d^{\star t}$ and $a^{\star r}\star c^{\star u} \star d^{\star t}$, where $r,s,t\in \mn$, $u>0$ and $a^{\star r}$ is the product of $a$'s $r$ times, etc.\ (see \cite[Lemma~1.3]{DC-L}).

\section[The loop algebra L\_\{0,1\}(g)]{The loop algebra $\boldsymbol{\Ll_{0,1}(\mathfrak g)}$}

 The algebra $\Ll_{0,1} = \Ll_{0,1}(\mathfrak g)$ first appeared in relation with the reflection equation (see \cite{KS} and the references therein), and as the {\it braided group} associated to $\Oq$, which is the algebra of automorphisms of the ribbon category $\mathcal{C}\otimes \mc\big(q^{1/D}\big)$ of $U_q$-modules of type~$1$ (see \cite[Examples~7.4.1 and~9.4.10]{Majid}). In this paper we use a third definition of $\Ll_{0,1}$, by means of a twist of the $U_q$-bimodule $\Oq$ for the left and right coregular actions, following~\cite{DM} (but using right modules instead of left ones). One merit of this definition is to be intrinsic, not given by generators and relations. That it is recovers the braided group of~$\Oq$ is shown in \cite[Section~5.2]{DKM}. The twist uses $R$-matrices, whence the need of categorical completions. The reflection equation will be recovered in Proposition~\ref{fusrel}.

Denote the left and right coregular actions of $U_q$ on $\Oq$ by
\[
x\rhd \alpha := \sum_{(\alpha)}\alpha_{(1)} \langle \alpha_{(2)},x\rangle,\qquad \alpha \lhd x := \sum_{(\alpha)} \langle \alpha_{(1)},x\rangle \alpha_{(2)}
\]
for all $x\in U_q$ and $\alpha\in \Oq$, where $\langle\, ,\, \rangle\colon \Oq \otimes U_q \to \mc(q)$ is the duality pairing defined in Section~\ref{OQ}, and $\textstyle \Delta(\alpha) = \sum_{(\alpha)}\alpha_{(1)} \otimes \alpha_{(2)}$ (Sweedler's coproduct notation). The actions~$\rhd$ and~$\lhd$ commute, and~$\Oq$ is a $U_q$-{\it module algebra} for both actions, i.e., for all $x\in U_q$ and $\alpha,\beta\in \Oq$ we have
\[x\rhd (\alpha\star \beta)=\sum_{(x)}(x_{(1)}\rhd \alpha)\star(x_{(2)}\rhd \beta), \qquad (\alpha\star\beta)\lhd x=\sum_{(x)}(\alpha\lhd x_{(1)} )\star(\beta\lhd x_{(2)}).\]
All this extends to define a structure of $\mUq$-module algebra on $\OqD$. Denote by $\mUq^{\rm cop}$ the Hopf algebra with the same algebra structure as $\mUq$ but the opposite coproduct ${\Delta}^{\rm cop}$ and the antipode ${S}^{-1}$. Consider the Hopf algebra (equipped with the standard Hopf algebra structure on tensor product)
\[{\mathbb D}_q=\mUq {\otimes} \mUq^{\rm cop}.\]
It has a right action on $\OqD$ defined for pure tensors by
\begin{equation*}
\alpha\cdot (x\otimes y) :={S}(y) \rhd \alpha \lhd x
\end{equation*}
for every $x\in \mUq$, $y\in \mUq^{\rm cop}$ and $\alpha\in \OqD$, where $\otimes$ denotes the algebraic tensor product.

Recall the universal $R$-matrix $R\in \mUq {\otimes} \mUq$. Let
\[F:= (R')_{23}(R')_{24} \in \big(\mUq {\otimes} \mUq^{\rm cop}\big)^{{\otimes} 2},\]
where $(R')_{kl} = \mathfrak{i}_{kl}(R')$, and $\mathfrak{i}_{kl}\colon \mUq {\otimes} \mUq \to \big(\mUq {\otimes} \mUq^{\rm cop}\big)^{{\otimes}2}$ identifies the subalgebra $\mUq \otimes 1$ (resp.\ $1\otimes \mUq$) of $\mUq {\otimes} \mUq$ with the $k$-th (resp.\ $l$-th) tensorand of $\big(\mUq {\otimes} \mUq^{\rm cop}\big)^{{\otimes} 2}$. The tensor $F$ is a {\it twist} of $\mUq {\otimes} \mUq^{\rm cop}$: by definition, this means that $F$ is invertible and satisfies
\begin{equation}\label{twist0}
(\varepsilon_{{\mathbb D}_q} \otimes {\rm id})(F) =({\rm id} \otimes \varepsilon_{{\mathbb D}_q} )(F) = 1
\end{equation}
and
\begin{equation}\label{twist}
F_{12}({\Delta}_{{\mathbb D}_q} {\otimes} {\rm id})(F) = F_{23}({\rm id} {\otimes} {\Delta}_{{\mathbb D}_q})(F).
\end{equation}
Put $\textstyle u:=\sum_{(F)}F_{(1)}S(F_{(2)})$. Denote by ${S}_{{\mathbb D}_q}$ the antipode of ${\mathbb D}_q$ and set
\begin{equation}\label{newstructure}
{\Delta}_{{\mathbb D}_q}^F(x) := F{\Delta}_{{\mathbb D}_q}(x) F^{-1} ,\qquad { S}_{{\mathbb D}_q}^F(x) := u{S}_{{\mathbb D}_q}(x)u^{-1},\qquad x\in {\mathbb D}_q.
\end{equation}
The maps ${\Delta}_{{\mathbb D}_q}^F$, ${ S}_{{\mathbb D}_q}^F$ and the counit of ${\mathbb D}_q$ define the comultiplication, the antipode and the counit of a new structure of Hopf algebra on the algebra ${\mathbb D}_q$. It is called {\it the twist of ${\mathbb D}_q$ by $F$}. Denote it by ${\mathbb A}_q$. Since ${\mathbb A}_q$ and ${\mathbb D}_q$ coincide as algebras, the right action of ${\mathbb D}_q$ on $\OqD$ is also a right action of ${\mathbb A}_q$. Define a new product on $\OqD$ by
\begin{equation}\label{mmtilde-1}
\alpha \beta = \sum_{(F)}(\alpha \cdot F_{(1)})\star(\beta\cdot F_{(2)}).
\end{equation}
Because of \eqref{twist0} and \eqref{twist}, this product defines on $\OqD$ a structure of associative unital algebra. Explicitly, by writing $\textstyle R=\sum_{(R)}R_{(1) }\otimes R_{(2)}=\sum_{(R)}R_{(1') }\otimes R_{(2')}$
we have \[F=\sum_{(R),(R)}1\otimes R_{(2)}R_{(2 ')}\otimes R_{(1)}\otimes R_{(1')}.\]
Then, by using $({S}\otimes {S})(R) =R$ we get
\begin{equation}\label{mmtilde}
\alpha\beta = \sum_{(R),(R)}(R_{(2')}{S}(R_{(2)}) \rhd \alpha) \star (R_{(1')}\rhd \beta \lhd R_{(1)})
\end{equation}
and conversely
\begin{equation}\label{relprod2}
\alpha \star\beta =\sum_{(R),(R)} (R_{(2)}R_{(2')} \rhd \alpha) ({S}\big(R_{(1')}\big)\rhd \beta \lhd R_{(1)}).
\end{equation}
Note that, by the expression of $R$ and the fact that the generators $E_i$, $F_i$ of $U_q$ act nilpotently on finite dimensional $U_q$-modules, there is only a finite number of non zero terms in the last sum. Therefore the expression (\ref{mmtilde-1}) is well defined. It is easy to check that the right ${\mathbb A}_q$-module $\OqD$ endowed with the product \eqref{mmtilde-1} is a right ${\mathbb A}_q$-module algebra. We denote it $\Ll_{0,1}\big(q^{1/D}\big)$, and call it the {\it twist} by $F$ of the ${\mathbb D}_q$-module algebra $\OqD$.

We claim that the product~(\ref{mmtilde-1}) restricted to $\Oq\subset \OqD$ is defined over the subfield $\mc(q)$. Recall that~$\Oq$ is generated by the matrix coefficients ${}_{V_{\lambda}}\phi{}^{e^j}_{e_i}$, where $\lambda\in {P_+}$ and $\{e_i\}$, $\{e^i\}$ are dual basis of weight vectors of~$V_{\lambda}$. By equation~(\ref{mmtilde}), and using $({\rm id}\otimes S)(R)=({\rm id}\otimes \ell) R^{-1} \big({\rm id}\otimes \ell^{-1}\big)$, we have
\begin{equation}
\label{mmmtilde}
\alpha\beta = \sum_{(R),(R^{-1})}\big(R_{(2')}\ell \big(R^{-1}\big)_{(2)}\ell^{-1} \rhd \alpha\big) \star \big(R_{(1')}\rhd \beta \lhd \big(R^{-1}\big)_{(1)}\big).
\end{equation}
Therefore, if $\alpha= {}_{V_{\lambda}}\phi{}^{e^j}_{e_i}$ and $\beta= {}_{V_{\mu}}\phi{}^{f^n}_{f_m}$, we obtain
\begin{gather}
{}_{V_{\lambda}}\phi{}^{e^j}_{e_i}{}_{V_{\mu}}\phi{}^{f^n}_{f_m}= \sum_{(R),(R^{-1}), j', n',n'' }\big(\pi_{V_\lambda}\big(R_{(2')}\ell \big(R^{-1}\big)_{(2)}\ell^{-1}\big)\big)^{e^{j'}}_{e_i}\nonumber\\
\hphantom{{}_{V_{\lambda}}\phi{}^{e^j}_{e_i}{}_{V_{\mu}}\phi{}^{f^n}_{f_m}= \sum_{(R),(R^{-1}), j', n',n'' }}{}\times
 \big(\pi_{V_\mu}\big(R_{(1')}\big)\big)^{f^{n'}}_{f_{m}} \big(\pi_{V_\mu}\big(R^{-1}\big)_{(1)}\big)^{f^{n}}_{f_{n''}}
{}_{V_{\lambda}}\phi{}^{e^j}_{e_{j'}} \star {}_{V_{\mu}}\phi{}^{f^{n''}}_{f_{n'}}.\label{devprodq}
\end{gather}
Recall that the matrix entries of the endomorphisms $R_{V_\lambda,V_\mu}^{\pm 1} \in \operatorname{End}_{\mathbb{C}(q^{1/D})}(V_\lambda \otimes V_\mu)$ belong to $q^{\pm (\lambda,\mu)}\mc(q)$ (see Theorem~\ref{integteo}(2) for the stronger integral statement). Then we see that the factors $q^{+(\lambda,\mu)}$ cancel the factors $q^{-(\lambda,\mu)}$ in the last expression. Noting moreover that the matrix entries of $\ell$ belong to~$A$, we finally obtain that the coefficients in the sum belong to $\mc(q)$.
\begin{prop} \label{L01surq} The subspace $\Oq\subset \OqD$ endowed with the product~\eqref{mmtilde-1} is a right module $\mc(q)$-algebra over $U_q\otimes U_q^{\rm cop}\subset {\mathbb A}_q$. We denote it $\Ll_{0,1}(\mathfrak g)$, or simply $\Ll_{0,1}$.
\end{prop}

The following result is due to Majid \cite{Majid} (see \cite{DM} for a simpler proof):
\begin{prop} \label{copHopf} The coproduct $\Delta\!\colon\! \mUq \!\ra\! \mUq {\otimes} \mUq$ yields a morphism of Hopf algebras \mbox{$ \Delta\!\colon\! \mUq \!\to\! \mathbb{A}_q$}.
\end{prop}

By using the morphism of Hopf algebras $U_q \hookrightarrow \mUq \xrightarrow{\Delta} \mathbb{A}_q$ we can pull-back the action of $\mathbb{A}_q$ on $\Ll_{0,1}$, which thus becomes a right $U_q$-module algebra. It is easily seen that the action of $U_q$ is the right coadjoint action, defined by
\[\operatorname{coad}^r(x)(\alpha) = \sum_{(x)}{S}(x_{(2)}) \rhd \alpha \lhd x_{(1)},\qquad
\forall x \in U_q, \quad \forall \alpha\in \Ll_{0,1}.\]
Next we are going to recall a fundamental relation between $\Ll_{0,1}$ and $\tilde U_q$. Recall that $U_q$ is a~right $U_q$-module algebra for the right adjoint action, defined by
 \[\operatorname{ad}^r(y)(x) = \sum_{(y)}{S}(y_{(1)}) x y_{(2)}, \qquad \forall x,y \in U_q.\]
This action extends to actions on the simply-connected quantum group $\tilde U_q$ and $\mUq$, and thus defines on them structures of right $U_q$-module algebras. Denote by $\mathcal{Z}(\tilde U_q)$ the center of $\tilde U_q$. We~have
\begin{equation}\label{Zprop}
\mathcal{Z}\big(\tilde U_q\big)= \big\{ z \in \tilde U_q, \,\forall x\in U_q, \,\operatorname{ad}^r(x)(z)=\epsilon(x) z\big\}.
\end{equation}
The set of \emph{locally finite} elements of $\tilde U_q$ is defined by
\[
\tilde U_q^{\rm lf} :=\big\{x\in \tilde U_q\, \vert \, \operatorname{rk}_{\mc(q)}(\operatorname{ad}^r(U_q)(x)) < \infty\big\}.
\]
It is a $U_q$-module subalgebra of $\tilde U_q$.
Finally, the set of $\operatorname{coad}^r$-invariant elements of $\Ll_{0,1}$ is defi\-ned~by
\[
\Ll_{0,1}^{U_q} := \big\{\alpha \in \Ll_{0,1}\, \vert \, \forall y \in U_q,\, \operatorname{coad}^r(y)(\alpha) = \varepsilon(y)\alpha\big\}.
\]
Since $\Ll_{0,1}$ is a $U_q$-module algebra, $\Ll_{0,1}^{U_q}$ is a subalgebra of $\Ll_{0,1}$.

\begin{teo}\label{drinfeldmap} Define $\Phi_1\colon \Ll_{0,1}\ra \mUq$, $\Phi_1(\alpha) = (\alpha \otimes {\rm id})(RR')$. We have:
\begin{enumerate}\itemsep=0pt
\item[$1.$] $\Phi_1$ is a morphism of algebra, equivariant and injective.
\item[$2.$] The image of $\Phi_1$ is $\tilde U_q^{\rm lf}$.
\item[$3.$] $\Phi_1$ induces an isomorphism from $\Ll_{0,1}^{U_q}$ to $\mathcal{Z}\big(\tilde U_q\big)$.
\end{enumerate}
\end{teo}
We call $\Phi_1$ the \emph{RSD map}, after Reshetikhin and Semenov-Tian-Shansky \cite{RSTS} and Drinfeld \cite{Dr,Dr1}, who considered it first.

A proof of Theorem \ref{drinfeldmap} can be found in \cite[Theorem~3]{Bau1}. To make a complete correspondence with that statement, note that~\eqref{Thetaeval} implies that $\Phi_1$ takes values in $\tilde U_q$. The difficult parts of Theorem~\ref{drinfeldmap} are the injectivity of $\Phi_1$ and the claim~(2). Note that the third claim follows from the first two and~\eqref{Zprop}.

An alternative proof that $\Phi_1$ is an equivariant morphism is given in \cite[Proposition~4.7]{DM}, based on the construction of a left $U_q^{\rm op}$-comodule structure on~$\Ll_{0,1}$. The equivariance of~$\Phi_1$ also follows from a simple computation shown in Proposition~\ref{Alekseevmap} below.

Note that because the elements of the restricted dual are necessarily $\operatorname{coad}^r$ finite, the equi\-va\-ri\-ance of~$\Phi_1$ implies that their images by~$\Phi_1$ are necessarily $\operatorname{ad}^r$-finite.

\begin{remk} We give a simple self-contained proof of Theorem~\ref{drinfeldmap} in the~${\mathfrak{sl}}(2)$ case in Proposition~\ref{imagephi}.
\end{remk}

\begin{remk}\label{remlf} It is a result of~\cite{JL} that the $\operatorname{ad}^r(U_q)$-module $\tilde{U}_q^{\rm lf}$ is generated by the elements $\ell_{-\lambda}$, $\lambda\in 2P_+$, where for every weight $\mu = \sum_{i=1}^m n_i\omega_i$ we set $\ell_\mu = \prod_{i=1}^m \ell_i^{n_i}$. Moreover, there is an Ore subset $S$ of $\tilde{U}_q^{\rm lf}$ such that $\tilde{U}_q$ is a free module of finite rank over the skew fraction ring $S^{-1}\tilde{U}_q^{\rm lf}$ (see~\cite{JL2}). The set $S$ is the Abelian group generated by the elements $\ell_{-\lambda}$, where $\lambda \in 2P_+ \cap Q_+$.
\end{remk}

Finally, we provide a definition of $\Ll_{0,1}$ by generators and relations. Though well-known, we include a proof for completeness. As $\Ll_{0,1}$ and $\Oq$ are isomorphic as linear spaces, the matrix coefficients ${}_V\phi{}^{e^i}_{e_j}$ generate $\Ll_{0,1}$. Let $V$ be an object of the category ${\mathcal C}$, $(e_i)$ a basis of $V$, $(e^i)$ the dual basis, and $E_i^j$ the corresponding basis of $\operatorname{End}(V)$, defined by $E_i^j (e_k)=\delta_{j,k}e_i$. Define
\begin{gather}\label{matdef}
\stackrel{V}{M}=\sum_{i,j }E_i^j\otimes {}_V\phi{}^{e^i}_{e_j} \in \operatorname{End}(V)\otimes \Ll_{0,1}.
\end{gather}
This expression defines $\stackrel{V}{M}$ independently of the choice of basis $(e_i)$, since we have the following naturality property: if $V$, $W$ are objects of ${\mathcal C}$, and $f\colon V\to W$ a morphism of $U_q$-modules, then
\begin{gather}\label{naturality}
\stackrel{W}{M} (f\otimes {\rm id})=(f\otimes {\rm id})\stackrel{V}{M} \!.
\end{gather}
Given two objects $V$ and $W$ of $\mathcal{C}$ we write
\[
\stackrel{V}{M}_1 = \sum_{i,j }E_i^j\otimes {\rm id} \otimes {}_V\phi{}^{e^i}_{e_j} \in \operatorname{End}(V)\otimes \operatorname{End}(W) \otimes \Ll_{0,1}
\]
and similarly $ \stackrel{W}{M}_2 = \sum_{i,j }{\rm id} \otimes E_i^j\otimes {}_W\phi{}^{e^i}_{e_j}$. We view $\stackrel{V\otimes W}{M}$ as an element of $\operatorname{End}(V)\otimes \operatorname{End}(W)$ by using the standard isomorphism
\[
\operatorname{End}(V)\otimes \operatorname{End}(W)\cong \operatorname{End}(V\otimes W).
\]
Finally we view ${R}_{V,W}, {R'}_{V,W}\in \operatorname{End}(V)\otimes \operatorname{End}(W)$ as elements of $\operatorname{End}(V)\otimes \operatorname{End}(W)\otimes 1 \subset \operatorname{End}(V)\otimes \operatorname{End}(W)\otimes \Ll_{0,1}$.

\begin{prop}\label{fusrel} The following fusion relation holds true in $\operatorname{End}(V)\otimes \operatorname{End}(W) \otimes \Ll_{0,1}{:}$
\begin{equation}\label{fusionrel}
\stackrel{V\otimes W}{M}=\stackrel{V}{M}_1 {R'}_{V,W}\stackrel{W}{M}_2(R_{V,W}')^{-1},
\end{equation}
where the product of $\Ll_{0,1}$ is used to multiply the matrix elements of $\stackrel{V}{M}$ and $\stackrel{W}{M}$. It implies the reflection equation:
\begin{equation}\label{reflexrel}
{R}_{V,W}\stackrel{V}{M}_1 {R'}_{V,W}\stackrel{W}{M}_2\, =\, \stackrel{W}{M}_2 {R}_{V,W}\stackrel{V}{M}_1{R'}_{V,W}.
\end{equation}
Conversely, the naturality relations \eqref{naturality} and the fusion relations~\eqref{fusionrel}, for every objects $V$ and $W$ of $\mathcal{C}$, are a defining set of relations for $\Ll_{0,1}$. That is, $\Ll_{0,1}$ can be viewed as the quotient of the algebra freely generated over $k=\mc(q)$ by the matrix coefficients ${}_V\phi^{e^i}_{e_j}$, for all objects $V$ of the category $\mathcal{C}$, by the ideal generated by the relations~\eqref{naturality} and~\eqref{fusionrel}.
\end{prop}

\begin{proof} Let us write $\textstyle R_{V,W}=\sum_{(R)}R_{(1) }\otimes R_{(2)}=\sum_{(R)}R_{(1') }\otimes R_{(2')}$. The product of $\Oq$ is defi\-ned~by
\[{}_{(V\otimes W)}\phi{}^{e^i\otimes e^k}_{e_j\otimes e_l} = {}_V\phi{}^{e^i}_{e_j}\star {}_W\phi{}^{e^k}_{e_l}.\]
Then the relation \eqref{fusionrel} is equivalent to
\begin{gather*}
\sum_{i,j,k,l} E_i^j\otimes E_k^l\otimes {}_V\phi{}^{e^i}_{e_j}\star {}_W\phi{}^{e^k}_{e_l} = \sum_{(R),(R),i,j,k,l} E_i^jR_{(2)}R_{(2')}\otimes R_{(1)}E_k^l S\big(R_{(1')}\big)\otimes {}_V\phi{}^{e^i}_{e_j}{}_W\phi{}^{e^k}_{e_l}.
\end{gather*}
The isomorphism $V \otimes V^* \rightarrow \operatorname{End}(V)$, $v\otimes f\mapsto (w\mapsto f(w)v)$, maps $R_{(1)}e_k \otimes S\big(R_{(1')}\big)^*e^l$ to $R_{(1)}E_k^l S\big(R_{(1')}\big)$ and $e_i \otimes (R_{(2)}R_{(2')})^*e^j$ to $E_i^jR_{(2)}R_{(2')}$. Hence the above relation can be writ\-ten~as
\[
\sum_{i,j,k,l} E_i^j\otimes E_k^l\otimes {}_V\phi{}^{e^i}_{e_j}\star {}_W\phi{}^{e^k}_{e_l} = \sum_{(R),(R),i,j,k,l} E_i^j\otimes E_k^l \otimes \big({}_V\phi^{e^i}_{R_{(2)}R_{(2')}e_j}\big)\Big({}_V\phi^{R_{(1)} ^*e^k}_{ S(R_{(1')})e_l}\Big).
\]
Now we have
\[
{}_V\phi^{R_{(1)} ^*e^k}_{ S(R_{(1')})e_l} = S\big(R_{(1')}\big) \rhd {}_V\phi^{e^k}_{e_l} \lhd R_{(1)} ,\qquad {}_V\phi^{e^i}_{R_{(2)}R_{(2')}e_j} = R_{(2)}R_{(2')}\rhd {}_W\phi^{e^i}_{e_j},
\]
where we use the coregular actions $\rhd$, $\lhd$ and we denote now by $R_{(1)}, R_{(1')}, R_{(2)}, R_{(2')}\in \mUq$ the components of the universal R-matrix, instead of $R_{V,W}$. Identifying the matrix coefficients in $i$, $j$, $k$, $l$ we recover the relation~\eqref{relprod2}. Hence it is equivalent to the fusion relation, which thus provides a defining set of relations for $\Ll_{0,1}$.

Finally, note that by \eqref{naturality} we have
\[
\sigma_{V,W}R_{V,W} \stackrel{V\otimes W}{M} = \stackrel{W\otimes V}{M}\sigma_{V,W}R_{V,W},
\]
where $\sigma_{V,W} \colon V\otimes W \ra W \otimes V$ is the flip map. Then
\begin{align*}
\sigma_{V,W}R_{V,W}\stackrel{V}{M}_1 {R'}_{V,W}\stackrel{W}{M}_2\big(R_{V,W}'\big)^{-1} & = \, \stackrel{W}{M}_1 {R'}_{W,V}\stackrel{V}{M}_2\big(R_{W,V}'\big)^{-1} \sigma_{V,W}R_{V,W}\\
& = \sigma_{V,W} \stackrel{W}{M}_2 {R}_{V,W}\stackrel{V}{M}_1,
\end{align*}
which implies the reflection equation.
\end{proof}

\begin{remk}
As usual, denote by $a_V$ or $\pi_V(a)\in \operatorname{End}(V)$ the component of an element $a=(a_V) \in \mUq$ associated to the object $V$ of $\mathcal{C}$. The right coadjoint action of $U_q$ on $\Ll_{0,1}$ can be written in matrix form as (see, e.g., the proof of~\eqref{actmat2} below for a similar computation)
\begin{equation}\label{actmat}
\operatorname{coad}^r(y)\big(\!\stackrel{V}{M}\!\big) = \sum_{(y)}\big((y_{(1)_V} \otimes {\rm id})\stackrel{V}{M}( S(y_{(2)})_V \otimes {\rm id})\big).
\end{equation}
\end{remk}

Let $\lambda\in P_+$, $V_\lambda$ the type $1$ simple $U_q$-module of highest weight $\lambda$, and $\pi_{V_{\lambda}}$ the associated representation. Denote by $\operatorname{Tr}_{V_{\lambda}}\colon \operatorname{End}(V_{\lambda})\ra k $ the trace on $\operatorname{End}(V_{\lambda}).$ Put
\[\operatorname{qTr}_{V_\lambda}\big(\!\stackrel{V_\lambda}{M}\big) := (\operatorname{Tr}_{V_\lambda}\otimes {\rm id})\big(\big(\pi_{V_{\lambda}}(\ell)\otimes {\rm id}\big)\stackrel{V_\lambda}{M} \big)\in \Ll_{0,1}.\]
\begin{prop} \label{linbasecentre2} \quad
\begin{enumerate}\itemsep=0pt
\item[$1.$] The elements $\operatorname{qTr}_{V_\lambda}\big(\!\stackrel{V_\lambda}{M}\big)$, $\lambda\in P_+$, form a basis of $\Ll_{0,1}^{U_q}$.
\item[$2.$] The elements $(\operatorname{Tr}_{V_\lambda}\otimes {\rm id})\big((\pi_{V_{\lambda}} \otimes {\rm id})\big((\ell\otimes 1)(RR')\big)\big)$, $\lambda\in P_+$, form a basis of $ \mathcal{Z}\big(\tilde{U}_q\big)$.
 \end{enumerate}
\end{prop}

\begin{proof} The first part is an immediate consequence of Proposition~\ref{linbaseinv} below (namely, it is the case $n=1$, so that $a_{\lambda}(x)\in \operatorname{End}_{U_q}(V_{\lambda})$ is a scalar for every~$\lambda\in P_+$). The second part is a~consequence of the first and Theorem~\ref{drinfeldmap}. \end{proof}

\begin{remk}\label{centermatrk} Let $V_{\omega_1}$ be the fundamental representation of~$\tilde{U}_q$. The center $\mathcal{Z}\big(\tilde{U}_q\big)$ contains the elements (see \cite{Bau2,FRT})
\begin{equation*}
(\operatorname{Tr}_{V_{\omega_1}}\otimes {\rm id})\big((\pi_{V_{\omega_1}}\otimes {\rm id})\big((\ell\otimes {\rm id})\big((RR')^k\big)\big)\big), \qquad k\in \{1,\dots,m\}.
\end{equation*}
These elements generate $\mathcal{Z}\big(\tilde{U}_q\big)$ when $\mathfrak{g}$ is of type $A_m$ or $C_m$; see \cite{Bau2} for a more precise description in the other cases. Using Theorem~\ref{drinfeldmap} and
\begin{equation*}
({\rm id} \otimes \Phi_1)\big(\!\stackrel{V}{M}\!\big) = (\pi_V \otimes {\rm id})(RR')
\end{equation*}
we deduce that the elements
\begin{equation*}
{}_k\omega := (\operatorname{Tr}_{V_{\omega_1}}\otimes {\rm id})\big((\pi_{V_{\omega_1}}(\ell)\otimes {\rm id})\stackrel{V_{\omega_1}}{M}{}^{\!\! k} \big),\qquad k\in \{1,\dots,m\},
\end{equation*}
belong to and generate $\Ll_{0,1}^{U_q}$ when $\mathfrak{g}$ is of type $A_m$ or $C_m$. In particular, for $\mathfrak{g} = {\mathfrak{sl}}(2)$ we have $\Ll_{0,1}^{U_q} = \mc(q)[\omega]$, with
\begin{equation}\label{omegadef}
\omega = \operatorname{qTr}_{V_2}\big(\!\stackrel{V_2}{M} \big) = \operatorname{Tr}_{V_2}\big(K_{V_2}\stackrel{V_2}{M} \big) ,\qquad K_{V_2} := \begin{pmatrix} q & 0\\ 0&q^{-1}\end{pmatrix}\!,
\end{equation}
where $V_{2}$ is the $2$-dimensional type~1 simple $U_q({\mathfrak{sl}}(2))$-module, and as usual $K_{V_2}$ is the endomorphism of $V_2$ given by the action of~$K$.
\end{remk}

Next we define the integral form $\Ll_{0,1}^A$ of $\Ll_{0,1}$. Recall that ${\mathcal O}_A \subset \Oq$ is a Hopf $A$-subalgebra, and that it is a free $A$-module, and a full $A$-sublattice of $\Oq$.

\begin{prop} \label{defL01A} The space $\Oo_A$ endowed with the product \eqref{mmtilde-1} is an $A$-algebra. We denote it~$\Ll_{0,1}^A$. Moreover, $\Ll_{0,1}^A$ is the $A$-algebra generated by the matrix coefficients ${}_X\phi^{e^i}_{e_j}$, for all objects~$X$ of the category~$\mathcal{C}_A$, with defining relations~\eqref{naturality} and~\eqref{fusionrel}.
\end{prop}

Note that, since $\Ll_{0,1}^A$ coincides with $\Oo_A$ as an $A$-module, it is a free $A
$-module and we have $\Ll_{0,1}^A\otimes_A {\mathbb C}(q)=\Ll_{0,1}$.
\begin{proof} We have to check that products of elements of an $A$-basis of ${\mathcal O}_A$ can be expressed as linear combinations in this basis with coefficients in $A$. Denote by $\{\alpha_i\}$ a basis of $\Oo_A$. By~Pro\-po\-sition~\ref{L01surq} we have
\begin{equation}\label{mijk}
\alpha_i\alpha_j=\sum_k m^{k}_{ij}\alpha_k,
\end{equation}
where each $m^{k}_{ij}\in \mc(q)$. The basis elements $\alpha_i$ are linear combinations over $A$ of matrix coef\-fi\-cients of $U_A^{\rm res}$-modules of type $1$. Let $X$ and $Y$ be two such modules, and $\{e_i\}$,~$\{e^i\}$ and~$\{f_n\}$,~$\{f^n\}$ dual basis of~$X$ and~$Y$ respectively. As in~\eqref{devprodq} we have
\begin{gather*}
{}_{X}\phi{}^{e^j}_{e_i}{}_{Y}\phi{}^{f^n}_{f_m}= \sum_{(R),(R^{-1}), j', n',n'' }\big(\pi_{X}\big(R_{(2')}\ell (R^{-1})_{(2)}\ell^{-1}\big)\big)^{e^{j'}}_{e_i}
 \big(\pi_{Y}\big(R_{(1')}\big)\big)^{f^{n'}}_{f_{m}} \\
 \hphantom{{}_{X}\phi{}^{e^j}_{e_i}{}_{Y}\phi{}^{f^n}_{f_m}= \sum_{(R),(R^{-1}), j', n',n'' }}{}\times
 \big(\pi_{Y}(R^{-1})_{(1)}\big)^{f^{n}}_{f_{n''}}
{}_{X}\phi{}^{e^j}_{e_{j'}} \star {}_{Y}\phi{}^{f^{n''}}_{f_{n'}}.
\end{gather*}
The elements ${}_{X}\phi{}^{e^j}_{e_{j'}} \star {}_{Y}\phi{}^{f^{n''}}_{f_{n'}}$ being in~$\Oo_A$, they can expressed as linear combinations over $A$ of the basis elements~$\alpha_i$. Also, by Theorem~\ref{integteo}(1) the coefficients in this sum belong to $\mc\big[q^{1/D},q^{-1/D}\big]$. Therefore each coefficient $m^{k}_{ij}$ in~\eqref{mijk} belongs to $\mc\big[q^{1/D},q^{-1/D}\big]$. Since $\mc(q) \cap \mc\big[q^{1/D},q^{-1/D}\big] = \mc\big[q,q^{-1}\big]$, this proves the first claim. The arguments of Proposition~\ref{fusrel} apply as well to $\Ll_{0,1}^A$, which implies the second claim.\end{proof}

\begin{lem} \label{intform} The action $\operatorname{coad}^r$ yields on $\Ll_{0,1}^A$ a structure of right $U_A^{\rm res}$-module algebra, whence of $U_A$-module algebra. Moreover, $\Phi_1\big(\Ll_{0,1}^A\big)$ is a full $A$-sublattice of $\tilde U_q^{\rm lf}$, that we denote by $\tilde{U}_A^{\rm lf}$.
\end{lem}

\begin{proof} The first claim is clear, for the action $\operatorname{coad}^r$ endows the $A$-module $\Ll_{0,1}^A$ with a structure of right $U_A^{\rm res}$-module algebra, dual by the (non degenerate) pairing $\Oo_A\otimes_A U_A^{\rm res}\ra A$ to the structure defined by $\operatorname{ad}^r$ on $U_A^{\rm res}$. The inclusion $U_A \subset U_A^{\rm res}$ yields the second claim. The third follows from Theorem~\ref{drinfeldmap} and the fact that $\Ll_{0,1}^A$ is a full $A$-sublattice of $\Ll_{0,1}$. \end{proof}

\section[The example of L\_\{0,1\}(sl(2))]{The example of $\boldsymbol{\Ll_{0,1}({\mathfrak{sl}}(2))}$}\label{sl2}

In this section we provide a presentation by generators and relations of $\Ll_{0,1}({\mathfrak{sl}}(2))$, and using it we give an elementary proof of Theorem~\ref{drinfeldmap} in this case.

Let $V_r$ be the $r$-dimensional type $1$ simple $U_q({\mathfrak{sl}}(2))$-module. Put on $V_2$ the basis vectors~$v_+$,~$v_-$ such that $K.v_+ = qv_+$, $E.v_+=0$, $v_- = F.v_+$, and define in this basis
\begin{equation}\label{notsl20}
\stackrel{V_2}{M}=\begin{pmatrix} a&b\\c&d\end{pmatrix} \in \operatorname{End}(V_2)\otimes \Ll_{0,1}({\mathfrak{sl}}(2)).
\end{equation}

\begin{lem}\label{genrel01} The algebra $\Ll_{0,1}({\mathfrak{sl}}(2))$ is generated by the matrix elements $a$, $b$, $c$, $d$ of $\stackrel{V_2}{M}$. These satisfy the relations
\begin{alignat}{3}
&ad=da,\qquad && ab-ba = -\big(1-q^{-2}\big)bd,&\nonumber \\
& db = q^2 bd ,\qquad && cb-bc = \big(1-q^{-2}\big)\big(da-d^2\big), & \label{rel01}\\
& cd = q^2dc ,\qquad && ac-ca=\big(1-q^{-2}\big)dc,&\nonumber
\end{alignat}
as well as $ad-q^2bc=1$. Moreover $\omega=qa+q^{-1}d$ is central.
\end{lem}

Note that $\omega$ coincides with the element~\eqref{omegadef}.

\begin{proof} The family of matrix elements of $\stackrel{V}{M}$ for all type $1$ $U_q({\mathfrak{sl}}(2))$-modules $V$ spans $\Ll_{0,1}({\mathfrak{sl}}(2))$ over $\mc(q)$, since as a vector space it is the same as $\Oq({\mathfrak{sl}}(2))$. Any finite dimensional $U_q({\mathfrak{sl}}(2))$-module is completely reducible, and any simple one is a direct summand of some tensor power of~$V_2$. Hence the fusion relation~\eqref{fusionrel} implies that $\Ll_{0,1}({\mathfrak{sl}}(2))$ is generated by the matrix elements $a$, $b$, $c$, $d$. The relations~\eqref{rel01} follow easily from the reflection equation~\eqref{reflexrel} associated to~$V_2$, using the expression
 \[
 R_{V_2 ,V_2}=q^{-1/2}\begin{pmatrix} q &0 &0& 0\\0& 1&q-q^{-1}&0\\0&0&1&0\\0&0&0&q \end{pmatrix}\!, \] and they imply that $qa+q^{-1}d$ is central.
Because $V_2\otimes V_2$ admits the trivial representations as a subrepresentation, there exist non zero intertwiners
\[
\phi\colon \ V_2\otimes V_2\rightarrow \mc(q) , \qquad \psi\colon \ \mc(q)\rightarrow V_2\otimes V_2.
\]
As a consequence $\phi \ \circ \stackrel{V_2\otimes V_2}{M}\circ \psi$ is proportional to the unit element of $\Ll_{0,1}({\mathfrak{sl}}(2))$. Using the fusion relation, an easy computation provides the additional relation $ad-q^2bc=1$.
\end{proof}

Consider the RSD map $\Phi_1\colon \Ll_{0,1}({\mathfrak{sl}}(2))\ra \mUq({\mathfrak{sl}}(2))$, $\alpha \mapsto (\alpha\otimes {\rm id} )(RR')$. A straightforward computation using the expression (\ref{Rsl2}) shows that
\begin{alignat}{3}
& \Phi_1(a) = K+q^{-1}\big(q-q^{-1}\big)^2 FE,\qquad && \Phi_1(b)= q^{-1}\big(q-q^{-1}\big) F,&\nonumber\\
& \Phi_1(c) = \big(q-q^{-1}\big) K^{-1}E,\qquad && \Phi_1(d)= K^{-1}.&\label{phigen}
\end{alignat}

Therefore $\operatorname{Im}(\Phi_1)$ is contained in $U_q({\mathfrak{sl}}(2))$. The image of the central element $\omega$ is
\begin{gather}\label{omega}
\Omega = \Phi_1(\omega) = qK+q^{-1}K^{-1}+\big(q-q^{-1}\big)^2FE,
\end{gather}
which is $(q-q^{-1})^2$ times the standard Casimir element of $U_q({\mathfrak{sl}}(2))$.

Next we show that the relations~\eqref{rel01} and $ad-q^2bc=1$ yield a presentation of ${\mathcal L}_{0,1}({\mathfrak{sl}}(2))$.
Let $\tilde{\mathcal L}_{0,1}$ be the algebra generated by elements $\tilde{a}$, $\tilde{b}$, $\tilde{c}$, $\tilde{d}$ satisfying all these relations. De\-note~by
\[
j\colon \ \tilde{\mathcal L}_{0,1}\to {\mathcal L}_{0,1}({\mathfrak{sl}}(2))
\]
the unique morphism of algebra sending $\tilde{x}$ to $x$ for $x\in \{a,b,c,d\}.$

\begin{prop}\label{drinfeldinj} The monomials $\tilde{a}^\alpha \tilde{b}^\beta \tilde{c}^\gamma$ and $\tilde{d}^\delta \tilde{b}^\beta \tilde{c}^\gamma$, where $\alpha, \beta, \gamma, \delta\in \mn$ and $\alpha \geq 1$, form a basis of $\tilde{\Ll}_{0,1}$ over $\mc(q)$. Moreover, $j$ is an isomorphism, $\Phi_1$ is injective, and the center of~$\Ll_{0,1}({\mathfrak{sl}}(2))$ is the polynomial algebra $\mc(q)[\omega]$.
\end{prop}
\begin{proof} By inspection of the relations of Lemma \ref{genrel01}, it is easily seen that the given monomials form a generating set of $\tilde{\Ll}_{0,1}$. As for linear independence, consider the Verma $U_q({\mathfrak{sl}}(2))$-module~$M_x$, $x\in {\mathbb C}^{\times}$, with basis vectors $v_n$, $n\in \mn$, and action
\[Ev_0=0,\qquad Kv_n=x q^{-2n}v_n,\qquad Fv_n=v_{n+1}.\]
Then
\[Ev_n = [n]_q \frac{xq^{1-n}-x^{-1}q^{n-1}}{q-q^{-1}}v_{n-1},\qquad \Omega v_n = \big(qx+q^{-1}x^{-1}\big)v_n,\]
where $\Omega$ is as above. By the formulas~\eqref{phigen} we have
\[\Phi_1\big(a^\alpha b^\beta c^\gamma\big) = q^{-\alpha-\beta}\big(q-q^{-1}\big)^{\beta+\gamma} \big(\Omega-q^{-1}K^{-1}\big)^\alpha F^\beta\big(K^{-1}E\big)^\gamma\]
and
\[\Phi_1\big(d^\delta b^\beta c^\gamma\big) = q^{-\beta}\big(q-q^{-1}\big)^{\beta+\gamma} K^{-\delta} F^\beta\big(K^{-1}E\big)^\gamma.\]
Consider a linear relation with coefficients $A_{\alpha\beta\gamma},B_{\delta\beta'\gamma'}\in \mc(q)$:
\[\sum_{\alpha,\beta,\gamma}A_{\alpha\beta\gamma} \tilde{a}^\alpha \tilde{b}^\beta \tilde{c}^\gamma+
\sum_{\delta,\beta',\gamma'} B_{\delta\beta'\gamma'}\tilde{d}^\delta \tilde{b}^{\beta'} \tilde{c}^{\gamma'}=0.\]
Applying $\Phi_1\circ j$ we get (keeping the same names for the resulting coefficients in $\mc(q)$):
\begin{equation}\label{linrel}
\sum_{\alpha,\beta,\gamma,\gamma'}\big(A_{\alpha\beta\gamma} \big(\Omega-q^{-1}K^{-1}\big)^\alpha F^\beta\big(K^{-1}E\big)^\gamma + B_{\delta\beta'\gamma'} K^{-\delta} F^{\beta'}\big(K^{-1}E\big)^{\gamma'} \big)=0.
\end{equation}
By acting on the highest weight vector $v_0$ all terms on the left hand side vanish, but those with $\gamma=\gamma'=0$. Hence
\[
\sum_{\alpha,\beta}A_{\alpha\beta 0} \big(\Omega-q^{-1}K^{-1}\big)^\alpha F^\beta v_0 + \sum_{\delta\beta'}B_{\delta\beta'0} K^{-\delta} F^{\beta'} v_0 =0
\]
implying for each $\beta$ the relation
\[\sum_{\alpha}A_{\alpha\beta 0} \big(qx+q^{-1}x^{-1}\big(1-q^{2\beta}\big)\big)^\alpha + \sum_{\delta}B_{\delta\beta 0} x^{-\delta}q^{2\beta \delta}=0.\]
This is a Laurent polynomial in $x$. Since $\delta\geq 1$, the highest degree term in $x$ has vanishing coefficient, $A_{\alpha\beta 0}=0$, and hence $B_{\delta\beta0}=0$. So~\eqref{linrel} has no terms with $\gamma=\gamma'=0$. Then, by acting on~$v_1$ it results again that~\eqref{linrel} has no term with $\gamma=\gamma'=1$. Iterating this argument, an obvious recurrence implies that all the coefficients vanish, $A_{\alpha\beta\gamma} = B_{\delta\beta'\gamma'}=0$, which implies the linear independence of the monomials $\tilde{a}^\alpha \tilde{b}^\beta \tilde{c}^\gamma$ and $\tilde{d}^\delta \tilde{b}^\beta \tilde{c}^\gamma$ and therefore proves that they form a~basis of $\tilde{\Ll}_{0,1}$.

As a by-product we see that $\Phi_1\circ j $ is injective, and therefore $j$ is injective too. It is also surjective because the monomials ${a}^\alpha {b}^\beta {c}^\gamma$ and ${d}^\delta {b}^\beta c^\gamma$ form a generating family of $\Ll_{0,1}({\mathfrak{sl}}(2))$. It~follows that~$\Phi_1$ is injective as well.

Finally, let $z$ be in the center of $\Ll_{0,1}$. Then $\Phi_1(z)$ commutes with $\Phi_1(\Ll_{0,1}({\mathfrak{sl}}(2)))$, and hence with $K$, $E$, $F$. Thus it belongs to the center of $U_q({\mathfrak{sl}}(2))$, which is $\mc(q)[\Omega]$. Therefore
$\Phi_1(z)=P(\Omega)=P(\Phi_1(\omega)) =\Phi_1(P(\omega)).$ The result follows from the injectivity of $\Phi_1$.\end{proof}

\begin{remk}A consequence of this proposition is that the algebras $\Ll_{0,1}({\mathfrak{sl}}(2))$ and $U_q({\mathfrak{sl}}(2))$ are not isomorphic. Indeed the former has the family of one dimensional representations~$\rho_{xy}$ ($x\in \mc$, $y\in \mc^\times$), defined by \[\rho_{xy}(d)=0,\qquad \rho_{xy}(a)=x, \qquad \rho_{xy}(b)=-q^{-1}y, \qquad \rho_{xy}(c)=q^{-1}y^{-1}.\]
Then $K\notin \operatorname{Im}(\Phi_1)$, for otherwise $K=\Phi_1(\lambda)$, and injectivity of $\Phi_1$ and $K^{-1}=\Phi_1(d)$ would imply $\lambda d=d\lambda=1$. But this is not possible since $\rho_{xy}(d)=0$. In fact the family $\{\rho_{xy}\}_{x,y}$ cannot be obtained by pull-back of representations of $U_q({\mathfrak{sl}}(2))$ and the category of finite dimensional modules of $\Ll_{0,1}({\mathfrak{sl}}(2))$ is not semisimple (see \cite[Proposition 9]{BRT}).
\end{remk}

From the above results we can now derive an easy proof of the particular case of Theorem~\ref{drinfeldmap} for $\mathfrak{g} = {\mathfrak{sl}}(2)$. Note that $\tilde U_q^{\rm lf}({\mathfrak{sl}}(2)) = U_q^{\rm lf}({\mathfrak{sl}}(2))$ by the result of \cite{JL} recalled in Remark~\ref{remlf}. Also, the inclusion $U_q^{\rm lf}({\mathfrak{sl}}(2))\subset U_q({\mathfrak{sl}}(2))$ is strict, for $K\notin U_q^{\rm lf}({\mathfrak{sl}}(2))$.

\begin{prop}\label{imagephi} The RSD map yields an isomorphism of $U_q$-module algebras
\[\Phi_1\colon \ \Ll_{0,1}({\mathfrak{sl}}(2))\ra U_q^{\rm lf}({\mathfrak{sl}}(2)).\]
\end{prop}
\begin{proof} That $\Phi_1$ is an equivariant morphism follows from the same arguments as for the Alekseev map (see Theorem \ref{Alekseevmap}). Injectivity was shown in Proposition~\ref{drinfeldinj}. We prove that $\operatorname{Im}(\Phi_1) = U_q^{\rm lf}({\mathfrak{sl}}(2))$ by following closely the arguments of Section~3.11 of Joseph--Letzter \cite{JL}. For every integer $m\geq 0$ we have
\begin{equation}\label{formad0}
\operatorname{ad}^r(K)\big(\big(EK^{-1}\big)^m\big) = q^{-2m} \big(EK^{-1}\big)^m ,\qquad \operatorname{ad}^r(E)\big(\big(EK^{-1}\big)^m\big) = 0.
\end{equation}
Moreover,
\begin{gather*}
\operatorname{ad}^r(F)\big(EK^{-1}\big) = \frac{q}{q-q^{-1}}\big(\Omega-(q+q^{-1})K^{-1} \big), \qquad \operatorname{ad}^r\big(F^2\big)\big(EK^{-1}\big) = -\big(q+q^{-1}\big)F, \\
\operatorname{ad}^r(F^3)\big(EK^{-1}\big)=0, \qquad \operatorname{ad}^r(E)\big(K^{-1}\big) = EK^{-1}\big(q^{-2}-1\big), \\
\operatorname{ad}^r(F)\big(K^{-1}\big) = F\big(1-q^{-2}\big), \qquad
\operatorname{ad}^r(E)(F) = - \operatorname{ad}^r(F)\big(EK^{-1}\big), \qquad \operatorname{ad}^r(F)(F) = 0.
\end{gather*}

These relations imply that $EK^{-1}$, $K^{-1}$ and $F$ belong to $U_q^{\rm lf}({\mathfrak{sl}}(2))$. Because $\Omega$ is invariant under the action $\operatorname{ad}^r$,
the formulas \eqref{phigen} imply that $\Phi_1(a)$, $\Phi_1(b)$, $\Phi_1(c)$ and $\Phi_1(d)$ belong to $U_q^{\rm lf}({\mathfrak{sl}}(2))$. Therefore $\operatorname{Im}(\Phi_1) \subset U_q^{\rm lf}({\mathfrak{sl}}(2))$.

These above relations imply also that $\Phi_1(c)$ is a highest weight vector for the action $\operatorname{ad}^r$, generating a copy of $V_3$, the $3$-dimensional simple $U_q({\mathfrak{sl}}(2))$-module of type $1$.

Let us show by induction that $\operatorname{ad}^r(U_q({\mathfrak{sl}}(2)))\big(\Phi_1(c)^m\big) = V_{2m+1}$, for an arbitrary positive integer~$m$. In view of \eqref{formad0}, it remains to prove that $\operatorname{ad}^r\big(F^k\big)\big(\big(EK^{-1}\big)^m\big)\ne 0$ for all integers $k\leq 2m$, and $\operatorname{ad}^r\big(F^k\big)\big(\big(EK^{-1}\big)^m\big)= 0$ for $k\geq 2m+1.$ This holds true for $m=1$.

If this is true for a given $m$, then using the formula of $\Delta(F)$, the $q$-binomial identity, and the fact that $U_q({\mathfrak{sl}}(2))$ is an $\operatorname{ad}^r$-module algebra, we get
\begin{align*}
\operatorname{ad}^r\big(F^k\big)\big(\big(EK^{-1}\big)^{m+1}\big) & = \sum_{i=0}^{k} \left[\begin{matrix} k \\ i \end{matrix}\right]_{q^{-2}} \operatorname{ad}^r\big(K^{-i}F^{k-i}\big)\big(\big(EK^{-1}\big)^{m}\big) \operatorname{ad}^r\big(F^i\big)\big(EK^{-1}\big)\\
 &= \sum_{i=0}^{2} \left[\begin{matrix} k \\ i \end{matrix}\right]_{q^{-2}} \operatorname{ad}^r\big(K^{-i}F^{k-i}\big)\big(\big(EK^{-1}\big)^{m}\big) \operatorname{ad}^r\big(F^i\big)\big(EK^{-1}\big),
 \end{align*}
since $\operatorname{ad}^r\big(F^i\big)\big(EK^{-1}\big)=0$ if and only if $i\geq 3$. By the induction hypothesis, for $k\geq 2m+3$ each of the three terms of the sum vanishes. For $k=2m+2$, only the term for $i=2$ is non zero, which is equal to
\begin{equation*}
\left[\begin{matrix} k \\ 2 \end{matrix}\right]_{q^{-2}} \operatorname{ad}^r\big(K^{-2}F^{2m}\big)\big(\big(EK^{-1}\big)^{m}\big) \big({-}\big(q+q^{-1}\big)F\big).
\end{equation*}
By induction this term is non zero. Therefore $\operatorname{ad}^r(F^k)\big(\big(EK^{-1}\big)^m\big)\ne 0$ for all integers $k\leq 2m+2.$ This proves our claim.

Finally the multiplication map $\left( \oplus_{m\geq 0} V_{2m+1}\right) \otimes \mc(q)[\Omega] \lra U_q^{\rm lf}({\mathfrak{sl}}(2))$ is easily shown to be an isomorphism, as in \cite[Section~3.11]{JL}. Since $\Omega=\Phi_1(\omega)$, and $V_{2m+1}$ is generated by~$\Phi_1(c^{2m})$, we deduce the inclusion $U_q^{\rm lf}({\mathfrak{sl}}(2))\subset \operatorname{Im}(\Phi_1)$.\end{proof}

\begin{remk}If we denote ${\mathbb H}=\oplus_{m\geq 0} V_{2m+1},$ the isomorphism of modules \[{\mathbb H} \otimes \mc(q)[\Omega] \lra U_q^{\rm lf}({\mathfrak{sl}}(2))\] is an explicit example of the theorem of separation of variables of \cite{Bau1, JL2} in the case of ${\mathfrak{sl}}(2)$ (this case was first described in~\cite{JL}). As shown in these works, the multiplicity $[{\mathbb H}, V_k]$ for $k\geq 1$ is the dimension ($0$ or $1$) of the set of zero weight vectors of $V_k$.
\end{remk}

Let us make explicit the result of \cite{JL2} mentionned in Remark \ref{remlf}. By the relations \eqref{rel01} and Proposition \ref{drinfeldinj}, it is immediate that for all $x\in \Ll_{0,1}({\mathfrak{sl}}(2))$ there exist elements $y$, $y'\in \Ll_{0,1}({\mathfrak{sl}}(2))$ such that $dy = xd$ and $y'd=dx$. The element $d$ is regular because $\Phi_1(d)$ is invertible and $\Phi_1$ is injective, so $\{d^n\}_{n\in \mn}$ is a {\it left and right multiplicative Ore set} in $\Ll_{0,1}({\mathfrak{sl}}(2))$. The localization of~$\Ll_{0,1}({\mathfrak{sl}}(2))$ over $\{d^n\}_{n\in \mn}$ is well-defined (see \cite[Theorem 1.9 and Corollary 6.4]{GW}). Let us denote it by ${}_{\rm loc}\Ll_{0,1}({\mathfrak{sl}}(2))$.

\begin{prop}\label{Drinfeldloc} The map $\Phi_1\colon {}_{\rm loc}\Ll_{0,1}({\mathfrak{sl}}(2))\to U_q({\mathfrak{sl}}(2))$ defined by $\Phi_1\big(d^{-1}\big)=K$ and the formulas~\eqref{phigen} is an isomorphism of $U_q({\mathfrak{sl}}(2))$-module algebras.
\end{prop}

\begin{proof} Since $\Phi_1(d)$ is invertible, $\Phi_1$ extends to a morphism ${}_{\rm loc}\Ll_{0,1}({\mathfrak{sl}}(2)) \to U_q({\mathfrak{sl}}(2))$ uniquely (see \cite[Proposition~6.3]{GW}). It is an isomorphism, because the monomials $d^\delta b^\beta c^\gamma$ with $\beta,\gamma\in \mn$ and $\delta\in \mz$ make a basis of ${}_{\rm loc}\Ll_{0,1}({\mathfrak{sl}}(2))$, and they are sent by $\Phi_1$ to a PBW basis of $U_q({\mathfrak{sl}}(2))$.
\end{proof}

Finally consider the integral form $\Ll^A_{0,1}({\mathfrak{sl}}(2))$.

\begin{lem} \label{intsl2} The following holds:
\begin{enumerate}\itemsep=0pt
\item[$1.$] The $A$-algebra $\Ll^A_{0,1}({\mathfrak{sl}}(2))$ is generated by $a$, $b$, $c$, $d$ with the defining relations \eqref{rel01}.

\item[$2.$] The localization of $\Ll^A_{0,1}({\mathfrak{sl}}(2))$ over the set $\{d^n\}_{n\in \mn}$, that we denote by ${}_{\rm loc}\Ll^A_{0,1}({\mathfrak{sl}}(2))$, is generated by $a$, $b$, $c$, $d^{\pm 1}$ with the defining relations \eqref{rel01}.

\item[$3.$] The RSD map yields embeddings of $U_A$-module algebras $\Phi_1 \colon \Ll^A_{0,1}({\mathfrak{sl}}(2)) \ra U_A^{\rm lf}({\mathfrak{sl}}(2))$ and $\Phi_1 \colon {}_{\rm loc}\Ll^A_{0,1}({\mathfrak{sl}}(2)) \ra U_A({\mathfrak{sl}}(2))$, and $\Phi_1\big({}_{\rm loc}\Ll^A_{0,1}({\mathfrak{sl}}(2))\big) = U_A'({\mathfrak{sl}}(2))$ is the $A$-subalgebra $U_A'({\mathfrak{sl}}(2))\subset U_A({\mathfrak{sl}}(2))$, generated by $\big(q-q^{-1}\big)E$, $\big(q-q^{-1}\big)F$ and~$K^{\pm 1}$.
\end{enumerate}
\end{lem}
\begin{proof} (1) It is shown in \cite{DC-L}, Proposition 1.3, that $\Oo_A({\mathfrak{sl}}(2))\subset \Oo_q({\mathfrak{sl}}(2))$ is the $A$-subalgebra generated by $a$, $b$, $c$, $d$. Since $\Ll^A_{0,1}({\mathfrak{sl}}(2))$ coincides with $\Oo_A({\mathfrak{sl}}(2))$ as an $A$-module, any element of~$\Ll^A_{0,1}({\mathfrak{sl}}(2))$ is a linear combination over $A$ of monomials in $a$, $b$, $c$, $d$ with respect to the product~$\star$ of~$\Oo_A({\mathfrak{sl}}(2))$. By the relation \eqref{relprod2} and Proposition \ref{defL01A}, any such monomial is a linear combination over $A$ of monomials in $a$, $b$, $c$, $d$ with respect to the product (\ref{mmtilde-1}) of~$\Ll^A_{0,1}({\mathfrak{sl}}(2))$. This implies the first claim. The second follows from the fact that the relations~\eqref{rel01} are defined over~$A$. The assertion (2) is an immediate consequence of~(1), and~(3) follows from Lemma~\ref{intform}. We use the formulas~\eqref{phigen} to determine $\Phi_1\big({}_{\rm loc}\Ll^A_{0,1}({\mathfrak{sl}}(2))\big)$.\end{proof}

\section[The ''daisy'' graph algebra L\_\{0,n\}(g)]{The ``daisy'' graph algebra $\boldsymbol{\Ll_{0,n}(\mathfrak g)}$}\label{Lgnalg}

\subsection{Definition and first properties} We define $\Ll_{0,n}$ as a {\it twisted product} of $n$ copies of the $U_q\otimes U_q^{\rm cop}$-module algebra $\Ll_{0,1}$ (see Proposition~\ref{L01surq}), following \cite{DKM} (but using right modules instead of left ones). Since the twist uses $R$-matrices, as for $\Ll_{0,1}$ the construction uses as well categorical completions. Eventually we will see that $\Ll_{0,n}$ is also the braided tensor product of $n$ copies of the $U_q$-module algebra~$\Ll_{0,1}$ (see, e.g.,~\cite{Majid} for the notion of braided tensor product).

We need to recall a few notions. Consider Hopf algebras $A$ and $B$, and a {\it bicharacter} $F\in B\otimes A$. By definition, $F$ is an invertible tensor and satisfies
\begin{equation}\label{bichar}
(\Delta_{B} \otimes {\rm id}_{A})(F) = F_{23}F_{13} , \qquad ({\rm id}_{B} \otimes \Delta_{A})(F) = F_{12}F_{13} .
\end{equation}
Viewing $F$ as an element of $(1\otimes B)\otimes (A\otimes 1) \subset (A\otimes B)^{\otimes 2}$, it is readily checked that
\[F_{12}(\Delta_{A\otimes B} \otimes {\rm id})(F) = F_{23}({\rm id} \otimes \Delta_{A\otimes B})(F)\]
and \[(\varepsilon_{A\otimes B} \otimes {\rm id})(F) =({\rm id} \otimes \varepsilon_{A\otimes B} )(F) = 1.\]
Therefore $F$ is a twist of $A\otimes B$, endowed with the standard Hopf algebra product structure (see~\eqref{newstructure} for the similar operation applied to ${\mathbb D}_q$). Denote by $A\otimes^F B$ the resulting structure of Hopf algebra. Given a right $A$-module algebra $M$, and a right $B$-module algebra $N$, one defines the twisted tensor product $M\otimes^F N$ as the space $M\otimes N$ endowed with the product
 \begin{equation*}
 (\alpha \otimes \beta)(\alpha'\otimes \beta') = \sum_{(F)} \alpha (\alpha'\cdot F_{(2)}) \otimes (\beta\cdot F_{(1)})\beta'
 \end{equation*}
for every $\alpha, \beta \in M$, $\alpha', \beta' \in N$. This product gives $M\otimes^F N$ a structure of associative and unital right module algebra over $A\otimes^F B$, containing $M\otimes^F 1$ and $1 \otimes^F N$ as $A$- and $B$-module algebras respectively.
\begin{remk}\label{twistbraided} In the case where $A=B$ is quasi-triangular with $R$-matrix $R$, and $F= R'$, this construction gives the braided tensor product of $A$-module algebras.
\end{remk}

These constructions generalizes straightforwardly to the categorical completions we consider. Taking $A=B=\mUq$, the identities \eqref{Rdef0} for the universal $R$-matrix imply that $R'$ is a bicharacter of $\mUq {\otimes} \mUq$, and that we have (compare with~\eqref{twist})
\[ \mathbb{A}_q = \mUq {\otimes} ^{R'}\mUq.\]

We are going to iterate the above twist constructions. To this aim, observe that given homomorphisms of Hopf algebras $f_A\colon \mUq \ra A$, $f_B\colon \mUq \ra B$, the element $(f_B \otimes f_A)(R')$ is a~bicharacter of $A\otimes B$, and by Proposition \ref{copHopf} the map
\[f_A \odot f_B := (f_{A} \otimes f_{B})\circ \Delta \colon \ \mUq \ra A \otimes^{(f_B \otimes f_A)(R')} B\]
 is a morphism of Hopf algebras. Because of~\eqref{bichar}, taking $H:=A=B$ the operation $\odot$ endows the set of all morphisms of Hopf algebras $\mUq\ra H$ with a structure of associative monoid. One checks easily that it has the counit as neutral element. Clearly, $\odot$ is natural in the sense that $(g_A\circ f_A) \odot (g_B\circ f_B) = (g_A\otimes g_B) \circ (f_A\odot f_B)$ for any morphisms of Hopf algebras $g_A\colon A\ra A'$, $g_B\colon B\ra B'$.

In particular, by using the family of morphisms ${\rm id}_{\mUq}^{\odot n}$, $n\in \mn$, we define Hopf algebras $\mUq^{\odot n}$ as follows:
 \[\mUq^{\odot 0} = \mc(q^{1/D}) , \qquad \mUq^{\odot 1} = \mUq\]
and for $n\geq 1$,
\[\mUq^{\odot (n+1)} := \mUq^{\odot n} {\otimes}^{J_n} \mUq,\qquad {\rm where} \quad J_n:= \big({\rm id}_{\mUq} \otimes {\rm id}_{\mUq}^{\odot n}\big)(R').\]
Therefore, an immediate induction shows that $\mUq^{\odot n}$ is the twist of $\mUq^{\otimes n}$ by $\prod_{k=n-1}^{1} J_k$. Because of \eqref{bichar}, for every $k$, $l\in \mn$ we have an equality of Hopf algebras
\begin{equation}\label{assodot}
\mUq^{\odot (k+l)} = \mUq^{\odot {k}} \odot \mUq^{\odot {l}}.
\end{equation}
For instance, since ${\mathbb A}_q = \mUq^{\odot 2}$ and ${\rm id}_{\mUq}^{\odot 2} = \Delta$, the Hopf algebra ${\mathbb A}_q^{\odot 2} = \mUq^{\odot 4}$ is the twist of ${\mathbb A}_q{\otimes} \mathbb{A}_q$ by
\begin{equation}\label{formulaF}
F_2:= \big({\rm id}_{\mUq}^{\odot 2} \odot {\rm id}_{\mUq}^{\odot 2}\big)(R') = ({\Delta} {\otimes} {\Delta})(R')=R'_{23}R'_{13}R'_{24}R'_{14}
 \in \mathbb{A}_q {\otimes} \mathbb{A}_q.
 \end{equation}
 We are mainly concerned with the Hopf algebras $\mathbb{A}_q^{\odot n} = \mUq^{\odot 2n}$. Denote by $F_n$ the twist from~$\mathbb{A}_q^{\otimes n}$ to~$\mathbb{A}_q^{\odot n}$.

 We define $\Ll_{0,n}\big(q^{1/D}\big)$ as the twist by $F_n$ of the right $\mathbb{A}_q^{\otimes n}$-module algebra $\Ll_{0,1}\big(q^{1/D}\big)^{\otimes n}$ (endowed with the componentwise action and product). By construction $\Ll_{0,n}\big(q^{1/D}\big)$ is a right $\mathbb{A}_q^{\odot n}$-module algebra, and it coincides with $\Ll_{0,1}\big(q^{1/D}\big)^{\otimes n}$ as a $\mathbb{A}_q^{\odot n}$-module.

The above definition of $\Ll_{0,n}\big(q^{1/D}\big)$ by induction can be formulated by means of exchange relations between its factors $\Ll_{0,1}$. Let us explain this. For every $k$, $l$ the right $\mathbb{A}_q^{\odot n}$-module algebra $\Ll_{0,k}\big(q^{1/D}\big) \odot \Ll_{0,l}\big(q^{1/D}\big)$ is the twist by $\textstyle F(l,k) = ({\Delta}^{\odot l} \otimes {\Delta}^{\odot k})(R')$ of the $\mathbb{A}_q^{\odot k} {\otimes} \mathbb{A}_q^{\odot l}$-module algebra $\Ll_{0,k}\big(q^{1/D}\big) \otimes \Ll_{0,l}\big(q^{1/D}\big)$. By associativity of ${\odot}$ (see~\eqref{assodot}) we have
 \begin{equation*}
 \Ll_{0,k}\big(q^{1/D}\big) \odot \Ll_{0,l}\big(q^{1/D}\big) = \Ll_{0,k+l}\big(q^{1/D}\big).
 \end{equation*}
 Under this factorisation, by setting $k+l=n$, for every $\alpha, \alpha' \in \Ll_{0,k}\big(q^{1/D}\big)$ and $\beta,\beta'\in \Ll_{0,l}\big(q^{1/D}\big)$ the product of $\Ll_{0,n}\big(q^{1/D}\big)$ takes the form
 \begin{equation}\label{braidprod}
 (\alpha \otimes \beta)(\alpha'\otimes \beta') = \sum_{(F(l,k))} \alpha (\alpha'\cdot F(l,k)_{(2)}) \otimes (\beta\cdot F(l,k)_{(1)})\beta',
 \end{equation}
where
 the products $\alpha (\alpha'\cdot F(l,k)_{(2)})$ and $(\beta\cdot F(l,k)_{(1)})\beta'$ are taken in $\Ll_{0,k}\big(q^{1/D}\big)$ and $\Ll_{0,l}\big(q^{1/D}\big)$ respectively. Moreover, for every $1\leq a\leq n$ the map $\mathfrak{i}_a\colon \Ll_{0,1}\big(q^{1/D}\big)\ra \Ll_{0,n}\big(q^{1/D}\big)$, identifying $\Ll_{0,1}\big(q^{1/D}\big)$ with the $a$-th factor of $\Ll_{0,1}\big(q^{1/D}\big)^{\otimes n}$ by putting $1$'s elsewhere \big($1$ being the unit of $\Ll_{0,1}\big(q^{1/D}\big)$\big), is an embedding of module algebras. We will use the notations
\[\Ll_{0,n}\big(q^{1/D}\big)^{(a)}:= \operatorname{Im}(\mathfrak{i}_a) ,\qquad (\alpha)^{(a)} := \mathfrak{i}_a(\alpha).\]
Taking $(\alpha)^{(a)} \in \Ll_{0,n}\big(q^{1/D}\big)^{(a)}$ and $(\beta)^{(b)} \in \Ll_{0,n}\big(q^{1/D}\big)^{(b)}$ with $a<b$, and $F=F(1,1)\in \mathbb{A}_q{\otimes}\mathbb{A}_q$ as in~\eqref{formulaF}, we deduce from~\eqref{braidprod} that \big(all products being taken in $\Ll_{0,n}\big(q^{1/D}\big)$\big)
\begin{equation}\label{commutrel}
(\beta)^{(b)}(\alpha)^{(a)} = \sum_{(F)} \big(\alpha\cdot F_{(2)}\big)^{(a)} \big(\beta\cdot F_{(1)}\big)^{(b)}.
\end{equation}
We have:
\begin{prop}\label{defiL0n} The subspace $\Ll_{0,1}^{\otimes n}\subset \Ll_{0,n}\big(q^{1/D}\big)$ endowed with the product defined by \eqref{mmtilde-1} on each component and by the exchange relation \eqref{commutrel} between distinct components is a right module $\mc(q)$-algebra over $\big(U_q\otimes U_q^{\rm cop}\big)^{\odot n}\subset \mathbb{A}_q^{\odot n}$. We denote it $\Ll_{0,n}$. Moreover, $\Ll_{0,n}$ is a $U_q$-module algebra, namely the braided tensor product of~$n$ copies of the $U_q$-module algebra $\Ll_{0,1}$, with $U_q$-action, that we still denote by $\operatorname{coad}^r$, given by
\begin{gather}
\operatorname{coad}^r(y)\big(\alpha^{(1)}\otimes \dots \otimes \alpha^{(n)}\big) = \sum_{(y)} \operatorname{coad}^r(y_{(1)})\big(\alpha^{(1)}\big) \otimes \dots \otimes \operatorname{coad}^r(y_{(n)})\big(\alpha^{(n)}\big)\label{actionprod}
\end{gather}
for all $y \in U_q$ and $\alpha^{(1)}\otimes \dots \otimes \alpha^{(n)} \in \Ll_{0,n}$.
\end{prop}
 \begin{proof} We have to show that the product \eqref{commutrel} restricted to $\Ll_{0,1}^{\otimes n}\subset \Ll_{0,n}\big(q^{1/D}\big)$ is defined over the subfield $\mc(q)$. Indeed, the injection $\mathfrak{i}_a\colon \Ll_{0,1}\ra \Ll_{0,n}\big(q^{1/D}\big)$ being a morphism of algebras, the product of elements $\alpha, \beta\in \mathfrak{i}_a(\Ll_{0,1})$ is defined over $\mc(q)$ (see Proposition \ref{L01surq}). Therefore it is enough to check that the exchange relation \eqref{commutrel} is defined over $\mc(q)$ as well. We proceed as in~\eqref{mmmtilde} and~\eqref{devprodq}. For matrix coefficients $\big({}_{V_\lambda}\phi^{e^i}_{e_j}\big)^{(a)}$ and $\big({}_{V_\mu}\phi^{e^k}_{e_l}\big)^{(b)}$, where $a<b$, \eqref{commutrel} can be written as
 \begin{gather*}
\big({}_{V_\mu}\phi^{e^k}_{e_l}\big)^{(b)} \big({}_{V_\lambda}\phi^{e^i}_{e_j}\big)^{(a)}
\\ \qquad
{}= \sum_{(R^1)\cdots (R^4)} \big(\big({}_{V_\lambda}\phi^{e^i}_{e_j}\big)^{(a)}\cdot \big(R^1_{(1)}R^2_{(1)} \otimes R^3_{(1)}R^4_{(1)}\big) \big)\big(\big({}_{V_\mu}\phi^{e^k}_{e_l}\big)^{(b)}\cdot \big(R^2_{(2)}R^4_{(2)} \otimes R^1_{(2)}R^3_{(2)}\big)\big)
\\ \qquad
{}= \sum_{(R^1)\cdots (R^4)}\!\!\! \big(S\big(R^3_{(1)}R^4_{(1)}\big)\!\rhd \big({}_{V_\lambda}\phi^{e^i}_{e_j}\big)^{(a)}\! \lhd R^1_{(1)}R^2_{(1)}\big)\big(S\big(R^1_{(2)}R^3_{(2)}\big)\!\rhd \big({}_{V_\mu}\phi^{e^k}_{e_l}\big)^{(b)}\!\lhd R^2_{(2)}R^4_{(2)}\big)
\\ \qquad
{} = \sum_{(R^1)\cdots (R^4), i',j',k',l'} \pi_{V_\lambda}\big(R^1_{(1)}R^2_{(1)}\big)^{e^i}_{e_{i'}} \ \pi_{V_\lambda}\big(S\big(R^3_{(1)}R^4_{(1)}\big)\big)^{e^{j'}}_{e_{j}}
 \\[-2ex] \qquad\hphantom{= \sum_{(R^1)\cdots (R^4), i',j',k',l'}}
{}\times \pi_{V_\mu}\big(R^2_{(2)}R^4_{(2)}\big)^{e^k}_{e_{k'}} \pi_{V_\mu}\big(S\big(R^1_{(2)}R^3_{(2)}\big)\big)^{e^{l'}}_{e_{l}}\ \big({}_{V_\lambda}\phi^{e^{i'}}_{e_{j'}}\big)^{(a)} \ \big({}_{V_\mu}\phi^{e^{k'}}_{e_{l'}}\big)^{(b)}
\\ \qquad
{} = \sum_{(R^1)\cdots (R^4), i',j',k',l'} \pi_{V_\lambda}\big(R^1_{(1)}\!\!{}^{-1}R^2_{(1)}\big)^{e^i}_{e_{i'}} \ \pi_{V_\lambda}\big(R^4_{(1)}\!\!{}^{-1} R^3_{(1)}\big)^{e^{j'}}_{e_{j}}
 \\[-2ex] \qquad\hphantom{=\sum_{(R^1)\cdots (R^4), i',j',k',l'}}
{}\times \pi_{V_\mu}\big(R^2_{(2)}R^4_{(2)}\!\!{}^{-1}\big)^{e^k}_{e_{k'}} \pi_{V_\mu}\big(R^3_{(2)}\ell R^1_{(2)}\!\!{}^{-1} \ell^{-1}\big)^{e^{l'}}_{e_{l}}\ \big({}_{V_\lambda}\phi^{e^{i'}}_{e_{j'}}\big)^{(a)} \big({}_{V_\mu}\phi^{e^{k'}}_{e_{l'}}\big)^{(b)},
\end{gather*}
where we denote by $R^1,\dots,R^4$ the four $R$-matrices in \eqref{formulaF}, $R^i = \sum_{(R)} R_{(1)}^i \otimes R_{(2)}^i$ as usual, and $\ell$ is the pivotal element.
Since the matrix entries of $R_{V_\lambda,V_\mu}^{\pm 1} \in \operatorname{End}_{\mathbb{C}(q^{1/D})}(V_\lambda \otimes V_\mu)$ belong to $q^{\pm (\lambda,\mu)}\mc(q)$ (see Theorem \ref{integteo}(2) for the stronger integral statement), the factors $q^{+(\lambda,\mu)}$ cancel the factors $q^{-(\lambda,\mu)}$ in the last expression. The matrix entries of~$\ell$ belong to~$A$, so we finally obtain that the coefficients in the sum belong to~$\mc(q)$.

Finally, as for the action of $U_q$ on $\Ll_{0,n}$, consider the structure of $U_q$-module algebra on $\Ll_{0,1}$ obtained by pulling-back the action of $\big(U_q\otimes U_q^{\rm cop}\big)$ via $\Delta$ (see the comment after Proposition~\ref{copHopf}). It follows from Remark~\ref{twistbraided} that $\Ll_{0,n}=\Ll_{0,1}^{\odot n}$ is the braided tensor of $n$ copies of $\Ll_{0,1}$. The restriction of the action of $U_q$ on each factor $\Ll_{0,1}$ being $\operatorname{coad}^r$, the formula~\eqref{actionprod} is a consequence of the properties of actions of module algebras.
\end{proof}

The commutation relation \eqref{commutrel} yields a presentation of $\Ll_{0,n}$ by generators and relations. This presentation is well-known (see \cite{AGS1,BR1, DKM,MW}), but for completeness we prefer to give a~proof. Let $V$ be an object of the category ${\mathcal C}$. Similarly to \eqref{matdef}, define
\[\stackrel{V}{M}{}^{\!\! (a)}=\sum_{i,j }E_i^j\otimes \big({}_V\phi^{e^i}_{e_j}\big)^{(a)} = \sum_{i,j }E_i^j\otimes \big(1^{\otimes (a-1)} \otimes ({}_V\phi{}^{e^i}_{e_j}) \otimes 1^{\otimes (n-a)}\big) \in \operatorname{End}(V)\otimes \Ll_{0,n}.\]
Note that the matrix coefficients of the set of matrices $\stackrel{V}{M}{}^{\!\! (a)}$, for every object $V$ of $\mathcal{C}$ and $1\leq a\leq n$, generate the algebra $\Ll_{0,n}$.
\begin{prop}\label{matrelL0n} For every $a<b$ the matrices $\stackrel{V}{M}{}^{\!\! (a)}$ satisfy the fusion equation \eqref{fusionrel}, the naturality relations \eqref{naturality}, and the following \emph{exchange relation} $($in $\operatorname{End}(V)\otimes \operatorname{End}(W) \otimes \Ll_{0,n})$
\begin{equation}\label{exrel2}
{R}_{V,W}\stackrel{V}{M_1}{}^{\!\!\! (a)} {R}_{V,W}^{-1}\stackrel{W}{M_2}{}^{\!\!\! (b)}= \stackrel{W}{M_2}{}^{\!\! (b)} {R}_{V,W}\stackrel{V}{M_1}{}^{\!\!\! (a)} R_{V,W}^{-1}.
\end{equation}
Moreover, all these relations determine the product of $\Ll_{0,n}$. Hence the algebra $\Ll_{0,n}$ can be viewed as the quotient of the algebra freely generated over $\mc(q)$ by the matrix coefficients $\big({}_V\phi^{e^i}_{e_j}\big)^{(a)}$, for all objects $V$ of $\mathcal{C}$ and every $1\leq a\leq n$, by the ideal generated by the relations \eqref{naturality}, \eqref{fusionrel} and~\eqref{exrel2}.
\end{prop}
\begin{proof} The fusion and naturality relations follow from the fact that $\mathfrak{i}_a\colon \Ll_{0,1}\ra \Ll_{0,n}$ is a morphism of algebras. The matrix coefficients $\big({}_V\phi^{e^i}_{e_j}\big)^{(a)}$ generate the subalgebras $\Ll_{0,n}^{(a)}$, whence $\Ll_{0,n}$ too. Conversely, by Proposition \ref{fusrel} the fusion and naturality relations determine the product of~$\Ll_{0,n}^{(a)}$. Hence it is enough to show that the exchange relations \eqref{exrel2} and the commutation relation~\eqref{commutrel} are equivalent. Let us write $\textstyle R_{V,W}=\sum_{(R)}R_{(1) }\otimes R_{(2)}=\sum_{(R)}R_{(1') }\otimes R_{(2')}$. Then~\eqref{exrel2} is equivalent to
\begin{gather*}
\sum_{(R),(R), i,j,k,l} R_{(1)}E_i^j S\big(R_{(1')}\big) \otimes R_{(2)}R_{(2')}E_k^l \otimes \big({}_V\phi^{e^i}_{e_j}\big)^{(a)} \big({}_W\phi^{e^k}_{e_l}\big)^{(b)}
\\ \qquad
{}= \sum_{(R),(R),i,j,k,l} R_{(1)}E_i^j S\big(R_{(1')}\big) \otimes E_k^l R_{(2)}R_{(2')} \otimes \big({}_W\phi^{e^k}_{e_l}\big)^{(b)}\big({}_V\phi^{e^i}_{e_j}\big)^{(a)}.
\end{gather*}
The isomorphism $V \otimes V^* \rightarrow \operatorname{End}(V)$, $v\otimes f\mapsto (w\mapsto f(w)v)$, maps $R_{(1)}e_i \otimes S\big(R_{(1')}\big)^*e^j$ to $R_{(1)}E_i^j S\big(R_{(1')}\big)$, $R_{(2)}R_{(2')}e_k \otimes e^l$ to $R_{(2)}R_{(2')}E_k^l$, and $e_k \otimes (R_{(2)}R_{(2')})^*e^l$ to $E_k^lR_{(2)}R_{(2')}$. Hence the above relation can be written as
\begin{gather*}
\sum_{(R),(R), i,j,k,l} E_i^j \otimes E_k^l \otimes \Big({}_V\phi^{R_{(1)} ^*e^i}_{ S(R_{(1')})e_j}\Big)^{(a)} \Big({}_W\phi^{(R_{(2)}R_{(2')})^*e^k}_{e_l}\Big)^{(b)}
\\ \qquad
{}= \sum_{(R),(R),i,j,k,l} E_i^j\otimes E_k^l \otimes \big({}_W\phi^{e^k}_{R_{(2)}R_{(2')}e_l}\big)^{(b)}\Big({}_V\phi^{R_{(1)} ^*e^i}_{ S(R_{(1')})e_j}\Big)^{(a)}.
\end{gather*}
Now we have
\begin{gather*}
 {}_V\phi^{R_{(1)} ^*e^i}_{S(R_{(1')})e_j} = S\big(R_{(1')}\big) \rhd {}_V\phi^{e^i}_{e_j} \lhd R_{(1)} = {}_V\phi^{e^i}_{e_j} \cdot \big(R_{(1)}\otimes R_{(1')}\big),
\\
{}_W\phi^{(R_{(2)}R_{(2')})^*e^k}_{e_l} = {}_W\phi^{e^k}_{e_l} \lhd R_{(2)}R_{(2')} = {}_W\phi^{e^k}_{e_l} \cdot \big(R_{(2)}R_{(2')} \otimes 1\big),
\\
{}_W\phi^{e^k}_{R_{(2)}R_{(2')}e_l} = R_{(2)}R_{(2')}\rhd {}_W\phi^{e^k}_{e_l} = {}_W\phi^{e^k}_{e_l} \cdot \big(1 \otimes S^{-1}\big(R_{(2)}R_{(2')}\big)\big),
\end{gather*}
where we use the coregular actions $\rhd, \lhd$ and on the right hand sides we denote by $R_{(1)}$, $R_{(1')}$, $R_{(2)}$, $R_{(2')}\in \mUq$ the components of the universal R-matrix (instead of $R_{V,W}$). Using that $( S\otimes S)(R) = R$, and denoting by $m$ the product map of $\Ll_{0,n}$, the above relation eventually becomes
 \begin{gather*}
 \sum_{i,j,k,l} E_i^j \otimes E_k^l \otimes m\big(\big(\big({}_V\phi^{e^i}_{e_j}\big)^{(a)} \otimes \big({}_W\phi^{e^k}_{e_l}\big)^{(b)}\big) \cdot (R_{13}R_{23})\big)
 \\ \qquad
 {}= \sum_{i,j,k,l} E_i^j\otimes E_k^l \otimes m\big(\big(\big({}_W\phi^{e^k}_{e_l}\big)^{(b)} \otimes \big({}_V\phi^{e^i}_{e_j}\big)^{(a)}\big) \cdot (R_{24}^-R_{23}^-)\big).
 \end{gather*}
Identifying terms we get the commutation relation
\begin{align}
\big({}_W\phi^{e^k}_{e_l}\big)^{(b)}\big({}_V\phi^{e^i}_{e_j}\big)^{(a)} & = m\big(\big(\big({}_V\phi^{e^i}_{e_j}\big)^{(a)} \otimes \big({}_W\phi^{e^k}_{e_l}\big)^{(b)}\big) \cdot (R_{14}R_{24}R_{13}R_{23})\big)
\nonumber
\\
& = \sum_{(F)} \big(\big({}_V\phi^{e^i}_{e_j}\big)^{(a)}\cdot F_{(2)}\big)\big(\big({}_W\phi^{e^k}_{e_l}\big)^{(b)}\cdot F_{(1)}\big),
\label{lastformMV}
\end{align}
where $\textstyle F = R'_{23}R'_{13}R'_{24}R'_{14} = \sigma_{12,34}(R_{14}R_{24}R_{13}R_{23})$. This is the same as \eqref{commutrel}, so it is equivalent to the exchange relation. \end{proof}

Let us denote $\stackrel{V}{M}_{ij} := {}_V\phi^{e^j}_{e_i}$. The following proposition gives a formulation of the action $\operatorname{coad}^r$ in \eqref{actionprod} on these matrix coefficients.

\begin{prop} 
For all $1\leq a\leq n$ we have:
\begin{equation}\label{coadmat}
\operatorname{coad}^r(y)\big(\!\stackrel{V}{M}_{ij}\!{}^{(a)}\big) = \sum_{(y),k,l}\pi_V(y_{(1)})_{ik}\stackrel{V}{M}_{kl}\!{}^{(a)}\pi_V\big(S(y_{(2)})_{lj}\big)
\end{equation}
also written as
\begin{equation*}
\operatorname{coad}^r(y)\big(\!\stackrel{V}{M}\!{}^{(a)}\big) = \sum_{(y)}\big(y_{(1)_V} \otimes {\rm id} \big)\stackrel{V}{M}\!{}^{(a)}\big( S(y_{(2)})_V \otimes {\rm id} \big).
\end{equation*}
\end{prop}
\begin{proof} The formula \eqref{coadmat} extends uniquely to any product of matrices $\stackrel{V}{M}\!{}^{(a)}$ by the defining property of actions of module algebras. That it is equivalent to \eqref{actionprod} is immediate from the definition of the coregular actions. \end{proof}
\begin{remk} The fact that $\Ll_{0,n}$ is a right $U_q$-module algebra (see Proposition \ref{defiL0n}) follows also by verifying that the ideal generated by its defining relations \eqref{naturality}, \eqref{fusionrel} and \eqref{exrel2} of Proposition \ref{matrelL0n} is stable under $\operatorname{coad}^r$, using the formula \eqref{coadmat}. In the original papers on qLGFTs, these defining relations were imposed in order that this $U_q$-module algebra structure holds.
\end{remk}

\subsection{The Alekseev map} We have defined $\Ll_{0,n}$ as a twist of $\Ll_{0,1}^{\otimes n}$ in Proposition \ref{defiL0n}, and obtained a presentation by generators and relations in Proposition \ref{matrelL0n}. Although this presentation seems complicated, the Alekseev map, defined below, identifies $\Ll_{0,n}$ as a module subalgebra of $\tilde U_q^{\otimes n}$.

Let $V$ be a $U_q$-module of type $1$. For every $1\leq a\leq n$ set
\[
R_{Va} = (\pi_V\otimes \mathfrak{i}_a)(R) \in \operatorname{End}(V)\otimes \tilde U_q^{{\otimes} n},
\]
where as usual $R\in \mUq{\otimes} \mUq$ is the universal $R$-matrix, $\pi_V\colon \mUq \ra \operatorname{End}(V)$ the canonical projection, and $\mathfrak{i}_a \colon \tilde U_q \ra \tilde U_q^{\otimes n}$ the map defined by $\mathfrak{i}_a(x) = 1^{\otimes (a-1)}\otimes x\otimes 1^{\otimes (n-a)}$. That $R_{Va}$ is a matrix with entries in $\tilde U_q^{{\otimes} n}$ follows from \eqref{Thetaeval}. Consider the linear map
\[
\Phi_n \colon\ \operatorname{Vect}_{\mc(q)}\big\{\big({}_V\phi^{e^i}_{e_j}\big)^{(a)}\, \big\vert \, V\in \operatorname{Ob}(\mathcal{C}),\, 1\leq a\leq n,\, 1\leq i,j \leq \dim (V) \big\}\longrightarrow \tilde U_q^{{\otimes} n}
\]
defined by
\begin{gather}
({\rm id} \otimes \Phi_n)\big(\!\stackrel{V}{M}{}^{\!\! (n)}\big) = R_{Vn}R_{Vn}', \nonumber
\\ 
({\rm id} \otimes \Phi_n)\big(\!\stackrel{V}{M}{}^{\!\! (a)}\big) = (R_{Vn}\cdots R_{Va+1}) R_{Va}R_{Va}' (R_{Vn}\cdots R_{Va+1})^{-1}, \qquad
1\leq a < n. \label{definitionofphin2}
\end{gather}
Since $\textstyle \stackrel{V}{M}{}^{\!\! (a)}=\sum_{i,j }E_i^j\otimes \big({}_V\phi^{e^i}_{e_j}\big)^{(a)}$, we have $({\rm id} \otimes \Phi_n)\big(\!\!\stackrel{V}{M}{}^{\!\! (a)}\big) \in \operatorname{End}(V)\otimes \tilde U_q^{{\otimes} n}$. When $n=1$, $\Phi_n$ coincides with the RSD map $\Phi_1\colon \Ll_{0,1} \ra \tilde U_q$. We call $\Phi_n$ the {\it Alekseev map}. It was first introduced by Alekseev in \cite{A}.

We can represent $({\rm id} \otimes \Phi_n)\big(\!\!\stackrel{V}{M}{}^{\!\! (a)}\big)$ by the oriented colored braid shown in the figure below, where we use the standard graphical encoding of invariant operators of ribbon categories (see~\cite{Tu}). In this figure, the vertical strand with label $a$ carries the $a$-th factor of $\tilde U_q^{{\otimes} n}$, and the one with label $V$ carries $\operatorname{End}(V)$. Positive crossings carry the operator $\check{R} = \sigma \circ R$, where $\sigma(x\otimes y) = y\otimes x$ as usual, and negative crossings carry the operator $\check{R}^{-1}$.

\vspace{2ex}
\def\svgwidth{0.25\textwidth}
\begin{figure}[h!]\centering
\begingroup%
 \makeatletter%
 \providecommand\color[2][]{%
 \errmessage{(Inkscape) Color is used for the text in Inkscape, but the package 'color.sty' is not loaded}%
 \renewcommand\color[2][]{}%
 }%
 \providecommand\transparent[1]{%
 \errmessage{(Inkscape) Transparency is used (non-zero) for the text in Inkscape, but the package 'transparent.sty' is not loaded}%
 \renewcommand\transparent[1]{}%
 }%
 \providecommand\rotatebox[2]{#2}%
 \ifx\svgwidth\undefined%
 \setlength{\unitlength}{595.27559055bp}%
 \ifx\svgscale\undefined%
 \relax%
 \else%
 \setlength{\unitlength}{\unitlength * \real{\svgscale}}%
 \fi%
 \else%
 \setlength{\unitlength}{\svgwidth}%
 \fi%
 \global\let\svgwidth\undefined%
 \global\let\svgscale\undefined%
 \makeatother%
 \begin{picture}(1,1.41428571)%
 \put(0,0){\includegraphics[width=\unitlength,page=1]{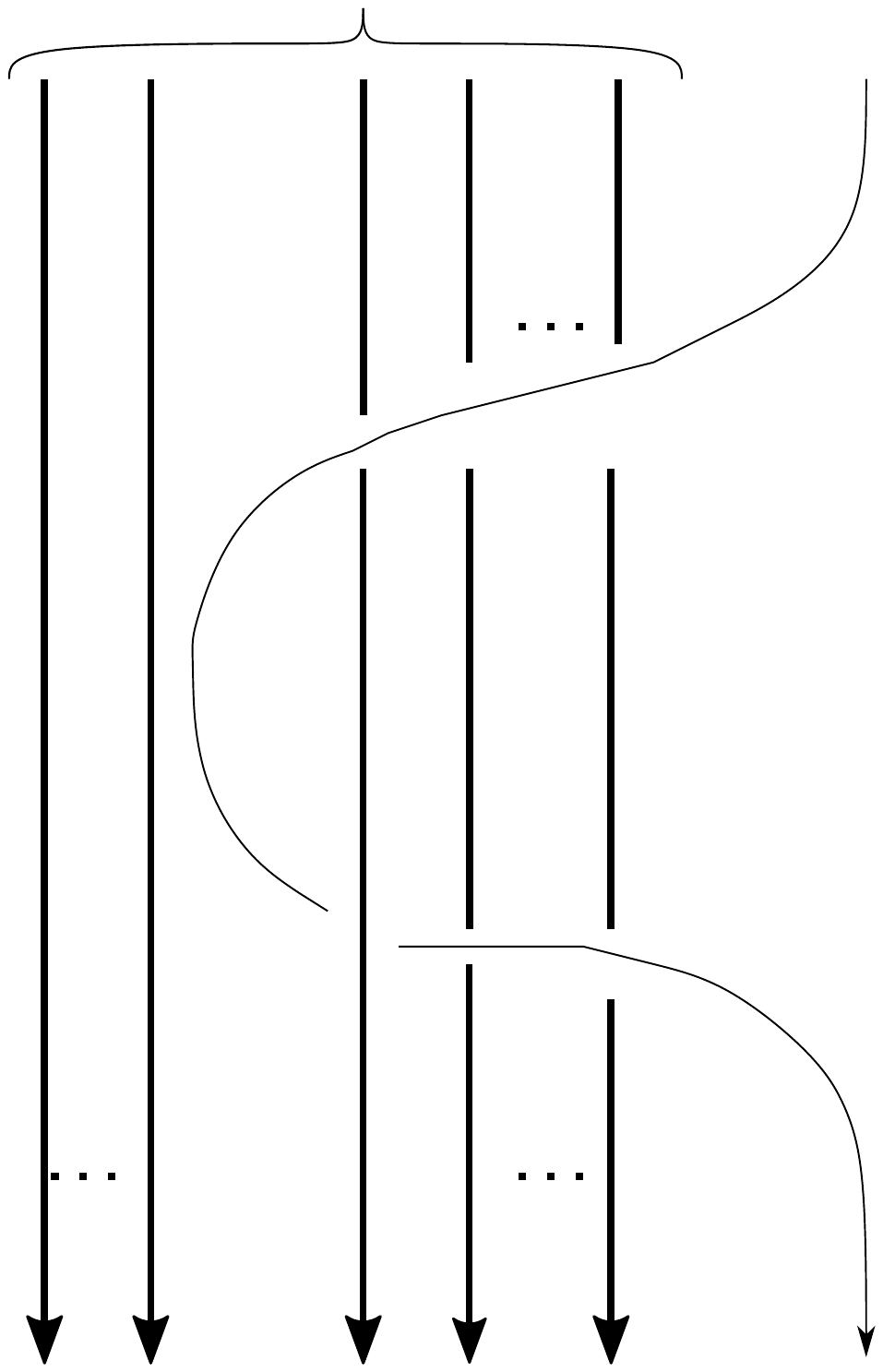}}%
 \put(0.2015873,0.73075394){\color[rgb]{0,0,0}\makebox(0,0)[lb]{\smash{}}}%
 \put(0.12599206,0.78115082){\color[rgb]{0,0,0}\makebox(0,0)[lb]{\smash{}}}%
 \put(0.0617654,0.07691941){\color[rgb]{0,0,0}\makebox(0,0)[lb]{\smash{\;\;\;$1$}}}%
 \put(0.39344072,0.08102897){\color[rgb]{0,0,0}\makebox(0,0)[lb]{\;\smash{$a$}}}%
 \put(0.62454827,0.08102897){\color[rgb]{0,0,0}\makebox(0,0)[lb]{\smash{$n$}}}%
 \put(0.83909643,0.08320443){\color[rgb]{0,0,0}\makebox(0,0)[lb]{\smash{$V$}}}%
 \put(0.42837314,1.37436254){\color[rgb]{0,0,0}\makebox(0,0)[lb]{\smash{$\tilde U_q^{{\otimes} n}$}}}%
 \put(0.19206485,0.07691912){\color[rgb]{0,0,0}\makebox(0,0)[lb]{\smash{}}}%
 \put(0.35277779,0.32757932){\color[rgb]{0,0,0}\makebox(0,0)[lb]{\smash{}}}%
 \end{picture}%
\endgroup%


\caption{The colored braid representing $({\rm id} \otimes \Phi_n)\big(\!\stackrel{V}{M}{}^{\!\! (a)}\big)$.}\label{fig6.1}
\end{figure}

Let us endow $\tilde U_q^{{\otimes} n}$ with the following action of $U_q$
\begin{equation}
\operatorname{ad}^r(y)(x) = \sum_{(y)}\Delta^{(n)}({S}(y_{(1)})) x \Delta^{(n)}(y_{(2)}) \label{adjointactiononn}
\end{equation}
for all $y \in U_q$, $x \in \tilde U_q^{{\otimes} n}$. Then $\tilde U_q^{{\otimes} n}$ becomes a right $U_q$-module algebra.

The next result is due to Alekseev \cite{A}.
\begin{teo}\label{Alekseevmap} The Alekseev map yields an embedding of module algebras $\Phi_n\colon \Ll_{0,n}\ra \tilde U_q^{{\otimes} n}$. Moreover it satisfies
\begin{equation}\label{phicommut}
({\rm id} \otimes \Phi_n)\big(\!\stackrel{V}{M}{}^{\!\! (1)}\cdots \stackrel{V}{M}{}^{\!\! (n)}\big) = \big(\pi_V\otimes \Delta^{(n-1)}\big)(RR').
\end{equation}
\end{teo}
\begin{proof} Let us extend $\Phi_n$ (by keeping the same notation) in the natural way to the algebra freely generated by the matrix coefficients $\big({}_V\phi^{e^i}_{e_j}\big)^{(a)}$.

By Proposition \ref{matrelL0n}, $\Phi_n$ induces a well-defined algebra morphism $\Ll_{0,n}\ra \tilde U_q^{{\otimes} n}$ if it preserves the fusion and exchange relations. Using the graphical encoding recalled in Figure~\ref{fig6.1}, this is shown in the next two figures in the case $n=2$, which generalizes immediately to any $n$. The symbol $\stackrel{\centerdot}{=}$ means equality up to isotopy. Similarly, the relation \eqref{phicommut} is proved by the third figure.
\begin{figure}[h!]
\vspace*{-2cm}
\def\svgwidth{0.35\textwidth}
\begin{center}
\begingroup%
 \makeatletter%
 \providecommand\color[2][]{%
 \errmessage{(Inkscape) Color is used for the text in Inkscape, but the package 'color.sty' is not loaded}%
 \renewcommand\color[2][]{}%
 }%
 \providecommand\transparent[1]{%
 \errmessage{(Inkscape) Transparency is used (non-zero) for the text in Inkscape, but the package 'transparent.sty' is not loaded}%
 \renewcommand\transparent[1]{}%
 }%
 \providecommand\rotatebox[2]{#2}%
 \ifx\svgwidth\undefined%
 \setlength{\unitlength}{595.27559055bp}%
 \ifx\svgscale\undefined%
 \relax%
 \else%
 \setlength{\unitlength}{\unitlength * \real{\svgscale}}%
 \fi%
 \else%
 \setlength{\unitlength}{\svgwidth}%
 \fi%
 \global\let\svgwidth\undefined%
 \global\let\svgscale\undefined%
 \makeatother%
 \begin{picture}(1,1.41428571)%
 \put(0,0){\includegraphics[width=\unitlength,page=1]{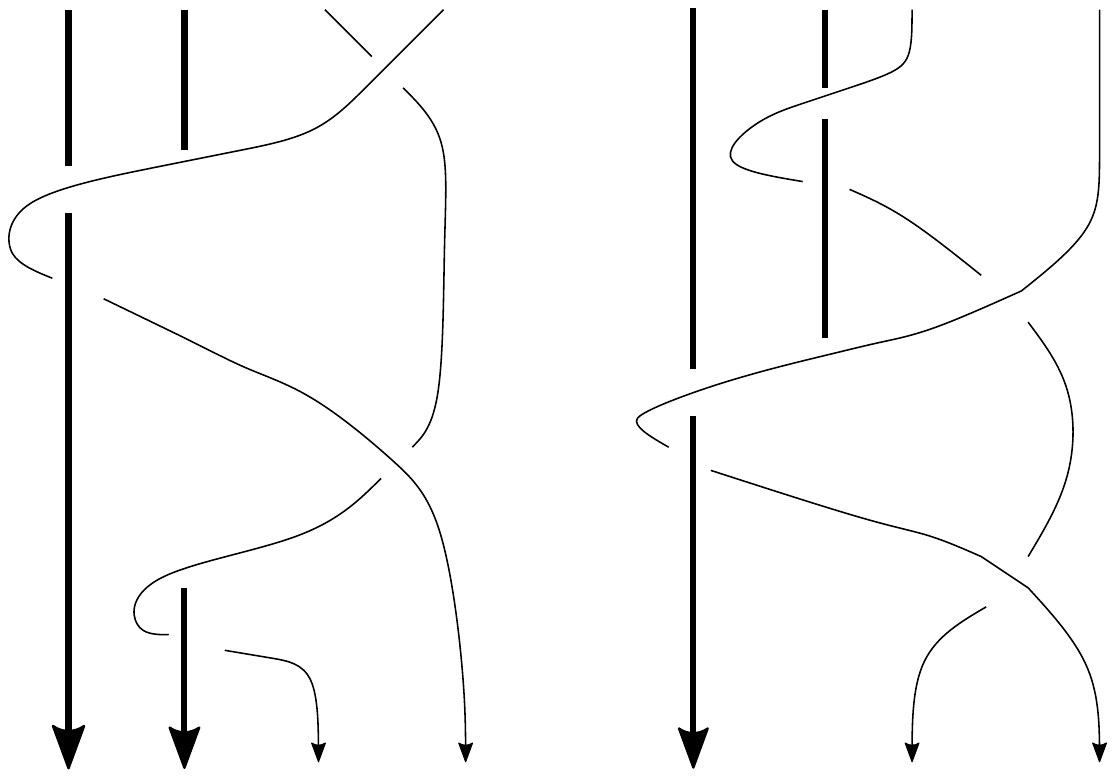}}%
 \put(0.03779762,0.36537694){\color[rgb]{0,0,0}\makebox(0,0)[lb]{\smash{}}}%
 \put(0,0){\includegraphics[width=\unitlength,page=2]{dessin1final.pdf}}%
 \put(0.61736109,0.36537694){\color[rgb]{0,0,0}\makebox(0,0)[lb]{\smash{}}}%
 \put(0.97013884,0.49136905){\color[rgb]{0,0,0}\makebox(0,0)[lb]{\smash{}}}%
 \put(0.45517938,0.68945334){\color[rgb]{0,0,0}\makebox(0,0)[lb]{\smash{$\stackrel{\centerdot}{=}$}}}%
 \put(0,0.25943008){\color[rgb]{0,0,0}\makebox(0,0)[lt]{\begin{minipage}{3.71428571\unitlength}\raggedright \end{minipage}}}%
 \put(0,0.25943008){\color[rgb]{0,0,0}\makebox(0,0)[lt]{\begin{minipage}{3.71428571\unitlength}\raggedright \end{minipage}}}%
 \put(0,0){\includegraphics[width=\unitlength,page=3]{dessin1final.pdf}}%
 \put(0.231309,0.36397632){\color[rgb]{0,0,0}\makebox(0,0)[lb]{\;\smash{$W$}}}%
 \put(0.37978552,0.36397632){\color[rgb]{0,0,0}\makebox(0,0)[lb]{\smash{$V$}}}%
 \put(0.73066617,0.36397632){\color[rgb]{0,0,0}\makebox(0,0)[lb]{\;\smash{$W$}}}%
 \put(0.87914269,0.36397632){\color[rgb]{0,0,0}\makebox(0,0)[lb]{\smash{$V$}}}%
 \end{picture}%
\endgroup%
\end{center}

\vspace*{-3cm}
\def\svgwidth{0.35\textwidth}
\begin{center}
\begingroup%
 \makeatletter%
 \providecommand\color[2][]{%
 \errmessage{(Inkscape) Color is used for the text in Inkscape, but the package 'color.sty' is not loaded}%
 \renewcommand\color[2][]{}%
 }%
 \providecommand\transparent[1]{%
 \errmessage{(Inkscape) Transparency is used (non-zero) for the text in Inkscape, but the package 'transparent.sty' is not loaded}%
 \renewcommand\transparent[1]{}%
 }%
 \providecommand\rotatebox[2]{#2}%
 \ifx\svgwidth\undefined%
 \setlength{\unitlength}{595.27559055bp}%
 \ifx\svgscale\undefined%
 \relax%
 \else%
 \setlength{\unitlength}{\unitlength * \real{\svgscale}}%
 \fi%
 \else%
 \setlength{\unitlength}{\svgwidth}%
 \fi%
 \global\let\svgwidth\undefined%
 \global\let\svgscale\undefined%
 \makeatother%
 \begin{picture}(1,1.41428571)%
 \put(0,0){\includegraphics[width=\unitlength,page=1]{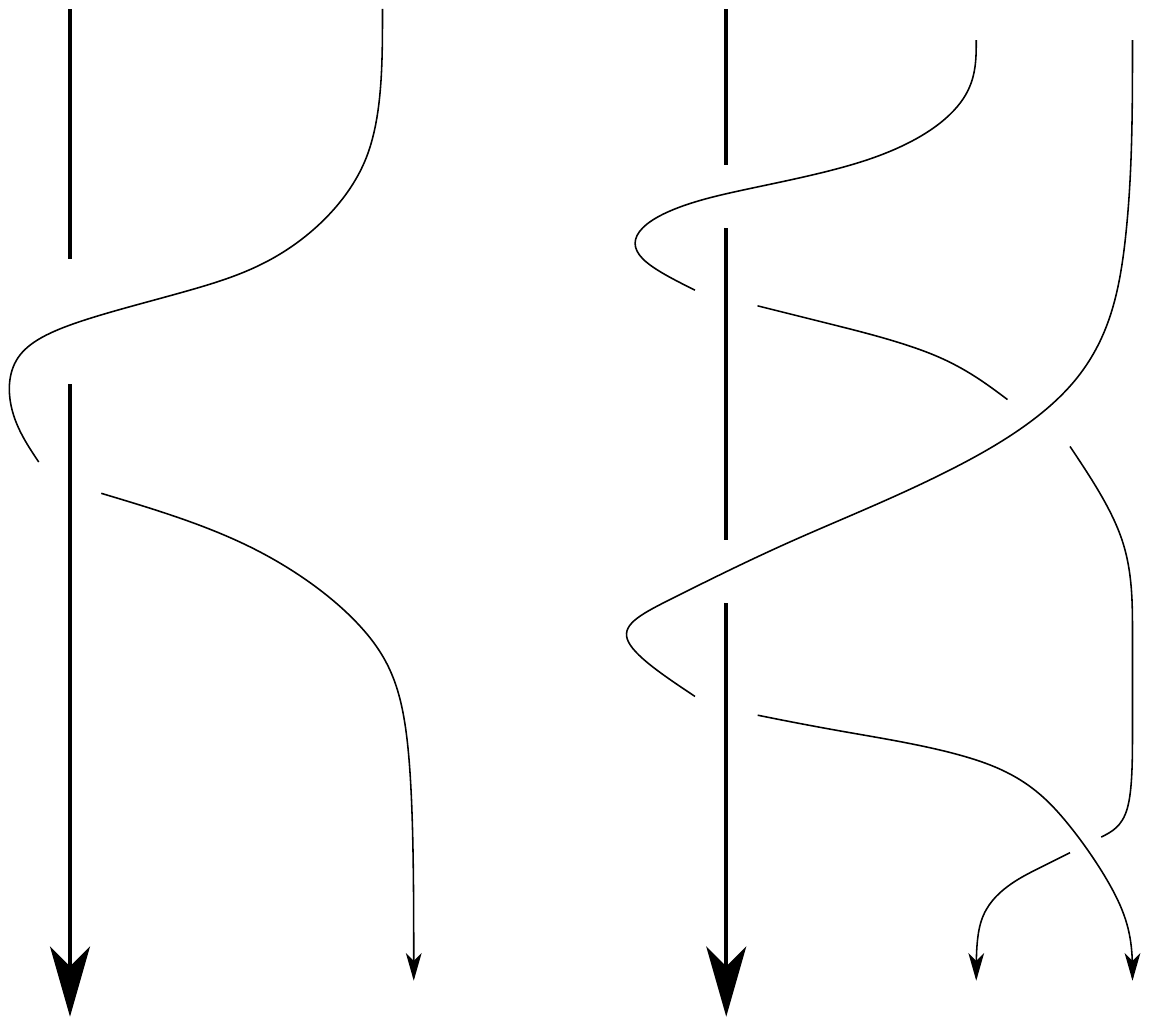}}%
 \put(0.3162529,0.36397632){\color[rgb]{0,0,0}\makebox(0,0)[lb]{\smash{$V\otimes W$}}}%
 \put(0.81131221,0.36397632){\color[rgb]{0,0,0}\makebox(0,0)[lb]{\smash{$V$}}}%
 \put(0.93348607,0.36397632){\color[rgb]{0,0,0}\makebox(0,0)[lb]{\smash{$W$}}}%
 \put(0.45357143,0.81894844){\color[rgb]{0,0,0}\makebox(0,0)[lb]{\smash{$\stackrel{\centerdot}{=}$}}}%
 \end{picture}%
\endgroup%
\end{center}

\vspace*{-3.5cm}
\def\svgwidth{0.35\textwidth}
\begin{center}
\begingroup%
 \makeatletter%
 \providecommand\color[2][]{%
 \errmessage{(Inkscape) Color is used for the text in Inkscape, but the package 'color.sty' is not loaded}%
 \renewcommand\color[2][]{}%
 }%
 \providecommand\transparent[1]{%
 \errmessage{(Inkscape) Transparency is used (non-zero) for the text in Inkscape, but the package 'transparent.sty' is not loaded}%
 \renewcommand\transparent[1]{}%
 }%
 \providecommand\rotatebox[2]{#2}%
 \ifx\svgwidth\undefined%
 \setlength{\unitlength}{595.27559055bp}%
 \ifx\svgscale\undefined%
 \relax%
 \else%
 \setlength{\unitlength}{\unitlength * \real{\svgscale}}%
 \fi%
 \else%
 \setlength{\unitlength}{\svgwidth}%
 \fi%
 \global\let\svgwidth\undefined%
 \global\let\svgscale\undefined%
 \makeatother%
 \begin{picture}(1,1.41428571)%
 \put(0,0){\includegraphics[width=\unitlength,page=1]{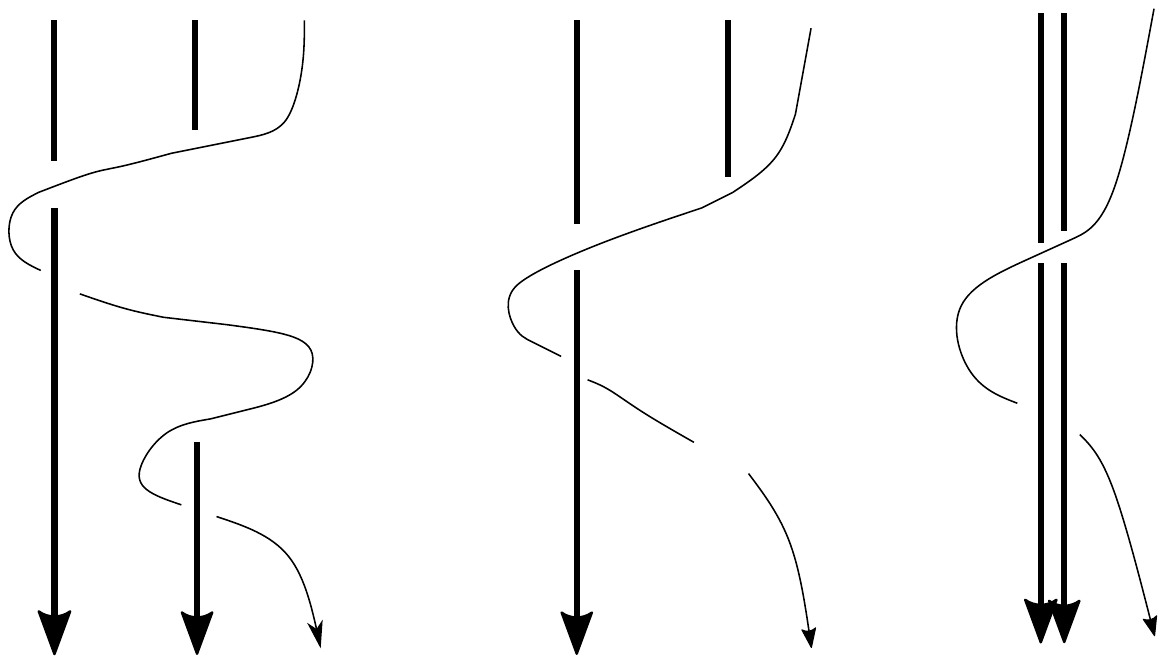}}%
 \put(0.34580125,0.76790234){\color[rgb]{0,0,0}\makebox(0,0)[lb]{\smash{$\stackrel{\centerdot}{=}$}}}%
 \put(0.73702596,0.77228168){\color[rgb]{0,0,0}\makebox(0,0)[lb]{\smash{$\stackrel{\centerdot}{=}$}}}%
 \put(0,0){\includegraphics[width=\unitlength,page=2]{dessin4final.pdf}}%
 \put(0.27718253,0.52916667){\color[rgb]{0,0,0}\makebox(0,0)[lb]{\smash{}}}%
 \put(0.30238095,0.48396823){\color[rgb]{0,0,0}\makebox(0,0)[lb]{\!\!\smash{V}}}%
 \put(0.70555558,0.52916667){\color[rgb]{0,0,0}\makebox(0,0)[lb]{\smash{}}}%
 \put(0.70555558,0.48396823){\color[rgb]{0,0,0}\makebox(0,0)[lb]{\!\!\smash{V}}}%
 \put(0.95753966,0.48396823){\color[rgb]{0,0,0}\makebox(0,0)[lb]{\smash{V}}}%
 \put(0,0){\includegraphics[width=\unitlength,page=3]{dessin4final.pdf}}%
 \end{picture}%
\endgroup%
\end{center}
\vspace*{-3cm}
\caption{The colored braid identities representing the exchange relation \eqref{exrel2} (top), the fusion relation \eqref{fusionrel} (middle), and the relation \eqref{phicommut} (bottom).}\label{fig6.2}
\end{figure}

The linear map defined by $\stackrel{V}{M}{}^{\!\! (a)} \mapsto R_{Va}R_{Va}'$ on the $\mc(q)$-vector space spanned by all the matrix coefficients $\big({}_V\phi^{e^i}_{e_j}\big)^{(a)}$, for every $1\leq a\leq n$ and object $V$ of $\mathcal{C}$, is injective by Theorem~\ref{drinfeldmap}. It differs from $\Phi_n$ by a linear isomorphism (induced on the $a$-th tensorand by conjugating with $R_{Vn}\cdots R_{Va+1}$), so $\Phi_n$ is injective.

Finally, let us show that $\Phi_n$ is a morphism of $U_q$-modules. Since both multiplications in $\Ll_{0,n}$ and $\tilde U_q^{{\otimes} n}$ commute with the respective actions of $U_q$, it is enough to check this on generators. For those given by the matrix coefficients we then have to show that
\begin{equation}\label{equivnonexplicit}
({\rm id} \otimes \Phi_n)\big( \operatorname{coad}^r(y)\big(\!\stackrel{V}{M}{}^{\!\! (a)}\big) \big) = ({\rm id} \otimes \operatorname{ad}^r(y))\big(({\rm id} \otimes \Phi_n)\big(\!\stackrel{V}{M}{}^{\!\! (a)}\big)\big)
\end{equation}
for every object $V$ in $\mathcal{C}$, $1\leq a\leq n$ and $y\in U_q$. In the case $a=n$, $\Phi_n$ has the same expression as the RSD map. Hence \eqref{equivnonexplicit} follows from the statement of equivariance in Proposition \ref{drinfeldmap}.
This can also be proved easily:
\begin{gather*}
({\rm id} \otimes \operatorname{ad}^r(y))\big(({\rm id} \otimes \Phi_n) \big(\!\stackrel{V}{M}{}^{\!\! (n)}\big)\big)
\\ \qquad
{}=({\rm id} \otimes \operatorname{ad}^r(y))(R_{Vn}R_{Vn}')
\\ \qquad
{}=\sum_{(RR'),(y)}\pi_{V}\big((RR')_{(1)}\big)S(y_{(1)}) (RR')_{(2)}y_{(2)}
\\ \qquad
{}=\sum_{(RR'),(y)} \pi_{V}\big(y_{(1)}S(y_{(2)})\big) \pi_{V}\big((RR')_{(1)}\big)S(y_{(3)} )(RR')_{(2)}y_{(4)}
\\ \qquad
{}=\sum_{(RR'),(y)} \pi_{V}(y_{(1)})\pi_{V}\big((RR')_{(1)}\big)\pi_V\big(S(y_{(2)})\big) (RR')_{(2)}S(y_{(3)} )y_{(4)}
\\ \qquad
{}= \sum_{(RR'),(y)} \pi_{V}(y_{(1)})\pi_{V}((RR')_{(1)})\pi_V\big(S(y_{(2)})\big) (RR')_{(2)}
\\ \qquad
{}= ({\rm id} \otimes \Phi_n)\big(\! \operatorname{coad}^r(y)\big(\!\stackrel{V}{M}{}^{\!\! (n)}\big) \big),
\end{gather*}
where we have used $[RR', (S\otimes S)(\Delta(y))]=0.$

More generally, by writing the actions explicitly, that result implies as well the relation
\begin{gather}
\sum_{(y)}\big(\big(\pi_V(y_{(1)}) \otimes 1\big) R_{Va}R_{Va}' \big(\pi_V( S(y_{(2)})) \otimes 1\big)\big)\nonumber
\\ \qquad
{} = \sum_{(y)} \big(1 \otimes \mathfrak{i}_{a}\big( S(y_{(1)})\big)\big) R_{Va}R_{Va}' \big(1 \otimes \mathfrak{i}_{a}(y_{(2)})\big).
\label{equivexplicit}
\end{gather}
Now set
\[R^{(a)} = R_{Vn}\cdots R_{Va+1}\]
and denote by $\mathfrak{i}_{an}\colon \tilde U_q^{{\otimes} n-a} \ra 1^{ \otimes a} \otimes\tilde U_q^{{\otimes} n-a} \subset \tilde U_q^{{\otimes} n}$ the identification map with the final $n-a$ tensorands. Then
\begin{gather*}
 (1 \otimes \Phi_n) \big(\! \operatorname{coad}^r(y) \big(\!\stackrel{V}{M}{}^{\!\! (a)}\big) \big) = \sum_{(y)}\big(\pi_V(y_{(1)}) \otimes 1\big) R^{(a)} R_{Va}R_{Va}' R^{(a)}{}^{\!-1} \big(\pi_V( S(y_{(2)})) \otimes 1\big)
 \\ \hphantom{ (1 \otimes \Phi_n) \big(\! \operatorname{coad}^r(y) \big(\!\stackrel{V}{M}{}^{\!\! (a)}\big) \big)}
{} = \sum_{(y)}\big(\pi_V(y_{(1)}) \otimes 1\big) R^{(a)} R_{Va}R_{Va}'
 \\ \hphantom{ (1 \otimes \Phi_n) \big(\! \operatorname{coad}^r(y) \big(\!\stackrel{V}{M}{}^{\!\! (a)}\big) \big)=}
{} \times\big(1 \otimes \mathfrak{i}_{an}\big(\Delta^{(n-a-1)}\big(\underbrace{ S(y_{(2)})y_{(3)}}_{=\varepsilon(y_{(2)})}\big)\big)\big) R^{(a)}{}^{\!-1} \big(\pi_V( S(y_{(4)})) \otimes 1\big)
 \\ \hphantom{ (1 \otimes \Phi_n) \big(\! \operatorname{coad}^r(y) \big(\!\stackrel{V}{M}{}^{\!\! (a)}\big) \big)}
{} = \sum_{(y)}\big(\pi_V(y_{(1)}) \otimes 1\big) R^{(a)} \big(1 \otimes \mathfrak{i}_{an}\big(\Delta^{(n-a-1)}\big( S(y_{(2)})\big)\big)\big)
 \\ \hphantom{ (1 \otimes \Phi_n) \big(\! \operatorname{coad}^r(y) \big(\!\stackrel{V}{M}{}^{\!\! (a)}\big) \big)=}
{}\times R_{Va}R_{Va}' \big(1\! \otimes \mathfrak{i}_{an}\big(\Delta^{(n-a-1)}(y_{(3)})\big)\big) R^{(a)}{}^{\!-1} \big(\pi_V\big( S(y_{(4)})\big)\! \otimes 1\big).
 \end{gather*}
Since $R \Delta = \Delta'R$ in $\mUq^{ \otimes 2}$, we have $\textstyle \sum_{(y)} R^{-1} (y_{(2)}\otimes y_{(1)}) = \sum_{(y)} (y_{(1)}\otimes y_{(2)}) R^{-1}$ and then
\begin{gather*}
\sum_{(y)} \big(1 \otimes S(y_{(1)})\big)R(y_{(2)} \otimes 1) = \sum_{(y)} (y_{(1)} \otimes 1) R \big(1 \otimes S(y_{(2)})\big),
\\
\sum_{(y)} \big( S(y_{(1)}) \otimes 1\big)R^{-1}(1 \otimes y_{(2)}) = \sum_{(y)} (1 \otimes y_{(1)}) R^{-1}\big(S(y_{(2)}) \otimes 1\big)
\end{gather*}
by applying ${\rm id} \otimes S$, and using $({\rm id} \otimes S)(R^{-1}) = R$. We deduce recursively
\begin{gather*}
\sum_{(y)}\big(\pi_V(y_{(1)}) \otimes 1\big) R^{(a)} \big(1 \otimes \mathfrak{i}_{an}\big(\Delta^{(n-a-1)}\big( S(y_{(2)})\big)\big)\big)
\\ \qquad
{}= \big(1 \otimes \mathfrak{i}_{an}\big(\Delta^{(n-a-1)}\big( S(y_{(1)})\big)\big) \big) R^{(a)} \big(\pi_V(y_{(2)}) \otimes 1\big)
\end{gather*}
and
\begin{gather*}
\sum_{(y)} \big(1 \otimes \mathfrak{i}_{an}\big(\Delta^{(n-a-1)}(y_{(1)})\big)\big) R^{(a)}{}^{\!-1} \big(\pi_V\big( S(y_{(2)})\big) \otimes 1\big)
\\ \qquad
{}= \big(\pi_V\big( S(y_{(1)})\big) \otimes 1\big) R^{(a)}{}^{\!-1} \big(1 \otimes \mathfrak{i}_{an}\big(\Delta^{(n-a-1)}(y_{(2)})\big)\big).
\end{gather*}
So
\begin{gather*}
(1 \otimes \Phi_n) \big(\!\operatorname{coad}^r(y) \big(\!\stackrel{V}{M}{}^{\!\! (a)}\big) \big)
\\ \qquad
{}= \sum_{(y)} \big(1 \otimes \mathfrak{i}_{an}\big(\Delta^{(n-a-1)}\big( S(y_{(1)})\big)\big)\big) R^{(a)} \big(\pi_V(y_{(2)}) \otimes 1\big)
\\ \qquad\hphantom{=\sum_{(y)}}
{} \times R_{Va}R_{Va}' \big(\pi_V\big( S(y_{(3)})\big) \otimes 1\big) R^{(a)}{}^{\!-1} \big(1 \otimes \mathfrak{i}_{an}\big(\Delta^{(n-a-1)}(y_{(4)})\big)\big)
\\ \qquad
{} = \sum_{(y)} \big(1 \otimes \mathfrak{i}_{an}\big(\Delta^{(n-a-1)}\big( S(y_{(1)})\big)\big)\big) R^{(a)} \big(1 \otimes \mathfrak{i}_{a}\big( S(y_{(2)})\big)\big)
\\ \qquad\hphantom{=\sum_{(y)}}
{} \times R_{Va}R_{Va}'\big(1 \otimes \mathfrak{i}_{a}(y_{(3)})\big) R^{(a)}{}^{\!-1} \big(1 \otimes \mathfrak{i}_{an}\big(\Delta^{(n-a-1)}(y_{(4)})\big)\big)
\\ \qquad
{} = \sum_{(y)} \big(1 \otimes \mathfrak{i}_{an}\big(\Delta^{(n-a-1)}\big( S(y_{(1)})\big)\big) \big) \big(1 \otimes \mathfrak{i}_{a}\big( S(y_{(2)})\big)\big) R^{(a)}
\\ \qquad\hphantom{=\sum_{(y)}}
{} \times R_{Va}R_{Va}' R^{(a)}{}^{\!-1} \big(1 \otimes \mathfrak{i}_{a}(y_{(3)})\big) \big(1 \otimes \mathfrak{i}_{an}\big(\Delta^{(n-a-1)}(y_{(4)})\big)\big)
\\ \qquad
{} = \sum_{(y)} \big(1 \otimes \mathfrak{i}_{a-1n}\big(\Delta^{(n-a)}\big( S(y_{(1)})\big)\big) \big) R^{(a)} R_{Va}R_{Va}' R^{(a)}{}^{\!-1} \big(1 \otimes \mathfrak{i}_{a-1n}\big(\Delta^{(n-a)}(y_{(2)})\big)\big)
\\ \qquad
{} = \sum_{(y)} \big(1 \otimes \Delta^{(n-1)}\big( S(y_{(1)})\big)\big) R^{(a)} R_{Va}R_{Va}' R^{(a)}{}^{\!-1} \big(1 \otimes \Delta^{(n-1)}(y_{(2)})\big),
 \end{gather*}
where we use \eqref{equivexplicit} in the second equality, and the others follow from trivial commutations between elements lying on different tensorands, and the property $( S \otimes S) \Delta' = \Delta S$. This proves~\eqref{equivnonexplicit}.
The relation (\ref{phicommut}) has been proved with a picture which encapsulates the following simple computation:
\begin{align*}
({\rm id} \otimes \Phi_n)\big(\!\stackrel{V}{M}{}^{\!\! (1)}\cdots \stackrel{V}{M}{}^{\!\! (n)}\big) & = \prod_{a=1}^{n} \big((R_{Vn}\cdots R_{Va+1}) R_{Va}R_{Va}' (R_{Vn}\cdots R_{Va+1})^{-1}\big)
\\
& =R_{Vn}\cdots R_{V1}R_{V1}'\cdots R_{Vn}'=\big(\pi_V\otimes \Delta^{(n-1)}\big)(RR').
\end{align*}
This concludes the proof.
\end{proof}

Denote by $\big({\tilde U}_q^{{\otimes} n}\big)^{\rm lf}$ the set of locally finite elements of ${\tilde U}_q^{{\otimes} n}$ with respect to the action \eqref{adjointactiononn}. Consider the subalgebras $({\mathcal L}_{0,n})^{U_q}$ and $\big(\tilde U_q^{\otimes n}\big)^{U_q}$ of invariant elements of the $U_q$-module algebras~${\mathcal L}_{0,n}$ and $\tilde U_q^{\otimes n}$ respectively. We have:

\begin{prop}\label{surjPhiinv}
The Alekseev map sends ${\mathcal L}_{0,n}$ isomorphically onto $\big({\tilde U}_q^{{\otimes} n}\big)^{\rm lf}$, and its restriction to invariant elements affords an isomorphism $\Phi_n\colon {\mathcal L}_{0,n}^{U_q}\rightarrow \big({\tilde U}_q^{{\otimes} n}\big)^{U_q}$.
\end{prop}
\begin{proof} Consider the first claim. By \eqref{coadmat} the action $\operatorname{coad}^r$ on ${\mathcal L}_{0,n}$ preserves the spaces of matrix coefficients of each object $V$ of $\mathcal{C}$. Since it gives ${\mathcal L}_{0,n}$ a structure of module algebra, it follows that $\operatorname{coad}^r$ is locally finite. The equivariance of $\Phi_n$ then implies that its image is contained in~$\big({\tilde U}_q^{{\otimes} n}\big)^{\rm lf}$. As $\Phi_n$ is injective, it remains to prove that it maps ${\mathcal L}_{0,n}$ surjectively onto $\big({\tilde U}_q^{{\otimes} n}\big)^{\rm lf}$. Let us denote by ${\rm ad}_n^r$ the action (\ref{adjointactiononn}). We can consider another action of $U_q$ on $\tilde U_q^{\otimes n}$ defined as the $n$-fold tensor product of the adjoint action of $U_q$ on $\tilde U_q^{\otimes n}$, which we denote $(\operatorname{ad}^r)^{\otimes n}$. It is thus defined as follows:
\[
(\operatorname{ad}^r)^{\otimes n}(y)\big(a^{(1)}\otimes \dots \otimes a^{(n)}\big) = \sum_{(y)} \operatorname{ad}^r(y_{(1)})\big(a^{(1)}\big) \otimes \dots \otimes \operatorname{ad}^r(y_{(n)})\big(a^{(n)}\big).
\]
We first show that these two actions on $\tilde U_q^{\otimes n}$ define isomorphic modules. For simplicity of notations let us consider the case $n=2$. Consider the map
\[
\psi\colon\quad {\mathbb U}_q^{\otimes 2}\rightarrow {\mathbb U}_q^{\otimes 2},\qquad
a\otimes b\mapsto R^{-1}_{12}(a\otimes 1)R_{12}(1\otimes b).
\]
We claim that $\psi$ intertwines $(\operatorname{ad}^r)^{\otimes 2}$ and ${\rm ad}_2^r$. This follows from a small variation of the previous proof. Indeed we have
\begin{align*}
{\rm ad}_2^r(y)(\psi(a\otimes b))&=\sum_{(y), (R), (R^{-1})}S(y_{(2)})(R^{-1})_{(1)}a R_{(1)} y_{(3)}\otimes S(y_{(1)})(R^{-1})_{(2)} R_{(2)} b y_{(4)}
\\
&=\sum_{(y), (R), (R^{-1})}(R^{-1})_{(1)}S(y_{(1)})a R_{(1)} y_{(3)}\otimes (R^{-1})_{(2)} S(y_{(2)}) R_{(2)} b y_{(4)}
\\
&=\sum_{(y), (R), (R^{-1})}(R^{-1})_{(1)}S(y_{(1)})a y_{(2)}R_{(1)}\otimes (R^{-1})_{(2)} R_{(2)} S(y_{(3)}) b y_{(4)}
\\
&=\psi((\operatorname{ad}^r)^{\otimes 2}(a\otimes b))) .
\end{align*}
It is easy to show that $\psi$ is an isomorphism, it descends to a map $\psi\colon \big(\tilde{U}_q^{\rm lf}\big)^{\otimes 2}\rightarrow \big({\tilde U}_q^{\otimes 2}\big)^{\rm lf}$, and that $\textstyle \psi(a\otimes b) = \sum_{(R)} \operatorname{ad}^r(R_{(1)})(a) \otimes R_{(2)}b$. Now, let $z\in \big({\tilde U}_q^{{\otimes} 2}\big)^{\rm lf}$. Then $\psi^{-1}(z)$ is locally finite for the action $(\operatorname{ad}^r)^{\otimes 2}$. By the main theorem of \cite{KLNY}, the set of locally finite elements of~${\tilde U}_q^{{\otimes} n}$ for $(\operatorname{ad}^r)^{\otimes n}$ is equal to $\big({\tilde U}_q^{\rm lf}\big)^{\otimes n}$. Therefore $\psi^{-1}(z)\in {\tilde U}_q^{\rm lf}\otimes {\tilde U}_q^{\rm lf}$, and from the surjectivity of the map $\Phi_1$ onto ${\tilde U}_q^{\rm lf}$, we deduce that $z=\psi\circ \Phi_1^{\otimes 2}(t)$ for some $t \in {\mathcal L}_{0,1}^{\otimes 2}$. But because of the identity $R_{12}^{-1}R_{01}R'_{01}R_{12}=R_{02}R_{01}R'_{01}R_{02}^{-1}$, which is a consequence of the Yang--Baxter equation, we obtain that $z=\psi\circ \Phi_1^{\otimes 2}(t)=\Phi_2(t)$, where we have identified the two vector spaces ${\mathcal L}_{0,2}$ and~${\mathcal L}_{0,1}^{\otimes 2}$. This shows $\Phi_2$ is a surjection onto $\big({\tilde U}_q^{{\otimes} 2}\big)^{\rm lf}$. The generalisation to any $n$ is straightforward.

By the previous theorem we have an inclusion $\Phi_n\big({\mathcal L}_{0,n}^{U_q}\big)\subset \big({\tilde U}_q^{{\otimes} n}\big)^{U_q}$. The above argument applies in particular to ${\rm ad}_n^r$-invariant elements $z\in \big({\tilde U}_q^{{\otimes} n}\big)^{U_q}$. Together with the injectivity of $\Phi_n$ it implies that this inclusion is an equality. Therefore the second claim follows.\end{proof}

\begin{remk}
This last result is a generalisation to $n\geq 1$ of Theorem \ref{drinfeldmap}(3).
\end{remk}

\subsection{Integral form} Let $\{\alpha_i\}_{i\in I}$ be a basis of the free $A$-module $\Ll_{0,1}^A$. Put $(\alpha_i)^{(a)} := \mathfrak{i}_a(\alpha_i)$, where as usual $\mathfrak{i}_a\colon \Ll_{0,1}^A\ra \big(\Ll_{0,1}^A\big)^{\otimes n}$ is the inclusion map, $1\leq a\leq n$. A basis of the free $A$-module $\big(\Ll_{0,1}^A\big)^{\otimes n}$ is given by the elements $\alpha_{i_1\cdots i_n} = (\alpha_{i_1})^{(1)} \otimes \dots \otimes(\alpha_{i_n})^{(n)}$, the products in $\Ll_{0,n}^A$ of the elements $(\alpha_{i_1})^{(1)},\dots, (\alpha_{i_n})^{(n)}$ (in this order from $1$ to $n$), indexed by the tuples $(i_1,\dots,i_n)\in I^n$. For every basis elements $\alpha_{i_1\cdots i_n}$, $\alpha_{j_1\cdots j_n}$ we have
\[
\alpha_{i_1\cdots i_n} \alpha_{j_1\cdots j_n} = \sum_{k_1\cdots k_n} m_{i_1\cdots i_n, j_1\cdots j_n}^{k_1\cdots k_n} \alpha_{k_1\cdots k_n}
\]
for some coefficients $m_{i_1\cdots i_n,j_1\cdots j_n}^{k_1\cdots k_n}\in \mc(q)$. By Proposition~\ref{defL01A} we know that $(\alpha_i)^{(a)}(\alpha_j)^{(a)}$ is a~linear combination over $A$ of basis elements $(\alpha_k)^{(a)}$, for every $i$, $j$ and $a$. Therefore, it remains to show that any product $(\alpha_i)^{(b)}(\alpha_j)^{(a)}$, for $a<b$, is a linear combination over $A$ of basis elements $(\alpha_k)^{(a)}(\alpha_l)^{(b)}$. Any $(\alpha_k)^{(a)}$ is a linear combination over $A$ of matrix cooefficients of $U_A^{\rm res}$-modules of type $1$. Let $X$ and $Y$ be two such modules. By \eqref{lastformMV} the exchange relation between matrix coefficients $\big({}_X\phi^{e^i}_{e_j}\big)^{(a)}$ and $\big({}_Y\phi^{e^k}_{e_l}\big)^{(b)}$, for $a<b$, can be written as
 \begin{align*}
\big({}_Y\phi^{e^k}_{e_l}\big)^{(b)}\big({}_X\phi^{e^i}_{e_j}\big)^{(a)} = \sum_{(F)} \big(\big({}_X\phi^{e^i}_{e_j}\big)^{(a)}\cdot F_{(2)}\big)\big(\big({}_Y\phi^{e^k}_{e_l}\big)^{(b)}\cdot F_{(1)}\big),
\end{align*}
where $\textstyle F = R'_{23}R'_{13}R'_{24}R'_{14}$. By using Theorem \ref{integteo} as in the proof of Proposition \ref{defL01A}, we see that the sum is a linear combination over $\mc\big[q^{1/D}, q^{-1/D}\big]$ of basis elements $(\alpha_k)^{(a)}(\alpha_l)^{(b)}$. Therefore each coefficient $m_{i_1\cdots i_n,j_1\cdots j_n}^{k_1\cdots k_n}\in \mc(q)\cap \mc\big[q^{1/D}, q^{-1/D}\big] = \mc\big[q, q^{-1}\big]$.

This proves the first claim of the following statement. The second follows from the arguments of Proposition \ref{matrelL0n}, which apply as well to $\Ll_{0,n}^A$.
\begin{prop}\label{defiL0nA} The $A$-submodule $\big(\Ll_{0,1}^A\big)^{\otimes n}$ of $\Ll_{0,n}$ is an $A$-subalgebra. We denote it $\Ll_{0,n}^A$. Moreover, $\Ll_{0,n}^A$ is the $A$-algebra generated by the matrix coefficients $\big({}_X\phi^{e^i}_{e_j}\big)^{(a)}$, for all objects $X$ of the category $\mathcal{C}_A$ and $1\leq a \leq n$, with defining relations given by the naturality relations \eqref{naturality}, the fusion relations \eqref{fusionrel}, and the exchange relations \eqref{exrel2}.
\end{prop}
Note that, by the properties of $\Ll_{0,1}^A$, $\Ll_{0,n}^A$ is a free $A$-module and we have $\Ll_{0,n}^A\otimes_A {\mathbb C}(q)=\Ll_{0,n}$.
\begin{lem}\label{intform-n} The action $\operatorname{coad}^r$ on $\Ll_{0,n}$ yields on $\Ll_{0,n}^A$ a structure of right $U_A$-module algebra, and the Alekseev map restricts to an embedding of $U_A$-module algebras $\Phi_n\colon \Ll_{0,n}^A\ra \tilde{U}_A^{\otimes n}$.
\end{lem}
\begin{proof} The first claim follows immediately from the case of $\Ll_{0,1}^A$ (see Lemma \ref{intform}), since $\Ll_{0,n}^A = \big(\Ll_{0,1}^A\big)^{\otimes n}$ as an $A$-module. For the second claim, one simply note that Theorem \ref{integteo} implies that $R_{Xn}\cdots R_{Xa+1} \in \operatorname{End}_A(X) \otimes \tilde{U}_A^{\otimes n}$, for every $U_A^{\rm res}$-module $X$ of type $1$, and use that ${\rm id} \otimes\Phi_n$ is defined on generators of $\operatorname{End}_A(X) \otimes \Ll_{0,n}^A$ as ${\rm id} \otimes \Phi_1^{\otimes n}$ followed by conjugations with matrices of the form $R_{Xn}\cdots R_{Xa+1}$.\end{proof}

\begin{prop}\label{nozeroq}
$\Ll_{0,n}$, and therefore its subalgebras $\Ll_{0,n}^A$ and $\big(\Ll_{0,n}^A\big)^{U_A}$, does not have non trivial zero divisors.
\end{prop}
\begin{proof} Because of the injectivity of $\Phi_n$ it is sufficient to show that $\tilde{U}_q^{\otimes n}$ and $\tilde{U}_A^{\otimes n}$ have no non trivial zero divisors. We note that $U_A(\mathfrak{g})^{\otimes n}=U_A(\mathfrak{g}^{\oplus n}).$ Then the result for ${U}_A^{\otimes n}$ is a~consequence of \cite[Corollary 1.8]{DC-K} applied to $\mathfrak{g}^{\oplus n}.$ Note that in that paper it is assumed that the Cartan matrix is indecomposable but their method, which consists in proving that an associated graded algebra is quasipolynomial, does not use this assumption. It applies as well for~$\tilde{U}_A^{\otimes n}$,~$\tilde{U}_q^{\otimes n}$. \end{proof}

In the case of $\mathfrak{g}={\mathfrak{sl}}(2)$ we can restrict the target:
\begin{prop}\label{Alekseevsl2} $\Phi_n\colon \Ll_{0,n}({\mathfrak{sl}}(2))\ra \tilde U_q({\mathfrak{sl}}(2))^{{\otimes} n}$ takes values in $U_q({\mathfrak{sl}}(2))^{{\otimes} n}$, and yields an embedding of $U_A$-module algebras $\Phi_n\colon \Ll_{0,n}^A({\mathfrak{sl}}(2))\ra U_A({\mathfrak{sl}}(2))^{{\otimes} n}$.
\end{prop}
\begin{proof}
The second claim follows from the first and Lemma \ref{intform-n}. Consider the first claim. For every $1\leq a\leq n$ and object $V$ of $\mathcal{C}$ we have $R_{Va}R_{Va}' \in \operatorname{End}(V)\otimes U_q^{\rm lf}({\mathfrak{sl}}(2))^{\otimes n}$ by Proposition \ref{imagephi}. We have to show that conjugating by $R_{Vn}\cdots R_{Va+1}\in \operatorname{End}(V)\otimes \tilde U_q({\mathfrak{sl}}(2))^{{\otimes} n}$ maps $\operatorname{End}(V)\otimes U_q^{\rm lf}({\mathfrak{sl}}(2))^{\otimes n}$ to $\operatorname{End}(V)\otimes U_q({\mathfrak{sl}}(2))^{\otimes n}$. It is enough to prove it for $V=V_2$, and for the conjugation by $R_{Va+1}$ only (the general case follows from this one by an easy induction). Recall the expression of $R$ in \eqref{Rsl2}. Using that $E$, $F$ act nilpotently on $V_2$ with order $2$, and $R^{-1} = (S\otimes {\rm id} )(R)$, we~get
\begin{gather}
R_{Va+1} = (\pi_{V_2}\otimes \mathfrak{i}_{a+1})(R) = \begin{pmatrix} 1\otimes q^{H/2} & (q-q^{-1}) 1\otimes q^{H/2}F \\[.5ex] 0 & 1\otimes q^{-H/2}\end{pmatrix}\!, \label{1formRV2}
\\[.5ex]
R_{Va+1}^{-1} = (\pi_{V_2}\otimes \mathfrak{i}_{a+1})(R^{-1}) = \begin{pmatrix} 1\otimes q^{-H/2} & -q(q-q^{-1}) 1\otimes q^{H/2}F \\[.5ex] 0 & 1\otimes q^{H/2}\end{pmatrix}\!, \label{2formRV2}
\end{gather}
where $q^{H/2}\in \tilde U_q({\mathfrak{sl}}(2))$ is defined in Section \ref{Uq}, and for each matrix entry we write only the components in the $a$-th and $a+1$-th tensorands of $\tilde U_q({\mathfrak{sl}}(2))^{{\otimes} n}$ (the others being $1$'s). With this convention, let
\begin{equation}\label{defAmat}
 N := \begin{pmatrix} u \otimes 1 & v \otimes 1 \\ w \otimes 1 & x \otimes 1 \end{pmatrix} \in \operatorname{End}(V)\otimes \tilde U_q({\mathfrak{sl}}(2))^{{\otimes} n}.
\end{equation}
Then
\[
R_{Va+1} N R_{Va+1}^{-1} = \begin{pmatrix} u \otimes 1 & -(q^2-1)u\otimes KF -q^2(q-q^{-1})w\otimes KF^2 \\ {} + q^{-1}(q-q^{-1})w\otimes F & + v \otimes K +q^{-1} (q-q^{-1})x\otimes KF \\[1.5ex] w\otimes K^{-1} & -q(q-q^{-1})w \otimes F + x\otimes 1\end{pmatrix}\!.
\]
This matrix has entries in $U_q({\mathfrak{sl}}(2))^{\otimes n}$ if the matrix entries of $N$ belong to $U_q({\mathfrak{sl}}(2))^{\otimes n}$. Conjugating recursively $R_{Va}R_{Va}'$ with $R_{Vi}$, for $i$ from $a+1$ to $n$, by the same computation we deduce that $\Phi_n$ takes values in $U_q({\mathfrak{sl}}(2))^{\otimes n}$. \end{proof}

\subsection[Localization when g=sl(2)]
{Localization when $\boldsymbol{\mathfrak{g}={\mathfrak{sl}}(2)}$}

We are now going to define a localization of $\Ll_{0,n}({\mathfrak{sl}}(2))$, which will satisfy a generalization of Proposition \ref{Drinfeldloc}. We need the following lemma. For every $u\in U_q({\mathfrak{sl}}(2))$, denote by $u^{(i)}$ the element of $U_q({\mathfrak{sl}}(2))^{\otimes n}$ with $u$ in the $i$-th tensorand and $1$'s elsewhere. Analogously to \eqref{notsl20}, for every $1\leq i \leq n$, put
\begin{equation}\label{notsl2}
\stackrel{V}{M}{}^{\!\! (i)} =\begin{pmatrix} a^{(i)}&b^{(i)}\\c^{(i)}&d^{(i)}\end{pmatrix} \in \operatorname{End}(V)\otimes \Ll_{0,n}({\mathfrak{sl}}(2))^{(i)},
\end{equation}
where $V=V_2$.

{\sloppy\begin{lem}\label{imOK}
For every $1\leq i \leq n$, $\Phi_1^{\otimes n}\big(\Ll_{0,n}^{(i)}({\mathfrak{sl}}(2))\big)$ is contained in the subalgebra of $U_q({\mathfrak{sl}}(2))^{\otimes n}$ generated by $\operatorname{Im}(\Phi_n)$ and the elements $K^{(i+1)},\dots, K^{(n)}$.
\end{lem}}

\begin{proof} The case $i=n$ is clear, as $\Phi_1^{\otimes n}\big(\Ll_{0,n}^{(n)}({\mathfrak{sl}}(2))\big) = \Phi_n\big(\Ll_{0,n}^{(n)}({\mathfrak{sl}}(2))\big)$ by definition. We argue by decreasing induction on $i\in \{1,\dots,n\}$. Take
\[
\tilde{M} = ({\rm id}\otimes \Phi_1^{\otimes n})\big(\!\stackrel{V}{M}{}^{\!\! (n-1)}\big)
\]
in \eqref{defAmat}. The entries of $\tilde{M}$ generate the algebra $\Phi_1^{\otimes n}\big(\Ll_{0,n}^{(n-1)}({\mathfrak{sl}}(2))\big)$. By the formula of $R_{Vn} \tilde{M} R_{0n}^{-1}$ we have $\Phi_n\big(c^{(n-1)}\big) = w\otimes K^{-1}$. Hence $w\otimes 1 = \big(w\otimes K^{-1}\big)(1\otimes K)$ belongs to the algebra generated by $\Phi_n(\Ll_{0,n}({\mathfrak{sl}}(2)))$ and $K^{(n)}$. Since $q^{-1}\big(q-q^{-1}\big)(1\otimes F) = \Phi_1^{\otimes n}\big(b^{(n)}\big) = \Phi_n\big(b^{(n)}\big)$, the same is true of $-q\big(q-q^{-1}\big)w \otimes F$, and hence eventually also of $x\otimes 1$, $u\otimes 1$ and $v\otimes 1$ by using again the formula of $R_{Vn} \tilde{M} R_{Vn}^{-1}$. This proves the statement for $i=n-1$. Inducting on~$i$, using the matrix $\tilde{M}$ of generators of $\Ll_{0,n}^{(i)}({\mathfrak{sl}}(2))$, the same reasoning proves the result for all values of~$i$. \end{proof}

\begin{lem} Define elements $\xi^{(i)}\in \Ll_{0,n}({\mathfrak{sl}}(2))$, $i=1,\dots,n$ by
\begin{equation*}
\xi^{(i)}=\big(\!\stackrel{V}{M}{}^{\!\! (i)}\cdots \stackrel{V}{M}{}^{\!\! (n)}\big)_{22},
\end{equation*}
where $22$ denotes the lower right matrix element. The elements $\xi^{(i)}$ are commuting and satisfy:
\begin{equation*}
\Phi_n\big(\xi^{(i)}\big)=\big(K^{-1}\big)^{(i)}\cdots \big(K^{-1}\big)^{(n)}.
\end{equation*}
\end{lem}
\begin{proof} As in \eqref{phicommut} we have
\begin{gather*}
({\rm id} \otimes \Phi_n)\big(\!\stackrel{V}{M}{}^{\!\! (i)}\cdots \stackrel{V}{M}{}^{\!\! (n)}\big) = R_{Vn}\cdots R_{Vi}R_{Vi}' \cdots R_{Vn}=\big(\pi_{V}\otimes 1^{\otimes (i-1)}\otimes \Delta^{(n-i)}\big)(RR').
\end{gather*}
The lower right matrix element of $(\pi_{V}\otimes {\rm id} )(RR')$ is equal to $K^{-1}$. As a result, by applying~$\Delta^{(n-i)}$ we obtain
 \begin{equation*}
\Phi_n\big(\big(\!\stackrel{V}{M}{}^{\!\! (i)}\cdots \stackrel{V}{M}{}^{\!\! (n)}\big)_{22}\big) = 1^{\otimes (i-1)}\otimes \Delta^{(n-i)}\big(K^{-1}\big)=\big(K^{-1}\big)^{(i)}\cdots \big(K^{-1}\big)^{(n)}.
\end{equation*}
By injectivity of $\Phi_n$ this proves that the elements $\xi^{(i)}$ are commuting. \end{proof}

The elements $\xi^{(i)}$ commute, and $\big\{\big(\xi^{(1)}{}^k\big)\big\}_{k\in {\mathbb N}}$ is an Ore set of $\Ll_{0,n}({\mathfrak{sl}}(2))$. In fact, it is easy to see that for every $i=1,\dots, n$ the element $\xi^{(i)}$ is an Ore element of $\Ll_{0,n}^{(i\leq)}({\mathfrak{sl}}(2))$, where $\Ll_{0,n}^{(i\leq)}({\mathfrak{sl}}(2))$ is the subalgebra of $\Ll_{0,n}({\mathfrak{sl}}(2))$ generated by the subalgebras $\Ll_{0,n}^{(a)}({\mathfrak{sl}}(2))$, $a\geq i$. Indeed, because~$\Ll_{0,n}$ has no non trivial zero divisors, $\xi^{(i)}$ is a regular element. The set $\big\{\big(\xi^{(i)}{}^k\big)\big\}_{k\in {\mathbb N}}$ is multiplicatively closed, and it is an Ore subset of $\Ll_{0,n}^{(i\leq)}({\mathfrak{sl}}(2))$, since for all $x\in \Ll_{0,n}^{(i\leq)}({\mathfrak{sl}}(2))$ there exist elements $y,y'\in \Ll_{0,n}^{(i\leq)}({\mathfrak{sl}}(2))$ such that
$x\xi^{(i)}=\xi^{(i)}y$ and $\xi^{(i)}x=y'\xi^{(i)}$. This is shown as follows: $\Phi_n(x)\Phi_n\big(\xi^{(i)}\big)=
\Phi_n(x)\big(K^{-1}\big)^{(i)}\cdots \big(K^{-1}\big)^{(n)} =\big(K^{-1}\big)^{(i)}\cdots \big(K^{-1}\big)^{(n)}\operatorname{ad}^r(K)(\Phi_n(x))$, where we have used the fact that $\Phi_n(x)\in 1^{\otimes i-1}\otimes U_q^{\otimes (n-i+1)}$. But $\operatorname{ad}^r(K)(\Phi_n(x))=\Phi_1(\operatorname{coad}^r(K)(x)),$ and therefore the Ore conditions are satisfied with $y=\operatorname{coad}^r(K)(x)$.

This argument does not permit to conclude that the multiplicative sets $\big\{\big(\xi^{(i)}{}^k\big)\big\}_{k\in {\mathbb N}}$ are Ore sets of $\Ll_{0,n}({\mathfrak{sl}}(2))$. We therefore cannot localize with respect to the elements $\xi^{(i)}$ as easily as for~$\xi^{(1)}$. We proceed in a different way.

We first explain the case $n=2$. The element $\xi^{(1)}=\big(\!\stackrel{V}{M}{}^{\!\! (1)}\stackrel{V}{M}{}^{\!\! (2)}\big)_{22}=d^{(1)}d^{(2)}+c^{(1)}b^{(2)}$ is a regular Ore element of $\Ll_{0,2}$, so we can define the localisation
$\Ll_{0,2}\big[\xi^{(1)}{}^{-1}\big].$ We want to define the inverse of the element $\xi^{(2)}$, and a new algebra $\Ll_{0,2}\big[\xi^{(1)}{}^{-1}\big]\big[\xi^{(2)}{}^{-1}\big]$ such that $\Ll_{0,2}\big[\xi^{(1)}{}^{-1}\big]\subset \Ll_{0,2}\big[\xi^{(1)}{}^{-1}\big]\big[\xi^{(2)}{}^{-1}\big]$ and $\Phi_2$ extends naturally to an homomorphism of algebras
\[
\Phi_2\colon\ \Ll_{0,2}\big[\xi^{(1)}{}^{-1}\big]\big[\xi^{(2)}{}^{-1}\big]\rightarrow U_q({\mathfrak{sl}}(2))^{\otimes 2}.
\]

This can be done by writing explicitely the exchange relations between $\stackrel{V}{M}{}^{\!\! (1)}$ and $\stackrel{V}{M}{}^{\!\! (2)}$ involving $d^{(2)}=\xi^{(2)}$:
\begin{gather*}
c^{(1)}d^{(2)}=d^{(2)}c^{(1)},\\
d^{(2)}a^{(1)}=a^{(1)}d^{(2)}+\big(1-q^{-2}\big) c^{(1)}b^{(2)},\\
d^{(2)}d^{(1)}=d^{(1)}d^{(2)}+\big(1-q^{2}\big)c^{(1)}b^{(2)},\\
d^{(2)}b^{(1)}=b^{(1)}d^{(2)}+\big(1-q^2\big)\big(a^{(1)}-d^{(1)}\big)b^{(2)}.
\end{gather*}

We define $\Ll_{0,2}\big[\xi^{(1)}{}^{-1}\big]\big[\xi^{(2)}{}^{-1}\big]$ to be the algebra generated by the elements
$a^{(1)}$, $b^{(1)}$, $c^{(1)}$, $d^{(1)}$,
$ a^{(2)}$, $b^{(2)}$, $c^{(2)}$, $d^{(2)}$, $\xi^{(1)}{}^{-1}$, $\xi^{(2)}{}^{-1}$,
where $a^{(1)},\dots, d^{(2)}$ satisfy the exchange relation (\ref{rel01}) and (\ref{exrel2}), $\xi^{(1)}{}^{-1}$ is the inverse in the sense of Ore of the element
$\xi^{(1)}=d^{(1)}d^{(2)}+c^{(1)}b^{(2)}$ and the following relations for $\xi^{(2)}{}^{-1}$, the inverse of $\xi^{(2)}=d^{(2)},$ hold true:
\begin{gather*}
\xi^{(2)} \xi^{(2)}{}^{-1}=\xi^{(2)}{}^{-1}\xi^{(2)}=1, \\
\xi^{(2)}{}^{-1}c^{(1)}=c^{(1)}\xi^{(2)}{}^{-1},\\
\xi^{(2)}{}^{-1}a^{(1)}=a^{(1)}\xi^{(2)}{}^{-1}-q^{-2}\big(1-q^{-2}\big) c^{(1)}b^{(2)}\xi^{(2)}{}^{-2},\\
\xi^{(2)}{}^{-1}d^{(1)}=d^{(1)}\xi^{(2)}{}^{-1}-q^{-2}\big(1-q^{2}\big)c^{(1)}b^{(2)}\xi^{(2)}{}^{-2},\\
\xi^{(2)}{}^{-1}b^{(1)}=b^{(1)}\xi^{(2)}{}^{-1}+q^{-1}\big(1-q^{-2}\big)\big(a^{(1)}-d^{(1)}\big)b^{(2)}\xi^{(2)}{}^{-2}
\\ \hphantom{\xi^{(2)}{}^{-1}b^{(1)}=}
{}+q^{-5}\big(1-q^{2}\big)\big(1-q^{-4}\big) c^{(1)}\big(b^{(2)}\big)^2\xi^{(2)}{}^{-3}.
\end{gather*}

The last relations are chosen in order that $\Ll_{0,2}\big[\xi^{(1)}{}^{-1}\big]\subset \Ll_{0,2}\big[\xi^{(1)}{}^{-1}\big]\big[\xi^{(2)}{}^{-1}\big]$ as an algebra and are a direct consequence of the exchange relations (\ref{exrel2}) and the invertibility of $\xi^{(2)}.$
As~a~result $\Phi_2$ extends to a morphism of algebra $\Phi_2\colon \Ll_{0,2}\big[\xi^{(1)}{}^{-1}\big]\big[\xi^{(2)}{}^{-1}\big]\rightarrow U_q({\mathfrak{sl}}(2))^{\otimes 2}.$ This morphism of algebra will be shown to be an isomorphism in Proposition \ref{Alekseevsl2loc}.

The construction of the localisation for $\Ll_{0,n}({\mathfrak{sl}}(2))$, $n>2$, is defined by the same procedure.
We want to define a localisation $\Ll_{0,n}({\mathfrak{sl}}(2))\big[\xi^{(n)}{}^{-1},\dots, \xi^{(1)}{}^{-1}\big]$. The set of elements $\big\{\xi^{(n)k}\big\}_k$ being an Ore subset of $\Ll_{0,n}({\mathfrak{sl}}(2))$, we can define the standard localisation $\Ll_{0,n}({\mathfrak{sl}}(2))\big[\xi^{(n)}{}^{-1}\big]$. Note that the matrix $\stackrel{V}{M}{}^{\!\! (2)}\cdots \stackrel{V}{M}{}^{\!\! (n)}$ has the same reflection equation with $\stackrel{V}{M}{}^{\!\! (1)}$ as
$\stackrel{V}{M}{}^{\!\! (2)}.$ Therefore we can define $\Ll_{0,n}({\mathfrak{sl}}(2))\big[\xi^{(n)}{}^{-1}\big]\big[\xi^{(n-1)}{}^{-1}\big]$ by the same method as when $n=2$, and so on.

\begin{defi}\label{Lonloc}
By iterating the preceeding construction we define:
\[
{}_{\rm loc}\Ll_{0,n}({\mathfrak{sl}}(2))=\Ll_{0,n}({\mathfrak{sl}}(2))\big[\xi^{(n)}{}^{-1}\big] \big[\xi^{(n-1)}{}^{-1}\big]\cdots \big[\xi^{(1)}{}^{-1}\big].
\]
In the sequel it will be convenient to define invertible elements $\delta^{(i)}\in {}_{\rm loc}\Ll_{0,n}({\mathfrak{sl}}(2))$, $i=1,\dots,n,$ satisfying
$\xi^{(i)}=\delta^{(i)}\cdots\delta^{(n)}$, i.e., $\delta^{(i)}=\xi^{(i)}\xi^{(i+1)}{}^{-1}$.
\end{defi}
The elements $\delta^{(i)}$, $i=1,\dots,n,$ are invertible, commute and satisfy $\Phi_n\big(\delta^{(i)}\big) = \big(K^{-1}\big)^{(i)}$ .

We can define the localization ${}_{\rm loc}\Ll_{0,n}^A({\mathfrak{sl}}(2))$ of the integral form $\Ll_{0,n}^A({\mathfrak{sl}}(2))$ in the very same way (recall Lemma \ref{intsl2} for the case $n=1$). Let
\[
\Phi_n\colon\quad {}_{\rm loc}\Ll_{0,n}({\mathfrak{sl}}(2))\ra U_q({\mathfrak{sl}}(2))^{{\otimes} n}
\]
be the unique morphism of module algebras extending the Alekseev map. Recall the subalgebra $U_A'({\mathfrak{sl}}(2))$ of $U_A({\mathfrak{sl}}(2))$, defined in Lemma~\ref{intsl2}.
\begin{prop}\label{Alekseevsl2loc} The extended Alekseev map $\Phi_n\colon {}_{\rm loc}\Ll_{0,n}({\mathfrak{sl}}(2))\ra U_q({\mathfrak{sl}}(2))^{{\otimes} n}$ is an isomorphism of $U_q$-module algebras, and it restricts to an isomorphism of $U_A$-module algebras $\Phi_n\colon {}_{\rm loc}\Ll_{0,n}^A({\mathfrak{sl}}(2))\ra U_A'({\mathfrak{sl}}(2))^{{\otimes} n}$.
\end{prop}
\begin{proof} All claims are clear by the previous results, except surjectivity. For the first claim, it follows from Lemma \ref{imOK}, since the algebra generated by $\Phi_1^{\otimes n}\big(\Ll_{0,n}^{(i)}({\mathfrak{sl}}(2))\big)$ and $K^{(i)}$ is $U_q({\mathfrak{sl}}(2))^{(i)}$ (the $i$-th tensorand). For the second claim, it follows from Lemma \ref{intsl2}(3).\end{proof}

\begin{remk} {\rm This proposition (first claim) justifies the localization of $\Ll_{0,n}({\mathfrak{sl}}(2))$ by the elements $\xi^{(1)},\dots,\xi^{(n)}$. The result of localization, if any, by the elements $d^{(1)},\dots,d^{(n)}$ is not clear to us; note that they do not commute.}
\end{remk}
\subsection{Invariant elements} We go back to the general situation of $\Ll_{0,n}^A = \Ll_{0,n}^A(\mathfrak{g})$ for an arbitrary $\mathfrak{g}$. We need the following fact, proved by Majid \cite{Majid} for $\Oq$. For completeness we recall the main ideas of the proof given in \cite{DKM}, which extends immediately to $\mathcal{O}_A$.
\begin{prop}\label{coprodprop} The iterated coproduct $\Delta^{(n-1)}\colon \mathcal{O}_A \ra \mathcal{O}_A^{\otimes n}$, considered as a linear map $\Ll_{0,1}^A\ra \Ll_{0,n}^A$, is an algebra morphism and satisfies the following commutative diagram:
\[\xymatrix{ \Ll_{0,1}^A \ar[r]^{\Delta^{(n-1)}} \ar[d]_{\Phi_1} & \Ll_{0,n}^A \ar[d]^{\Phi_n} \\ \tilde{U}_A\ar[r]^{\Delta^{(n-1)}} & \tilde U_A^{\otimes n}.}\]
\end{prop}
\begin{proof} The main point of the proof is that, for every $U_A^{\rm res}$-module $X$ of type $1$, we have
\begin{equation}\label{prodmatrix}
\big({\rm id} \otimes \Delta^{(n-1)}\big)\big(\!\stackrel{X}{M}\!\big)
= \,\stackrel{X}{M}{}^{\!\! (1)}\cdots \stackrel{X}{M}{}^{\!\! (n)}.
\end{equation}
This follows by a direct application of the definitions. But Proposition \ref{fusrel} implies that the map $\Delta^{(n-1)}\colon \Ll_{0,1}^A\ra \Ll_{0,n}^A$ is an algebra morphism if and only if the matrix on the left hand side of~\eqref{prodmatrix} satisfies the fusion relation. This is easily shown to be true for the right hand side by using the exchange relation \eqref{exrel2} recursively. This proves the proposition. The commutativity diagram is a reformulation of \eqref{phicommut}.\end{proof}

Consider now the algebras of invariant elements,
\begin{gather*}
\Ll_{0,n}^{U_q} :=\big\{\alpha \in \Ll_{0,n}\mid \forall y \in U_q,\ \operatorname{coad}^r(y)(\alpha) = \varepsilon(y)\alpha\big\},
\\
\big(\Ll_{0,n}^A\big)^{U_A} :=\big\{\alpha \in \Ll_{0,n}^A\mid \forall y \in U_A,\ \operatorname{coad}^r(y)(\alpha) = \varepsilon(y)\alpha\big\}.
\end{gather*}
Note that $\big(\Ll_{0,n}^A\big)^{U_A} = \Ll_{0,n}^{U_q} \cap \Ll_{0,n}^A$, and therefore $\Ll_{0,n}^{U_q} = \big(\Ll_{0,n}^A\big)^{U_A} \otimes_A \mc(q)$. Since $\Ll_{0,n}^A$ is a free $A$-module, and $A$ is a principal ideal domain, the $A$-submodule $\big(\Ll_{0,n}^A\big)^{U_A}$ is free.
\begin{prop} \label{invL0n} The algebra $\big(\Ll_{0,n}^A\big)^{U_A}$ is the centralizer of $\Delta^{(n-1)}\big(\Ll_{0,1}^A\big)$ in $\Ll_{0,n}^A$. As a corollary we have
$\mathcal{Z}\big(\Ll_{0,n}^A\big)^{U_A} = \mathcal{Z}\big(\Ll_{0,n}^A\big)$.
\end{prop}
\begin{proof} Clearing denominators it is enough to give the details for $\Ll_{0,n}^{U_q}$. First note that an element $z$ of $\tilde U_q^{ \otimes n}$ is invariant under the right adjoint action of $U_q$ if and only if it commutes with $ \Delta^{(n-1)}(x)$ for every $x\in U_q$. Indeed, we have
\begin{align*}
z\Delta^{(n-1)}(x) & = \sum_{(x)} z \varepsilon(x_{(1)})\Delta^{(n-1)}(x_{(2)}) = \sum_{(x)} \varepsilon(x_{(1)}) z \Delta^{(n-1)}(x_{(2)})
\\
& =\! \sum_{(x)} \Delta^{(n-1)}(x_{(1)})\Delta^{(n-1)}( S(x_{(2)})) z \Delta^{(n-1)}(x_{(3)})
=\! \sum_{(x)}\!\Delta^{(n-1)}(x_{(1)}) \operatorname{ad}^r(x_{(2)})(z).
\end{align*}
Hence $\textstyle z\Delta^{(n-1)}(x) = \sum_{(x)} \Delta^{(n-1)}(x_{1}) \varepsilon(x_{(2)}) z = \Delta^{(n-1)}(x)z$ if $z$ is an invariant element. Conversely, if $z$ commutes with $\Delta^{(n-1)}(x)$ for every $x\in U_q$, then
\[
\operatorname{ad}^r(x)(z) = \sum_{(x)} \Delta^{(n-1)}( S(x_{(1)})) \Delta^{(n-1)}(x_{(2)}) z = \Delta^{(n-1)}(\varepsilon(x)) z = \varepsilon(x) z.
\]
This proves our claim. Now, let $\alpha\in \Ll_{0,n}^{U_q}$. Then $\Phi_n(\alpha)$ is an $\operatorname{ad}^r(U_q)$-invariant element of $\tilde U_q^{\otimes n}$, by Theorem \ref{Alekseevmap}. By the claim above, and the fact that commuting with elements of $\tilde U_q$ or~$U_q$ is the same, $\Phi_n(\alpha)$ commutes with the matrix coefficients of $\big(\pi_V \otimes \Delta^{(n-1)}\big)(RR')$, for every object~$V$ of $\mathcal{C}$. By the injectivity of $\Phi_n$ and the relation \eqref{phicommut}, $\alpha$ commutes with the matrix coefficients of \eqref{prodmatrix} for all objects $V$ of $\mathcal{C}$. These generate the algebra $\Delta^{(n-1)}(\Ll_{0,1})$, so $\Ll_{0,n}^{U_q}$ lies in the centralizer of the latter.

{\sloppy
Conversely, the same reasoning shows that if $\alpha$ lies in the centralizer of $\Delta^{(n-1)}(\Ll_{0,1})$, then~$\Phi_n(\alpha)$ commutes with the matrix coefficients of $\big(\pi_V \otimes \Delta^{(n-1)}\big)(RR')$, for every object~$V$ of~$\mathcal{C}$. From Theorem \ref{drinfeldmap} we deduce that $\Phi_n(\alpha)$ lies in the centralizer of $\Delta^{(n-1)}\big(\tilde U_q^{\rm lf}\big)$. This is the same as the centralizer of $\Delta^{(n-1)}(\tilde U_q)$ by Remark \ref{remlf}(1), so as above we deduce that $\Phi_n(\alpha)$ is an invariant element of $\tilde U_q^{\otimes n}$, and by injectivity and equivariance of $\Phi_n$ that $\alpha$ is an invariant element of~$\Ll_{0,n}$.

}

The corollary is immediate, because a central element of $\Ll_{0,n}$ necessarily commutes with $\Delta^{(n-1)}(\Ll_{0,1})$, so it is invariant.\end{proof}

We now give an explicit basis of the algebra of invariant elements $\Ll_{0,n}^{U_q}$. Let $n$ be an integer greater than $1$. For every $2\leq k\leq n$ we denote by $S(k)\in {\mUq}^{{\otimes} n}$ the element defined by
\[
S(k)={\rm id}^{\otimes (k-2)} \otimes \big({\rm id} \otimes \Delta^{(n-k)}\big)(R).
\]
Let $\lambda_1,\dots,\lambda_n\in P_+$, and $V_{\lambda_1},\dots,V_{\lambda_n}$ the type $1$ simple $U_q$-modules of highest weights $\lambda_1,\dots,\lambda_n$ respectively. Put $[\lambda]=(\lambda_1,\dots,\lambda_n)$, and consider the $U_q$-module
\[
V_{[\lambda]}=\bigotimes_{j=1}^n V_{\lambda_j}.
\]
Define $\stackrel{[\lambda]}{\mathbb M}\,\in \operatorname{End}\big(V_{[\lambda]}\big)\otimes\Ll_{0,n}$ by $\stackrel{[\lambda]}{\mathbb M}\,=\,\stackrel{V_\lambda}{ M}$ if $n=1$, and if $n\geq 2$ by
\begin{equation}\label{Matrixlambda}
\stackrel{[\lambda]}{\mathbb M}\,=\, \stackrel{V_{\lambda_1}}{ M}_1{}^{\!\!\! (1)} \prod_{k=2}^{n}
\big(S(k)^{-1}_{V_{\lambda_1},\dots,V_{\lambda_n}}\stackrel{V_{\lambda_k}}{ M}_k{}^{\!\!\! (k)} \big)\prod_{k=n}^{2} S(k)_{V_{\lambda_1},\dots,V_{\lambda_n}} .
\end{equation}
For example for $n=2$ we have
\[
\stackrel{[\lambda]}{\mathbb M}\,=\,
\stackrel{V_{\lambda_1}}{M}_1{}^{\!\!\! (1)}R_{12}^{-1}\stackrel{V_{\lambda_2}}{M}_2{}^{\!\!\! (2)}R_{12}
\]
 and for $n=3$,
\[
\stackrel{[\lambda]}{\mathbb M}\,=\, \stackrel{V_{\lambda_1}}{M}_1{}^{\!\!\! (1)}({\rm id}\otimes \Delta )(R^{-1})\stackrel{V_{\lambda_2}}{M}_2{}^{\!\!\! (2)}
(R^{-1})_{23}\stackrel{V_{\lambda_3}}{M}_3{}^{(3)}R_{23}({\rm id}\otimes \Delta )(R).
\]
The reason for considering the matrices $\stackrel{[\lambda]}{\mathbb M}$ comes from the two following propositions. The first generalizes the identity \eqref{actmat}.
\begin{prop}\label{coadmat2} The right coadjoint action of $U_q$ on the matrix elements of $\stackrel{[\lambda]}{\mathbb M}$ can be written in matrix form as
\begin{equation}\label{actmat2}
\operatorname{coad}^r(y)\big(\!\stackrel{[\lambda]}{\mathbb M}\!\big) = \sum_{(y)}\big((\pi_{V_{[\lambda]}}(y_{(1)}) \otimes {\rm id} )
\stackrel{[\lambda]}{{\mathbb M}}(\pi_{V_{[\lambda]}}( S(y_{(2)})) \otimes {\rm id} )\big).
\end{equation}
\end{prop}
\begin{proof} We show it for $n=2$, since the general proof follows from it by an easy induction on~$n$. We have
\begin{align*}
\operatorname{coad}^r(y)\big(\!\stackrel{[\lambda]}{\mathbb M}\!\big)&=
\operatorname{coad}^r(y)\big(\!\!\stackrel{V_{\lambda_1}}{M}_1{}^{\!\!\! (1)}R_{12}^{-1}\stackrel{V_{\lambda_2}}{M}_2{}^{\!\!\! (2)}R_{12}\big)
\\
&=\sum_{(y)}(y_{(1)})_1\stackrel{V_{\lambda_1}}{M}_1{}^{\!\!\! (1)} S(y_{(2)})_1R_{12}^{-1}(y_{(3)})_2\stackrel{V_{\lambda_2}}{M}_2{}^{\!\!\! (2)} S(y_{(4)})_2R_{12}
\\
&=\sum_{(y)}(y_{(1)})_1\stackrel{V_{\lambda_1}}{M}_1{}^{\!\!\! (1)}(y_{(2)})_2 R_{12}^{-1}( S(y_{(3)}))_1\stackrel{V_{\lambda_2}}{M}_2{}^{\!\!\! (2)}
S(y_{(4)})_2R_{12}
\\
&=\sum_{(y)}(y_{(1)}\otimes y_{(2)}) \stackrel{V_{\lambda_1}}{M}_1{}^{\!\!\! (1)}R_{12}^{-1}\stackrel{V_{\lambda_2}}{M}_2{}^{\!\!\! (2)}( S(y_{(3)}) \otimes S(y_{(4)}))R_{12}
\\
&=\sum_{(y)}(y_{(1)}\otimes y_{(2)})\stackrel{V_{\lambda_1}}{M}_1{}^{\!\!\! (1)}R_{12}^{-1}\stackrel{V_{\lambda_2}}{M}_2{}^{\!\!\! (2)}R_{12}
\big( S(y_{(4)}) \otimes S(y_{(3)})\big)
\\
&=\sum_{(y)}\big(\big(\pi_{V_{[\lambda]}}(y_{(1)}) \otimes {\rm id} \big)
\stackrel{[\lambda]}{{\mathbb M}}\big(\pi_{V_{[\lambda]}}( S(y_{(2)})) \otimes {\rm id} \big)\big).
\end{align*}
All equalities are clear but the third and fifth, which follow from $\big({\rm id} \otimes S^{-1}\big)(R) = R^{-1}$, $R \Delta R^{-1} = \Delta^{\rm cop}$, $ \Delta \circ S = ( S \otimes S) \circ \Delta$, and the fact that $S$ is an algebra antimorphism. The result follows. \end{proof}

Denote by $\big(\!\stackrel{\lambda}{e}_i\!\big)$ a basis of $V_{\lambda}$, by $\big(\!\stackrel{\lambda}{e}{}^{\!\! i}\big)$ the dual basis, and put $\stackrel{[\lambda]}{e}_{[i]} := \stackrel{\lambda_1}{e}_{i_1}\otimes \dots \otimes \stackrel{\lambda_n}{e}_{i_n}$. We~thus get a basis $\big(\!\stackrel{[\lambda]}{e}{}^{\!\! [i]}\big)$ of $V_{[\lambda]}$. Denote by $\big(\!\stackrel{[\lambda]}{e}{}^{\!\! [i]}\big)$ the dual basis.

\begin{prop}\label{canonicform}
The elements $ \big(\!\stackrel{[\lambda]}{e}{}^{\!\!\! [j]}\otimes {\rm id} \big)\stackrel{[\lambda]}{\mathbb M}\big(\!\stackrel{[\lambda]}{e}{}_{\!\!\! [i]} \otimes {\rm id} \big)$ of $\Ll_{0,n}$, where $[\lambda] \in P_+^n$, $[i]$ labels the basis of $V_{[\lambda]}$, and $[j]$ the dual basis, form a basis of $\Ll_{0,n}$. Equivalently, for every $x\in \Ll_{0,n}$ there exists a unique family of endomorphisms $a_{[\lambda]}(x)\in \operatorname{End}_{\mc(q)}\big(V_{[\lambda]}\big)$, $[\lambda]\in P_+^n$, which are zero except possibly for a finite number of terms, such that
\begin{equation}\label{linbase}
x=\sum_{[\lambda]}\big(\!\operatorname{Tr}_{{V_{[\lambda]}}}\otimes {\rm id} \big)\big(\big(\pi_{V_{\lambda}}(\ell)a_{[\lambda]}(x)\otimes {\rm id} \big)\stackrel{[\lambda]}{\mathbb M}\!\big)=\sum_{[\lambda]}\operatorname{qTr}_{{V_{[\lambda]}}}\big(\big(a_{[\lambda]}(x)\otimes {\rm id} \big)\stackrel{[\lambda]}{\mathbb M}\!\big),
\end{equation}
where as usual $\ell$ is the pivotal element and $\pi_{V_{\lambda}}(\ell)$ is the endomorphism of $V_{[\lambda]}$ given by the action of $\ell$.
\end{prop}
Note that $\pi_{V_{\lambda}}(\ell)$ is introduced in \eqref{linbase} in order to simplify the statement of the next proposition.
\begin{proof} We show the first claim of the proposition for $n=2$; the general proof is similar. In~this case it is sufficient to prove that the set of matrix elements of the tensors
\[
\stackrel{V_{\lambda_1}}{M}_1{}^{\!\!\! (1)}R_{12}^{-1}\stackrel{V_{\lambda_2}}{M}_2{}^{\!\!\! (2)}=\sum_{(R)}(R_{(2)})_2\stackrel{V_{\lambda_1}}{M}_1{}^{\!\!\! (1)}\stackrel{V_{\lambda_2}}{M}_2{}^{\!\!\! (2)}S(R_{(1)})_1,\qquad [\lambda]\in P_+^n,
\]
is a basis of $\Ll_{0,2}$. Note that we use $\textstyle R^{-1}=(S \otimes {\rm id} )(R)=\sum_{(R)}S(R_{(1)}) \otimes R_{(2)}$ in the above equality. Let $T=( S\otimes S^{-1})\big(R^{-1}\big)$, that we denote as usual by $\textstyle T=\sum_{(T)}T_{(1)} \otimes T_{(2)}$. From the identity $1 \otimes 1 = ( S \otimes {\rm id} )\big(R^{-1}R\big) = \sum_{(R),(T)} S(R_{(1)})T_{(1')}\otimes S(T_{(2') })R_{(2)}$, we obtain
\[
\stackrel{V_{\lambda_1}}{M}_1{}^{\!\!\! (1)}\stackrel{V_{\lambda_2}}{M}_2{}^{\!\!\! (2)}=\sum_{(T)} S(T_{(2)})_2\stackrel{V_{\lambda_1}}{M}_1{}^{\!\!\! (1)}R_{12}^{-1}\stackrel{V_{\lambda_2}}{M}_2{}^{\!\!\! (2)}(T_{(1)})_1.
\]
The conclusion follows from this, since the set of matrix elements of $\stackrel{V_{\lambda_1}}{M}_1{}^{\!\!\! (1)}\stackrel{V_{\lambda_2}}{M}_2{}^{\!\!\! (2)}$, $[\lambda]\in P_+^2$, is a basis of $\Ll_{0,2}$. The second claim is a trivial consequence of the first one. \end{proof}

By combining the two previous propositions we obtain the following characterisation of the invariants elements.

\begin{prop} \label{linbaseinv} An element $x\in \Ll_{0,n}$ belongs to $\Ll_{0,n}^{U_q}$ if and only if for every $[\lambda] \in P_+^n$ we have $a_{[\lambda]}(x)\in \operatorname{End}_{U_q}\big(V_{[\lambda]}\big)$. Therefore, the elements
\begin{equation}\label{defvectinv}
v_{[\lambda]}\big(a_{[\lambda]}^{(k)}\big) := \operatorname{qTr}_{V_{[\lambda]}} \big( a_{[\lambda]} ^{(k)}\stackrel{[\lambda]}{\mathbb M}\!\big),
\end{equation}
where $\big\{a_{[\lambda]}^{(k)}\big\}_{k}$ is any basis of $\operatorname{End}_{U_q}\big(V_{[\lambda]}\big)$, make a basis of $\Ll_{0,n}^{U_q}$.
\end{prop}
\begin{proof}
We have $x\in \Ll_{0,n}^{U_q}$ if and only if $\operatorname{coad}^r(y)(x)=\epsilon(y)x$ for every $y\in U_q$, or equivalently for every $y\in U_q$. From Proposition \ref{coadmat2}, the basis provided by Proposition \ref{canonicform}, and the fact that $S^2(y) = \ell y\ell^{-1}$, this is equivalent to $\textstyle \sum_{(y)}S^{-1}(y_{(2)})_{V_{[\lambda]}}a_{[\lambda]}(x)(y_{(1)})_{V_{[\lambda]}}=\epsilon(y) a_{[\lambda]}(x)$ for every $[\lambda] \in P_+^n$, which is also equivalent to $y_{V_{[\lambda]}}a_{[\lambda]}(x)=a_{[\lambda]}(x)y_{V_{[\lambda]}}$.\end{proof}

Finally, let us consider the integral form $\big(\Ll_{0,n}^A\big)^{U_A}$. Matrices $\stackrel{[X]}{\mathbb M}$ can be defined by the formula~\eqref{Matrixlambda}, replacing $V_{\lambda_1}, \dots, V_{\lambda_n}$ with arbitrary $U_A^{\rm res}$-modules $X_1,\dots,X_n$ of type $1$. Clearly these matrices still satisfy the equivariance property of Proposition \ref{coadmat2}, and the arguments of Proposition \ref{canonicform} imply that their matrix elements form a generating family of the $A$-mo\-du\-le~$\Ll_{0,n}^A$. When the $U_A^{\rm res}$-modules $X_i$ span the set of full $A$-sublattices ${}_AV_{\lambda_i}$ of the $U_q$-modules~$V_{\lambda_i}$, such elements form a free family. One can still define invariant elements $v_{X}(a_{X})\in \big(\Ll_{0,n}^A\big)^{U_A}$ by
\[
v_{X}(a_{X}) := \operatorname{qTr}_{X} \big( a_{X}\stackrel{X}{\mathbb M}\!\big),
\]
where $X:=X_1\otimes \dots \otimes X_n$ for arbitrary $U_A^{\rm res}$-modules $X_1,\dots,X_n$ of type $1$, and $a_{X} \in \operatorname{End}_{U_A^{\rm res}}(X)$. The $A$-module $\operatorname{End}_{U_A^{\rm res}}(X)$ is free, as it is a submodule of the free $A$-mo\-du\-le~$\operatorname{End}_{U_A^{\rm res}}(X)$ and $A$ is a principal ideal domain. Then, taking $X_i := {}_AV_{\lambda_i}$ and $a_{X}^{(k)}$ basis elements of $\operatorname{End}_{U_A^{\rm res}}(X)$, the elements $v_{X}\big(a_{X}^{(k)}\big)$ form a free family over $A$, and by Proposition \ref{linbaseinv} a basis of $\Ll_{0,n}^{U_q}$.

In the case of $\mathfrak{g} = {\mathfrak{sl}}(2)$, let us state the following result. It is a direct consequence of Theo\-rem~\ref{teoOBS} and the fact that the skein algebra $K_\zeta(\Sigma)$ is finitely generated and Noetherian, which is proved in \cite[Theorem 3]{PS}, by topological means. The method relies on the fact that~$K_\zeta(\Sigma)$ has a natural filtration, and consists in proving that the associated graded algebra is finitely generated and Noetherian.
\begin{teo}
The algebra $\big(\Ll_{0,n}^A({\mathfrak{sl}}(2))\big)^{U_A({\mathfrak{sl}}(2))}$ is finitely generated and Noetherian.
\end{teo}

\subsection{Centers}
Next we turn to the center $\mathcal{Z}(\Ll_{0,n})$ of $\Ll_{0,n}$ and $\mathcal{Z}\big(\Ll_{0,n}^{U_q}\big)$ of $\Ll_{0,n}^{U_q}$.
\begin{prop} \label{center0n}
The Alekseev map affords an isomorphism from $\mathcal{Z}(\Ll_{0,n})$ to $\mathcal{Z}\big(\tilde U_q^{\otimes n}\big)$. In~particular this implies $\mathcal{Z}(\Ll_{0,n}) = \mathcal{Z}(\Ll_{0,1})^{\otimes n}$. Moreover, the elements
\[
{}_{\lambda}\omega^{(i)}:= \operatorname{qTr}_{V_{\lambda}}\big(\!\stackrel{V_{\lambda}}{M}{}^{\!\! (i)}\big), \qquad
\lambda\in P_+, i\in \{1,\dots,n\},
\]
belong to the center $\mathcal{Z}(\Ll_{0,n})$ and the family of elements
$\textstyle \prod_{i=1}^n {}_{\lambda_i}\omega^{(i)}$, where $\lambda_1,\dots,\lambda_n\in P_+$, form a basis of $\mathcal{Z}(\Ll_{0,n})$. The elements ${}_{\lambda}\omega^{(i)}$ belong to $\Ll_{0,n}^A$, and they form a free family.
\end{prop}
\begin{proof}
We have $\Phi_n(\mathcal{Z}(\Ll_{0,n})) = \mathcal{Z}(\Phi_n(\Ll_{0,n}))$ by Theorem \ref{Alekseevmap}. We claim that if $x\in \tilde U_q^{\otimes n}$ commutes with $\Phi_n(\Ll_{0,n})$, it is central. Indeed, $x$ commutes with $\Phi_n\big(\Ll_{0,n}^{(n)}\big)=1^{\otimes (n-1)} \otimes \tilde U_q^{\rm lf}$, and hence with $1^{\otimes (n-1)} \otimes \tilde U_q$, and the matrix coefficients of $R_{Vn}^{\pm 1}$ by Remark \ref{remlf}. (For instance, when $\mathfrak{g}={\mathfrak{sl}}(2)$, $K=q^H$ and $q^{H/2}$ commute and are diagonalizable on objects of $\mathcal{C}$, and since $x$ commutes with $1^{\otimes (n-1)} \otimes K$, it commutes with $1^{\otimes (n-1)} \otimes q^{H/2}$ too, and the conclusion follows from \eqref{1formRV2} and \eqref{2formRV2}). Using that $x$ commutes also with $\Phi_n\big(\Ll_{0,n}^{(n-1)}\big)$, and hence with the matrix coefficients of $R_{Vn} R_{Vn-1}R_{Vn-1}'R_{Vn}^{-1}$, we deduce that it commutes with the matrix coefficients of $R_{Vn-1}R_{Vn-1}'$, and hence with every element of $1^{\otimes (n-2)} \otimes \tilde U_q \otimes 1$. Continuing in this way recursively, we get that $x\in \mathcal{Z}\big(\tilde U_q^{\otimes n}\big)$, which proves our claim. It implies that $\Phi_n(\mathcal{Z}(\Ll_{0,n}))\subset \mathcal{Z}\big(\tilde U_q^{\otimes n}\big)$.

Let us prove the converse inclusion. Put $R^{(a)}:=R_{Vn}\cdots R_{Va+1}$. For every $1\leq a\leq n$ we have
\begin{align}
\Phi_n\big(\!\operatorname{qTr}_{V_\lambda}\big(\!\stackrel{V_\lambda}{M}{}^{\!\! (a)} \big)\big) & = (\operatorname{Tr}_{V_{\lambda}}\otimes {\rm id} )\big(\big(\pi_{V_{\lambda}}(\ell)\otimes {\rm id} \big)R^{(a)}R_{Va}R_{Va}'R^{(a)}{}^{\!-1}\big) \nonumber
\\
& = (\operatorname{Tr}_{V_{\lambda}}\otimes {\rm id} )\big(\big(\pi_{V_{\lambda}}(\ell)\otimes {\rm id} \big)R_{Va}R_{Va}'\big)\nonumber
\\
& = \mathfrak{i}_a\big(\Phi_1\big(\!\operatorname{qTr}_{V_{\lambda}} \big(\!\stackrel{V_\lambda}{M}\!\big)\big)\big),\label{alekseevonz}
 \end{align}
where $\Phi_1$ in the last equality is the RSD map, and the second equality follows from $R^{-1} = (S\otimes {\rm id} )(R)$ and $S^2(x) = \ell x \ell^{-1}$ for every $x\in \mUq$. By Theorem \ref{drinfeldmap} and Proposition \ref{linbasecentre2}, the family of these elements, when $\lambda$ spans $P_+$, forms a basis of the center of the $a$-th tensorand of~$\tilde U_q^{\otimes n}$. Therefore $\Phi_n(\mathcal{Z}(\Ll_{0,n}))= \mathcal{Z}\big(\tilde U_q^{\otimes n}\big)$. By injectivity of $\Phi_n$ this proves the first claim. The computation \eqref{alekseevonz} concludes the description of $\mathcal{Z}(\Ll_{0,n})$. The last claim is clear.\end{proof}

In particular, recall from \eqref{notsl2} that for $\mathfrak{g}={\mathfrak{sl}}(2)$ and every $1\leq i \leq n$ we put
\[
\stackrel{V_2}{M}{}^{\!\! (i)} =\begin{pmatrix} a^{(i)}&b^{(i)}\\c^{(i)}&d^{(i)}\end{pmatrix}\in \operatorname{End}(V_2)\otimes \Ll_{0,n}({\mathfrak{sl}}(2))^{(i)}.
\]
Set
\[
\omega^{(i)} = qa^{(i)} + q^{-1}d^{(i)} = \operatorname{qTr}_{V_2}\big(\!\stackrel{V_2}{M}{}^{\!\! (i)} \big).
\]
By (\ref{alekseevonz}) we have \[\Phi_n\big(\omega^{(i)}\big)=\Omega^{(i)},\] that is, $(q-q^{-1})^2$ times the Casimir element of the $i$-th tensorand of $U_q({\mathfrak{sl}}(2))^{\otimes n}$, and then $\mathcal{Z}(U_q({\mathfrak{sl}}(2))^{\otimes n}) = \mc(q)\big[\Omega^{(1)},\dots,\Omega^{(n)}\big]$. The proposition implies
\begin{equation*}
\mathcal{Z}(\Ll_{0,n}({\mathfrak{sl}}(2))) = \mc(q)\big[\omega^{(1)},\dots,\omega^{(n)}\big].
\end{equation*}
More generally, by the same arguments and the results recalled in Remark \ref{centermatrk}, there is an ana\-lo\-gous description of $\mathcal{Z}(\Ll_{0,n}(\mathfrak{g}))$ for an arbitrary finite dimensional complex simple Lie algebra~$\mathfrak{g}$ of type $A$, $B$, $C$ or $D$.

\begin{remk}\label{centersl2A}{\rm When $\mathfrak{g}={\mathfrak{sl}}(2)$, by using the $A$-basis of the $A$-module $\Ll_{0,1}^A({\mathfrak{sl}}(2))=\Oo_A({\mathfrak{sl}}(2))$ given in \cite[Proposition~1.3]{DC-L}, the arguments of Proposition \ref{drinfeldinj} show that $\mathcal{Z}\big(\Ll_{0,1}^A({\mathfrak{sl}}(2))\big)=A[\omega]$. Then it is easy to see as above that the family of elements $\textstyle \prod_{i=1}^n {}_{\lambda_i}\omega^{(i)}$ form an $A$-basis of $\mathcal{Z}\big(\Ll_{0,n}^A({\mathfrak{sl}}(2))\big)$, and that $\mathcal{Z}\big(\Ll_{0,n}^A({\mathfrak{sl}}(2))\big) = A\big[\omega^{(1)},\dots,\omega^{(n)}\big]$.}
\end{remk}

For every $\lambda\in P_+$ denote
\[
{}_\lambda\eta= \operatorname{qTr}_{V_\lambda}\big(\!\stackrel{V_\lambda}{M}{}^{\!\! (1)}\cdots \stackrel{V_\lambda}{M}{}^{\!\! (n)}\big).
\]
In the specific case of $\mathfrak{g}={\mathfrak{sl}}(2)$, let us put
\begin{equation}\label{eta}
\eta:=\operatorname{qTr}_{V_2}\big(\!\stackrel{V_2}{M}{}^{\!\! (1)}\cdots \stackrel{V_2}{M}{}^{\!\! (n)}\big).
\end{equation}

 \begin{lem} \label{center1lem} $\mathcal{Z}\big(\Ll_{0,n}^{U_q}\big)$ contains the commutative algebra generated over $\mc(q)$ by the elements ${}_{\lambda_1}\omega^{(1)},\dots, {}_{\lambda_n}\omega^{(n)}$ and ${}_\lambda\eta$ for all $\lambda,\lambda_1,\dots,\lambda_n\in P_+$.
\end{lem}

\begin{proof} Clearly $\mathcal{Z}\big(\Ll_{0,n}^{U_q}\big)$ contains $\mathcal{Z}(\Ll_{0,n})^{U_q}$, whence the elements ${}_{\lambda_1}\omega^{(1)},\dots, {}_{\lambda_n}\omega^{(n)}$ for every $\lambda_1,\dots,\lambda_n\in P_+$, by Propositions~\ref{center0n} and~\ref{invL0n}. Moreover, by this result $\Ll_{0,n}^{U_q}$ is the centralizer of $\Delta^{(n-1)}(\Ll_{0,1})$ in $\Ll_{0,n}$. Therefore $\mathcal{Z}\big(\Ll_{0,n}^{U_q}\big)\supset \Delta^{(n-1)}(\mathcal{Z}(\Ll_{0,1}))$. By Proposition \ref{linbasecentre2} and the relation~\eqref{prodmatrix}, the elements ${}_\lambda\eta$, $\lambda\in P_+$, form a basis of $\Delta^{(n-1)}(\mathcal{Z}(\Ll_{0,1}))$. The result follows. \end{proof}

Note that the $\mc(q)$-algebra generated by the elements ${}_{\lambda_1}\omega^{(1)},\dots, {}_{\lambda_n}\omega^{(n)}$ and ${}_\lambda\eta$ for all weights $\lambda,\lambda_1,\dots,\lambda_n\in P_+$ is the polynomial algebra generated by these elements when $\lambda,\lambda_1,\dots,\lambda_n$ are fundamental weights. We have even a much better result:
\begin{teo} \label{CenterLinv}
We have an isomorphism
\[
\mathcal{Z}\big(\Ll_{0,n}^{U_q}\big) \cong \Delta^{(n-1)}(\mathcal{Z}(\Ll_{0,1}))\otimes_{\mc(q)} \mathcal{Z}(\Ll_{0,1})^{\otimes n}.
\]
Therefore $\mathcal{Z}\big(\Ll_{0,n}^{U_q}\big)$ is the polynomial algebra generated over $\mc(q)$ by the elements ${}_{\lambda_1}\omega^{(1)}$, $\dots$, ${}_{\lambda_n}\omega^{(n)}$ and ${}_\lambda\eta$ for all fundamental weights $\lambda,\lambda_1,\dots,\lambda_n$. In particular, when $\mathfrak{g}={\mathfrak{sl}}(2)$ we have
\[
\mathcal{Z}\big(\Ll_{0,n}({\mathfrak{sl}}(2))^{U_q({\mathfrak{sl}}(2))}\big)= \mc(q)\big[\omega^{(1)},\dots,\omega^{(n)}, \eta\big]
\]
and
\[
\mathcal{Z}\big(\Ll_{0,n}^A({\mathfrak{sl}}(2))^{U_A({\mathfrak{sl}}(2))}\big)= A\big[\omega^{(1)},\dots,\omega^{(n)}, \eta\big].
\]
\end{teo}
This is a consequence of the following lemmas, which are interesting by themselves. Denote by $U_h := U_{h}(\mathfrak{g})$ the Hopf algebra over ${\mathbb C}[[h]]$ generated topologically (in the $h$-adic sense) by~$E_i$, $F_i$,~$H_i$, where $i=1,\dots,m$, satisfying the relations $[H_i, H_j]=0$, $[H_i, E_j]= a_{ij} E _j$, $[H_i, F_j]=- a_{ij}F_j$ and \eqref{EFK}--\eqref{Serre2}, where $K_i$ is replaced by $q_i^{H_i}$ and $q_i$ by $e^{d_ih}$. By the same formula as~\eqref{adjointactiononn} with $U_{h}$ in place of $U_q$, $U_{h}$ acts on $U_{h}^{\otimes n}$ by the right adjoint action. Let $\mathcal{Z}(U_{h})$ be the center of~$U_{h}$.

\begin{lem} \label{centerh-n}
The $h$-adic completion of $\Delta^{(n-1)}(\mathcal{Z}(U_h))\otimes_{{\mathbb C}[[h]]}\mathcal{Z}(U_h)^{\otimes n}$, considered as a subalgebra of $U_{h}^{\otimes n}$, is the center of $\big(U_{h}^{\otimes n}\big)^{U_h}$.
\end{lem}

\begin{proof}
Denote by $U = U(\mathfrak g)$ the envelopping algebra of ${\mathfrak g}$, and by $\Delta_0$ its canonical coproduct. We have $U_h / h U_h=U$.
As shown in \cite[Lemma 3.10]{RY}, a direct application of Theorem~10.1 of~\cite{Knop} proves the result for $h =0$, i.e., for $U$. Because ${\mathfrak g}$ is finite dimensional and semisimple, Drinfeld's results in \cite{Dr,Dr1} show that there exists an isomorphism of algebras $\phi\colon U_h \rightarrow U[[h]]$, equal to the identity on the quotient $U_h / h U_h$, and there exists an invertible element $J \in U^{\otimes 2}[[h]]$ such that $\forall x\in U_h, (\phi\otimes \phi)(\Delta(x))=J \Delta_0 (\phi(x)) J^{-1}$. By using $J$ one can easily define an invertible element $J_n \in U^{\otimes n}[[h]]$ such that $\forall x\in U_h$, $(\phi^{\otimes n})(\Delta^{(n-1)}(x))=J_n\Delta_0^{(n-1)} (\phi(x)) J_n^{-1}$. As a result the map $\psi\colon \big(U_{h}^{\otimes n}\big)^{U_h}\rightarrow (U^{\otimes n})^{U}[[h]]$, $\psi(x) = J_n^{-1}(\phi^{\otimes n}(x)) J_n$, is an isomorphism of algebras. We obviously have $\phi(\mathcal{Z}(U_h))=\mathcal{Z}(U)[[h]]$. The lemma follows after having checked that
$\psi^{-1}\big(\mathcal{Z}(U)^{\otimes n}[[h]]\big)=\mathcal{Z}(U_h)^{\otimes n}$ and $\psi^{-1}\big(\Delta_0^{(n-1)}(\mathcal{Z}(U)[[h]])\big)=\Delta^{(n-1)}(\mathcal{Z}(U_h))$.
\end{proof}

\begin{lem} 
The center of $(U_{q}^{\otimes n})^{U_q}$ is $\mathcal{Z}\big((U_{q}^{\otimes n})^{U_q}\big)=\Delta^{(n-1)}(\mathcal{Z}(U_q))\otimes_{{\mathbb C}(q)}\mathcal{Z}(U_q)^{\otimes n}$.
\end{lem}

\begin{proof}
By Lemma \ref{center1lem} it only remains to prove the inclusion ``$\subset$''. Let $\varphi\colon U_A\otimes_A {\mathbb C}[[h]] \rightarrow U_h$ be the morphism of algebras defined by $\varphi(K_i)=e^{h d_i H_i}$, $\varphi(E_i)=E_i$, $\varphi(F_i)=F_i$. Let $z\in\mathcal{Z}\big((U_{q}^{\otimes n})^{U_q}\big)$. Up to multiplication by an element of $A$ we can assume that $z\in \mathcal{Z}\big((U_{A}^{\otimes n})^{U_A}\big)$. Because $\varphi^{\otimes n}(z)$ commutes with $\varphi^{\otimes n}\big((U_A^{\otimes n})^{U_A}\otimes_A {\mathbb C}[[h]]\big)$, it commutes with its $h$-adic completion. Consider the family of elements $v_{[\lambda]}(a_{[\lambda]})$ defined by \eqref{defvectinv} with ${}_{A}V_{[\lambda]} := {}_AV_{\lambda_1}\otimes \dots \otimes {}_AV_{\lambda_1}$ and~$a_{[\lambda]}$ in $\operatorname{End}_{U_A^{\rm res}}\big({}_AV_{\lambda_1}\otimes \dots \otimes {}_AV_{\lambda_n}\big)$. By Proposition \ref{linbaseinv} (see the comments therefater) they form a $\mc(q)$-basis of $\Ll_{0,n}^{U_q}$. Proposition \ref{surjPhiinv} then implies that the elements $\Phi_n(v_{[\lambda]}\big(a_{[\lambda]}^{(k)}\big)$ form a basis of~$\big(U_q^{\otimes n}\big)^{U_q}$. These results still hold by working with $U_h$ and the topological version of $\Ll_{0,n}$ defined over $\mc[[h]]$, so the elements $\varphi^{\otimes n}\big(\Phi_n\big(v_{[\lambda]}\big(a_{[\lambda]}^{(k)}\big)\big)\big)$ form a ${\mathbb C}[[h]]$-topological basis of $\big(U_{h}^{\otimes n}\big)^{U_h}$. Therefore $\varphi^{\otimes n}(z)$ centralizes $\big(U_{h}^{\otimes n}\big)^{U_{h}}$, whence $\varphi^{\otimes n}(z)\in \Delta^{(n-1)}(\mathcal{Z}(U_h))\otimes_{{\mathbb C}[[h]]}\mathcal{Z}(U_h)^{\otimes n}$ by Lemma \ref{centerh-n}. This proves $z\in \Delta^{(n-1)}(\mathcal{Z}(U_q))\otimes_{{\mathbb C}(q)}\mathcal{Z}(U_q)^{\otimes n}$. \end{proof}

\begin{proof}[Proof of Theorem \ref{CenterLinv}]
We have an isomorphism ${\mathcal L}_{0,n}^{U_q}\cong \big(\tilde{U}_q^{\otimes n}\big)^{U_q}$ by Proposition \ref{surjPhiinv}, and $\mathcal{Z}\big(\big(\tilde{U}_q^{\otimes n}\big)^{U_q}\big)=\Delta^{(n-1)}\big(\mathcal{Z}(\tilde{U}_q)\big)\otimes_{{\mathbb C}(q)}\mathcal{Z}(\tilde{U}_q)^{\otimes n}$ by the last lemma applied to $\tilde U_q$ instead of~$U_q$. But $\Phi_n^{-1}\big(\mathcal{Z}(\tilde{U}_q)^{\otimes n}\big) = \mathcal{Z}({\mathcal L}_{0,1})^{\otimes n}$ by Proposition \ref{center0n}, which is the vector space generated by ${}_{\lambda_1}\omega^{(1)},\dots, {}_{\lambda_n}\omega^{(n)}$ for $\lambda_1,\dots,\lambda_n\in P_+$ (by Proposition \ref{linbasecentre2}), and $\Phi_n^{-1} \big(\Delta^{(n-1)}\big(\mathcal{Z}(\tilde{U}_q)\big)\big)$ is the vector space generated by ${}_\lambda\eta$ for all $\lambda\in P_+$ (by Proposition \ref{linbasecentre2} and the relation~\eqref{prodmatrix}). This proves all claims but the last one, which in turn follows from $\mathcal{Z}\big(\big(\Ll_{0,n}^A\big)^{U_A}\big) = \mathcal{Z}\big(\Ll_{0,n}^{U_q}\big) \cap \Ll_{0,n}^A$, and the fact that $\omega^{(1)},\dots,\omega^{(n)}$, $\eta\in \Ll_{0,n}^A$.
\end{proof}

\subsection{Specializations}
Let $\epsilon\in {\mathbb C}^\times$. We defined the unrestricted specialization $U_\epsilon := U_A \otimes_A \mc_\e$ of $U_q$ in \eqref{restrictU}. Recall that ${\mathbb C}_{\e}$ is the $A$-module ${\mathbb C}$, where $q$ acts by multiplication by~$\e$. Similarly, the unrestricted specialisation of $\Ll_{0,n}^A$ at $\e$ is the $U_\epsilon$-module algebra
\begin{equation*}
 \Ll_{0,n}^\e =\Ll_{0,n}^A\otimes_A {\mathbb C}_{\e}.
 \end{equation*}
We need to consider the specialization at $q=\e$ of the Alekseev map. By Lemma \ref{intform-n} and~Pro\-position~\ref{Alekseevsl2} we know that $\Phi_n\colon \Ll_{0,n}^A\ra \tilde{U}_A^{\otimes n}$ is an embedding of $U_A$-module algebras, and in the case $\mathfrak{g}={\mathfrak{sl}}(2)$ it maps into $U_A({\mathfrak{sl}}(2))^{{\otimes} n}$. Moreover, by Proposition~\ref{Alekseevsl2loc} the latter embedding extends to an isomorphism $\Phi_n\colon {}_{\rm loc}\Ll_{0,n}^A({\mathfrak{sl}}(2))\ra U_A'({\mathfrak{sl}}(2))^{{\otimes} n}$.

When $q=\e$ the formulas \eqref{phigen} show that $\Phi_1\colon \Ll_{0,1}^\e({\mathfrak{sl}}(2))\ra U_\e({\mathfrak{sl}}(2))$ is an embedding, and $\Phi_1\colon {}_{\rm loc}\Ll_{0,1}^\e({\mathfrak{sl}}(2))\ra U_\e({\mathfrak{sl}}(2))$ an isomorphism. For an arbitrary $\mathfrak{g}$, $\Phi_1\colon \Ll_{0,1}^\e\ra \tilde{U}_\e$ is also an embedding; this follows from the following facts. In \cite[ Sections 4 and 6]{DC-L}, De Concini--Lyubashenko introduced an embedding of algebras $\mu''\colon \Oo_A \ra \tilde U_A(\mathfrak{b}_-)\otimes \tilde U_A(\mathfrak{b}_+)$ (where $\tilde U_A(\mathfrak{b}_\pm)$ is the subalgebra of $\tilde U_A$ associated to the Borel subalgebra $\mathfrak{b}_\pm$ of $\mathfrak{g}$), and they proved that it affords an embedding $\mu''_\e\colon \Oo_\e \ra \tilde U_\e(\mathfrak{b}_-)\otimes \tilde U_\e(\mathfrak{b}_+)$. We have $\Phi_1 = m \circ ({\rm id} \otimes S^{-1})\circ \mu''$. Moreover, $\operatorname{Im}(\mu'')$ is contained in the subalgebra of $\tilde U_\e(\mathfrak{b}_-)\otimes \tilde U_\e(\mathfrak{b}_+)$ with basis elements
\[
\bar{F}_{\beta_1}^{n_1}\cdots \bar{F}_{\beta_N}^{n_N}K_{n_1\beta_1+\dots+ n_N \beta_N}K_{\lambda}\otimes K_{-\lambda}K_{-p_1\beta_1-\dots- p_N \beta_N}\bar{E}_{\beta_1}^{p_1}\cdots \bar{E}_{\beta_N}^{p_N},
\]
where $\lambda\in P$ and $n_1,\dots,n_N, p_1,\dots,p_N \in {\mathbb N}$. The map $m \circ \big({\rm id} \otimes S^{-1}\big)$ sends this basis to a~free family of $\tilde U_\e$. As a result $\Phi_1\colon \Ll_{0,1}^\e\ra \tilde{U}_\e$ is injective. Since $\Phi_n$ differs from $\Phi_1^{\otimes n}$ by a linear isomorphism (see \eqref{definitionofphin2}), it follows that{\samepage
\[
\Phi_n\colon\ \Ll_{0,n}^\e\to \tilde U_\epsilon^{\otimes n}
\]
is an embedding of $U_\epsilon$-module algebras.}

Moreover, in the ${\mathfrak{sl}}(2)$ case, by Proposition \ref{Alekseevsl2loc} and the fact that $\Phi_1\colon {}_{\rm loc}\Ll_{0,1}^\e({\mathfrak{sl}}(2))\ra U_\e({\mathfrak{sl}}(2))$ is an isomorphism, it follows that
\begin{equation}\label{Phinlociso}
\Phi_n\colon\ {}_{\rm loc}\Ll_{0,n}^\e({\mathfrak{sl}}(2))\ra U_\e({\mathfrak{sl}}(2))^{\otimes n}
\end{equation}
is an isomorphism. Note that, by Lemma \ref{intform}, when $n=1$ we have an isomorphism
\begin{equation*}
\Phi_1\colon \ \Ll_{0,1}^\epsilon\to \tilde{U}_A^{\rm lf}\otimes_A {\mathbb C}_{\epsilon}.
\end{equation*}
It is important to note that taking the specialization at a root of unity $q=\e$ and taking locally finite elements are non commuting operations. For instance $\tilde U_\epsilon$ is a free module of finite rank over its center (see, e.g., \cite[Section~9.2]{CP}). Hence it coincides with its subalgebra of locally finite elements $(\tilde U_\epsilon)^{\rm lf}$. On another hand, $\tilde{U}_A^{\rm lf}\otimes_A {\mathbb C}_{\epsilon}$ is strictly contained in $\tilde U_\epsilon$; for instance it does not contain the elements~$K_i$.

The algebra $\tilde U_q^{\rm lf}$ is very complicated~\cite{Jos}, we know neither generators nor basis for arbitrary~$\mathfrak{g}$. This hurdle prevents us to give a precise description of $\tilde U_A^{\rm lf}$, $\tilde{U}_A^{\rm lf}\otimes_A {\mathbb C}_{\epsilon}$, and $\big(\tilde U_A^{\otimes n}\big)^{\rm lf}\otimes_A {\mathbb C}_{\epsilon}$ for $n>1$.

Finally, recall that $\big(\Ll_{0,n}^A\big)^{U_A}$ is the centralizer of $\Delta^{(n-1)}\big(\Ll_{0,1}^A\big)$, see Proposition~\ref{invL0n}. By~the same arguments the algebra $(\Ll_{0,n}^\e)^{U_\e}$ of $U_\e$-invariant elements of $\Ll_{0,n}^\e$ is the centralizer of $\Delta^{(n-1)}(\Ll_{0,1}^\e)$. Therefore we have an inclusion $\mathcal{Z}(\Ll_{0,n}^\epsilon) \subset (\Ll_{0,n}^\e)^{U_\e}$, and since clearly $\big(\Ll_{0,n}^A\big)^{U_A}\allowbreak\otimes_A \mc_\e$ is a subset of $(\Ll_{0,n}^\e)^{U_\e}$, multiplication defines a morphism of algebras
\[
\big(\big(\Ll_{0,n}^A\big)^{U_A}\otimes_A \mc_\e\big) \otimes \mathcal{Z}(\Ll_{0,n}^\epsilon) \ra (\Ll_{0,n}^\e)^{U_\e}.
\]
In the rest of this paper we will simplify notations by setting
\begin{equation}\label{AepsL0n}
\big(\Ll_{0,n}^A\big)^{U_A}_\e := \big(\Ll_{0,n}^A\big)^{U_A}\otimes_A \mc_\e.
\end{equation}
The arguments of Proposition \ref{nozeroq} apply as well to the specialization $q=\epsilon$ (using that $\tilde U_\e$ has no non trivial zero divisors), so we have:

\begin{prop}
The algebra $\Ll_{0,n}^\epsilon$ does not have non trivial zero divisors, and therefore the subalgebras $(\Ll_{0,n}^\e)^{U_\e}$ and $\big(\Ll_{0,n}^A\big)^{U_A}_\e$ too.
\end{prop}

\section[Center of L\_\{0,n\}\textasciicircum{}{e}(sl(2)) and quantum coadjoint action]
{Center of $\boldsymbol{\Ll_{0,n}^\e({\mathfrak{sl}}(2))}$ and quantum coadjoint action}
\label{specialization}

\subsection[Center of L\_\{0,n\}\textasciicircum{}{e}(sl(2))]
{Center of $\boldsymbol{\Ll_{0,n}^\e({\mathfrak{sl}}(2))}$}\label{centersl2}

From now on $\mathfrak{g} = {\mathfrak{sl}}(2)$ and $\e$ is a primitive $l$-th root of unity. We assume that $l\geq 3$ and $l$ is odd. We make this latter choice to simplify the exposition; the case of $l$ even can be treated in a similar way, as all our constructions below rely on the description of the center of $U_\e({\mathfrak{sl}}(2))$, which is done for all primitive roots of unity $\e$ in \cite[Section~3]{DC-K}. We omit $\mathfrak{g}$ from the notations of the various algebras, and denote $U_\epsilon({\mathfrak{sl}}(2))$ by $U_\epsilon$, etc.

By the relations \eqref{phigen} and \eqref{omega}, $\Phi_1(\Ll_{0,1}^\e) = \tilde{U}_A^{\rm lf}\otimes_A {\mathbb C}_{\epsilon} = U_A^{\rm lf}\otimes_A {\mathbb C}_{\epsilon}$ is the algebra generated over $\mc$ by $\Omega$, $EK^{-1}$, $F$ and $K^{-1}$. Adding the generator $K$ gives $U_\epsilon$. By results of \cite{DC-K}, the center $\mathcal{Z}(U_\epsilon)$ of $U_\epsilon$ is the $\mc$-algebra generated by $E^l$, $F^l$, $K^{\pm l}$ and $\Omega$ satisfying the relation
 \begin{equation}\label{relcenter}
 \prod_{j=1}^l (\Omega- c_j) = \big(\epsilon-\epsilon^{-1}\big)^{2l}E^lF^l+K^l+K^{-l}- 2,
 \end{equation}
 where $ c_j = \epsilon^j+\epsilon^{-j}$. Let $\mathcal{Z}_0(U_\epsilon)$ be the subalgebra of $\mathcal{Z}(U_\epsilon)$ generated by $E^l$, $F^l$ and $K^{\pm l}$. It is a sub-Hopf algebra of $U_\epsilon$, with
\begin{gather*}
\Delta\big(K^{\pm l}\big)=K^{\pm l}\otimes K^{\pm l},\qquad
\Delta\big(E^l\big)=E^l\otimes K^l+1\otimes E^l,\qquad
\Delta\big(F^l\big)=K^{-l}\otimes F^l + F^l\otimes 1,
\\
S\big(E^l\big) = -E^lK^{-l},\qquad
S\big(F^l\big) = -K^lF^l,\qquad
S\big(K^{\pm l}\big) = K^{\mp l},
\\
\varepsilon\big(E^l\big) = \varepsilon\big(F^l\big)=0,\qquad
\varepsilon\big(K^l\big)=1.
\end{gather*}
Consider the sequence of polynomials $T_k$, $k\in \mathbb{N}$, defined recursively by
\begin{equation}\label{ChPol}
T_0(x)=2,\qquad
T_1(x) =x,\qquad
T_k(x) = xT_{k-1}(x)-T_{k-2}(x)\qquad
{\rm for}\quad k\geq 2.
\end{equation}
Note that $T_k(x)/2$ is the $k$-th Chebyshev polynomial of the first type in the variable $x/2$. One has $T_k\big(u+u^{-1}\big) = u^k+u^{-k}$, from which one derives easily that $\textstyle T_l(x) -2 = \prod_{j=1}^l (x-c_j)$.
Therefore, the relation \eqref{relcenter} can be written as:
\begin{equation}\label{relcenter3}
T_l(\Omega) = \big(\epsilon-\epsilon^{-1}\big)^{2l}E^lF^l+K^l+K^{-l}.
 \end{equation}
Since $l$ is odd, it is also equivalent to
\begin{equation}\label{relcenter2}
 \prod_{j=1}^l (\Omega+ c_j) = \big(\epsilon-\epsilon^{-1}\big)^{2l}E^lF^l+K^l+K^{-l}+2.
 \end{equation}
By \eqref{phigen} we have
 \begin{gather}
 \Phi_1(\omega)=\Omega, \qquad
 \Phi_1\big(b^l\big) = \big(\epsilon-\epsilon^{-1}\big)^l F^{l}, \nonumber\\
 \Phi_1\big(c^{l}\big) = \big(\epsilon-\epsilon^{-1}\big)^l (EK^{-1})^l, \qquad
 \Phi_1\big(d^{l}\big) = K^{-l}.\label{genroot}
 \end{gather}
Hence
\begin{equation}\label{Tcheb0}
T_l(\omega) = b^lc^ld^{-l} + d^{-l}+d^l.
\end{equation}
Using that $\Phi_1$ is equivariant, injective, and surjective when extended to the localization, we deduce
\begin{gather*}
\mathcal{Z}({}_{\rm loc}\Ll_{0,1}^\epsilon) = \mc\big[\omega,b^l,c^l,d^{\pm l}\big]/\mathcal{I},
\\
\mathcal{Z}(\Ll_{0,1}^\epsilon) = \mc\big[\omega,b^l,c^l,d^{l}\big]/\mathcal{I},
\end{gather*}
where $\mathcal{I}$ is the ideal of $\mathcal{Z}(\Ll_{0,1}^\epsilon)$ generated by $\big(T_l(\omega) - d^l\big)d^l - b^lc^l - 1$. By the presentation of $\Ll_{0,1}^A$ in Lemma \ref{intsl2}, we have
\[
\mathcal{Z}(\Ll_{0,1}^\epsilon) = (\Ll_{0,1}^\epsilon)^{U_\epsilon}.
\]
Alternatively, we have $\big(U_A^{\rm lf}\otimes_A {\mathbb C}_{\epsilon}\big)^{U_\epsilon} = \mathcal{Z}\big(U_A^{\rm lf}\otimes_A {\mathbb C}_{\epsilon}\big)$ by the arguments of Corollary \ref{invL0n} for $n=1$. Therefore $\mathcal{Z}(\Ll_{0,1}^\epsilon) = \Phi_1^{-1}\big(\big(U_A^{\rm lf}\otimes_A {\mathbb C}_{\epsilon}\big)^{U_\epsilon}\big) = (\Ll_{0,1}^\epsilon)^{U_\epsilon}$.

We can now define a notion of quantum Frobenius homomorphism for $\Ll_{0,1}^\epsilon,$ similar in spirit to the one defined for $\mathcal{O}_\epsilon$ in \cite{PW}, which is a map $\mathcal{O}_1\ra \mathcal{Z}(\mathcal{O}_\epsilon)$. Consider the specialization $\Ll_{0,1}^1$ of $\Ll_{0,1}^A$ at $q=1$. We have $\Ll_{0,1}^1=\mathcal{O}_1=\mathcal{O}(G)$ as commutative algebras (with $G={\rm SL}(2,\mc))$). Denote by $\underline{a}$, $\underline{b}$, $\underline{c}$, $\underline{d}$ the images of the generators $a$, $b$, $c$, $d$ of $\Ll_{0,1}^A$ under the specialization map $\Ll_{0,1}^A\ra \Ll_{0,1}^1$. They satisfy $\underline{a}\underline{d}-\underline{b}\underline{c}=1.$
Let us define $Q_l\in {\mathbb C}[X,Y]$ by $Q_l(X,Y)=T_l\big(\epsilon X+\epsilon^{-1}Y\big)-X^l-Y^l$.
Recall that $a$ and $d$ commute. So we put:
\begin{defi} \label{defFr1}
The Frobenius map $\operatorname{Fr}\colon \Ll_{0,1}^1\rightarrow \mathcal{Z}(\Ll_{0,1}^\epsilon)$ is the homomorphism of algebras given by
\[
\operatorname{Fr}(\underline{a})=a^l+Q_l(a,d)=T_l(\omega)-d^l, \qquad
\operatorname{Fr}(\underline{b})=b^l, \qquad
\operatorname{Fr}(\underline{c})=c^l,\qquad
\operatorname{Fr}(\underline{d})=d^l.
\]
\end{defi}
We shall denote
\begin{equation*}
\underline{\stackrel{V_2}{M}}=\begin{pmatrix} \underline{a} & \underline{b}\\\underline{c}&\underline{d}\end{pmatrix}\qquad \text{and}\qquad
\operatorname{Fr}\underline{\stackrel{V_2}{M}}=\begin{pmatrix} a^l+Q_l(a,d)& b^l\\c^l&d^l
\end{pmatrix}\!.
\end{equation*}
Note that:
\begin{itemize}\itemsep=0pt
\item
$\det\big(\!\operatorname{Fr}\underline{\stackrel{V_2}{M}}\big)-1 = (T_l(\omega) - d^l)d^l - b^lc^l - 1$, i.e., the generator of the ideal $\mathcal{I}$.

\item
$T_l\big(\!\operatorname{qTr}\big(\!\stackrel{V_2}{M}\big)\big)
=\operatorname{Tr}\big(\!\operatorname{Fr}\underline{\stackrel{V_2}{M}}\big)$.
\end{itemize}

The notions above can be developed similarly for every $\Ll_{0,n}^\e$, $n\geq 1$. First, recall the additional generators $\delta^{(i)-1}$ of the localization ${}_{\rm loc}\Ll_{0,n}^\e$ introduced in Definition \ref{Lonloc}.
\begin{prop} \label{center0nroot}
We have
\begin{gather*}
\mathcal{Z}(\Ll_{0,n}^\e) =\mc\big[\omega^{(i)},b^{(i)l},c^{(i)l},d^{(i)l},\, i=1,\dots, n\big]/\big(\mathcal{I}^{(i)},\, i=1,\dots, n\big),
\\
\mathcal{Z}({}_{\rm loc}\Ll_{0,n}^\e) = \mc\big[\omega^{(i)},b^{(i)l},c^{(i)l},d^{(i) l}, \big(\delta^{(i)l}\big)^{-1},\, i=1,\dots, n\big]/\big(\mathcal{I}^{(i)},\, i=1,\dots, n\big),
\end{gather*}
 where $\mathcal{I}^{(i)}$ is the ideal generated by the element $\big(T_l\big(\omega^{(i)}\big) - d^{(i)l}\big)d^{(i)l} - b^{(i)l}c^{(i)l} - 1$.
\end{prop}

\begin{proof} First we prove that $\Phi_n\big(b^{(i)l}\big), \Phi_n\big(c^{(i)l}\big), \Phi_n\big(d^{(i)l}\big)\in \mathcal{Z}(U_\epsilon^{\otimes n})$ for every $i=1,\dots,n$. By~injectivity of $\Phi_n$ it will follow that $b^{(i)l}, c^{(i)l}, d^{(i)l}\in \mathcal{Z}(\Ll_{0,n}^\e)$. By~\eqref{genroot} the claim is true for $i=n$. Let $1\leq i\leq n-1$, and denote $R^{(i)}:=R_{V_2n}\cdots R_{V_2i+1} \in \operatorname{End}(V_2)\otimes \tilde U_q^{\otimes n}$, with $q$ an indeterminate, as in the previous sections. Define $r_{11}^{(i)}, r_{12}^{(i)}, r_{21}^{(i)}, r_{22}^{(i)}\in \tilde U_q^{\otimes n}$ by
\[
R^{(i)} =\begin{pmatrix} r_{11}^{(i)}&r_{12}^{(i)}\\[1ex] r_{21}^{(i)}&r_{22}^{(i)}\end{pmatrix}
\]
and $m_{11}^{(i)}, m_{12}^{(i)}, m_{21}^{(i)}, m_{22}^{(i)} \in 1^{\otimes (i-1)} \otimes U_q^{\rm lf} \otimes 1^{\otimes (n-i)}$ by
\[
{\rm id}_{V_2}\otimes \big(1^{\otimes (i-1)} \otimes \Phi_1 \otimes 1^{\otimes (n-i)}\big)\begin{pmatrix} a^{(i)}&b^{(i)}\\ c^{(i)}&d^{(i)}\end{pmatrix} = \begin{pmatrix} m_{11}^{(i)} &m_{12}^{(i)}\\m_{21}^{(i)}&m_{22}^{(i)}\end{pmatrix}\!,
\]
where as in \eqref{notsl2} we put
\[
\stackrel{V_2}{M}{}^{\!\! (i)} =\begin{pmatrix} a^{(i)}&b^{(i)}\\c^{(i)}&d^{(i)}\end{pmatrix}\in \operatorname{End}(V_2)\otimes \Ll_{0,n}^{(i)}.
\]
By Proposition \ref{imagephi}, $m_{11}^{(i)}$, $m_{12}^{(i)}$, $m_{21}^{(i)}$, $m_{22}^{(i)}$ generate the subalgebra $1^{\otimes (i-1)} \otimes U_q^{\rm lf} \otimes 1^{\otimes (n-i)}$ of~$U_q^{\otimes n}$. These elements satisfy the relations \eqref{rel01}, with $m_{11}^{(i)}$, $m_{12}^{(i)}$, $m_{21}^{(i)}$, $m_{22}^{(i)}$ replacing $a$, $b$, $c$, $d$ respectively, and they commute with $r_{11}^{(i)}$, $r_{12}^{(i)}$, $r_{21}^{(i)}$, $r_{22}^{(i)}$, since the non trivial tensor components of the latters do not lie on the $i$-th tensorand of $U_q^{\otimes n}$. Denote $\mathfrak{i}_{in}\colon \tilde U_q^{\otimes (n-i)} \ra \tilde U_q^{\otimes n}$ the identification map with the last $n-i$ tensorands. We have $R^{(i)} = \big(\pi_{V_2}\otimes \mathfrak{i}_{in}\circ \Delta^{(n-i-1)}\big)(R)$ by the relation~\eqref{Rdef0}. So the formulas \eqref{1formRV2} and~\eqref{2formRV2} yield
\[
r_{21}^{(i)}=0,\qquad
r_{22}^{(i)} = r_{11}^{(i)-1},\qquad
r_{12}^{(i)}r_{11}^{(i)}= q r_{11}^{(i)}r_{12}^{(i)},\qquad
r_{12}^{(i)} = r_{11}^{(i)}f_{12}^{(i)}
\]
and
\[
R^{(i)-1} =\begin{pmatrix} r_{11}^{(i)-1}&-qr_{12}^{(i)}\\[.5ex] 0 &r_{11}^{(i)}\end{pmatrix}\!,
\]
where
\begin{equation*}
r_{11}^{(i)} = \mathfrak{i}_{in}\circ \Delta^{(n-i-1)}\big(q^{H/2}\big),\qquad f_{12}^{(i)} = \mathfrak{i}_{in}\circ \Delta^{(n-i-1)}\big(\big(q-q^{-1}\big)F\big).
\end{equation*}
Therefore
\begin{gather*}
R^{(i)}\big({\rm id}_{V_2}\otimes \big(1^{\otimes (i-1)} \otimes \Phi_1 \otimes 1^{\otimes (n-i)}\big)\big)\begin{pmatrix} a^{(i)}&b^{(i)}\\ c^{(i)}&d^{(i)}\end{pmatrix}R^{(i)-1}
\\ \qquad
{} = \begin{pmatrix} m_{11}^{(i)}+r_{12}^{(i)}r_{11}^{(i)-1}m_{21}^{(i)} & * \\[.5ex] r_{11}^{(i)-2}m_{21}^{(i)} & -qr_{11}^{(i)-1}r_{12}^{(i)}m_{21}^{(i)} + m_{22}^{(i)} \end{pmatrix}\!.
\end{gather*}
Note also that, by definition,
\[
\big({\rm id}_{V_2} \otimes \Phi_n\big)\begin{pmatrix} a^{(i)}&b^{(i)}\\ c^{(i)}&d^{(i)}\end{pmatrix} = R^{(i)}\big( {\rm id}_{V_2}\otimes \big(1^{\otimes (i-1)} \otimes \Phi_1 \otimes 1^{\otimes (n-i)}\big)\big)\begin{pmatrix} a^{(i)}&b^{(i)}\\ c^{(i)}&d^{(i)}\end{pmatrix}R^{(i)-1}.
\]
Hence
\begin{gather*}
\Phi_n\big(c^{(i)}\big) = r_{11}^{(i)-2}m_{21}^{(i)},
\\
\Phi_n\big(d^{(i)}\big) = -qr_{11}^{(i)-1}r_{12}^{(i)}m_{21}^{(i)} + m_{22}^{(i)} .
\end{gather*}
Let now take the specialization $q=\epsilon$ as above, a primitive $l$-th root of unity where $l \geq 3$ is odd. We have
\begin{gather}
\Phi_n\big(c^{(i)l}\big) = \big(r_{11}^{(i)-2}m_{21}^{(i)}\big)^l = r_{11}^{(i)-2l}m_{21}^{(i)l},\label{Phincl}
\\
\Phi_n\big(d^{(i)l}\big) = \big({-}\epsilon r_{11}^{(i)-1}r_{12}^{(i)}m_{21}^{(i)} + m_{22}^{(i)}\big)^l = -r_{11}^{(i)-l}r_{12}^{(i)l}m_{21}^{(i)l} + m_{22}^{(i)l},\label{Phindl}
\end{gather}
by using $m_{21}^{(i)}m_{22}^{(i)} = \epsilon^2 m_{22}^{(i)}m_{21}^{(i)}$, the $q$-binomial formula (see, e.g., \cite[Proposition IV.2.2]{Kas}), and the vanishing at $q=\epsilon$ of the $q$-Gauss binomial coefficients $[l]_q/[k]_q[l-k]_q$, $0< k< l$. Now
\[
m_{21}^{(i)l}= \big(\epsilon-\epsilon^{-1}\big)^l1^{\otimes (i-1)} \otimes \big(K^{-1}E\big) ^l \otimes 1^{\otimes (n-i)},\qquad
r_{11}^{(i)-2l}= \mathfrak{i}_{in}\circ \Delta^{(n-i-1)}\big(K^{-l}\big).
\]
These are central elements of $U_\epsilon^{\otimes n}$, so $\Phi_n\big(c^{(i)l}\big)$ is central. As
\[
r_{11}^{(i)-l}r_{12}^{(i)l}=f_{12}^{(i)l }= \mathfrak{i}_{in}\circ \Delta^{(n-i-1)}\big(\big(\epsilon-\epsilon^{-1}\big)^l F^l\big)
\]
and $\mathcal{Z}_0(U_\epsilon)$ is a Hopf algebra, $r_{11}^{(i)-l}
r_{12}^{(i)l}$ is a central element of
$U_\epsilon^{\otimes n}$. Again, $m_{21}^{(i)l}$ and $m_{22}^{(i)l}$ being central in $U_\epsilon^{\otimes n}$, $\Phi_n\big(d^{(i)l}\big)$ is central.

Finally, recalling that $\Ll_{0,n}^{(i)}$ is isomorphic to $\Ll_{0,1}$, by specializing $q$ to $\epsilon$ we get
\begin{equation}\label{relcentreomega}
d^{(i)l} T_l\big(\omega^{(i)}\big) - d^{(i)2l} - 1 = b^{(i)l}c^{(i)l}.
\end{equation}
We know that $\omega^{(i)}\in \mathcal{Z}(\Ll_{0,n}^\epsilon)$ and we just proved that $d^{(i)l}, c^{(i)l}\in \mathcal{Z}(\Ll_{0,n}^\epsilon)$. Therefore, for every $x\in \Ll_{0,n}$ we get $b^{(i)l}c^{(i)l} x = x b^{(i)l}c^{(i)l}=b^{(i)l} x c^{(i)l}$, i.e., $\big(b^{(i)l} x - x b^{(i)l}\big) c^{(i)l} = 0, $ $c^{(i)l}$ is not a zero divisor, we deduce $b^{(i)l} x - x b^{(i)l}=0$. Hence $b^{(i)l}$ is central in $\Ll_{0,n}^\epsilon$.

 A formula of $\Phi_n\big(b^{(i)l}\big)$ can be obtained as follows. Recall that \eqref{alekseevonz} implies
 \begin{equation}\label{omegaOmega}
\Phi_n\big(\omega^{(i)}\big) = \Omega^{(i)}.
\end{equation}
{\samepage
This can also be checked by using the above formulas:
\begin{gather*}
\Phi_n\big(\omega^{(i)}\big) = q\Phi_n\big(a^{(i)}\big)+q^{-1}\Phi_n\big(d^{(i)}\big)
= q\big(m_{11}^{(i)} + q^{-1}f_{12}^{(i)}m_{21}^{(i)}\big) + q^{-1}\big(m_{22}^{(i)} - qf_{12}^{(i)}m_{21}^{(i)}\big)
\\ \hphantom{\Phi_n\big(\omega^{(i)}\big)}
{} = qm_{11}^{(i)}+q^{-1}m_{22}^{(i)}.
\end{gather*}} \noindent
Then, by applying $\Phi_n$ to the relation \eqref{relcentreomega} and using \eqref{Phincl}, \eqref{Phindl} and \eqref{omegaOmega}, one finds
\begin{equation}\label{phinbl}
\Phi_n\big(b^{(i)l}\big) = -r_{11}^{(i)l}r_{12}^{(i)l}\big(T_l\big(\Omega^{(i)}\big) - 2 m_{22}^{(i)l}\big)+ r_{11}^{(i)2l}m_{12}^{(i)l}-r_{12}^{(i)2l}m_{21}^{(i)l}.
\end{equation}
We can now achieve the proof. Note that $\mc\big[\omega^{(1)},\dots,\omega^{(n)}\big] = \mathcal{Z}\big(\Ll_{0,n}^A\big) \otimes_A {\mathbb C}_{\e}$ by Remark \ref{centersl2A}. The inclusion $\mathcal{Z}\big(\Ll_{0,n}^A\big) \otimes_A {\mathbb C}_{\e}\subset \mathcal{Z}(\Ll_{0,n}^\e)$ is clear, and the natural embedding $\Ll_{0,n}^\e \rightarrow {}_{\rm loc}\Ll_{0,n}^\e$ maps $\mathcal{Z}(\Ll_{0,n}^\e)$ into $\mathcal{Z}({}_{\rm loc}\Ll_{0,n}^\e)$. By the case $n=1$ the elements $\omega^{(i)}$, $b^{(i)l}$, $c^{(i)l}$, $d^{(i) l}$ generate $\mathcal{Z}\big(\Ll_{0,n}^\e{}^{(i)}\big)$, with ideal of relations $\mathcal{I}^{(i)}$. The set of these ideals for $i=1,\dots,n$ provide all the relations in~$\mathcal{Z}(\Ll_{0,n}^\e)$, for there are no others in $\operatorname{Im}(\Phi_n)$ (as shows, e.g., Corollary~\ref{dressingcenter} below). Therefore one has inclusions
\[
\mc\big[\omega^{(i)},b^{(i)l},c^{(i)l},d^{(i)l},\, i=1,\dots, n\big]/\big(\mathcal{I}^{(i)},\, i=1,\dots, n\big) \subset \mathcal{Z}\big(\Ll_{0,n}^\e\big) \subset \mathcal{Z}\big({}_{\rm loc}\Ll_{0,n}^\e\big).
\]
The conclusion follows at once, since by their very definition the elements $\delta^{(i)\pm l}$ are central in~${}_{\rm loc}\Ll_{0,n}^\e{}$, and $\Phi_n$ maps the algebra generated by them and the left-hand side isomorphically to~$U_\epsilon^{\otimes n}$.\end{proof}

Analogously to the case $n=1$, recalling that $\Ll_{0,1}^1 = \mathcal{O}(G)$ we put:
\begin{defi} \label{defFrn} The Frobenius map ${\rm Fr}\colon \big(\Ll_{0,1}^1\big)^{\otimes n}\rightarrow \mathcal{Z}\big(\Ll_{0,n}^\epsilon\big)$ is the homomorphism of algebras given by
\[
\operatorname{Fr}\big(\underline{a}^{(i)}\big)=a^{(i)l}+Q_l\big(a^{(i)},d^{(i)}\big), \qquad \operatorname{Fr}\big(\underline{b}^{(i)}\big)=b^{(i)l},\qquad
\operatorname{Fr}\big(\underline{c}^{(i)}\big)=c^{(i)l}, \qquad \operatorname{Fr}\big(\underline{d}^{(i)}\big)=d^{(i)l}.
\]
\end{defi}
We shall denote
\begin{equation}\label{notFrmat}
\underline{\stackrel{V_2}{M}}{}^{(i)}=\begin{pmatrix} \underline{a}^{(i)} & \underline{b}^{(i)}\\ \underline{c}^{(i)}&\underline{d}^{(i)}\end{pmatrix}
\qquad {\rm and}\qquad
\operatorname{Fr}\underline{\stackrel{V_2}{M}}{}^{(i)}=\begin{pmatrix} a^{(i)l}+Q_l\big(a^{(i)},d^{(i)}\big)& b^{(i)l}\\[.5ex] c^{(i)l}&d^{(i)l}
\end{pmatrix}\!,
\end{equation}
where we recall that $a^{(i)l}+Q_l\big(a^{(i)},d^{(i)}\big) = T_l\big(\omega^{(i)}\big)- d^{(i)l}$. We can express $\Phi_n \circ \big(\Phi_1^{\otimes n}\big)^{-1}$ on the center as follows. Set
\begin{gather*}
\mathcal{R}^{(i)} =\begin{pmatrix} r_{11}^{(i)l}& r_{12}^{(i)l}\\[.5ex] 0 &r_{11}^{(i)-l}\end{pmatrix}\!,
\\
\mathcal{M}^{(i)} = \big(1^{\otimes (i-1)} \otimes \Phi_1 \otimes 1^{\otimes (n-i)}\big)\big(\!\operatorname{Fr}\underline{\stackrel{V_2}{M}}{}^{(i)}\big)
 = \begin{pmatrix} T_l(\Omega^{(i)})- m_{22}^{(i)l} & m_{12}^{(i)l}\\[.5ex] m_{21}^{(i)l} & m_{22}^{(i)l} \end{pmatrix}\!,
\end{gather*}
where we use the notations in the proof of Proposition \ref{center0nroot}. These matrices belong to $\operatorname{End}(V_2) \otimes \mathcal{Z}_0(\tilde U_\epsilon)^{\otimes n}$ and $\operatorname{End}(V_2) \otimes \mathcal{Z}_0(U_\epsilon)^{\otimes n}$ respectively. Here we note that $\mathcal{Z}(\tilde U_\epsilon)$ is generated by $\Omega$, $E^l$, $F^l$ and $K^{\pm \frac{l}{2}}$ satisfying the relation \eqref{relcenter}, and we define $\mathcal{Z}_0(\tilde U_\epsilon)$ as the subalgebra generated by~$E^l$,~$F^l$ and~$K^{\pm \frac{l}{2}}$.

By using \eqref{Phincl}, \eqref{Phindl}, \eqref{omegaOmega} and \eqref{phinbl} it is easy to check that:
\begin{cor} \label{dressingcenter} The map $\Phi_n \circ \big(\Phi_1^{\otimes n}\big)^{-1} \colon \mathcal{Z}_0(U_\epsilon)^{\otimes n} \ra \mathcal{Z}_0(U_\epsilon)^{\otimes n}$ is given by
 \[
 \big({\rm id}_{V_2} \otimes \big(\Phi_n \circ (\Phi_1^{\otimes n})^{-1}\big)\big)\big(\mathcal{M}^{(i)}\big) = \mathcal{R}^{(i)}\mathcal{M}^{(i)}\mathcal{R}^{(i)-1}.
 \]
\end{cor}
We will find useful later to have explicit formulas. Let us introduce the following generators of $\mathcal{Z}_0(U_\epsilon)$:
\begin{equation}\label{pref}
x= -\big(\e-\e^{-1}\big)^{l}E^lK^{-l}, \qquad
y=\big(\e-\e^{-1}\big)^{l}F^l,\qquad
z^{\pm 1}= K^{\pm l}.
\end{equation}
Similarly, denote by $\mathcal{Z}_0(\tilde U_\epsilon)\subset \tilde U_\epsilon$ the subalgebra generated by $x$, $y$, $z$ and
\[
z'{}^{\pm 1}= K^{\pm\frac{ l}{2}}.
\]
For every $a\in \big\{x,y,z^{\pm 1}\big\}$ put $a^{(i)} = 1^{\otimes (i-1)} \otimes a \otimes 1^{\otimes (n-i)}$. We can view $\mathcal{Z}_0(U_\epsilon)^{\otimes n}$ as a polynomial algebra in the variables $x^{(i)}$, $y^{(i)}$, $\big(z^{\pm 1}\big)^{(i)}$, and $\mathcal{Z}_0(\tilde U_\epsilon)^{\otimes n}$ as a polynomial algebra in~$x^{(i)}$, $y^{(i)}$, $\big(z'{}^{\pm 1}\big)^{(i)}$.
Then
\begin{equation}\label{Mcal}
\mathcal{M}^{(i)} = \begin{pmatrix} T_l\big(\Omega^{(i)}\big)- (z^{-1})^{(i)} & y^{(i)}\\[.5ex] -x^{(i)} & (z^{-1})^{(i)} \end{pmatrix} = \begin{pmatrix} z^{(i)}\big(1-x^{(i)}y^{(i)}\big) & y^{(i)}\\[.5ex] -x^{(i)} & (z^{-1})^{(i)} \end{pmatrix}\!.
 \end{equation}
Also,
\begin{gather}
r_{11}^{(i)l} = \big(z'{}^{-1}\big)^{(i+1)}\cdots \big(z'{}^{-1}\big)^{(n)}, \label{g11l}
\\
r_{12}^{(i)l} = \big(z'{}^{-1}\big)^{(i+1)}\cdots \big(z'{}^{-1}\big)^{(n)} \bigg(y^{(i+1)} + \sum_{j=1}^{n-i-1} \big(z^{-1}\big)^{(i+1)}\cdots \big(z^{-1}\big)^{(i+j)}y^{(i+j+1)}\bigg)\label{g12l}
\end{gather}
yields
\begin{gather*}
\Phi_n\big(c^{(i)l}\big) = r_{11}^{(i)-2l}m_{21}^{(i)l} =-x^{(i)} \big(z^{-1}\big)^{(i+1)}\cdots \big(z^{-1}\big)^{(n)},
\label{phinclbis}
\\
\Phi_n\big(d^{(i)l}\big) = m_{22}^{(i)l} - m_{21}^{(i)l} r_{11}^{(i)-l}r_{12}^{(i)l} \notag
\\ \hphantom{\Phi_n\big(d^{(i)l}\big)}
{} = \big(z^{-1}\big)^{(i)} + x^{(i)}\bigg(y^{(i+1)} + \sum_{j=1}^{n-i-1} \big(z^{-1}\big)^{(i+1)}\cdots \big(z^{-1}\big)^{(i+j)} y^{(i+j+1)}\bigg)\label{phindlbis} .
\end{gather*}
One can readily express $\Phi_n\big(b^{(i)l}\big) = -r_{11}^{(i)l}r_{12}^{(i)l}\big((T_l(\Omega^{(i)}) - 2 m_{22}^{(i)l}\big)+ r_{11}^{(i)2l}m_{12}^{(i)l}-r_{12}^{(i)2l}m_{21}^{(i)l}$ as well as a polynomial in the variables $x^{(j)}$, $y^{(j)}$, $\big(z^{\pm 1}\big)^{(j)}$.

In \cite{PW}, Parshall--Wang showed that the quantum Frobenius homomorphism $\mathcal{O}_1 \ra \mathcal{Z}(\mathcal{O}_\e)$ they defined is a morphism of coalgebra. As we now explain, our quantum Frobenius homomorphism ${\rm Fr}\colon \Ll_{0,1}^1\rightarrow \mathcal{Z}(\Ll_{0,1}^\epsilon)$ satisfies a similar property.

Recall the algebra morphism $\Delta^{(n-1)}\colon \Ll_{0,1}^A\rightarrow \Ll_{0,n}^A$ (see Proposition \ref{coprodprop}). Denote again by $\Delta^{(n-1)}\colon \Ll_{0,1}^\epsilon\rightarrow \Ll_{0,n}^\epsilon$ its evaluation at $\epsilon$. For every $1\leq i_1< \dots < i_k\leq n$ we define $\iota_{i_1,\dots,i_k}\colon \Ll_{0,k}^\e \ra \Ll_{0,n}^\e$ to be the identification map of the $j$-th tensorand of $\Ll_{0,k}^\e$ with the $i_j$-th tensorand of $\Ll_{0,n}^\e$, $1\leq j\leq k$.

\begin{lem} For every $1\leq i_1< \dots < i_k\leq n$ we have
\begin{equation}\label{DeltaFr0}
\iota_{i_1,\dots,i_k}\circ \Delta^{(k-1)}\big(\!\operatorname{Fr}\underline{\stackrel{V_2}{M}}\big) =\operatorname{Fr}\underline{\stackrel{V_2}{M}}{}^{(i_1)}\cdots \operatorname{Fr}\underline{\stackrel{V_2}{M}}{}^{(i_k)}.
\end{equation}
\end{lem}

\begin{proof} First we show \eqref{DeltaFr0} in the case $k=n=2$, that is
\begin{equation}\label{DeltaFr}
\Delta\big(\!\operatorname{Fr}\underline{\stackrel{V_2}{M}}\big)
=\operatorname{Fr}\underline{\stackrel{V_2}{M}}{}^{(1)} \operatorname{Fr}\underline{\stackrel{V_2}{M}}{}^{(2)}.
\end{equation}
The commutation relations of $\stackrel{V_2}{M}{}^{\!\!\!(1)}$ and $\stackrel{V_2}{M}{}^{\!\!\!(2)}$ being complicated, we cannot compute the matrix components directly (which is the way used in \cite{PW} for their quantum Frobenius homomorphism $\mathcal{O}_1 \ra \mathcal{Z}(\mathcal{O}_\e)$, where such computations reduce to the $q$-binomial identity). Instead, we first use the Alekseev map. Indeed, because $\Phi_2$ is an algebra embedding, it is sufficient to show that~\eqref{DeltaFr} holds after having been composed with $\Phi_2$. We have
\[
\Phi_1\big(\!\operatorname{Fr}\underline{\stackrel{V_2}{M}}\big)
= \begin{pmatrix} T_l(\Omega)- K^{-l} & \big(\epsilon-\epsilon^{-1}\big)^l F^l\\[.5ex]
 \big(\epsilon-\epsilon^{-1}\big)^l (EK^{-1})^l& K^{-l}\end{pmatrix}\!.
 \]
Therefore \eqref{DeltaFr} is a consequence of the following four equations in $U_\epsilon^{\otimes 2}$:
\begin{gather*}
\Delta\big(K^{-l}\big)=\Phi_2\big(c^{(1)l}b^{(2)l}+d^{(1)l}d^{(2)l}\big),
\\
\Delta\big(\big(\epsilon-\epsilon^{-1}\big)^l \big(EK^{-1}\big)^l\big)=\Phi_2\big(c^{(1)l}\big(a^{(2)l}+Q_l\big(a^{(2)},d^{(2)} \big)\big)+d^{(1)l}c^{(2)l}\big),
\\
\Delta\big(\big(\epsilon-\epsilon^{-1}\big)^l F^l\big)=\Phi_2\big(\big(a^{(1)l}+Q_l\big(a^{(1)},d^{(1)} \big)\big) b^{(2)l}+b^{(1)l}d^{(2)l}\big),
\\
\Delta(T_l(\Omega)- K^{-l})=\Phi_2\big(\big(a^{(1)l}+Q_l\big(a^{(1)},d^{(1)} \big)\big) \big(a^{(2)l}+Q_l\big(a^{(2)},d^{(2)}\big)\big)+b^{(1)l}c^{(2)l}\big).
\end{gather*}
We listed them in order of complexity. The first equation, using the explicit expression of $\Phi_2$ on components, is rewritten as
\begin{gather*}
\Delta\big(K^{-l}\big)=\big(\epsilon-\epsilon^{-1}\big)^{2l}\big(K^{-1}E\otimes K^{-1} \big)^l \big(1\otimes F^l\big)
\\ \hphantom{\Delta\big(K^{-l}\big)=}
{}+ \big(K^{-1}\otimes 1- \epsilon\big(\epsilon-\epsilon^{-1}\big)^2 K^{-1}E\otimes F\big)^l \big(1\otimes K^{-l}\big).
\end{gather*}
This relation holds thanks to the $q$-binomial identity. The other three relations, although more complicated, are shown similarly by a direct computation and using the expressions \eqref{Phincl}, \eqref{Phindl}, \eqref{omegaOmega} and \eqref{phinbl}.

Because the relative commutations relations of ${\stackrel{V_2}{M}}{}^{(1)},\dots, {\stackrel{V_2}{M}}{}^{(k)}$ and ${\stackrel{V_2}{M}}{}^{(i_1)},\dots ,{\stackrel{V_2}{M}}{}^{(i_k)}$ for general sequences $i_1< \cdots<i_k$ are the same, the proof of \eqref{DeltaFr0} follows immediately from the case where $i_j=j$, for $1\leq j\leq k$. This, in turn, follows from \eqref{DeltaFr} by induction on $k$. This concludes the proof.
\end{proof}

The following consequence of the Lemma will be a key tool in Section \ref{appskein}.
\begin{prop}\label{Thcentral} For every $1\leq i_1< \dots < i_k\leq n$ we have
\[
T_l\big(\!\operatorname{qTr}\big(\!\stackrel{V_2}{M}{}^{\!\! (i_1)}\cdots \stackrel{V_2}{M}{}^{\!\! (i_k)}\big)\big)=\operatorname{Tr}\big(\!\operatorname{Fr}\underline{\stackrel{V_2}{M}}{}^{(i_1)}\cdots \operatorname{Fr}\underline{\stackrel{V_2}{M}}{}^{(i_k)}\big).
\]
In particular, this element is central in $\Ll_{0,n}^\e$.
\end{prop}

\begin{proof}
Let $\iota_{i_1,\dots,i_k}\colon \Ll_{0,k}^\e \ra \Ll_{0,n}^\e$ be the identification map of the $j$-th tensorand of $\Ll_{0,k}^\e$ with the $i_j$-th tensorand of $\Ll_{0,n}^\e$, $1\leq j\leq k$. We have
\begin{align*}
T_l\big(\!\operatorname{qTr}\big(\!\stackrel{V_2}{M}{}^{\!\! (i_1)}\cdots \stackrel{V_2}{M}{}^{\!\! (i_k)}\big)\big) &= T_l\big(\!\operatorname{qTr}\big(\big({\rm id}\otimes \big(\iota_{i_1,\dots,i_k}\circ \Delta^{(k-1)}\big)\big)\big(\!\stackrel{V_2}{M}\big)\big)\big)
\\
& = T_l\big(\big(\iota_{i_1,\dots,i_k}\circ \Delta^{(k-1)}\big)\big(\!\operatorname{qTr}\big(\!\stackrel{V_2}{M}\big)\big)\big)
 \\
& = \big(\iota_{i_1,\dots,i_k}\circ \Delta^{(k-1)}\big)\big(T_l\big(\!\operatorname{qTr}\big(\!\stackrel{V_2}{M}\big)\big)\big)
 \\
& = \big(\iota_{i_1,\dots,i_k}\circ \Delta^{(k-1)}\big)\big(T_l(\omega)\big).
\end{align*}
In the first equality we used the formula \eqref{prodmatrix}, the second and fourth equalities follow from definitions, and the third comes from the fact that $\Delta^{(k-1)}\colon \Ll_{0,1}^A\rightarrow \Ll_{0,k}^A$ is a homomorphism of algebras. By the identity \eqref{Tcheb0} and the fact that $\mathcal{Z}_0(U_\epsilon) \cong \mc\big[b^l,c^l,d^{\pm l}\big]$ is a Hopf subalgebra of $U_\e$, we have $\Delta^{(k-1)}(T_l(\omega)) \in \mathcal{Z}({}_{\rm loc}\Ll_{0,1}^\e)^{\otimes k}$. This and Proposition \ref{center0nroot} imply that the above element is central in $\Ll_{0,n}^\e$. Moreover
\vspace{-1ex}
\begin{align*}
\iota_{i_1,\dots,i_k}\circ \Delta^{(k-1)}(T_l(\omega)) & = \iota_{i_1,\dots,i_k}\circ \Delta^{(k-1)}\big(\!\operatorname{Tr}\big(\!\operatorname{Fr}\underline{\stackrel{V_2}{M}}\big)\big)
\\
& = \operatorname{Tr}\big(\big({\rm id}\otimes \iota_{i_1,\dots,i_k}\circ \Delta^{(k-1)}\big)\big(\!\operatorname{Fr}\underline{\stackrel{V_2}{M}}\big)\big)
\\
& = \operatorname{Tr}\big(\!\operatorname{Fr}\underline{\stackrel{V_2}{M}}{}^{(i_1)}\cdots \operatorname{Fr}\underline{\stackrel{V_2}{M}}{}^{(i_k)}\big),
\end{align*}
where the first equality follows from the observations we made before Proposition \ref{center0nroot}, and the others from the previous lemma.
\end{proof}

\subsection[Gg-invariant central elements and SL(2,C)-characters]
{$\boldsymbol{\Gg}$-invariant central elements and $\boldsymbol{{\rm SL}(2,\mc)}$-characters}\label{QCAL0nsec}
We are going to relate $\mathcal{Z}(\Ll_{0,n}^\e)$ with the algebra of regular functions on the variety of ${\rm SL}(2,\mc)$-characters of the sphere with $n+1$ punctures, endowed with the Atiyah--Bott--Goldman Poisson structure. This is achieved in Section \ref{EXTENSION}. To this aim we recall a few preliminary results in the next section.

\subsubsection[The quantum coadjoint action for U\_e]{The quantum coadjoint action for $\boldsymbol{U_\e}$}
We refer to \cite{DC-K,DC-K-P1,DC-K-P2} for details about the material discussed in this section. It can be formulated for any of the quantum groups $U_\e(\mathfrak{g})$, but we restrict to $U_\e({\mathfrak{sl}}(2))$ as we shall need Proposition~\ref{center0nroot} in Section~\ref{EXTENSION}.

Consider the sets $\operatorname{Spec}(\mathcal{Z}(U_\epsilon))$ and $\operatorname{Spec}(\mathcal{Z}_0(U_\epsilon))$ of algebra homomorphisms from $\mathcal{Z}(U_\epsilon)$ and~$\mathcal{Z}_0(U_\epsilon)$ to $\mc$, respectively. They are affine algebraic sets. An element of $\operatorname{Spec}(\mathcal{Z}(U_\epsilon))$ is called a~\textit{central character} of $\Ue$. The inclusion $\mathcal{Z}_0(U_\epsilon)\subset \mathcal{Z}(U_\epsilon)$ induces a regular (restriction) map\vspace{-1ex}
\begin{equation}\label{tau}
\tau\colon\ \operatorname{Spec}(\mathcal{Z}(U_\epsilon)) \lra \operatorname{Spec}(\mathcal{Z}_0(U_\epsilon)).
\end{equation}
Since $\mathcal{Z}_0(U_\epsilon)$ is a polynomial algebra, any $\chi\in \operatorname{Spec}(\mathcal{Z}_0(U_\epsilon))$ is entirely defined by its values $(x_\chi,y_\chi,z_\chi)\in \mc^2\times\mc^*$ on the tuple $(x,y,z)$ of generators of $\mathcal{Z}_0(U_\epsilon)$ defined in \eqref{pref}. By \eqref{relcenter}--\eqref{relcenter2}, any $\chi\in \operatorname{Spec}(\mathcal{Z}(U_\epsilon))$ is entirely defined by its values $(x_\chi,y_\chi,z_\chi,\Omega_\chi)\in \mc^2\times\mc^*\times\mc$ on the tuple $(x,y,z,\Omega)$ of generators of $\mathcal{Z}(U_\epsilon)$, solutions to one of the equivalent equations\vspace{-1ex}
\[
\prod_{j=1}^l (\Omega_\chi\mp c_j) = -x_\chi y_\chi z_\chi+z_\chi+z_\chi^{-1}\mp 2,
\]
where $c_j = \e^{j}+\e^{-j}$. Hence $\tau$ has degree $l$, and $\operatorname{Spec}(\mathcal{Z}(U_\epsilon))$ is a branched covering of $\operatorname{Spec}(\mathcal{Z}_0(U_\epsilon)) = \mc^2\times \mc^*$ of degree $l$, a hypersurface in $\mc^2\times\mc^*\times\mc$ with quadratic singularities at the points $(0,0,\pm 1,\pm c_j)$, $j=1,\dots,(l-1)/2$.

Because $\mathcal{Z}_0(U_\epsilon)$ is a commutative Hopf algebra, $\operatorname{Spec}(\mathcal{Z}_0(U_\epsilon))$ has a canonical group structure defined dually by
\begin{equation}
\chi_1\chi_2(u) := (\chi_1\otimes \chi_2)\Delta(u),\qquad
\chi_1^{-1}(u) := \chi_1(S(u)),\qquad
e(u) := \varepsilon(u)
 \end{equation}
for any $u\in \mathcal{Z}_0(U_\epsilon)$ and $\chi_1, \chi_2\in \operatorname{Spec}(\mathcal{Z}_0(U_\epsilon))$, where $e\in \operatorname{Spec}(\mathcal{Z}_0(U_\epsilon))$ is the identity element. In formulas:
\begin{gather*}
x_{\chi_1\chi_2} = x_{\chi_1}+z_{\chi_1}^{-1}x_{\chi_2},\qquad
y_{\chi_1\chi_2} = y_{\chi_1}+y_{\chi_2}z_{\chi_1}^{-1},\qquad
z_{\chi_1\chi_2} = z_{\chi_1}z_{\chi_2},
\\
x_{\chi^{-1}} =-z_{\chi}x_{\chi},\qquad
y_{\chi^{-1}} = -y_{\chi}z_{\chi},\qquad
z_{\chi^{-1}} = z_{\chi}^{-1},\\
x_e = 0,\qquad y_e = 0,\qquad z_e = 1.
\end{gather*}
This can be formulated as follows. Put $G = {\rm SL}(2,\mc)$, and let 
$G^*$ be the group formed by the pairs of matrices
\[
\left(\begin{pmatrix} a & b \\ 0 & a^{-1}\end{pmatrix}\!, \begin{pmatrix} a^{-1} & 0 \\ c & a \end{pmatrix}\right) \in {\rm SL}(2,\mc)^{\rm op}\times {\rm SL}(2,\mc)^{\rm op},
\]
where ${\rm SL}(2,\mc)^{\rm op}$ is ${\rm SL}(2,\mc)$ endowed with the opposite multiplication. Set
\[
\psi\left( \begin{pmatrix} z' & z'y \\ 0 & z'{}^{-1} \end{pmatrix}\!, \begin{pmatrix} z'{}^{-1} & 0 \\ z'x & z' \end{pmatrix} \right) = \big(x,y,z'{}^2\big).
\]
Identifying $\operatorname{Spec}(\mathcal{Z}_0(U_\epsilon)))$ with $\mc^2\times \mc^*$ by mapping $\chi$ to $(x_\chi,y_\chi,z_\chi)$ defined as above, it is readily checked that this defines a surjective morphism of algebraic groups
\[
\psi\colon\ G^* \ra \operatorname{Spec}(\mathcal{Z}_0(U_\epsilon)))
\]
with kernel the subgroup generated by $-(I, I)$, where $I$ is the $2$-by-$2$ identity matrix. Put
\[
\bar G^* = G^*/\{\pm(I,I)\}.
\]
We will denote the quotient isomorphism by
\[
\bar \psi\colon\ \bar G^* \ra \operatorname{Spec}(\mathcal{Z}_0(U_\epsilon))).
\]

Let us endow the ring of regular functions $\mathcal{O}(G)$ with the Sklyanin--Drinfeld Poisson bracket $\{\, ,\, \}$, associated to the classical $r$-matrix
\[
\mathfrak{r} = \frac{1}{4} H\otimes H + E \otimes F\in \mathfrak{g} \otimes \mathfrak{g}.
\]
Recall that it can be given the following expression (see, e.g., \cite{STS}, or \cite{Au,CP} for a setup close to ours). First note that it is entirely determined by its values on the matrix coefficients (coordinate functions) $l_{11}$, $l_{12}$, $l_{21}$, $l_{22}$ of the fundamental representation of $G$ on $\mc^2$. Put
\begin{equation}\label{Ldef}
L = \begin{pmatrix} l_{11} & l_{12} \\ l_{21} & l_{22} \end{pmatrix} = \sum_{r,s=1}^2 E_r^s\otimes l_{rs}\in \operatorname{End}\big(\mc^2\big) \otimes \mathcal{O}(G).
\end{equation}
Denote by $\big\{\!\stackrel{1}{L},\stackrel{2}{L}\!\big\}$ and $\stackrel{1}{L}\stackrel{2}{L}$ the tensors in $\operatorname{End}\big(\mc^2\big) \otimes \operatorname{End}\big(\mc^2\big) \otimes \mathcal{O}(G)$ defined by
\[
\big\{\!\stackrel{1}{L},\stackrel{2}{L}\!\big\} = \sum_{r,s,t,u=1}^2 E_r^s\otimes E_t^u \otimes \{l_{rs},l_{tu}\},\qquad \stackrel{1}{L}\stackrel{2}{L} = \sum_{r,s,t,u=1}^2 E_r^s\otimes E_t^u \otimes l_{rs}l_{tu}
\]
and as usual put $\stackrel{1}{L} = \sum_{r,s=1}^2 E_r^s\otimes 1 \otimes l_{rs}$, $\stackrel{2}{L} = \sum_{t,u=1}^2 1 \otimes E_t^u \otimes l_{tu}$. The $r$-matrix in the fundamental representation is
\[
\mathfrak{r} = \frac{1}{4} \begin{pmatrix} 1 & 0 \\ 0 & -1 \end{pmatrix} \otimes \begin{pmatrix} 1 & 0 \\ 0 & -1 \end{pmatrix} + \begin{pmatrix} 0 & 1 \\ 0 & 0 \end{pmatrix} \otimes \begin{pmatrix} 0 & 0 \\ 1 & 0 \end{pmatrix}\in \operatorname{End}\big(\mc^2\big) \otimes \operatorname{End}\big(\mc^2\big),
\]
which we identify with $\mathfrak{r}\otimes 1 \in \operatorname{End}\big(\mc^2\big) \otimes \operatorname{End}\big(\mc^2\big) \otimes \mathcal{O}(G)$. Then the bracket $\{\, ,\, \}$ on $\mathcal{O}(G)$ is defined by the identity
\begin{equation*}
\big\{\!\stackrel{1}{L},\stackrel{2}{L}\!\big\} = \big[\mathfrak{r},\stackrel{1}{L}\stackrel{2}{L}\!\big] = \mathfrak{r}\stackrel{1}{L}\stackrel{2}{L} - \stackrel{1}{L}\stackrel{2}{L}\mathfrak{r}.
\end{equation*}
The group $G$ with the bracket $\{\, ,\, \}$ is a Poisson--Lie group. As such, it admits a dual Poisson--Lie group, which is precisely the group $G^*$ defined as above. Its Poisson bracket $\{\,,\, \}_{G^*}$ is defi\-ned~by
\begin{gather}\label{bracketdual}
\big\{\!\stackrel{1}{L}_\pm,\stackrel{2}{L}_\pm \big\}_{G^*} = \big[\mathfrak{r},\stackrel{1}{L}_\pm\stackrel{2}{L}_\pm\!\big],\qquad
\big\{\!\stackrel{1}{L}_+,\stackrel{2}{L}_- \big\}_{G^*} = \big[\mathfrak{r},\stackrel{1}{L}_+\stackrel{2}{L}_-\!\big],
\end{gather}
where $(L_+,L_-)$ is the pair of matrices of coordinate functions on $G^*\subset {\rm SL}(2,\mc)^{\rm op}\times {\rm SL}(2,\mc)^{\rm op}$. We will denote by $\{\,,\,\}_{\bar G^*}$ the induced bracket on $\bar G^*$. The push-forward $\psi_*\{\, ,\, \}_{G^*}$ is a Poisson bracket on $\operatorname{Spec}(\mathcal{Z}_0(U_\epsilon))$ that can be defined directly via the quantum coadjoint action, that we now recall.

The specialization morphism $ev_\epsilon\colon U_A \otimes \mc\big[q,q^{-1}\big]\ra U_\e$, $q\mapsto \epsilon$, is surjective, with kernel the subalgebra $\big(q^l-q^{-l}\big)U_A\otimes \mc\big[q,q^{-1}\big]$. Given an element $x\in U_\e$, let us denote by $\tilde x \in ev_\epsilon^{-1}(x)$ any preimage of $x$. For every $a\in \mathcal{Z}(U_\epsilon)$, $u \in U_\e$ we have $[\tilde{a},\tilde u] = \tilde{a}\tilde u - \tilde{u}\tilde a\in \big(q^l-q^{-l}\big)U_A\otimes \mc\big[q,q^{-1}\big]$. So, let us put
\begin{equation}\label{defder}
D_a(u) = - \lim_{q\rightarrow \e}\frac{\big[\tilde{a},\tilde u\big]}{l\big(q^l-q^{-l}\big)},\qquad
a\in \mathcal{Z}_0(U_\epsilon).
\end{equation}
It is easy to check that the maps $D_a\colon \Ue \ra \Ue$ are well-defined (they do not depend on the choices of $\tilde a$ and $\tilde u$), and that they are derivations of $U_\e$ preserving $\mathcal{Z}_0(U_\epsilon)$ and $\mathcal{Z}(U_\epsilon)$. Hence they define algebraic vector fields on $\operatorname{Spec}(\mathcal{Z}_0(U_\epsilon))$. Since $\Omega$ is central in $U_q({\mathfrak{sl}}(2))$, $D_a$ is the zero map when $a\in \mc[\Omega]$.

Let us introduce the elements $e,f$ defined in \cite{DC-K} by
\begin{equation}\label{efdef}
e=-xz=\big(\e-\e^{-1}\big)^{l}E^l,\qquad
f=-yz = -\big(\e-\e^{-1}\big)^{l}F^lK^l.
\end{equation}
Direct computations from the definition and the fact that $D_e \Omega =0$ imply the formulas:
\begin{gather}\label{derform2}
D_z(K)= 0, \qquad D_z(E) = -\frac{1}{l}zE,\qquad D_z(F) = \frac{1}{l}zF,
\\
\label{derform3}
D_e(K)=\frac{1}{l}eK,\qquad D_e(E) =0, D_e(F) =- \frac{1}{l}\big(\e-\e^{-1}\big)^{l-1}[K;1]E^{l-1},
\\
\label{derform4}
D_y(K)=-\frac{1}{l}yK,\qquad D_y(F) =0, D_y(E) = \frac{1}{l}\big(\e-\e^{-1}\big)^{l-1}[K;-1]F^{l-1},
\end{gather}
where $[K;r] = \frac{K\e^r-K^{-1}\e^{-r}}{\e-\e^{-1}}$, $r\in {\mathbb Z}$.

\smallskip
A Poisson bracket $\{\, ,\, \}_{\rm QCA}$ is defined on $\mathcal{Z}(U_\epsilon)$ by
\[
\{a,b\}_{\rm QCA}= D_a(b).
\]
With this Poisson bracket $\mathcal{Z}(U_\epsilon)$ is a Poisson algebra, $\mathcal{Z}_0(U_\epsilon)$ being a Poisson ideal.
The Poisson structure on $\mathcal{Z}_0(U_\epsilon)$ is completely defined by the formulas
\begin{equation}\label{bracketQCA}
\{y,x\} _{\rm QCA}= -1+xy+z^{-2},\qquad \{z,x\}_{\rm QCA} = -zx,\qquad \{z,y\}_{\rm QCA} = yz
\end{equation}
which can be straithforwardy checked using the definition. From them it is easy to show that the Poisson center of $\mathcal{Z}_0(U_\epsilon)$ is the commutative algebra generated by $-xyz+z+z^{-1}.$

Let us introduce the family of automorphisms of $U_q({\mathfrak{sl}}(2))$ defined as follows. Let $r\in \mathbb Z$, and~$\tau_r$,~$T_r$ the automorphisms of $U_q({\mathfrak{sl}}(2))$ defined by
\begin{gather*}
\begin{split}
&\tau_r(K)=K,\qquad \tau_r(E)=K^r E,\qquad \tau_r(F)=FK^{-r},\\
&T_r(K)=K^{-1},\qquad T_r(E)=-FK^{-r},\qquad T_r(F)=-K^r E.
\end{split}
\end{gather*}
Note that $T_r=\tau_r\circ T_0$ and $T_1$ is the braid group automorphism. We have
\[
T_1(x)=y,\qquad T_1(y)=z^2 x.
\]
These automorphisms of $U_q({\mathfrak{sl}}(2))$ define automorphisms of the specialisation $U_\epsilon$ and we will keep the same notation for them. If $T$ is any automorphism of the type $T_r$ or $\tau_r$, it is easy to show from the definition of $D_a$ for $a \in {\mathcal Z}(U_\epsilon)$ that
\[
D_{T(a)}\circ T=T\circ D_a.
\]

Denote by $G^0$ the {\it big cell} of $G={\rm SL}(2,\mc)$. It consists of the matrices with non vanishing lower right entry, and satisfies $G^0 = B_+B_-$, where $B_+$ and $B_-$ are the subgroups of $G$ of upper and lower triangular matrices. We have an unramified $2$-fold covering
\[
\fonc{\sigma}{G^*}{G^0,}{(b_+,b_-)}{b_+ b_-^{-1}.}
\]
It induces a diffeomorphism
\[
\bar \sigma \colon\ \bar G^* \ra G^0.
\]
Setting $z:=z'{}^2$ and
\begin{equation}\label{formmatG}
\mathcal{M} =\sigma \left(\begin{pmatrix} z' & z'y \\ 0 & z'{}^{-1} \end{pmatrix}\!, \begin{pmatrix} z'{}^{-1} & 0 \\ z'x & z' \end{pmatrix} \right)
\end{equation}
we have
\[
\mathcal{M} =\begin{pmatrix} z-z x y&y \\-x&z^{-1}\end{pmatrix}\!.
\]
Consider the Poisson bracket $\{\, ,\, \}_{\rm FR}$ on $\mathcal{O}(G)$ defined by
\[
\big\{\!\stackrel{1}{L},\stackrel{2}{L}\!\big\}_{\rm FR} = \mathfrak{r}\stackrel{1}{L}\stackrel{2}{L} -\stackrel{1}{L}\stackrel{2}{L}\mathfrak{r}'+\stackrel{2}{L}\mathfrak{r}'\stackrel{1}{L} -\stackrel{1}{L}\mathfrak{r}\stackrel{2}{L},
\]
where $\mathfrak{r}'$ is $\mathfrak{r}$ post-composed with the flip map $a\otimes b \mapsto b\otimes a$. The bracket $\{\, ,\, \}_{\rm FR}$ has been introduced in \cite{STS} and generalized in the work of Fock--Rosly \cite{FR} (which explains our notation).
Note that $\mathcal{O}\big(G^0\big)$ is the localization of the algebra $\mathcal{O}(G)$ with respect to the matrix coefficient~$l_{22}$ in~\eqref{Ldef}. As a consequence, the Poisson bracket $\{ ,\}_{\rm FR}$ being quadratic, it can be extended to~$\mathcal{O}(G^0)$.

The next result sums up the relationships between the brackets $\{\, ,\, \}_{\rm FR}$, $\{\, ,\, \}_{G^*}$, and $\{\, ,\, \}_{\rm QCA}$:

\begin{teo}\label{relbracket}
\quad

\begin{enumerate}\itemsep=0pt
\item[$1.$] The map $\bar \psi\colon (\bar G^*,\{\, ,\, \}_{\bar G^*})\ra (\operatorname{Spec}((\mathcal{Z}_0(U_\epsilon)),\{\, ,\, \}_{\rm QCA})$ is an isomorphism of complex Pois\-son--Lie groups {\rm \cite{DC-K-P1,DC-K-P2}}.

\item[$2.$] The map $\bar \sigma\colon (\bar G^*,\{\, ,\, \}_{\bar G^*})\ra (G^0,\{\, ,\, \}_{\rm FR})$ is a diffeomorphism of complex Poisson mani\-folds~{\rm \cite{STS}}.
 \end{enumerate}
\end{teo}
The first claim is proved in \cite{DC-K-P1,DC-K-P2} for the simply-connected quantum groups $\tilde U_\epsilon(\mathfrak{g})$, where $\mathfrak{g}$ is a complex finite dimensional simple Lie algebra (see Remark \ref{convdif} below). We state it rather in the case of $U_\epsilon({\mathfrak{sl}}(2))$, the adjoint quantum group for $\mathfrak{g} = {\mathfrak{sl}}(2)$, where it follows from straighforward computations using the formulas \eqref{bracketdual} and \eqref{bracketQCA}. The second statement follows as well from straighforward computations. Namely, by using the identities \eqref{bracketdual} and $\sigma^*(L) = L_+L_-^{-1}$, and the fact that $\mathfrak{r}+\mathfrak{r}'$ is ${\rm ad}$-invariant, one can identify the formulas of $\{\, ,\, \}_{\rm FR}$ and
\[
\sigma_*(\{\, ,\, \}_{G^*})(L\otimes L) = \big\{\sigma^*\big(\!\stackrel{1}{L}\!\big),\sigma^*\big(\!\stackrel{2}{L}\!\big)\big\}_{G^*}.
\]
We leave the verifications to the reader.

The next statement summarizes the results of De Concini--Kac--Procesi on the quantum coadjoint action. Recall the elements $e= -xz, f=-yz$ of $\mathcal{Z}_0(U_\epsilon)$ (see~\eqref{efdef}). Let us identify $\operatorname{Spec}(\mathcal{Z}_0(U_\epsilon))$ with $G^0$ using the map $\bar{\sigma} \circ \bar{\psi}^{-1}$ of Theorem \ref{relbracket}, and hence the derivations $D_e$, $D_f$, $D_z$ of $\mathcal{Z}_0(U_\epsilon)$ with algebraic vector fields on $G^0$ (and hence on $G$). Denote by $\underline{H}$, $\underline{X}$, $\underline{Y}$ the left-invariant vector fields on $G$ associated to the generators $H$, $X$, $Y$ of ${\mathfrak{sl}}(2)$, where
\[
H = \begin{pmatrix} 1 & 0\\ 0& -1 \end{pmatrix}\!,\qquad
X = \begin{pmatrix} 0 & 1\\ 0& 0 \end{pmatrix}\!,\qquad
Y=\begin{pmatrix} 0 & 0\\ 1&0 \end{pmatrix}\!.
\]
Denote by $\hat{\mathcal{Z}}_0(U_\epsilon)$ the vector space of power series in the generators $x$, $y$, and $z^{\pm 1}$ of $\mathcal{Z}_0(U_\epsilon)$ whose sums converge when evaluated on any $\chi\in \operatorname{Spec}(\mathcal{Z}_0(U_\e))\cong \mc^2\times \mc^*$, $\chi =(x_\chi,y_\chi,z_\chi)$ (thus defining holomorphic functions). Set
\begin{equation*}
\hat{U}_\e = \Ue \otimes_{\mathcal{Z}_0(U_\epsilon)} \hat{\mathcal{Z}}_0(U_\epsilon),\qquad \hat{\mathcal{Z}}(\Ue) = \mathcal{Z}(U_\e) \otimes_{\mathcal{Z}_0(U_\epsilon)} \hat{\mathcal{Z}}_0(U_\epsilon).
\end{equation*}

\begin{teo} [{\cite{DC-K-P1,DC-K-P2}}]\qquad\label{DCKPteo}

\begin{enumerate}\itemsep=0pt
\item[$(a)$] We have $D_{z} = z \underline{H}/2$, $D_e= -z\underline{Y}$, $D_f=z\underline{X}$.

\item[$(b)$] For every $t\in \mc$ the power series $\exp(tD_e)$, $\exp(tD_f)$ converge to automorphisms of $\hat{U}_\e$ preserving $\hat{\mathcal{Z}}_0(U_\epsilon)$ and $\hat{\mathcal{Z}}(U_\epsilon)$, and fixing $\mc[\Omega]$.
\end{enumerate}
 Denote by $\mathcal{G}_{\rm DCK}$ the subgroup of ${\rm Aut}(\hat{U}_\e)$ generated by the $1$-parameter groups $(\exp(tD_e))_{t\in \mathbb{C}}$ and $(\exp(tD_f))_{t\in \mathbb{C}}$. It acts dually by holomorphic transformations on $\operatorname{Spec}(\mathcal{Z}_0(U_\epsilon))$ and $\operatorname{Spec}(\mathcal{Z}(U_\epsilon))$ by defining $g.\chi$ on $u\in \hat{\mathcal{Z}}_0(U_\epsilon)$, for every $g\in \Gg$ and $\chi\in \operatorname{Spec}(\mathcal{Z}(U_\epsilon))$, by $u(g.\chi) = (g^{-1}.u)(\chi)$. Then:

\begin{enumerate}\itemsep=0pt

\item[$(c)$] The diffeomorphism $\bar \sigma \circ \bar \psi^{-1}\colon \operatorname{Spec}(\mathcal{Z}_0(U_\epsilon))\ra G^0$ maps the action of $\Gg_{\rm DCK}$ on the tangent spaces of $(0,0,\pm 1)\in \operatorname{Spec}(\mathcal{Z}_0(U_\epsilon))$ to the coadjoint action of $G$ on ${\mathfrak{sl}}(2,\mc)^*$, the tangent spaces of ${\pm I} \in G^0$.

\item[$(d)$] For any conjugacy class $\Gamma$ in ${\rm SL}(2,\mc)$, $\big(\bar \psi\circ \bar \sigma^{-1}\big)(\Gamma \cap G^0)$ is a $($non empty$)$ $\Gg_{\rm DCK}$-orbit in $\operatorname{Spec}(\mathcal{Z}_0(U_\epsilon))$, and $\tau^{-1}$ of this orbit has $l$ connected components, all of whose are $\Gg_{\rm DCK}$-orbits in $\operatorname{Spec}(\mathcal{Z}(U_\epsilon))$.

\item[$(e)$] An element $a\in {\mathcal{Z}}(\Ue)$ is $\mathcal{G}_{\rm DCK}$-invariant if and only if $a\in \mc[\Omega]$. Dually, the sets of fixed points of the action of $\mathcal{G}_{\rm DCK}$ on $\operatorname{Spec}(\mathcal{Z}_0(U_\epsilon))$ and $\operatorname{Spec}(\mathcal{Z}(U_\epsilon))$ are respectively $(0,0,\pm 1) := \big(\bar \psi \circ \bar \sigma^{-1}\big)(\{{\rm \pm I}\})$ and
\[
\mathcal{D} := \tau^{-1}((0,0,\pm 1)) = \{(0,0,\pm 1,\pm c_j), \, j=1,\dots,(l-1)/2\}.
\]

\item[$(f)$] The $\Gg_{\rm DCK}$-orbits on $\operatorname{Spec}(\mathcal{Z}_0(U_\epsilon))$ and $\operatorname{Spec}(\mathcal{Z}(U_\epsilon))$ are the symplectic leaves of $\{\, ,\, \}_{\rm QCA}$.

\end{enumerate}
\end{teo}

\begin{remk} \label{convdif}
In the statements of Theorems \ref{relbracket} and~\ref{DCKPteo} we use our conventions, which differ from those in \cite{DC-K-P1,DC-K-P2} in the following ways:

\begin{enumerate}\itemsep=0pt

\item[$(i)$] Theorem \ref{relbracket}(1) for the simply connected quantum group $\tilde U_\epsilon({\mathfrak{sl}}(2))$ instead of $U_\epsilon\big({\mathfrak{sl}}(2)\big)$, as in \cite{DC-K-P1,DC-K-P2}, states an isomorphism of $(G^*,\{\, ,\, \}_{G^*})$ with $\big(\!\operatorname{Spec} (\mathcal{Z}_0(\tilde U_\epsilon)),\{\, ,\, \}_{\rm QCA}\big)$; as we take opposite comultiplications on $\tilde U_\epsilon({\mathfrak{sl}}(2))$, we get opposite multiplications of $G^*$ \big(whence~$\bar G^*$ for $U_\epsilon({\mathfrak{sl}}(2))\big)$.

\item[$(ii)$] Our derivations $D_a$ differ by a sign, which we introduce in order to get the equality of $\bar \sigma_*\{\, ,\, \}_{\bar G^*}$ with $\{\, ,\, \}_{\rm FR}$ in Theorem~\ref{relbracket}$(2)$.

\item[$(iii)$] In \cite{DC-K,DC-K-P1,DC-K-P2}, the matrix ${\mathcal M}$ in \eqref{formmatG} is different because they use the opposite coproduct on $\tilde U_\epsilon({\mathfrak{sl}}(2))$; their braid group automorphism $T$ is in our notation $T_{-1}$, which satisfies $T_{-1}(e)=f$.
\end{enumerate}
\end{remk}
It is easily checked that the identities \eqref{bracketQCA} imply
\begin{equation}\label{DonM}
D_z\mathcal{M} = \bigg[\frac{z}{2}\underline{H},\mathcal{M}\bigg], \qquad
D_x\mathcal{M} = \big[{-}z^{-1}\underline{X},\mathcal{M}\big],\qquad
D_y\mathcal{M} = \big[z^{-1}\underline{Y},\mathcal{M}\big],
\end{equation}
where $\mathcal{M}$ is defined by \eqref{formmatG}. Hence we get
\begin{equation*}
D_z = \frac{z}{2}\underline{H}, \qquad
D_x = -z^{-1}\underline{X},\qquad
D_y = z^{-1}\underline{Y}.
\end{equation*}
This is the content of Theorem \ref{DCKPteo}$(a)$ expressed in terms of the derivations $D_z$, $D_x$, $D_y$. In~Theo\-rem \ref{DCKPteo}$(b)$, $U_\epsilon$ is enlarged to $\hat{U}_\e$ in order to define the automorphisms $\exp(tD_e)$, $\exp(tD_f)$ because the derivations $D_e$, $D_f$ are not nilpotent.

For reasons that will be explained in Remark \ref{explanation}, we will need a result analogous to Theorem~\ref{DCKPteo} but based on different derivations. This leads to substantial differences in the details, so we give all proofs.

Consider the following derivations on $U_{\epsilon}$:
\begin{equation}
\mathcal {E}=zD_x,\qquad
\mathcal{F}=-zD_y,\qquad
\mathcal{H}=-2z^{-1}D_z.
\end{equation}

\begin{lem} \label{Liecommutationslem}
The derivations $\mathcal {E}$, $\mathcal {F}$ and $\mathcal {H}$ satisfy the following commutation relations, gene\-ra\-ting therefore the Lie algebra ${\mathfrak{sl}}(2)$:
\begin{equation}\label{Liecommutations}
[\mathcal{H},\mathcal{E}]=2\mathcal{E},\qquad
[\mathcal{H},\mathcal{F}]=-2\mathcal{F}, \qquad
[\mathcal{E},\mathcal{F}]=\mathcal{H}.
\end{equation}
\end{lem}

\begin{proof}
Straightforward computation using (\ref{bracketQCA}).
\end{proof}

For any ${\mathbb C}$-algebra $A$ and endomorphism ${\mathcal V}\in {\operatorname{End}}(A)$ we can define
\[
\fonc{\exp(t{\mathcal V})}{A}{A[[t]],}{a}{\sum_{n\geq 0}\frac{t^n}{n!}{\mathcal V}^n(a).}
\]
If ${\mathcal V}$ is a derivation, $\exp(t{\mathcal V})$ is a morphism of algebras. It admits a unique extension (by continuity for the $t$-adic topology)
$\exp(t{\mathcal V})\colon A[[t]]\rightarrow A[[t]]$ which is an automorphism with inverse $\exp(-t{\mathcal V}).$

We now give expressions of $\exp(t{\mathcal F})$, $\exp(t{\mathcal E})\colon U_\e\ra U_\e[[t]]$. For every $\alpha\in {\mathbb C}$ consider the following elements of ${\mathbb C}[[t]]$:
\begin{gather*}
(1+t)^{\alpha} =\sum_{n\geq 0}\frac{1}{n!}t^n \alpha(\alpha-1)\cdots (\alpha-n+1),
\\
\psi_{\alpha}(t) =\sum_{n\geq 1}\frac{(-1)^n}{n!}t^{n-1} \alpha(\alpha-1)\cdots (\alpha-n+1)= \frac{(1-t)^{\alpha}-1 }{t} .
\end{gather*}
For any $a\in \mathcal{Z}_0(U_\e)$, we can similarly define elements $(1+ta)^{\alpha}$, $\psi_{\alpha}(ta)\in \mathcal{Z}_0(U_\e)[[t]]$.

\begin{lem}
The action of $\exp(t{\mathcal F})$ is given by:
\begin{gather*}
\exp(t{\mathcal F})(K) =(1-tyz)^{-1/l}K, \qquad
\exp(t{\mathcal F})\big(K^{-1}\big)=(1-tyz)^{1/l}K^{-1},
\\
\exp(t{\mathcal F})(F) =F,
\\
\exp(t{\mathcal F})(E) =E-\big(\epsilon-\epsilon^{-1}\big)^{l-2}\big(K\epsilon^{-1} tz\psi_{-1/l}(tyz)+K^{-1}\epsilon tz\psi_{1/l}(tyz)\big)F^{l-1}.
\end{gather*}
\end{lem}

\begin{proof}
From (\ref{derform4}) we have
\begin{equation}\label{actionF}
{\mathcal F}(K)=\frac{zy}{l}K,\qquad
{\mathcal F}(F)=0, {\qquad
\mathcal F}(E)=-\big(\epsilon-\epsilon^{-1}\big)^{l-1}\frac{z}{l}[K,-1]F^{l-1}.
\end{equation}
 Iterating this, a straightforward computation proves the lemma.\end{proof}

\begin{lem} 
The action of $\exp(t{\mathcal E})$ is given by
\begin{gather*}
\exp(t{\mathcal E})(K) =(1-txz)^{-1/l}K, \qquad
\exp(t{\mathcal E})\big(K^{-1}\big)=(1-txz)^{1/l}K^{-1},
\\
\exp(t{\mathcal E})\big(K^{-1}E\big) =K^{-1}E,
\\
\exp(t{\mathcal E})(FK) =FK+\big(\epsilon-\epsilon^{-1}\big)^{l-2}\big(K\epsilon tz\psi_{-1/l}(txz)+K^{-1}\epsilon^{-1} tz\psi_{1/l}(txz)\big)\big(K^{-1}E\big)^{l-1}.
\end{gather*}

\end{lem}
\begin{proof} This is a little trickier than the previous proof, because it is $D_y$ which appears in (\ref{derform4}) and ${\mathcal E}$ involves $D_x.$
We therefore use the automorphism $T_{-1}$ which is such that $T_{-1}(y)=-x$. As a result, by applying it to (\ref{actionF}), we obtain
\begin{equation*}
D_x\big(K^{-1}\big)=-\frac{x}{l}K^{-1}, \qquad
D_x\big(K^{-1}E\big)=0, \qquad
D_x(FK)=\big(\epsilon-\epsilon^{-1}\big)^{l-1}\frac{1}{l}[K,1]\big(K^{-1}E\big)^{l-1}
\end{equation*}
which imply
\begin{equation*}
{\mathcal E}(K)=\frac{zx}{l}K,\qquad
{\mathcal E}\big(K^{-1}E\big)=0, \qquad
{\mathcal E}(FK)=\big(\epsilon-\epsilon^{-1}\big)^{l-1}\frac{z}{l}[K,1]\big(K^{-1}E\big)^{l-1}.
\end{equation*}
 Iterating this, a straightforward computation proves the lemma.
 \end{proof}

Because ${\mathcal E}$ and ${\mathcal F}$ leave invariant $\mathcal{Z}(U_\e)$ and $\mathcal{Z}_0(U_\e)$ they define maps $\exp(t{\mathcal E})$, $\exp(t{\mathcal F})$: ${\mathcal Z}(U_\e)\rightarrow {\mathcal Z}(U_\e)[[t]]$ and ${\mathcal Z}_0(U_\e)\rightarrow {\mathcal Z}_0(U_\e)[[t]]$.
\begin{prop} \label{actioncenter}
We have
\begin{gather*}
\exp(t{\mathcal E})(x) =x,\qquad \exp(t{\mathcal F})(y) =y,\qquad \exp(t{\mathcal E})(\Omega) =\exp(t{\mathcal F})(\Omega)=\Omega,
\\
\exp(t{\mathcal E})(z)=(1-tzx)^{-1}z,\qquad \exp(t{\mathcal E})(y)=y+t\big({-}xyz+z-z^{-1}\big)+t^2x,
\\
\exp(t{\mathcal F})(z)=(1-tzy)^{-1}z,\qquad \exp(t{\mathcal F})(x)=x+t\big({-}xyz+z-z^{-1}\big)+t^2y.
\end{gather*}
\end{prop}

\begin{proof}
Rather than using the explicit action of $\exp(t{\mathcal E})$and $\exp(t{\mathcal F})$ on $U_\epsilon$, we remark that from (\ref{DonM}) we have
\begin{equation}\label{coadjointactionproof}
\exp(t{\mathcal E})({\mathcal M})=\exp(-tX){\mathcal M}\exp(tX),\qquad
\exp(t{\mathcal F})({\mathcal M})=\exp(-tY){\mathcal M}\exp(tY).
\end{equation}
The announced expressions follow by writing the matrix elements of these equations.\end{proof}

Let ${\mathcal G}$ be the free product of $({\mathbb C}, +)$ with itself; it can be equivalently seen as the group generated by elements $\phi_s$, $\psi_s$, where $s\in {\mathbb C}$, with the relations $\phi_s\phi_{s'}=\phi_{s+s'}$, $\psi_s\psi_{s'}=\psi_{s+s'}$ for every $s, s'\in {\mathbb C}$.

We can define a partial action of ${\mathcal G}$ on $\operatorname{Spec}(\mathcal{Z}(U_\e))$ and on $\operatorname{Spec}(\mathcal{Z}_0(U_\e))$, in the sense of~\cite{Ex1}. Let $s\in {\mathbb C}$, and denote by $D(\phi_s)$ the set of $\chi\in \operatorname{Spec}(\mathcal{Z}_0(U_\epsilon))$ such that for all $u\in \mathcal{Z}_0(U_\epsilon)$ the series $(\exp(s{\mathcal E})(u))(\chi)$ is normally convergent in a small neighborhood of $s$. Equivalently $D(\phi_s)=\{\chi\in \operatorname{Spec}(\mathcal{Z}_0(U_\epsilon)), \,\vert s x_\chi z_\chi\vert <1\}$. We define an action of $\phi_s$ on $D(\phi_s)$ by
\begin{gather}\label{defcoadjointaction1}
u(\phi_s.\chi) = (\exp(-s{\mathcal E})(u))(\chi).
\end{gather}
Similarly, put $D(\psi_s)=\{\chi\in \operatorname{Spec}(\mathcal{Z}_0(U_\epsilon)), \vert s y_\chi z_\chi\vert <1\}$ and define an action of $\psi_s$ on $D(\psi_s)$ by
\begin{gather}\label{defcoadjointaction2}
u(\psi_s.\chi) = (\exp(-s{\mathcal F})(u))(\chi).
\end{gather}
The domains $D(\phi_s)$ (resp.~$D(\psi_s)$) cover $\operatorname{Spec}(\mathcal{Z}_0(U_\epsilon))$ as $s$ varies in ${\mathbb C}$. By results of Exel (\cite{Ex1} and \cite[Examples~1 and~4]{Ex2}), the set formed by the local actions (\ref{defcoadjointaction1}) (resp.~(\ref{defcoadjointaction2})) on the domains $D(\phi_s)$ (resp.~$D(\psi_s)$) defines a partial action of the one-parameter group $(\phi_s)_{s\in \mathbb C}$ (resp.~$(\psi_s)_{s\in \mathbb C}$) on $ \operatorname{Spec}(\mathcal{Z}(U_\epsilon))$, and combining the two we get a partial action of their free product ${\mathcal G}$ on $ \operatorname{Spec}(\mathcal{Z}(U_\epsilon))$. Similarly we get a partial action of ${\mathcal G}$ on $ \operatorname{Spec}(\mathcal{Z}(U_\epsilon))$ by~repla\-cing~$\mathcal{Z}_0(U_\epsilon)$ by~$\mathcal{Z}(U_\epsilon)$.

We call them the {\it partial quantum coadjoint actions} of $\mathcal{G}$. Orbits (called {\it partial orbits} in \cite{Ex1}) are defined as for group actions: the $\Gg$-orbit of a point $\chi \in\operatorname{Spec}(\mathcal{Z}(U_\epsilon))$ is the set of points $g.\chi$, for every possible $g\in \mathcal{G}$ such that $g.\chi$ is defined.

The next result states the analogs of Theorem~\ref{DCKPteo}$(c)$--$(f)$ obtained by replacing $\mathcal{G}_{\rm DCK}$ by $\mathcal{G}$. In particular, it describes the $\mathcal{G}$-orbits.

\begin{teo}\label{QCAteo}\qquad

\begin{enumerate}\itemsep=0pt
\item[$(a)$] The diffeomorphism $\bar \sigma \circ \bar \psi^{-1}$ maps the action of $\Gg$ on the tangent spaces of $(0,0,\pm 1)\in \operatorname{Spec}(\mathcal{Z}_0(U_\epsilon))$ to the coadjoint action of $G$ on ${\mathfrak{sl}}(2,\mc)^* = T_{\pm I}(G^0)$.

\item[$(b)$] For any conjugacy class $\Gamma$ in $G$, $\big(\bar \psi\circ \bar \sigma^{-1}\big)\big(\Gamma \cap G^0\big)$ is a $($non empty$)$ $\Gg$-orbit in $\operatorname{Spec}(\mathcal{Z}_0(U_\epsilon))$, and $\tau^{-1}$ of this orbit has $l$ connected components, all of whose are $\Gg$-orbits in $\operatorname{Spec}(\mathcal{Z}(U_\epsilon))$.

\item[$(c)$] An element $a\in {\mathcal{Z}}(\Ue)$ is $\mathcal{G}$-invariant if and only if $a\in \mc[\Omega]$. Dually, the sets of fixed points of the partial $\Gg$-action on $\operatorname{Spec}(\mathcal{Z}_0(U_\epsilon))$ and $\operatorname{Spec}(\mathcal{Z}(U_\epsilon))$ are respectively $(0,0,\pm 1):= \big(\bar \psi \circ \bar \sigma^{-1}\big)(\{{\rm \pm I}\})$ and
\[
\mathcal{D} := \tau^{-1}((0,0,\pm 1)) = \{(0,0,\pm 1,\pm c_j), \, j=1,\dots,(l-1)/2\}.
\]

\item[$(d)$] The $\Gg$-orbits on $\operatorname{Spec}(\mathcal{Z}_0(U_\epsilon))$ and $\operatorname{Spec}(\mathcal{Z}(U_\epsilon))$ are the symplectic leaves of $\{\, ,\, \}_{\rm QCA}$.

 \end{enumerate}
\end{teo}

\begin{proof} $(a)$ Using $\bar \sigma^{-1*}\circ \bar{\psi}^*$ to identify $\mathcal{Z}_0(U_\epsilon)$ with $\mathcal{O}\big(G^0\big)$, this is a direct consequence of (\ref{coadjointactionproof}) and the fact that ${\mathcal M}=\pm I $ at the fixed point.

$(b)$ That $\Gamma \cap G^0$ is non-empty and connected is classical. Then so is $\big(\bar \psi\circ \bar \sigma^{-1}\big)\big(\Gamma \cap G^0\big)$. Take a point $\chi\in \big(\bar \psi\circ \bar \sigma^{-1}\big)\big(\Gamma \cap G^0\big)$. Evaluating (\ref{coadjointactionproof}) at $\chi$ shows that $\mathcal{G}\cdot \chi$ is contained in $\big(\bar \psi\circ \bar \sigma^{-1}\big)\big(\Gamma \cap G^0\big)$. We claim that $\mathcal{G}\cdot \chi$ is an open and closed subset, so by connectedness of $\big(\bar \psi\circ \bar \sigma^{-1}\big)\big(\Gamma \cap G^0\big)$ it coincides with it. Indeed, it is an open subset because any point of $\operatorname{Spec}(\mathcal{Z}_0(U_\epsilon))$ (whence of $\mathcal{G}\cdot \chi$) belongs to the domains $D(\phi_s)$ and $D(\psi_s)$ for $s$ small enough. Since the one-parameter groups $(\phi_s)_{s\in \mathbb C}$ and $(\psi_s)_{s\in \mathbb C}$ are obtained by integrating the derivations~${\mathcal E}$ and~${\mathcal F}$, it follows from Lemma \ref{Liecommutationslem} that $\mathcal{G}\cdot \chi$ contains a neighborhood in $\big(\bar \psi\circ \bar \sigma^{-1}\big)\big(\Gamma \cap G^0\big)$ of each of its points. By the same reason any limit point of $\mathcal{G}\cdot \chi$ in $\operatorname{Spec}(\mathcal{Z}_0(U_\epsilon))$ has a neighborhood where the partial $\Gg$-action is defined and which intersects $\mathcal{G}\cdot \chi$. Therefore it belongs to $\mathcal{G}\cdot \chi$, which shows that $\mathcal{G}\cdot \chi$ is also closed in $\big(\bar \psi\circ \bar \sigma^{-1}\big)\big(\Gamma \cap G^0\big)$. Because the covering map $\tau$ is unramified of degree $l$, the result for $\operatorname{Spec}(\mathcal{Z}(U_\epsilon))$ follows at once.

$(c)$ Any $\chi\in \operatorname{Spec}(\mathcal{Z}_0(U_\epsilon))$ belongs to the domains $D(\phi_s)$, $D(\psi_s)$ for $s$ small enough. Solving the equations $\phi_s(\chi)=\chi$ and $\psi_s(\chi)=\chi$ by using the formulas in Proposition~\ref{actioncenter} imposes $\chi=(0,0,\pm 1)$. This gives the fixed point set of $\operatorname{Spec}(\mathcal{Z}_0(U_\epsilon))$; the result for $\operatorname{Spec}(\mathcal{Z}(U_\epsilon))$ follows immediately. As for the first claim, note that for any conjugacy class $\Gamma$ of maximal dimension the set $\Gamma \cap G^0$ contains a diagonal matrix, and that the union of such sets forms a Zariski open and dense subset of $G$. Hence (b) above implies that a central invariant element is completely determined by its value at points $(0,0,z_\chi,\Omega_\chi)$. Therefore it is non zero if and only if it belongs to $\mc[\Omega]$.

$(d)$ The groups $\mathcal{G}$ and $\mathcal{G}_{\rm DCK}$ are obtained by integrating the derivations ${\mathcal E}$, ${\mathcal F}$, and~$D_e$,~$D_f$ respectively, and the Lie algebras generated by these two pairs of derivations have the same span at every point of $\operatorname{Spec}(\mathcal{Z}_0(U_\epsilon))$ or $\operatorname{Spec}(\mathcal{Z}(U_\epsilon))$. Then the conclusion follows from Theorem~\ref{DCKPteo}$(f)$.\end{proof}

\begin{remk}\label{explanation}
A difficulty with the group $\mathcal{G}$ is that its elements act only on subsets of $\operatorname{Spec}(\mathcal{Z}(U_\epsilon))$: in the formulas (\ref{coadjointactionproof}), the parameter $t\in \mc$ must be such that the lower right entries of the computed matrices are non zero. However it has various merits as compared to the group $\mathcal{G}_{\rm DCK}$ of Theorem \ref{DCKPteo}:
\begin{enumerate}\itemsep=0pt

\item[$(i)$] $\mathcal{G}$ is finite dimensional, associated to the Lie algebra ${\mathfrak{sl}}(2)$ (by Lemma \ref{Liecommutationslem}), whereas $\mathcal{G}_{\rm DCK}$ is infinite-dimensional.

\item[$(ii)$] The partial action of $\mathcal{G}$ on $\operatorname{Spec}(\mathcal{Z}(U_\epsilon))$ is by birational transformations (by Proposition~\ref{actioncenter}), whereas $\mathcal{G}_{\rm DCK}$ is a subgroup of ${\rm Aut}(\hat U_\e)$, acting on $\operatorname{Spec}(\mathcal{Z}(U_\epsilon))$ by holomorphic (entire) transformations.

\item[$(iii)$] $\mathcal{G}$ can be generalized straightforwardly to $\Ll_{0,n}^\e$, whereas $\mathcal{G}_{\rm DCK}$ does not. This is our main motivation for developing this construction. We do it in the next section.

\end{enumerate}
\end{remk}

\begin{remk}\label{tildeG}{\rm We think useful to have in mind the following description of $\operatorname{Spec}(\mathcal{Z}(U_\epsilon))$ (see~\cite{DC-K}). Denote by $G/\!/G$ the affine variety with coordinate ring $\mathcal{O}(G)^G$, the ring of regular functions on $G$ invariant under the coadjoint action of $G$. We have an isomorphism
\begin{equation*}
G/\!/G \cong \mc^*/\big(t\sim t^{-1}\big).
\end{equation*}
In fact, denote by $T$ the torus of $G$ formed by the diagonal matrices, and by $\mathcal{O}(T)$ its coordinate ring. Then $T\cong \mc^*$ and $\mathcal{O}(T)\cong \mc[t,t^{-1}]$, where $t$ is the coordinate function of the upper left entry of elements of $T$. The Weyl group $W$ of $G$ acts on $T$ by inversion, and
\begin{equation*}
\mathcal{O}(G)^G \cong \mathcal{O}(T)^W = \mc\big[t+t^{-1}\big].
\end{equation*}
 Consider the maps
\begin{equation*}
 p\colon\ G\lra G/\!/G,\qquad p_k\colon\ G/\!/G \lra G/\!/G,
\end{equation*}
where $p\colon G\ra G/\!/G$ is the quotient map, and $p_k$ is induced by the $k$-th power map $g\mapsto g^k$, $g\in G$, $k\in \mz$. Note here that on coordinate functions we have $p_k(t+t^{-1}) = t^k+t^{-k}$, so $p_k$ is just realized by the $k$-th Chebyshev polynomial $T_k$. Consider the fibered product of $p$ and $p_l$, that is, the affine variety
\begin{equation*}
G\times_{G/\!/G} G/\!/G= \{ (g,[t])\in G\times G/\!/G\ |\ p(g) = p_l([t])\}.
\end{equation*}
Set
\begin{equation*}
\tilde G = G\times_{G/\!/G} G/\!/G,\qquad
\tilde G^0 = G^0\times_{G/\!/G} G/\!/G.
\end{equation*}
Then $\operatorname{Spec}(\mathcal{Z}(U_\epsilon))$ is isomorphic to $\tilde G^0$. In fact, by Theorem \ref{relbracket} and the defining relation~\eqref{relcenter} we know that $\mathcal{Z}_0(U_\epsilon)$ is isomorphic to $\mathcal{O}(G^0)$, and that $\mathcal{Z}(U_\epsilon) = \mathcal{Z}_0(U_\epsilon)\otimes_{\mathcal{Z}_0(U_\epsilon)\cap \mc[\Omega]} \mc[\Omega]$. By~the quantum Harish-Chandra homomorphism, see \cite{DC-K} for details, $\mc[\Omega]\cong \mc\big[K+K^{-1}\big]$ and $\mathcal{Z}_0(U_\epsilon) \cap \mc[\Omega] \cong \mc\big[K^l+K^{-l}\big]$. Let us identify $\mc\big[K,K^{- 1}\big]$ with $\mathcal{O}(T)$. Then $\mc[\Omega]\cong \mathcal{O}(T)^W$ and $\mc[K^l+K^{-l}]\cong \mathcal{O}(T/\mu_l)^W$, where $\mu_l$ is the subgroup of $T$ corresponding to the $l$th-roots of unity under the isomorphism $T\cong \mc^*$. Hence
\[\mathcal{Z}(U_\epsilon)\cong \mathcal{O}\big(G^0\big) \otimes_{\mathcal{O}(T/\mu_l)^W} \mathcal{O}(T)^W.\]
The isomorphism of $\operatorname{Spec}(\mathcal{Z}(U_\e))$ with $G^0 \times_{G/\!/G} G/\!/G$ follows by duality.}
\end{remk}

\subsubsection[Extension to L\_\{{0,n}\}\textasciicircum{}e]
{Extension to $\boldsymbol{\Ll_{0,n}^\e}$}\label{EXTENSION}

The results of the previous section extend naturally to $\Ll_{0,n}^\epsilon$. First we consider the generalization of Theorem~\ref{relbracket}.

The bracket $\{\, ,\, \}_{\rm FR}$ on $G$ has been extended to $G^n$ by Fock--Rosly \cite{FR}. From its very definition, it is readily checked to be defined by
\begin{align}
\big\{\!\stackrel{1}{L}\!{}^{(i)},\stackrel{2}{L}\!{}^{(i)}\big\}_{\rm FR} = \mathfrak{r}\stackrel{1}{L}\!{}^{(i)}\!\!\stackrel{2}{L}\!{}^{(i)} -\stackrel{1}{L}\!{}^{(i)}\!\!\stackrel{2}{L}\!{}^{(i)}\mathfrak{r}'+\stackrel{2}{L}\!{}^{(i)}\mathfrak{r}'\!\stackrel{1}{L}\!{}^{(i)} -\stackrel{1}{L}\!{}^{(i)}\mathfrak{r}\!\stackrel{2}{L}\!{}^{(i)},\label{bra1n}
\\
\big\{\!\stackrel{1}{L}\!{}^{(i)},\stackrel{2}{L}\!{}^{(j)} \big\}_{\rm FR} = \mathfrak{r}\stackrel{1}{L}\!{}^{(i)}\!\!\stackrel{2}{L}\!{}^{(j)} +\stackrel{1}{L}\!{}^{(i)}\!\!\stackrel{2}{L}\!{}^{(j)}\mathfrak{r}-\stackrel{2}{L}\!{}^{(j)}\mathfrak{r}\!\stackrel{1}{L}\!{}^{(i)} -\stackrel{1}{L}\!{}^{(i)}\mathfrak{r}\!\stackrel{2}{L}\!{}^{(j)},\label{bra2n}
 \end{align}
where $L^{(i)} = \underline{\stackrel{V_2}{M}}{}^{(i)}$ (see \eqref{notFrmat}), $i,j\in \{1,\dots,n\}$, and $i<j$.

The formula \eqref{defder} extends naturally to define derivations $D_a\colon {}_{\rm loc}\Ll_{0,n}^\epsilon \ra {}_{\rm loc}\Ll_{0,n}^\epsilon$ preserving $\mathcal{Z}({}_{\rm loc}\Ll_{0,n}^\epsilon)$, for every $a\in \mathcal{Z}({}_{\rm loc}\Ll_{0,n}^\epsilon)$. Hence we can define a Poisson bracket $\{\, ,\, \}_{\rm QCA}$ on $\mathcal{Z}({}_{\rm loc}\Ll_{0,n}^\epsilon)$ (keeping voluntarily the same notation as on $\mathcal{Z}(U_\e)$) by \[\{a,b\}_{\rm QCA}= D_a(b)\]
for every $a,b\in\mathcal{Z}({}_{\rm loc}\Ll_{0,n}^\epsilon).$ As $\Phi_n$ is an isomorphism of algebras, both derivations and Poisson bracket can be defined and computed by means of their pushforwards on $U_\e^{\otimes n}$. Indeed, for any $a\in\mathcal{Z}({}_{\rm loc}\Ll_{0,n}^\e)$, $\Phi_n(a) \in \mathcal{Z}(U_\e^{\otimes n})$, so the derivation $D_{\Phi_n(a)}\colon U_\e^{\otimes n} \ra U_\e^{\otimes n}$ is defined and we can put $D_a=\Phi_n^{-1} D_{\Phi_n(a)} \Phi_n$. For every $a\in \mathcal{Z}(\Ll_{0,n}^\epsilon)$, by the inclusion $\Ll_{0,n}^\epsilon \subset {}_{\rm loc}\Ll_{0,n}^\epsilon $ it yields a~derivation $D_a\colon \Ll_{0,n}^\epsilon \ra \Ll_{0,n}^\epsilon$ preserving $\mathcal{Z}(\Ll_{0,n}^\epsilon)$.

Denote by $\mathcal{Z}_0(\Ll_{0,n}^\epsilon)$ (resp.~$\mathcal{Z}_0({}_{\rm loc}\Ll_{0,n}^\epsilon)$) the subalgebras of $\mathcal{Z}(\Ll_{0,n}^\epsilon)$ \big(resp.~$\mathcal{Z}({}_{\rm loc}\Ll_{0,n}^\epsilon)$\big) generated by $T_l\big(\omega^{(i)}\big)$, $b^{(i)l}$, $c^{(i)l}$, $d^{(i) l}$ \big(resp.~$T_l\big(\omega^{(i)}\big)$, $b^{(i)l}$, $c^{(i)l}$, $d^{(i) l}$ and $\delta^{(i) \pm l}$\big) for $i =1,\dots, n$.
Recall the Frobenius morphism ${\rm Fr}\colon \big(\Ll_{0,1}^1\big)^{\otimes n}\rightarrow \mathcal{Z}_0(\Ll_{0,n}^\epsilon)$ (see \eqref{notFrmat}), and the identification $\Ll_{0,1}^1=\mathcal{O}(G)$.

\begin{teo}\label{scteo} The map
\[
\Phi := \big(\bar \sigma^{-1*}\circ \bar{\psi}^*{}\big)^{\otimes n}\circ \Phi_1^{\otimes n}\colon\ \big(\mathcal{Z}_0({}_{\rm loc}\Ll_{0,n}^\epsilon),\{\ ,\}_{\rm QCA}\big)\ra \big(\mathcal{O}(G^0)^{\otimes n},\{\, ,\, \}_{\rm FR}\big)
\]
is an isomorphism of Poisson algebras, and $\Phi^{-1}$ restricted to $\mathcal{O}(G)^{\otimes n}\subset \mathcal{O}(G^0)^{\otimes n}$ coincides with the Frobenius homomorphism. Hence ${\rm Fr} \colon \big(\mathcal{O}(G)^{\otimes n},\{\, ,\, \}_{\rm FR}\big) \ra \big(\mathcal{Z}_0(\Ll_{0,n}^\epsilon),\{\ ,\}_{\rm QCA}\big)$ is an isomorphism of Poisson algebras.
\end{teo}

\begin{proof}
The map $\bar \sigma^{-1*}\circ \bar{\psi}^*{}\colon \big(\mathcal{Z}_0(U_\epsilon),\{\, ,\, \}_{\rm QCA}\big) \ra \big(\mathcal{O}(G^0),\{\, ,\, \}_{\rm FR}\big)$ is an isomorphism of Poisson algebras by the results recalled in the previous section. Since $\Phi_1 \colon {}_{\rm loc}\Ll_{0,1}^A\ra U_A'$ is an algebra isomorphism by Lemma \ref{intsl2}, when specializing at $q=\e$ it yields a Poisson isomorphism between $\big(\mathcal{Z}_0({}_{\rm loc}\Ll_{0,1}^\epsilon),\{\ ,\}_{\rm QCA}\big)$ and $\big(\mathcal{Z}_0(U_\epsilon),\{\, ,\, \}_{\rm QCA}\big)$. This proves the case $n=1$.

When $n\geq 2$, the map is well-defined because $\mathcal{Z}_0({}_{\rm loc}\Ll_{0,n}^\epsilon) = {}_{\rm loc}\mathcal{Z}_0(\Ll_{0,n}^\epsilon)$ (the localization by the powers $\delta^{(i) lk}$, $k\in \mathbb{Z}$), $\mathcal{Z}_0(\Ll_{0,n}^\epsilon) = \mathcal{Z}_0(\Ll_{0,1}^\epsilon)^{\otimes n}$ by Proposition \ref{center0nroot}, and $\Phi_1^{\otimes n}$ extends from $\mathcal{Z}_0(\Ll_{0,n}^\epsilon)$ to ${}_{\rm loc}\mathcal{Z}_0(\Ll_{0,n}^\epsilon)$ for obvious reasons (the algebra being commutative). On another hand, since $\Phi_n \colon {}_{\rm loc}\Ll_{0,n}^\e\ra U_\e^{\otimes n}$ is an equivariant algebra isomorphism (see \eqref{Phinlociso}) and $\Phi_n\big(\omega^{(i)}\big) = \Omega^{(i)}$, it yields a Poisson isomorphism between $\big(\mathcal{Z}_0({}_{\rm loc}\Ll_{0,n}^\epsilon),\{\ ,\}_{\rm QCA}\big)$ and $(\mathcal{Z}_0(U_\epsilon),\{\, ,\, \}_{\rm QCA})^{\otimes n}$, i.e., the algebra $\mathcal{Z}_0(U_\epsilon)^{\otimes n}$ with the product Poisson structure. Hence it is enough to prove that
\[
\Phi \circ \Phi_n^{-1}\colon\ \big(\mathcal{Z}_0(U_\epsilon\big),\{\, ,\, \}_{\rm QCA})^{\otimes n} \ra \big(\mathcal{O}\big(G^0\big)^{\otimes n},\{\, ,\, \}_{\rm FR}\big)
\]
is an isomorphism of Poisson algebras. This can be checked on generators, which is most easily done by using the inverse map. Put $\Pi = {\rm id}_{V_2} \otimes (\Phi_n\circ (\Phi_1^{\otimes n})^{-1})$. By Corollary~\ref{dressingcenter} we have
\begin{equation}\label{defdresscentre}
\Pi\big(\mathcal{M}^{(i)}\big)=\mathcal{R}^{(i)}\mathcal{M}^{(i)}\mathcal{R}^{(i)-1}.
\end{equation}
Moreover, \eqref{Mcal} and \eqref{formmatG} show that $\bar \sigma^{-1*}\circ \bar{\psi}^*$ identifies the matrix coefficients of $\Phi_1(\mathcal{Z}_0(\Ll_{0,1}^\epsilon))$ in the fundamental representation $V_2$ of $U_\e$, and the matrix coefficients of the fundamental representation of $G$ on $\mc^2$. Therefore it is enough to check that
\begin{equation}\label{finalformP}
\big\{\Pi^*\big(\!\stackrel{1}{L}\!{}^{(i)}\big),\Pi^*\big(\!\stackrel{2}{L}\!{}^{(j)}\big)\big\}_{\rm QCA}^{\otimes n} = \big\{\!\stackrel{1}{L}\!{}^{(i)},\stackrel{2}{L}\!{}^{(j)}\big\}_{\rm FR}
\end{equation}
for every $i,j\in \{1,\dots,n\}$, where $\Pi^*\big(L^{(i)}\big)$ is the pull-back of $L^{(i)}$ via $\Pi$, i.e., its expression as a~function of the matrices $\mathcal{M}^{(i)}$. Now we have
\begin{equation}\label{splitmat}
\mathcal{M}^{(i)} = \mathcal{M}^{(i)}_+\mathcal{M}^{(i)-1}_-,
\end{equation}
where
\[
\mathcal{M}^{(i)}_+ = \begin{pmatrix} z'{}^{(i)} & z'{}^{(i)}y{}^{(i)} \\ 0 & z'{}^{(i)-1} \end{pmatrix}\!, \qquad
\mathcal{M}^{(i)}_- = \begin{pmatrix} z'{}^{(i)-1} & 0 \\ z'{}^{(i)}x{}^{(i)} & z'{}^{(i)} \end{pmatrix}\!.
\]
Note also that \eqref{g11l} and~\eqref{g12l} imply
\begin{equation}\label{Gi+-}
\mathcal{R}^{(i)} = \mathcal{M}^{(n)}_+\cdots \mathcal{M}^{(i+1)}_+.
\end{equation}
Using the map $\sigma^{-1*}\circ \psi^*$ to identify $\mathcal{M}^{(i)}_+$, $\mathcal{M}^{(i)}_-$ with matrices $L_+$, $L_-$ of coordinate functions in $G^*$, we see that they satisfy the bracket identities \eqref{bracketdual}. Since $G^*$ is a Poisson--Lie group, products of matrices $L_+$ have the same bracket as $L_+$. In particular this applies to $\mathcal{R}^{(i)}$, and we can write
\[
\Pi^*(L^{(i)}) =\mathcal{R}^{(i)}\!\!\mathcal{M}^{(i)}_+\!\!\mathcal{M}^{(i)-1}_-\mathcal{R}^{(i)-1},
\]
where all matrices in the product have known brackets. Then \eqref{finalformP} can be straightforwardly compared with \eqref{bra1n} and~\eqref{bra2n}, using the Leibniz rule for simplifications, and the fact that $r+r'$ is ad-invariant. We leave the verifications to the reader.

The equality of $\Phi^{-1}$ with Fr on $\mathcal{O}(G)^{\otimes n}$ follows immediatly from the fact that $\bar \sigma^{-1*}\circ \bar{\psi}^*$ identifies matrix coefficients, as discussed above. The image of $Fr$ is $\mathcal{Z}_0(\Ll_{0,n}^\epsilon)$, by Proposition~\ref{center0nroot}. This achieves the proof.
\end{proof}

Next we turn to our generalization of Theorem \ref{QCAteo}. First we define a partial action of the group $\mathcal{G}$ on $\operatorname{Spec}({\mathcal{Z}}_0({}_{\rm loc}\Ll_{0,n}^\epsilon))$ and $\operatorname{Spec}({\mathcal{Z}}({}_{\rm loc}\Ll_{0,n}^\epsilon))$ by generalizing the method we used in the $n=1$ case.

Recall the derivations $D_a\colon {}_{\rm loc}\Ll_{0,n}^\epsilon \ra {}_{\rm loc}\Ll_{0,n}^\epsilon$, defined for every $a\in \mathcal{Z}({}_{\rm loc}\Ll_{0,n}^\e)$. Denote by $\hat{x}^{(i)}, \hat{y}^{(i)}, \hat{z}^{(i)}\in \mathcal{Z}({}_{\rm loc}\Ll_{0,n}^\e)$ the inverse images by $\Phi_n$ of the elements $x^{(i)}$, $y^{(i)}$, $z^{(i)}$, and let ${\mathcal E}^{(i)}$, ${\mathcal F}^{(i)}$, ${\mathcal H}^{(i)}$ be the derivations of ${}_{\rm loc}\Ll_{0,n}^\epsilon$ defined by
\begin{equation*}
{\mathcal E}^{(i)}= \hat{z}^{(i)} D_{\hat{x}^{(i)}},\qquad
{ \mathcal F}^{(i)}=-\hat{z}^{(i)} D_{\hat{y}^{(i)}},\qquad
{\mathcal H}^{(i)}=-2\hat{z}^{(i)-1} D_{\hat{z}^{(i)}}.
\end{equation*}
Note that because of the relations (\ref{phinclbis}) and the definition of $\delta^{(i)}$ we can obtain simple formulas for $\hat{z}^{(i)}$ and $\hat{x}^{(i)},$ namely:
\begin{equation*}
\hat{z}^{(i)\pm 1}=\delta^{(i)\mp l},\qquad
\hat{x}^{(i)}=-c^{(i)l}\prod_{k=i+1}^n \delta^{(k)-l}.
\end{equation*}
Set
\begin{equation*}
{\mathcal E}^{\Delta}= \sum_{i=1}^n{\mathcal E}^{(i)},\qquad
{\mathcal F}^{\Delta}= \sum_{i=1}^n{\mathcal F}^{(i)},\qquad
{\mathcal H}^{\Delta}= \sum_{i=1}^n{\mathcal H}^{(i)}.
\end{equation*}

\begin{prop}\label{convexpn} \quad

\begin{enumerate}\itemsep=0pt

\item[$(1)$] The derivations ${\mathcal E}^{\Delta}$, ${\mathcal F}^{\Delta}$ and ${\mathcal H}^{\Delta}$ satisfy the following commutation relations, generating the Lie algebra ${\mathfrak{sl}}(2)$:
\begin{equation*}
\big[\mathcal{H}^{\Delta},\mathcal{E}^{\Delta}\big]=2\mathcal{E}^{\Delta}, \qquad \big[\mathcal{H}^{\Delta},\mathcal{F}^{\Delta}\big]=-2\mathcal{F}^{\Delta}, \qquad \big[\mathcal{E}^{\Delta},\mathcal{F}^{\Delta}\big]=\mathcal{H}^{\Delta}.
\end{equation*}

\item[$(2)$] The power series $\exp\big(t{\mathcal E}^{\Delta}\big)$, $ \exp\big(t{\mathcal F}^{\Delta}\big)$ define morphisms of algebras from ${}_{\rm loc}{\Ll}_{0,n}^\e$ to \allowbreak ${}_{\rm loc}{\Ll}_{0,n}^\e[[t]]$, sending ${\mathcal{Z}}({}_{\rm loc}\Ll_{0,n}^\epsilon)$ to ${\mathcal{Z}}({}_{\rm loc}\Ll_{0,n}^\epsilon)[[t]]$ and fixing the elements $\omega^{(i)}$. Moreover
\begin{gather*}
\exp\big(t\Phi_n\big({\mathcal E}^{\Delta}\big)\big)\big({\mathcal M}^{(i)}\big)=\exp(-tX){\mathcal M}^{(i)}\exp(tX),
\\
\exp\big(t\Phi_n\big({\mathcal F}^\Delta\big)\big)\big({\mathcal M}^{(i)}\big)=\exp(-tY){\mathcal M}^{(i)}\exp(tY).
\end{gather*}

\end{enumerate}

As a result $\exp\big(t{\mathcal E}^{\Delta}\big)$ and $\exp\big(t{\mathcal F}^{\Delta}\big)$ are sending ${\mathcal{Z}}_0({}_{\rm loc}\Ll_{0,n}^\epsilon)$ to ${\mathcal{Z}}_0({}_{\rm loc}\Ll_{0,n}^\epsilon)[[t]]$.
\end{prop}

\begin{proof} (1) Straighforward from (\ref{Liecommutations}) and from the local Poisson-commutativity, i.e., the pro\-perty that $\{a,b\}=0$ when $a\in \big\{x^{(i)}, y^{(i)}, z^{(i)}\big\}$, $b\in \big\{x^{(j)}, y^{(j)}, z^{(j)}\big\}$, $i\not=j$.

(2) The first claim is straightforward, the second is a direct application of (\ref{coadjointactionproof}) and the local Poisson-commutativity.
\end{proof}

As in \eqref{defcoadjointaction1} let $s\in {\mathbb C}$, and denote by $D\big(\phi_s^{(i)}\big)$ the set of points $\chi\in \operatorname{Spec}(\!{\mathcal{Z}}_0({}_{\rm loc}\Ll_{0,n}^\epsilon))$ such that for all $u\in {\mathcal{Z}}_0({}_{\rm loc}\Ll_{0,n}^\epsilon)$ the series $\big(\!\exp(s{\mathcal E}^{(i)})(u)\big)(\chi)$ is normally convergent in a small neighborhood of $s$. Thus $D\big(\phi_s^{(i)}\big)=\big\{\chi\in \operatorname{Spec}\big({\mathcal{Z}}_0({}_{\rm loc}\Ll_{0,n}^\epsilon)\big), \big\vert s \hat{x}^{(i)}_\chi \hat{z}^{(i)}_\chi\big\vert <1\big\}$. Define an action of the element $\phi_s^{(i)}$ on $D\big(\phi_s^{(i)}\big)$ by
\begin{equation}\label{defcoadjointaction1bis}
u\big(\phi_s^{(i)}.\chi\big) = \big(\!\exp\big({-}s{\mathcal E}^{(i)}\big)(u)\big)(\chi).
\end{equation}
Similarly, put $D\big(\psi_s^{(i)}\big)=\big\{\chi\in \operatorname{Spec}\big({\mathcal{Z}}_0({}_{\rm loc}\Ll_{0,n}^\epsilon)\big), \big\vert s \hat{y}^{(i)}_\chi \hat{z}^{(i)}_\chi\big\vert <1\big\}$, and define an action of
$\psi_s^{(i)}$ on $D\big(\psi_s^{(i)}\big)$ by
\begin{equation}\label{defcoadjointaction2bis}
u\big(\psi_s^{(i)}.\chi\big) = \big(\!\exp\big({-}s{\mathcal F}^{(i)}\big)(u)\big)(\chi).
\end{equation}
Denote by ${\mathcal G}^{(i)}$ the group generated by $\phi_s^{(i)}, \psi_s^{(i)}, s\in {\mathbb C}$. It is isomorphic to ${\mathcal G}$. Denote by~${\mathcal G}^{\rm tot}$ the direct product of the groups ${\mathcal G}^{(i)}$, and by ${\mathcal G}^{\Delta}$ the subgroup of ${\mathcal G}^{\rm tot}$ generated by the dia\-go\-nal elements $\phi_s^{\Delta}=\big(\phi_s^{(1)},\dots,\phi_s^{(n)}\big)$ \big(resp.~$\psi_s^{\Delta}=\big(\psi_s^{(1)},\dots,\psi_s^{(n)}\big)$\big). These elements act on $D\big(\phi_s^{\Delta}\big)=\cap_{i=1}^n D\big(\phi_s^{(i)}\big), $ \big(resp.~$D(\psi_s^{\Delta})=\cap_{i=1}^n D\big(\psi_s^{(i)}\big) \big)$, by acting dually with the series $\exp\big(s{\mathcal E}^{\Delta}\big)$
\big(resp.~$\exp\big(s{\mathcal F}^{\Delta}\big)$\big).

As in the case $n=1$, the results of Exel \cite{Ex1, Ex2} imply that (\ref{defcoadjointaction1bis}) and (\ref{defcoadjointaction2bis}) define partial action of the group ${\mathcal G}^{\rm tot}\cong {\mathcal G}^n$ and ${\mathcal G}^{\Delta}\cong {\mathcal G}$ on $\operatorname{Spec}({\mathcal{Z}}_0({}_{\rm loc}\Ll_{0,n}^\epsilon)).$ We call the first one the {\it total partial quantum coadjoint action} and the second one the {\it diagonal partial quantum coadjoint action}.

Next we need the following result of Fock--Rosly. Recall that we denote by $X_G(\Sigma)$ the variety of $G$-characters of the sphere $\Sigma$ with $n+1$ punctures, and by $\{\, ,\, \}_{\rm Gold}$ the Goldman Poisson bracket on $X_G(\Sigma)$. Denote by $G^n/\!/G$ the algebraic quotient of $G^n$ by the adjoint action of $G$. It~is the affine variety with coordinate ring $\mathcal{O}(G^n)^{G}$, the ring of regular functions on $G^n$ invariant under the coadjoint action of $G$. The points of $X_G(\Sigma)$ are in one-to-one correspondence with the trace equivalence classes of representations $\pi_1(\Sigma)\ra G$. Therefore, choosing a basepoint and generators of the fundamental group of $\Sigma$ affords an isomorphism of algebraic sets 
\[
\mathfrak{c}\colon\ X_G(\Sigma)\ra G^n/\!/G.
\]

\begin{teo}[\cite{FR}] \label{FRteo} The adjoint action of the Poisson--Lie group $(G,\{\, ,\, \})$ on the Poisson manifold $(G^n,\{\, ,\, \}_{\rm FR})$ is a Poisson map. Hence $\{\, ,\, \}_{\rm FR}$ defines a Poisson bracket on $G^n/\!/G$. Moreover the map $\mathfrak{c}$ is a Poisson isomorphism: $\mathfrak{c}_*\{\, ,\, \}_{\rm Gold} = \{\, ,\, \}_{\rm FR}$.
\end{teo}

The Theorems \ref{scteo} and \ref{FRteo} relate $\mathcal{Z}_0(\Ll_{0,n}^\epsilon)$ with $\mathcal{O}(G^n)$ and $\mathcal{O}(G^n)^G$ with $\mathcal{O}(X_G(\Sigma))$. We~need to ``lift'' these results to the whole center $\mathcal{Z}(\Ll_{0,n}^\epsilon)$ and corresponding rings of regular functions.

At first, recall the isomorphism $\mathcal{O}(\tilde G^0)\cong \mathcal{Z}(U_\epsilon)$ of Remark \ref{tildeG}. It maps $\mathcal{O}(T)^W$ to $\mc[\Omega]$, the $\mathcal{G}$-invariant subalgebra of $\mathcal{Z}(U_\e)$. Composing it with $\Phi_1^{-1}$, we get an isomorphism of $\mathcal{O}(\tilde G^0)$ with $\mathcal{Z}({}_{\rm loc}\Ll_{0,1}^\epsilon)$, mapping the subalgebras $\mathcal{O}(\tilde G)$ and $\mathcal{O}(G)$ to $\mathcal{Z}(\Ll_{0,1}^\epsilon)$ and $\mathcal{Z}_0(\Ll_{0,1}^\epsilon)$ respectively. It follows straightforwardly from the arguments of Theorem \ref{scteo} that $Fr$ extends to an isomorphism
\[
\widetilde{\rm Fr} \colon\ \big(\mathcal{O}\big(\tilde{G}\big)^{\otimes n},\{\, ,\, \}_{\rm FR}\big) \ra \big(\mathcal{Z}(\Ll_{0,n}^\epsilon),\{\ ,\}_{\rm QCA}\big)
\]
mapping the $n$ copies of $\mathcal{O}(T)^W$ associated to the factors of $\mathcal{O}\big(\tilde G^n\big)$ to $\mc\big[\omega^{(1)}\big],\dots,\mc\big[\omega^{(n)}\big]$. Note in particular that $\{\, ,\, \}_{\rm FR}$ extends trivially from $\mathcal{O}(G)^{\otimes n}$ to $\mathcal{O}(\tilde{G})^{\otimes n}$. Also, the restriction map (keeping voluntarily the notation of \eqref{tau})
\begin{equation*}\label{tau-n}
\tau\colon\ \operatorname{Spec}(\mathcal{Z}(\Ll_{0,n}^\epsilon)) \lra \operatorname{Spec}(\mathcal{Z}_0(\Ll_{0,n}^\epsilon)),
\end{equation*}
which by Proposition \ref{center0nroot} is a regular map of degree $l^n$, corresponds under $\widetilde{\rm Fr}$ to the projection map $\tilde G^n\ra G^n$.

Now, let us identify as above $G^n$ with the space of representation $R_G(\Sigma) = \operatorname{Hom}(\pi_1(\Sigma),G)$ (fixing a basepoint and generators of $\pi_1(\Sigma)$). The conjugation action extends trivially from $G^n$ to~$\tilde G^n$. Then the projection map $\tilde G^n\ra G^n$ provides an identification of $\tilde G^n$ with a branched covering space $\tilde R_G'(\Sigma)$ of $R_G(\Sigma)$, endowed with the conjugation action of $G$. The points of~$\tilde R_G'(\Sigma)$ are given by representations $\rho\colon \pi_1(\Sigma)\ra G$ endowed with a choice of solution $x\in \mc$ of the equation $T_l(x) = \operatorname{Tr}(\rho(\gamma_i))$, $1\leq i \leq n$, where $\gamma_1,\dots,\gamma_n$ are the $n$ chosen generators of~$\pi_1(\Sigma)$. Taking algebraic quotients yields a branched covering map $\tilde X_G'(\Sigma)\ra X_G(\Sigma)$ of the same degree~$l^n$, and one can lift Theorem \ref{FRteo} to an isomorphism (with, again, $\{\, ,\, \}_{\rm Gold}$ trivially lifted from $X_G(\Sigma)$):
\[
\tilde{\mathfrak{c}} \colon\ \big(\tilde X_G'(\Sigma),\{\, ,\, \}_{\rm Gold}\big) \ra
\big(\tilde G^n/\!/G,\{\, ,\, \}_{\rm FR}\big).
\]

We can now state and prove our generalization of Theorem \ref{QCAteo}. When considering the (partial) action of $\mathcal G$ on $\operatorname{Spec}(\mathcal{Z}_0(\Ll_{0,n}^\epsilon))$ it will always be meant to be the diagonal action, by means of $\mathcal G^{\Delta}$.

\begin{cor} \label{equivcor} \quad

\begin{enumerate}\itemsep=0pt
\item[$(1)$] The dual diffeomorphism ${\rm Fr}^*\colon \big(\!\operatorname{Spec}(\mathcal{Z}_0(\Ll_{0,n}^\epsilon)),\{\ ,\}_{\rm QCA}\big)\ra (G^{n},\{\, ,\, \}_{\rm FR})$ maps the action of $\Gg$ on the tangent spaces of $(0,0,\pm 1)^n$ onto the coadjoint action of $G$ on ${\mathfrak{sl}}(2,\mc)^{*n}$.

\item[$(2)$] {\sloppy For any conjugacy class $\Gamma$ in $G^n$, $({\rm Fr}^*)^{-1}(\Gamma)$ is a $($non empty$)$ $\Gg$-orbit in $\operatorname{Spec}(\mathcal{Z}_0(\Ll_{0,n}^\epsilon))$, and $\tau^{-1}$ of this orbit has $l^n$ connected components, all of whose are $\Gg$-orbits in $\operatorname{Spec}(\mathcal{Z}(\Ll_{0,n}^\epsilon))$.

 }

\item[$(3)$] The map $\widetilde{\rm Fr} \circ \tilde{\mathfrak{c}}^{-1*}$ takes values in $\mathcal{Z}(\Ll_{0,n}^\epsilon)^{\mathcal{G}}$, and therefore affords an isomorphism of Poisson algebras
\[
\widetilde{\rm Fr} \circ \tilde{\mathfrak{c}}^{-1*}\colon\ \big(\mathcal{O}(\tilde X_G'(\Sigma)),\{\, ,\, \}_{\rm Gold}\big) \ra \big(\mathcal{Z}(\Ll_{0,n}^\epsilon)^{\mathcal{G}},\{\, ,\, \}_{\rm QCA}\big).
\]

\item[$(4)$] The orbits of the group $\Gg^{\rm tot}$ in $\operatorname{Spec}(\mathcal{Z}(\Ll_{0,n}^\e))$ are the symplectic leaves of $\{\, ,\, \}_{\rm QCA}$. These project onto the symplectic leaves of $\{\, ,\, \}_{\rm QCA}$ in $\operatorname{Spec}(\mathcal{Z}(\Ll_{0,n}^\e)^{\mathcal{G}})$.

\end{enumerate}
\end{cor}

\begin{proof} (1) Using as usual the isomorphism $\bar \sigma^{-1*}\circ \bar{\psi}^*$ to identify $\mathcal{Z}_0(U_\epsilon)$ with $\mathcal{O}(G^0)$, it follows easily from \eqref{defdresscentre}, \eqref{splitmat} and \eqref{Gi+-} that the automorphism $\Phi_n\circ \big(\Phi_1^{-1}\big)^{\otimes n}$ of $\mathcal{O}(G^0)^{\otimes n}$ is equi\-va\-riant with respect to the coadjoint action of $G$. Since $\Phi_n^{-1*}((0,0,\pm 1)^n) = (\pm I,\dots,\pm I)$, the formulas in Proposition \ref{convexpn}(2) show that $d\Phi_n^{-1*}$ maps the action of ${\mathcal E}^{\Delta}$, ${\mathcal F}^{\Delta}$ on the tangent spaces of $(0,0,\pm 1)^n$ to the coadjoint action of $X$, $Y$ on ${\mathfrak{sl}}(2,\mc)^{*n}$. Post-composing $\Phi_n^{-1*}$ with $\big(\Phi_1^{-1*}\big)^{\otimes n}\circ \Phi_n^*$, Proposition \ref{convexpn}(1) proves that $\Phi^{-1*}$ maps the action of $\Gg$ on the tangent spaces of $(0,0,\pm 1)^n \in \operatorname{Spec}(\mathcal{Z}_0({}_{\rm loc}\Ll_{0,n}^\epsilon))$ onto the coadjoint action of $G$ on ${\mathfrak{sl}}(2,\mc)^{*n} \cong T_{(\pm I,\dots,\pm I)}^*\big(G^0\big)^n$. The result follows, ${\rm Fr}^*$ being an extension of $\Phi^{-1*}$.

(2) This is an integrated version of (1) above. It follows from the arguments used to prove Theorem \ref{QCAteo}$(c)$, by replacing \eqref{coadjointactionproof} with the formulas in Proposition~\ref{convexpn}(2), and Lemma~\ref{Liecommutationslem} with Proposition \ref{convexpn}(1).

(3) It is enough to prove that ${\rm Fr}\circ \mathfrak{c}^{-1*}$ establishes a Poisson isomorphism between $\mathcal{O}(X_G(\Sigma))$ and $\mathcal{Z}_0(\Ll_{0,n}^\epsilon)^{\mathcal{G}}$. By Theorem \ref{scteo}, ${\rm Fr} = \Phi^{-1}_{\vert \mathcal{O}(G^n)} \colon \big(\mathcal{O}(G)^{\otimes n},\{\, ,\, \}_{\rm FR}\big) \ra \big(\mathcal{Z}_0(\Ll_{0,n}^\epsilon),\{\ ,\}_{\rm QCA}\big)$ is an isomorphism of Poisson algebras. By (2) above $Fr$ maps invariants functions to invariant functions. That the bracket $\{\, ,\, \}_{\rm QCA}$ is well-defined on $\mathcal{Z}_0(\Ll_{0,n}^\epsilon)^{\mathcal{G}}$ is an immediate consequence of its definition and the structure of module algebra of $\Ll_{0,n}^\epsilon$. Then the conclusion follows from the last claim of Theorem \ref{FRteo}.

(4) The first claim follows from Theorem \ref{QCAteo}$(d)$, and the facts that $\Phi_n$ is a Poisson isomorphism from $\big(\mathcal{Z}_0({}_{\rm loc}\Ll_{0,n}^\epsilon),\{\ ,\}_{\rm QCA}\big)$ to $\big(\mathcal{Z}_0(U_\epsilon),\{\, ,\, \}_{\rm QCA}\big)^{\otimes n}$ (see the proof of Theorem~\ref{scteo}), and that $\Phi_n^{-1*}$ maps the partial action of $\Gg^{\rm tot}$ on $\operatorname{Spec}(\mathcal{Z}_0(\Ll_{0,n}^\e))$ to the partial action of $\Gg^{n}$ (the $n$-fold direct product) on $\operatorname{Spec}(\mathcal{Z}_0(U_\epsilon))^n$ (by the definition of $\Gg^{\rm tot}$). The second claim follows from the first one in Theorem~\ref{FRteo}.
\end{proof}

\section{Topological formulation}\label{SKEIN}

\subsection{The Wilson loop functor}\label{Wfunctor} Recall that we denote by $\mathcal{C}_A$ the category of $U_A^{\rm res}$-modules of type $1$, and that $\mathcal{C}_A \otimes {\mathbb C}\big[q^{1/D},q^{-1/D}\big]$ is a ribbon category.

Recall also the following notions (see \cite{Tu}). Denote by ${\rm Rib}_{\mathcal{C}_A}$ the category whose morphisms are the isotopy classes rel$(\partial)$ of oriented ribbon graphs in $[0,1]^3$ {\it colored} over $\mathcal{C}_A$ (ie. with each component labelled by an object of $\mathcal{C}_A$), with boundary segments (if any) in $]0,1[\times \{1/2\}\times \{i\}$, where $i\in \{0,1\}$. The objects of ${\rm Rib}_{\mathcal{C}_A}$ are the tuples $((V_1,\varepsilon_1),\dots,(V_k,\varepsilon_k))$, where $\varepsilon_1,\dots,\varepsilon_k=\pm$ and $V_1,\dots,V_k$ are objects of $\mathcal{C}_A$. The source and target objects of a morphism of ${\rm Rib}_{\mathcal{C}_A}$ thus correspond to tuples of segments in $]0,1[\times \{1/2\}\times \{i\}$, $i\in \{0,1\}$, endowed with normal co-orientations specifying the associated signs $\pm$. We denote by $\mathbb{RT}$ the Reshetikhin--Turaev functor
\begin{gather*}
\mathbb{RT}\colon\ {\rm Rib}_{\mathcal{C}_A}\to \mathcal{C}_A \otimes {\mathbb C}\big[q^{1/D},q^{-1/D}\big].
\end{gather*}
Now, fix points $p_1< \cdots <p_n$ in $]0,1/2[$, and define ${\rm Rib}_{n,\mathcal{C}_A}$ as the category with the same objects as ${\rm Rib}_{\mathcal{C}_A}$ but morphisms the isotopy classes rel$(\partial)$ of oriented ribbon graphs in $[0,1]^3\setminus (\{p_1,\dots,p_n\}\times \{1/2\}\times [0,1])$, colored over $\mathcal{C}_A$ and with boundary segments (if any) in $]1/2,1[\allowbreak\times \{1/2\}\times \{i\}$, $i\in\{0,1\}$. Figure~\ref{fig8.1} shows an example.

\def\svgwidth{0.4\textwidth}
\begin{figure}[h!]\centering
\vspace*{-10mm}
\begingroup%
 \makeatletter%
 \providecommand\color[2][]{%
 \errmessage{(Inkscape) Color is used for the text in Inkscape, but the package 'color.sty' is not loaded}%
 \renewcommand\color[2][]{}%
 }%
 \providecommand\transparent[1]{%
 \errmessage{(Inkscape) Transparency is used (non-zero) for the text in Inkscape, but the package 'transparent.sty' is not loaded}%
 \renewcommand\transparent[1]{}%
 }%
 \providecommand\rotatebox[2]{#2}%
 \ifx\svgwidth\undefined%
 \setlength{\unitlength}{595.27559055bp}%
 \ifx\svgscale\undefined%
 \relax%
 \else%
 \setlength{\unitlength}{\unitlength * \real{\svgscale}}%
 \fi%
 \else%
 \setlength{\unitlength}{\svgwidth}%
 \fi%
 \global\let\svgwidth\undefined%
 \global\let\svgscale\undefined%
 \makeatother%
 \begin{picture}(1,1.41428571)%
 \put(0,0){\includegraphics[width=\unitlength,page=1]{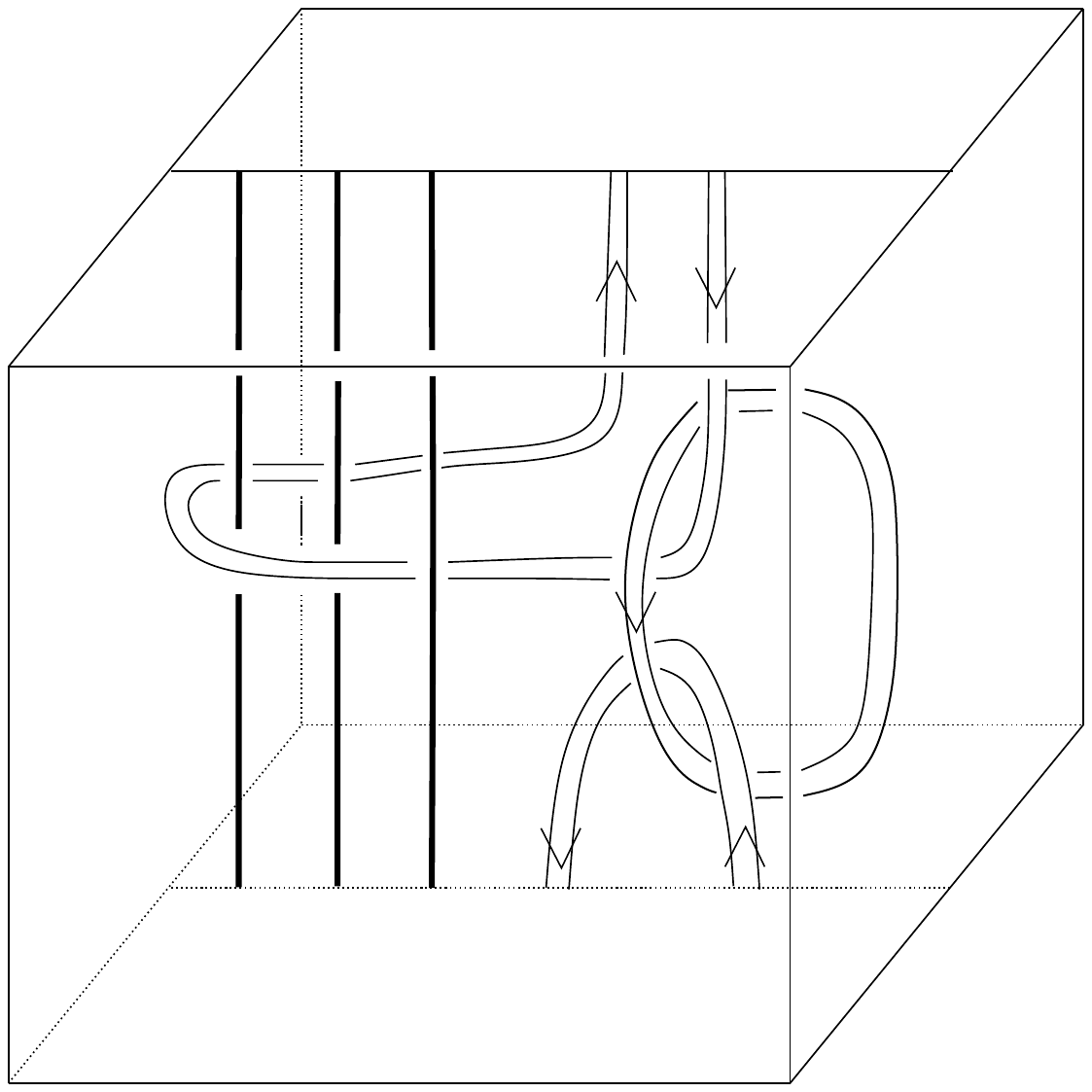}}%
 \put(0.13013142,0.22461044){\color[rgb]{0,0,0}
 \makebox(0,0)[lt]{\begin{minipage}{0.2014363\unitlength}\raggedright \end{minipage}}}%
 \put(1.15870441,0.4064379){\color[rgb]{0,0,0}\makebox(0,0)[lb]{\smash{}}}%
 \put(0.14941896,0.40186688){\color[rgb]{0,0,0}\makebox(0,0)[lb]{\smash{\quad $p_1$}}}%
 \put(0.39259349,0.40186702){\color[rgb]{0,0,0}\makebox(0,0)[lb]{\smash{$p_3$}}}%
 \put(0,0){\includegraphics[width=\unitlength,page=2]{dessin-10.pdf}}%
 \put(0.27100618,0.40186688){\color[rgb]{0,0,0}\makebox(0,0)[lb]{\smash{\ $p_2$}}}%
 \end{picture}%
\endgroup%
\vspace{-18mm}
\caption{A morphism of ${\rm Rib}_{n,\mathcal{C}_A}$.}\label{fig8.1}
\end{figure}
The composition of morphisms of ${\rm Rib}_{n,\mathcal{C}_A}$ is defined as for ${\rm Rib}_{\mathcal{C}_A}$. That is, given morphisms~$T_1$,~$T_2$ with a same pattern of co-oriented boundary segments on the bottom of~$T_1$ and top of~$T_2$, $T_1\circ T_2$ is obtained by placing $T_1$ atop~$T_2$, gluing the corresponding boundary segments, and deforming the result by isotopy into $[0,1]^3\setminus (\{p_1,\dots,p_n\}\times \{1/2\}\times [0,1])$. Identifying $[0,1]^3$ with the ``right half'' cube $[1/2,1]\times [0,1]^2$ in $[0,1]^3\setminus (\{p_1,\dots,p_n\}\times \{1/2\}\times [0,1])$ yields an obvious faithful functor
\[
\iota\colon\ {\rm Rib}_{\mathcal{C}_A}\ra {\rm Rib}_{n,\mathcal{C}_A}.
\]
The morphisms of ${\rm Rib}_{n,\mathcal{C}_A}$ are obtained by composing morphisms of ${\rm Rib}_{\mathcal{C}_A}$ and elementary morphisms as shown in the following picture, for $a=1,\dots,n$.

\def\svgwidth{0.3\textwidth}
\begin{figure}[h!]\centering
\begingroup%
 \makeatletter%
 \providecommand\color[2][]{%
 \errmessage{(Inkscape) Color is used for the text in Inkscape, but the package 'color.sty' is not loaded}%
 \renewcommand\color[2][]{}%
 }%
 \providecommand\transparent[1]{%
 \errmessage{(Inkscape) Transparency is used (non-zero) for the text in Inkscape, but the package 'transparent.sty' is not loaded}%
 \renewcommand\transparent[1]{}%
 }%
 \providecommand\rotatebox[2]{#2}%
 \ifx\svgwidth\undefined%
 \setlength{\unitlength}{595.27559055bp}%
 \ifx\svgscale\undefined%
 \relax%
 \else%
 \setlength{\unitlength}{\unitlength * \real{\svgscale}}%
 \fi%
 \else%
 \setlength{\unitlength}{\svgwidth}%
 \fi%
 \global\let\svgwidth\undefined%
 \global\let\svgscale\undefined%
 \makeatother%
 \begin{picture}(1,1.41428571)%
 \put(0.2015873,0.73075394){\color[rgb]{0,0,0}\makebox(0,0)[lb]{\smash{}}}%
 \put(0.12599206,0.78115082){\color[rgb]{0,0,0}\makebox(0,0)[lb]{\smash{}}}%
 \put(0,0){\includegraphics[width=\unitlength,page=1]{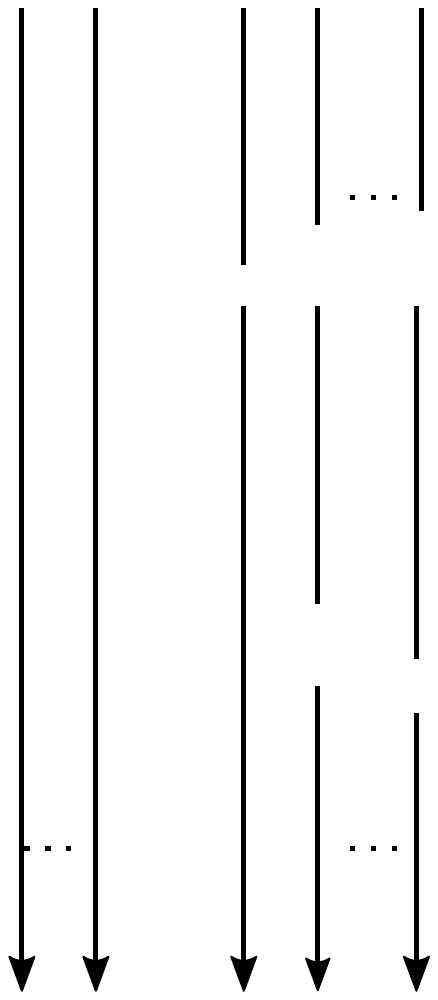}}%
 \put(0.35277779,0.32757932){\color[rgb]{0,0,0}\makebox(0,0)[lb]{\smash{}}}%
 \put(0,0){\includegraphics[width=\unitlength,page=2]{dessin14.pdf}}%
 \put(0.09755224,0.52542946){\color[rgb]{0,0,0}\makebox(0,0)[lb]{\ \smash{$1$}}}%
 \put(0.31066514,0.52833291){\color[rgb]{0,0,0}\makebox(0,0)[lb]{\smash{$a$}}}%
 \put(0.45916005,0.52833291){\color[rgb]{0,0,0}\makebox(0,0)[lb]{\smash{$n$}}}%
 \put(0.5441472,0.52841795){\color[rgb]{0,0,0}\makebox(0,0)[lb]{\smash{$V_1\otimes \cdots \otimes V_k$}}}%
 \put(0.18127401,0.5254293){\color[rgb]{0,0,0}\makebox(0,0)[lb]{\smash{$$}}}%
 \end{picture}%
\endgroup%
\vspace{-25mm}
\caption{A generating morphism of ${\rm Rib}_{n,\mathcal{C}_A}$ not in ${\rm Rib}_{\mathcal{C}_A}$.}\label{fig8.2}
\end{figure}

Recall that we denote by $\Sigma$ the sphere with $n+1$ open disks removed. Denote by $\mathcal{L}_{\mathcal{C}_{A}}(\Sigma)$ the $A$-module freely generated by the empty set and the isotopy classes of closed oriented ribbon graphs in $\Sigma \times [0,1]$ colored over $\mathcal{C}_{A}$. The stacking product $L . L'$ of elements $L, L'\in\mathcal{L}_{\mathcal{C}_{A}}(\Sigma)$ is defined as the isotopy class of the disjoint union of representatives of $L$ and $L'$, pushed in $\Sigma \times [0,1/2]$ and $\Sigma \times [1/2,1]$ respectively. The stacking product makes $\mathcal{L}_{\mathcal{C}_{A}}(\Sigma)$ an associative $A$-algebra. Since $\Sigma$ is diffeomorphic to the closure of $[0,1]^2\setminus (\{p_1,\dots,p_n\}\times \{1/2\})$, $\mathcal{L}_{\mathcal{C}_{A}}(\Sigma)$ can be identified with the $A$-algebra generated by the morphisms $\varnothing \ra \varnothing$ of ${\rm Rib}_{n,\mathcal{C}_A}$, the stacking product being the $A$-linear extension of the composition of morphisms.

Define a new category $\mathcal{C}_{A}\big(\Ll_{0,n}^A\big)$ with objects the couples $(n,V)$, where $V$ is an object of $\mathcal{C}_A\otimes {\mathbb C}[q^{1/D},q^{-1/D}]$, and with spaces of morphisms
\[
\operatorname{Hom}( (n,V), (n,W))=\big(\!\operatorname{Hom}_ {U_A^{\rm res}}(V,W)\otimes {\mathbb C}\big[q^{1/D},q^{-1/D}\big]\big)\otimes \Ll_{0,n}^A,
\]
where the composition is the tensor product of the composition of $U_A^{\rm res}$-module homomorphisms and the product in $\Ll_{0,n}^A$.

The following result follows from the arguments in Sections~6.1 and~6.2 of Faitg's PhD The\-sis~\cite{Faitg}, to which we refer for full details (see also \cite{Faitg4}). These arguments generalize and simplify those of \cite{BR2}. In these works only the Wilson loop map $W$ was considered, but the extension to~$\mathbb{W}$ is straightforward.

\begin{teo} \label{Winv}
There is unique functor ${\mathbb W}\colon {\rm Rib}_{n,\mathcal{C}_A}\rightarrow \mathcal{C}_{A}\big(\Ll_{0,n}^A\big)$ such that on objects we have ${\mathbb W}((V_1,\varepsilon_1),\dots,(V_k,\varepsilon_k)) = \big(n,V_1^{\varepsilon_1}\otimes \dots\otimes V_k^{\varepsilon_k}\big)$, where $V_j^{+}=V_j$ and $V_j^{-} = V_j^*$ $($the dual of $V_j)$, and on morphisms:
\begin{itemize}\itemsep=0pt
\item ${\mathbb W}(\iota(T))=\mathbb{RT}(T) \otimes 1$ for every morphism $T$ of ${\rm Rib}_{\mathcal{C}_A}$,
\item ${\mathbb W}$ gives 
 the value $\stackrel{V}{M}{}^{\!\!(a)}$ to the morphism shown in Figure~$\ref{fig8.2}$, in the case where there is a~single ribbon colored by $V:=V_1$.
\end{itemize}
Moreover the morphism of algebras $W\colon \mathcal{L}_{\mathcal{C}_{A}}(\Sigma) \ra \Ll_{0,n}^A\otimes {\mathbb C}\big[q^{1/D},q^{-1/D}\big]$ obtained by restric\-ting~$\mathbb{W}$ to $\mathcal{L}_{\mathcal{C}_{A}}(\Sigma)$ takes values in the invariant subalgebra $\big(\Ll_{0,n}^A\big)^{U_A}\otimes {\mathbb C}\big[q^{1/D},q^{-1/D}\big]$.
\end{teo}
We call ${\mathbb W}$ the Wilson loop functor, and $W$ the Wilson loop map. By construction, for every element $\hat{T}\in \mathcal{L}_{\mathcal{C}_{A}}(\Sigma)$, and any morphism $T$ of ${\rm Rib}_{n,\mathcal{C}_A}$ obtained from $\hat{T}$ by cutting open ribbons with colors $V_1,\dots,V_k$, setting $V=V_1\otimes \dots \otimes V_k$ we have
\begin{equation}\label{Trhol}
W(\hat{T}) = \operatorname{qTr}_V\left( {\mathbb W}(T)\right).
\end{equation}
Here is an alternative way of defining ${\mathbb W}$, purely in terms of the Reshetikhin--Turaev functor $\mathbb{RT}$. Define a category $\mathcal{C}_{A}\big(\tilde U_A^{\otimes n}\big)$ with same objects $(n,V)$ as $\mathcal{C}_{A}\big(\Ll_{0,n}^A\big)$ but with spaces of morphisms
\[
\operatorname{Hom}( (n,V), (n,W))=\big( \operatorname{Hom}_ {U_A^{\rm res}}(V,W)\otimes {\mathbb C}\big[q^{1/D},q^{-1/D}\big]\big) \otimes \tilde U_A^{\otimes n},
\]
where the composition is the tensor product of the composition of $U_A^{\rm res}$-module homomorphisms and the product in $\tilde U_A^{\otimes n}$.

The Alekseev map being a morphism of algebra, it defines a functor
\[
\Phi_n\colon\ \mathcal{C}_{A}\big(\Ll_{0,n}^A\big)\rightarrow {\mathcal{C}}_{A}\big(\tilde U_A^{\otimes n}\big)
\] by setting $\Phi_n(f\otimes a)=f\otimes \Phi_n(a)$ for every $a\in \Ll_{0,n}^A$, $f \in \operatorname{Hom}_ {U_A^{\rm res}}(V,W)$. To any morphism $T$ of ${\rm Rib}_{n,\mathcal{C}_A}$ we can associate a colored oriented ribbon graph in $[0,1]^3$,
\[
T^{\sharp}= \bigg(\bigcup_{i=1}^n [p_i-\delta,p_i+\delta]\times \{1/2\} \times [0,1]\bigg)_{\tilde U_A} \cup T,
\]
where $\delta>0$ is small and the cores of the ribbons $[p_i-\delta,p_i+\delta]\times \{1/2\} \times [0,1]$ are oriented from $\{p_i\}\times \{1/2\} \times \{1\}$ to $\{p_i\}\times \{1/2\} \times \{0\} $, and the subscript $\tilde U_A$ means that these ribbons are colored by the regular representation of $\tilde U_A$. Denote by $\overline{\rm Rib}_{\mathcal{C}_A}$ the category with morphisms given by the ribbon graphs $T^\sharp$; clearly the map $T\mapsto T^\sharp$ yields a functor
\[
\sharp\colon\ {\rm Rib}_{n,\mathcal{C}_A}\ra \overline{\rm Rib}_{n,\mathcal{C}_A}.
\]
By associating $\sigma \circ ({\rm id}\otimes \pi_{V})\big(R^{\pm}\big)\colon \operatorname{End}(V)\otimes \tilde U_A \ra \tilde U_A \otimes \operatorname{End}(V)$ to crossings of index $\pm 1$ colored over $\mathcal{C}_{A}$ and $\tilde{U}_A$, in the same way as $\sigma \circ (\pi_{W}\otimes \pi_{V})(R^{\pm 1})$ is associated by $\mathbb{RT}$ to crossings colored by objects~$V$,~$W$ of~$\mathcal{C}_{A}$, we extend~$\mathbb{RT}$ to a functor
\[
\mathbb{RT}\colon\ \overline{\rm Rib}_{n,\mathcal{C}_A} \rightarrow \mathcal{C}_{A}\big(\tilde U_A^{\otimes n}\big).
\]
The next fact is a direct consequence of the definition of $\Phi_n$, and its representation by Figure~\ref{fig6.1}. Note that $\Phi_n$ being injective, it defines ${\mathbb W}$ uniquely from $\mathbb{RT}$.

\begin{prop} \label{commutRT} We have a commutative diagram of functors:
\[
\xymatrix{{\rm Rib}_{n,\mathcal{C}_A} \ar[r]^{{}^\sharp} \ar[d]_{{\mathbb W}} & \overline{\rm Rib}_{n,\mathcal{C}_A} \ar[d]^{\mathbb{RT}} \\ \mathcal{C}_{A}\big(\Ll_{0,n}^A\big) \ar[r]^{\Phi_n} & {\mathcal{C}}_{A}\big(\tilde U_A^{\otimes n}\big).}
\]
\end{prop}

\subsection{The Wilson loop isomorphism}\label{Wiso} One can adjoin to $\mathcal{L}_{\mathcal{C}_{A}}(\Sigma)$ the $A$-span of closed oriented ribbon graphs in $\Sigma \times [0,1]$ colored over~$\mathcal{C}_{A}$ {\it and possibly having coupons} (see \cite{Tu} for this notion). Denote the resulting set by $\mathcal{L}_{\mathcal{C}_{A}}(\Sigma)^c$. It is an $A$-algebra for the stacking product, containing $\mathcal{L}_{\mathcal{C}_{A}}(\Sigma)$ as a subalgebra. Define a {\it one-coupon multicurve} as an element of $\mathcal{L}_{\mathcal{C}_{A}}(\Sigma)^c$ that can be represented by a (oriented, $\mathcal{C}_A$-colored) ribbon graph embedded in $\Sigma\times \{0\}$ (whence edges are unlinked and unknotted) and having a single coupon. An example of one-coupon multicurve with the coupon colored by $a_X$ is shown in Figure~\ref{fig8.3}, where $\Sigma$ is the four-holed sphere (i.e., $n=3$).

\def\svgwidth{0.3\textwidth}
\begin{figure}[h!]\centering
\begingroup%
 \makeatletter%
 \providecommand\color[2][]{%
 \errmessage{(Inkscape) Color is used for the text in Inkscape, but the package 'color.sty' is not loaded}%
 \renewcommand\color[2][]{}%
 }%
 \providecommand\transparent[1]{%
 \errmessage{(Inkscape) Transparency is used (non-zero) for the text in Inkscape, but the package 'transparent.sty' is not loaded}%
 \renewcommand\transparent[1]{}%
 }%
 \providecommand\rotatebox[2]{#2}%
 \ifx\svgwidth\undefined%
 \setlength{\unitlength}{595.27559055bp}%
 \ifx\svgscale\undefined%
 \relax%
 \else%
 \setlength{\unitlength}{\unitlength * \real{\svgscale}}%
 \fi%
 \else%
 \setlength{\unitlength}{\svgwidth}%
 \fi%
 \global\let\svgwidth\undefined%
 \global\let\svgscale\undefined%
 \makeatother%
 \begin{picture}(1,1.41428571)%
 \put(0,0){\includegraphics[width=\unitlength,page=1]{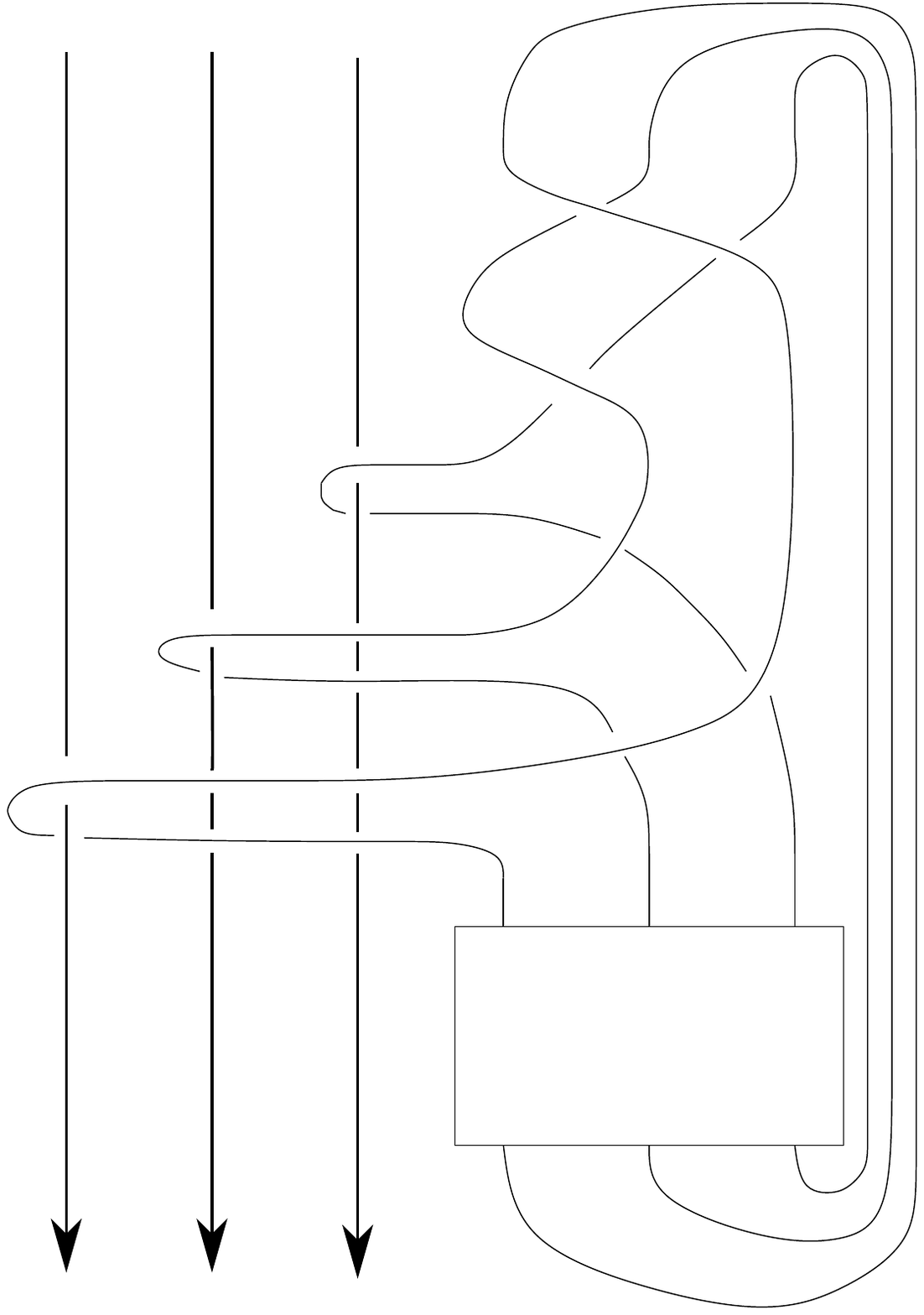}}%
 \put(0.54673843,0.39054486){\color[rgb]{0,0,0}\makebox(0,0)[lb]{\ \ \smash{$a_X$}}}%
 \end{picture}%
\endgroup%
\vspace{-7mm}
\caption{A one-coupon multicurve in $\Sigma_{0,4}\times [0,1]$.}\label{fig8.3}
\end{figure}
Recall the matrices $\stackrel{X}{\mathbb{M}}$ and the invariant elements $v_{X}(a_{X})$ defined for any $U_A^{\rm res}$-module \mbox{$X=X_1\otimes \dots \otimes X_n$} of type $1$ and any $a_{X} \in \operatorname{End}_{U_A^{\rm res}}(X)$ (see the comments after Proposition~\ref{linbaseinv}). Denote
by $L(a_X)$ the one-coupon multicurve in Figure~\ref{fig8.3}; so $n=3$. By \eqref{Trhol} it is clear that
\begin{equation}\label{onecoup}
W(L(a_{X})) = v_{X}(a_{X}).
\end{equation}
This generalizes immediately to any $n\geq 1$. By the comments after Proposition \ref{linbaseinv} elements of the form $v_{X}\big(a_{X}^{(k)}\big)$ form a basis of $\Ll_{0,n}^{U_q}$. So Theorem \ref{Winv} gives \big(with coefficients in $q^{1/D}$\big):
\begin{lem}\label{Wsurj}
The Wilson loop map $W\colon \mathcal{L}_{\mathcal{C}_{A}}(\Sigma)^c \otimes_A \mc\big(q^{1/D}\big) \ra \Ll_{0,n}^{U_q}\otimes_{\mc(q)} \mc\big(q^{1/D}\big)$ is surjective.
\end{lem}

In the sequel we assume that $\mathfrak{g}={\mathfrak{sl}}(2)$. We are going to see that in this case we can strengthen Lemma \ref{Wsurj}.

Set $\zeta := {\rm i}q^{1/2}$, and denote by $\mathcal{L}_\zeta(\Sigma)\subset \mathcal{L}_{\mathcal{C}_{A}}(\Sigma)\otimes {\mathbb C}\big[\zeta,\zeta^{-1}\big]$ the subalgebra freely generated as a $\mz\big[\zeta,\zeta^{-1}\big]$-module by the empty set and the isotopy classes of oriented ribbon links in~$\Sigma \times [0,1]$ colored by the fundamental representation~$V_2$. Recall that the Kauffman bracket skein algebra~$K_\zeta(\Sigma)$ is the $\mz[\zeta,\zeta^{-1}]$-algebra obtained from $\mathcal{L}_\zeta(\Sigma)$ by forgetting the link orientations, and taking the quotient by the ideal generated by the relations
\begin{gather}
L=\zeta L_+ +\zeta^{-1}L_-,\label{relsk1}
\\
L\sqcup \bigcirc=-\big(\zeta^2+\zeta^{-2}\big)L,
\label{relsk2}
\end{gather}
where in the second identity $ \bigcirc$ is the trivial ribbon link in a ball disjoint from $L$, and in the first identity $L,L_+, L_-\in \mathcal{L}_\zeta(\Sigma)$ are identical up to isotopy except in a ball in which they look like (the strands representing flat ribbons on a meridional projection disk):
\def\svgwidth{0.3\textwidth}
\begin{figure}[h!]\centering
\vspace*{1mm}
\begingroup%
 \makeatletter%
 \providecommand\color[2][]{%
 \errmessage{(Inkscape) Color is used for the text in Inkscape, but the package 'color.sty' is not loaded}%
 \renewcommand\color[2][]{}%
 }%
 \providecommand\transparent[1]{%
 \errmessage{(Inkscape) Transparency is used (non-zero) for the text in Inkscape, but the package 'transparent.sty' is not loaded}%
 \renewcommand\transparent[1]{}%
 }%
 \providecommand\rotatebox[2]{#2}%
 \ifx\svgwidth\undefined%
 \setlength{\unitlength}{255.97274719bp}%
 \ifx\svgscale\undefined%
 \relax%
 \else%
 \setlength{\unitlength}{\unitlength * \real{\svgscale}}%
 \fi%
 \else%
 \setlength{\unitlength}{\svgwidth}%
 \fi%
 \global\let\svgwidth\undefined%
 \global\let\svgscale\undefined%
 \makeatother%
 \begin{picture}(1,0.41302497)%
 \put(0,0){\includegraphics[width=\unitlength,page=1]{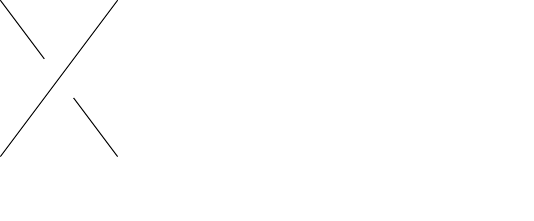}}%
 \put(0.11060748,0.00960071){\color[rgb]{0,0,0}\makebox(0,0)[lb]{\smash{$L$}}}%
 \put(0.45069672,0.02006506){\color[rgb]{0,0,0}\makebox(0,0)[lb]{\smash{$L_+$}}}%
 \put(0.81694657,0.02006506){\color[rgb]{0,0,0}\makebox(0,0)[lb]{\smash{$L_-$}}}%
 \put(0,0){\includegraphics[width=\unitlength,page=2]{dessin13.pdf}}%
 \end{picture}%
\endgroup%
\vspace{-2 mm}
\caption{Skein related ribbon links.}
\end{figure}

The following result is \cite[Theorem~1]{BFK2} (see also~\cite[Theorem~10]{BFK}). For completeness we give a proof by using the notions we have introduced.
\begin{teo}\label{teoCharlie}
The linear map defined by $\mathcal{W}(L) = {\rm i}^{{\rm lk}(L)}W(L)$, where $L\in \mathcal{L}_\zeta(\Sigma)$ and ${\rm lk}(L)$ is the linking number of $L$, descends to an isomorphism of $\mc(\zeta)$-algebras $\big($where $\zeta := {\rm i}q^{1/2}\big)$:
\[
\mathcal{W}\colon\ K_\zeta(\Sigma)\otimes \mc(\zeta)\ra \Ll_{0,n}^{U_q}\otimes \mc(\zeta).
\]
\end{teo}

\begin{proof}
The facts that the values of $\mathcal{W}$ do not depend on the link orientations, and that $\mathcal{W}$ maps to $0$ the ideal of $\mathcal{L}_\zeta(\Sigma)$ generated by the skein relations \eqref{relsk1} and~\eqref{relsk2}, follow from Proposition~\ref{commutRT}, the injectivity of $\Phi_n$, and the properties of $\mathbb{RT}$ and the $R$-matrix $R_{V_2,V_2}$ proved, e.g., in \cite[Lemma 3.18, Theorem 4.3 and Corollary 4.13]{KM}. We stress that, using as usual the standard pivotal element $K$ in the definition of $\mathbb{RT}$, the choice of variable $\zeta := {\rm i}q^{1/2}$ and the normalisation of the map $W$ by $i^{{\rm lk}(L)}$ are necessary for the skein relations to be in the kernel of~$\mathcal{W}$. Therefore $\mathcal{W}$ descends to a module map defined on $K_\zeta(\Sigma)$. Since the stacking product is induced from $\mathcal{L}_\zeta(\Sigma)$, $\mathcal{W}$ is a morphism of algebras.

To show that $\mathcal{W}$ is an isomorphism we use the following basis of $K_\zeta(\Sigma)\otimes \mc(\zeta)$ (see, e.g., \cite{KLins,Lick}). Let $\Gamma$ be a trivalent ribbon graph onto which $\Sigma$ retracts by deformation. Recall that an admissible coloring of $\Gamma$ is an assignment of a nonnegative integer to each edge, called the color of the edge, so that the colors adjacent to each vertex form admissible triples. A triple of colors $(a,b,c)$ is admissible if $a\leq b+c$, $b \leq a+c$, $c \leq a+b$ and $a+b+c$ is even. Admissible colorings $\gamma$ of $\Gamma$ parametrize multicurves carried by $\Gamma$, a color on an edge $e$ giving the number of components running parallel to $e$. Given any admissible coloring $\gamma$ of $\Gamma$, there is a~skein~$\Gamma_\gamma$ in $K_\zeta(\Sigma)\otimes \mc(\zeta)$ corresponding to $\gamma$, obtained by replacing an edge with color $m$ by the $m$-th Jones--Wenzl idempotent and vertices with Kauffman triads (these being defined in, e.g., \cite[p.~136 and Figure~14.7]{Lick}). The set of multicurves on $\Sigma$ forms a basis of $K_\zeta(\Sigma)$; since $\Gamma$ is a~spine of~$\Sigma$, the set of skeins $\Gamma_\gamma$ forms a basis of $K_\zeta(\Sigma)\otimes \mc(\zeta)$ as well. Note that, because the Jones--Wenzl idempotents have denominators, the skeins $\Gamma_\gamma$ must be defined over $\mc(\zeta)$.

By an isotopy of $\Sigma\times [0,1]\cong [0,1]^3\setminus (\{p_1,\dots,p_n\}\times \{1/2\}\times [0,1])$ we can deform $\Gamma$ (keeping the same notation) so that a neighborhood of the vertices lies inside a coupon embedded in $]1/2,1[\times \{1/2\} \times [0,1]$, and the $n$ ribbons forming the portion $\Gamma'$ of $\Gamma$ outside the coupon are attached to its left side, encircling $\{p_i\}\times \{1/2\} \times [0,1]$, for $i=1,\dots,n$. Let $\gamma$ be an admissible coloring of $\Gamma$. Put $\Gamma'' = \Gamma \setminus \Gamma'$, and denote by $\Gamma'_\gamma$, $\Gamma''_\gamma$ the graphs $\Gamma'$, $\Gamma''$ with edges colored by~$\gamma$. Denote by $e_1,\dots,e_n$ the edges of $\Gamma'$ ordered by increasing height, by $\gamma_i$ the color of $e_i$ plus $1$, and fix some orientation of the edges of $\Gamma$ so that each $e_i$ is oriented from top to bottom. In~order to fit with the standard framework for computations with $\mathbb{RT}$, let us rotate the cube $[0,1]^3$ clockwise by an angle of $\pi/2$ around the axis $\{1/2\}\times [0,1]\times \{1/2\}$. Then, using Theorem \ref{Winv} and the fact that the $m$-th Jones--Wenzl idempotent projects $V_2^{\otimes m}$ onto a subspace isomorphic to $V_{m+1}$, we see that
\[
\mathbb{W}(\Gamma'_\gamma) = \stackrel{V_{\gamma_1}}{M}{}^{\!\!(1)}\otimes \dots \otimes \stackrel{V_{\gamma_n}}{M}{}^{\!\!(n)},
\]
which is an element of $\otimes_{i=1}^n \big(V_{\gamma_i}\otimes V_{\gamma_i}^*\big)\otimes \Ll_{0,n}$ by identifying $\operatorname{End}(V)$ with $V\otimes V^*$, and
\[
\mathbb{W}(\Gamma''_\gamma) \in \operatorname{Hom}_{U_q}\big(\otimes_{i=1}^n \big(V_{\gamma_i}\otimes V_{\gamma_i}^*\big),1\big) \otimes_{\mc(q)} \mc(\zeta).
\]
By definitions we have $\mathcal{W}(\Gamma_\gamma) = \mathbb{W}(\Gamma''_\gamma)(\mathbb{W}(\Gamma'_\gamma))$, and the map $\gamma \mapsto \mathbb{W}(\Gamma''_\gamma)$ maps the admissible colorings of $\Gamma$ to a basis of the space of invariant elements of $\otimes_{i=n}^1 \operatorname{End}(V_{\gamma_i})^*\otimes_{\mc(q)} \mc(\zeta)$ for the action $(\operatorname{ad}^r)^{\otimes n}$. Therefore the set $\{\mathcal{W}(\Gamma_\gamma),\gamma\ {\rm admissible}\}$ is a basis of the space of invariant elements of $\Ll_{0,1}^{\otimes n}\otimes_{\mc(q)} \mc(\zeta)$ for the action $(\operatorname{coad}^r)^{\otimes n}$. By applying to it the linear isomorphism $\Phi_n^{-1}\circ \psi\circ \Phi_1^{\otimes n}\colon \Ll_{0,1}^{\otimes n} \ra \Ll_{0,n}$, with $\psi$ the intertwiner of $(\operatorname{ad}^r)^{\otimes n}$ and $\operatorname{ad}^r_n$ defined in the proof of~Pro\-position~\ref{surjPhiinv}, we get a basis of $\Ll_{0,n}^{U_q}\otimes_{\mc(q)} \mc(\zeta)$. This concludes the proof.
\end{proof}

\begin{remk}{\rm Above we can replace $\mathbb{W}(\Gamma'_\gamma)$ with $\stackrel{[\lambda]}{\mathbb{M}}$, where $[\lambda] = (\gamma_1,\dots,\gamma_n)$, by ``parsing'' the ends of $e_1,\dots,e_n$ as do the strands just above the coupon $a_X:=a_{[\lambda]}$ in Figure~\ref{fig8.3}. Cor\-res\-pon\-dingly $\mathbb{W}(\Gamma''_\gamma)$ becomes an element of $\operatorname{Hom}_{U_q}\big(1, (\otimes_{i=1}^n V_{\gamma_i})\otimes (\otimes_{i=1}^n V_{\gamma_i})^*\big) \otimes_{\mc(q)} \mc(\zeta)$. With these choices, the pairing of $\mathbb{W}(\Gamma'_\gamma)$ and $\mathbb{W}(\Gamma''_\gamma)$ is equivariant with respect to the actions $\operatorname{coad}^r_n$ and~$\operatorname{ad}^r_n$, and it puts in duality the product of $\Ll_{0,n}$ and the comultiplication map constructed in \cite{BFK}.}
\end{remk}

The following result is a version of Theorem \ref{teoCharlie} over the ring $\mc\big[\zeta,\zeta^{-1}\big]$:
\begin{teo}\label{teoOBS} We have an isomorphism of algebras $\mathcal{W}\colon K_\zeta(\Sigma)\ra \big(\Ll_{0,n}^A\big)^{U_A}\otimes_A \mc\big[\zeta,\zeta^{-1}\big]$.
\end{teo}

\begin{proof} We are going to use stated skein algebras and a result of Faitg (see \cite{Faitg4}). Denote by $\Sigma_{0,n}^{\circ,\bullet}$ our surface $\Sigma$ with one point removed on its boundary, by $t^{1/2}$ an indeterminate, and by $\mathcal{S}_t^s\big(\Sigma_{0,n}^{\circ,\bullet}\big)$ the stated skein algebra of $\Sigma_{0,n}^{\circ,\bullet}$, as defined in~\cite{Le2}. It is an algebra over $\mc[t^{1/2},t^{-1/2}]$, which contains $K_t(\Sigma)$ as a $\mz\big[t,t^{-1}\big]$-subalgebra. By a result of~\cite{CL}, $\mathcal{S}_t^s\big(\Sigma_{0,n}^{\circ,\bullet}\big)$ is a $\Oo_{t^2}$-comodule algebra, isomorphic to the braided tensor product of $n$ copies of $\mathcal{S}_t^s\big(\Sigma_{0,1}^{\circ,\bullet}\big)$. Here $\Oo_{t^2}$ is $\Oo_A$ with variable $q=t^2$; by the duality between comodule algebras and module algebras we can canonically regard $\mathcal{S}_t^s\big(\Sigma_{0,n}^{\circ,\bullet}\big)$ as a $U_A^{\rm res}$-module algebra, whence a $U_A$-module algebra. As coefficients contain~$\mc$ we can replace $t^{1/2}$ by $\zeta^{1/2}$ so that $\zeta= {\rm i}t$. On the other side, by the integral form of the last claim of Proposition \ref{defiL0n}, which follows immediately from Proposition \ref{defiL0nA} and Lemma \ref{intform-n}, we know that $\Ll_{0,n}^A$ is isomorphic to the braided tensor product of $n$ copies of the $U_A$-module algebra $\Ll_{0,1}^A$. In~\cite[Theorem 5.3]{Faitg4} it is constructed an isomorphism of $U_A$-module algebras $\tilde{\mathcal{W}}\colon \mathcal{S}_\zeta^s\big(\Sigma_{0,n}^{\circ,\bullet}\big)\ra \Ll_{0,n}^A\otimes_A \mc\big[\zeta^{1/2},\zeta^{-1/2}\big]$. The restriction of $\tilde{\mathcal{W}}$ on the subalgebra $K_\zeta(\Sigma)$ is just the Wilson loop map $\mathcal{W}$. To obtain $\tilde{\mathcal{W}}$, an explicit isomorphism is constructed between $\mathcal{S}_\zeta^s\big(\Sigma_{0,1}^{\circ,\bullet}\big)$ and $\Ll_{0,1}^A\otimes_A \mc\big[\zeta^{1/2},\zeta^{-1/2}\big]$ (see \cite[Lemmas~5.7 and~5.8]{Faitg4}).

Let us consider the isomorphism $\mathcal{S}_\zeta^s\big(\Sigma_{0,n}^{\circ,\bullet}\big)^{U_A} \ra \big(\Ll_{0,n}^A\big)^{U_A} \otimes_A \mc\big[\zeta^{1/2},\zeta^{-1/2}\big]$ induced by $\tilde{\mathcal{W}}$. Since $\big(\Ll_{0,n}^A\big)^{U_A} \otimes_A \mc\big[\zeta^{1/2},\zeta^{-1/2}\big] = \big(\Ll_{0,n}^{U_q}\otimes \mc\big(\zeta^{1/2}\big)\big)\cap \big(\Ll_{0,n}^A\otimes_A \mc\big[\zeta^{1/2},\zeta^{-1/2}\big]\big)$, we have
\begin{align*}
\mathcal{S}_\zeta^s\big(\Sigma_{0,n}^{\circ,\bullet}\big)^{U_A} & = \tilde{\mathcal{W}}^{-1}\big(\big(\Ll_{0,n}^A\big)^{U_A} \otimes_A \mc\big[\zeta^{1/2},\zeta^{-1/2}\big]\big) \\ & = \tilde{\mathcal{W}}^{-1}\big(\Ll_{0,n}^{U_q}\otimes \mc\big(\zeta^{1/2}\big)\big) \cap \mathcal{S}_\zeta^s\big(\Sigma_{0,n}^{\circ,\bullet}\big) = \big(K_\zeta(\Sigma)\otimes \mc\big(\zeta^{1/2}\big)\big) \cap \mathcal{S}_\zeta^s\big(\Sigma_{0,n}^{\circ,\bullet}\big)
\end{align*}
by Theorem \ref{teoCharlie}, which is just $K_\zeta(\Sigma)\otimes \mc\big[\zeta^{1/2},\zeta^{-1/2}\big]$. This proves the extension of $\mathcal{W}$ to the scalars $\mc\big[\zeta^{1/2},\zeta^{-1/2}\big]$ is an isomorphism $K_\zeta(\Sigma) \otimes \mc\big[\zeta^{1/2},\zeta^{-1/2}\big]\ra \big(\Ll_{0,n}^A\big)^{U_A}\otimes_A \otimes \mc\big[\zeta^{1/2},\zeta^{-1/2}\big]$. Therefore $\mathcal{W}$ is an isomorphism as well. This concludes the proof. \end{proof}

\begin{remk}\label{teoOBSrem} \quad
\begin{enumerate}\itemsep=0pt

\item[$(1)$] Faitg results discussed above hold for surfaces of arbitrary genus. Also, in order to have at hands a calculus based on non oriented diagrams he uses a pivotal element different from~$K$ to define the map $W$ (whence $\mathbb{RT}$) (see~\cite[Remark~5.2]{Faitg4}). This leads to the choice of $t$ in his constructions, instead of $\zeta={\rm i}q^{1/2}$ as we obtained in the proof of Theorem \ref{teoCharlie}.

\item[$(2)$] As a by-product of the proof we have $K_\zeta(\Sigma) \otimes \mc\big[\zeta^{1/2},\zeta^{-1/2}\big] = \mathcal{S}_\zeta^s\big(\Sigma_{0,n}^{\circ,\bullet}\big)^{U_A}$.

\item[$(3)$] The module algebra $\mathcal{S}_\zeta^s\big(\Sigma_{0,n}^{\circ,\bullet}\big)$ is defined over $\mc\big[\zeta^{1/2},\zeta^{-1/2}\big]$, and $K_\zeta(\Sigma)$ over $\mz\big[\zeta,\zeta^{-1}\big]$. One can observe that $\Ll_{0,n}^A$ is in fact defined over $\mz\big[q,q^{-1}\big]$, for arbitrary $\mathfrak{g}$, and Faitg's isomorphism $\tilde{\mathcal{W}}$ holds true over $\mz\big[\zeta^{1/2},\zeta^{-1/2}\big]$. Therefore $\mathcal{W}$ in Theorem \ref{teoOBS} holds true over $\mz\big[\zeta,\zeta^{-1}\big]$.
\end{enumerate}
\end{remk}

\subsection[The threading map and G-invariant central elements]{The threading map and $\boldsymbol{\Gg}$-invariant central elements}

A {\it multicurve} on $\Sigma$ is a union of disjoint simple non trivial (ie. not bounding a disk) closed curves considered up to isotopy. Any multicurve $\gamma$ is a stacking product $\textstyle \prod_{i=1}^k \gamma_i^{c_i}$, where $\gamma_1,\dots,\gamma_k$ are disjoint (hence commuting) non-isotopic simple non trivial closed curves on $\Sigma$, and $\gamma_i^{c_i}$, $c_i\in \mathbb{N}$, consists of $c_i$ parallel copies of $\gamma_i$.

Denote by $\mathcal{S}(\Sigma)$ the set of multicurves on $\Sigma$. Recall the normalized Chebyshev polynomials~$T_k$, $k\in \mathbb{N}$, defined in \eqref{ChPol}. For $\gamma\in \mathcal{S}(\Sigma)$, $\gamma = \prod_{i=1}^k \gamma_i^{c_i}$, set
\[
T(\gamma):= \prod_{i=1}^k T_{c_i}(\gamma_i).
\]
It is standard that $\mathcal{S}(\Sigma)$ is a $\mz\big[\zeta,\zeta^{-1}\big]$-basis of $K_\zeta(\Sigma)$ (see, e.g., \cite[Theorem~7]{PS}) and that~$\{T_k\}_k$ is a basis of $\mz[X]$. Therefore $\{T(\gamma), \gamma\in \mathcal{S}(\Sigma)\}$ is a~$\mz\big[\zeta,\zeta^{-1}\big]$-basis of $K_\zeta(\Sigma)$. It is called the {\it Chebyshev basis} (see~\cite{FKL2} and the references therein). Theorem~\ref{teoOBS} implies that $\{\mathcal{W}(T(\gamma)), \gamma\in\mathcal{S}(\Sigma)\}$ is a~$\mz\big[\zeta,\zeta^{-1}\big]$-basis of $\big(\Ll_{0,n}^A\big)^{U_A}\otimes_A \mz\big[\zeta,\zeta^{-1}\big]$.

Now let as usual $l\geq 3$ be an odd integer. Define $\mathcal{S}_l(\Sigma)\subset \mathcal{S}(\Sigma)$ as the set of multicurves of the form $\gamma_\partial\prod_{i=1}^k \gamma_i^{c_i}$, where $\gamma_\partial$ is peripheral, i.e., a monomial in the skein classes of the boundary components of $\Sigma$, $\gamma_i$ is a non peripheral curve, and $l$ divides $c_i$ for every $i\in \{1,\dots,k\}$.

Let $\epsilon$ be a primitive root of unity of odd order $l$. Recall the element $\eta$ in \eqref{eta}, the specialization $\big(\Ll_{0,n}^A\big)^{U_A}_\e$ in \eqref{AepsL0n}, and the algebra $\mathcal{Z}(\Ll_{0,n}^\e)^\Gg$ and the derivations $D_a\colon \Ll_{0,n}^\e \ra \Ll_{0,n}^\e$, $a\in \mathcal{Z}(\Ll_{0,n}^\e)$, defined in Section~\ref{EXTENSION}.

\begin{teo}\label{centerprop1}\quad\samepage

\begin{enumerate}\itemsep=0pt

\item[$(1)$] The set $\{\mathcal{W}(T(\gamma)), \gamma\in\mathcal{S}_l(\Sigma)\}$ is a $\mc$-basis of the algebra generated by $\mathcal{Z}(\Ll_{0,n}^\e)^\Gg$ and $\eta$. In~par\-ti\-cu\-lar, it is a central subalgebra of $\big(\Ll_{0,n}^A\big)^{U_A}_\e$.

\item[$(2)$] The derivations $D_a$, $a\in \mathcal{Z}(\Ll_{0,n}^\e)^\Gg$, act on $\big(\Ll_{0,n}^A\big)^{U_A}_\e$.
 \end{enumerate}
\end{teo}

\begin{proof} (1) First note that for every $\gamma\in\mathcal{S}(\Sigma)$, $\textstyle \gamma = \gamma_\partial\prod_{i=1}^k \gamma_i^{c_i}$ with $\gamma_\partial$ peripheral, we have
\begin{equation}\label{exprWT}
\mathcal{W}(T(\gamma)) = \mathcal{W}(T(\gamma_\partial))\prod_{i=1}^k T_{c_i}(\mathcal{W}(\gamma_i)).
\end{equation}
Denote by $\partial_1,\dots,\partial_{n+1}$ the boundary components of $\Sigma$, ordered so that the diffeomorphism identifying $\Sigma$ with the closure of $[0,1]^2\setminus (\{p_1,\dots,p_n\}\times \{1/2\})$ maps $\partial_{n+1}$ to $\partial\big( [0,1]^2\big)$, and $\partial_{i}$ to a small loop encircling $(p_i,1/2)$, $1\leq i\leq n$. By~\eqref{Trhol} we have
\begin{equation}\label{formcentq}
\mathcal{W}(\partial_{n+1}) = \operatorname{qTr}\big(\!\stackrel{V_2}{M}{}^{\!\! (1)}\cdots \stackrel{V_2}{M}{}^{\!\! (n)} \big) = \eta,\ \mathcal{W}(\partial_i) = \operatorname{qTr}\big(\!\stackrel{V_2}{M}{}^{\!\! (i)}\big) = \omega^{(i)}, \qquad 1\leq i\leq n.
\end{equation}
Therefore the elements $\mathcal{W}(T(\gamma_\partial))$ form a basis of $\mc\big[\omega^{(1)},\dots,\omega^{(n)},\eta\big]$, which is a central subalgebra of $\big(\Ll_{0,n}^A\big)^{U_A}_\e$ (see Theorem \ref{CenterLinv}). We can give formulas of the other terms as follows. Deform a ribbon neighborhood of $\gamma$ in $\Sigma\times [0,1]$ so that for every $i\in \{1,\dots,k\}$, $\gamma_i$ is represented by a~one-coupon multicurve like in Figure~\ref{fig8.3}, the coupon being filled with~$r$ simple arcs (possibly pairwise intersecting), for some $r\in \mn$. Since $\gamma_i$ is colored by $V_2$, we see that \eqref{onecoup} takes the form
\begin{gather*}
\mathcal{W}(\gamma_i) = \operatorname{qTr}_{V_2^{\otimes r}} \big( a(\gamma_i) \mathbb M(\gamma_i)\big),
\end{gather*}
where $a(\gamma_i)\in \operatorname{End}_{U_A^{\rm res}}\big(V_2^{\otimes r}\big)$, and $\mathbb M(\gamma_i)$ is defined as $\stackrel{[\lambda]}{\mathbb{M}}$, replacing the modules $V_{\lambda_1},\dots,V_{\lambda_n}$ associated to the sequence $[\lambda]=(\lambda_1,\dots,\lambda_n)$ by $V_2$ or the trivial module $V_1$ (this latter case happens when the module labels a strand which has to be removed to get the one-coupon multicurve representing $\gamma_i$). In order to get simultaneously a simple expression of all the elements $\mathcal{W}(\gamma_1),\dots, \mathcal{W}(\gamma_k)$, note that each $\gamma_i$, being a simple closed curve, separates $\Sigma$ in two disks with punctures. We can choose the above diffeomorphism of $\Sigma$ with the closure of $[0,1]^2\setminus (\{p_1,\dots,p_n\}\times \{1/2\})$ so that the ordering of $p_1,\dots,p_n$ makes $\gamma_i$ bounding a punctured disk $D_i\subset [0,1]^2$ with successive punctures $p_{j_i}, p_{j_i +1},\dots, p_{j_i + n_i-1}$. Then, for every $i\in \{1,\dots,k\}$ we have
\begin{gather*}
\mathcal{W}(\gamma_i) = \operatorname{qTr}_{V_2} \big(\! \stackrel{V_2}{M}{}^{\!\! (j_i)}\cdots \stackrel{V_2}{M}{}^{\!\! (j_i + n_i-1)}\big).
\end{gather*}
The algebras $\big(\Ll_{0,n}^A\big)^{U_A}\otimes_A \mz[\zeta,\zeta^{-1}]$, and therefore the expressions of $\mathcal{W}(\gamma_i)$, associated to different presentations of $\Sigma$ as above are related by isomorphisms induced by the mapping class group of~$\Sigma$.

Now assume that $\gamma\in\mathcal{S}_l(\Sigma)$, so $l$ divides $c_i$, $1\leq i\leq k$, in \eqref{exprWT}. Put $c'_i:=c_i/l$. Then
\begin{gather}
T_{c_i}(\mathcal{W}(\gamma_i)) = T_{c_i'}(T_l(\mathcal{W}(\gamma_i)))
 = T_{c_i'}\big(\!\operatorname{Tr}\big(\!\operatorname{Fr}\underline{\stackrel{V_2}{M}}{}^{(j_i)}\cdots \operatorname{Fr}\underline{\stackrel{V_2}{M}}{}^{(j_i+n_i-1)}\big)\big), \label{exprWT2}
 \end{gather}
where we use the standard identity $T_{c_i} = T_{c_i'}\circ T_l$ in the first equality, and Proposition \ref{Thcentral} in the second.

By the first fundamental theorem of classical invariant theory for ${\rm SL}_2$ (see, e.g., \cite{KP}) and the Cayley--Hamilton identity, the set of trace functions
\[
t_{j_i,\dots,j_i+n_i-1}\colon\ \big(\underline{\stackrel{V_2}{M}}{}^{(1)},\dots,\underline{\stackrel{V_2}{M}}{}^{(n)}\big)\mapsto \operatorname{Tr}\big(\underline{\stackrel{V_2}{M}}{}^{(j_i)}\cdots \underline{\stackrel{V_2}{M}}{}^{(j_i+n_i-1)}\big)
\]
for all possible tuples $(j_i,\dots,j_i+n_i-1)$ are generating functions of $\mathcal{O}(G^n)^G$. By \eqref{exprWT2} each of these functions is sent by $Fr$ to $T_l\left(\mathcal{W}(\gamma_i)\right)$ for some $\gamma_i$. By Corollary\ref{equivcor}(3), $\widetilde{\rm Fr}$ is an isomorphism from $\mathcal{O}(\tilde G^n)^G$ to $\mathcal{Z}(\Ll_{0,n}^\e)^\Gg$, mapping $\mathcal{O}(G^n)^G$ to $\mathcal{Z}_0(\Ll_{0,n}^\e)^\Gg$, and the $n$ copies of~$\mathcal{O}(T)^W$ embedded in the $n$ factors of $\mathcal{O}(\tilde G)^{\otimes n} = \mathcal{O}(\tilde G^n)$ to $\mc\big[\omega^{(1)}\big],\dots,\mc\big[\omega^{(n)}\big]\subset \big(\Ll_{0,n}^A\big)^{U_A}_\e$. Therefore $\mathcal{Z}_(\Ll_{0,n}^\e)^\Gg$ is spanned over $\mc$ by the elements $\mathcal{W}(T(\gamma))$ with $\gamma\in\mathcal{S}_l(\Sigma)$ and $\gamma_\partial$ a~monomial in $\partial_1,\dots,\partial_n$. In particular $\mathcal{Z}(\Ll_{0,n}^\e)^\Gg$ is a (central) subalgebra of $\big(\Ll_{0,n}^A\big)^{U_A}_\e$.

(2) It is enough to show that the derivations $D_a$, $a\in \mathcal{Z}(\Ll_{0,n}^\e)^\Gg$, restrict to endomorphisms of $\big(\Ll_{0,n}^A\big)^{U_A}_\e$. This follows from the inclusion $\mathcal{Z}(\Ll_{0,n}^\e)^\Gg\subset \big(\Ll_{0,n}^A\big)^{U_A}_\e$, since $[\tilde{a},u] \in \big(\Ll_{0,n}^A\big)^{U_A}$ for every $\tilde{a}, u\in \big(\Ll_{0,n}^A\big)^{U_A}$ ($\Ll_{0,n}^A$ being an $U_A$-module algebra), which by the definition of $D_a$ implies that $D_a\big(\big(\Ll_{0,n}^A\big)^{U_A}_\e\big) \subset \big(\Ll_{0,n}^A\big)^{U_A}_\e$ when $a\in \mathcal{Z}(\Ll_{0,n}^\e)^\Gg$.
\end{proof}

\begin{remk}[{threading $T_l(V_2)$}] 
For every $1\leq i\leq n$, $0\leq k\leq n-i$ we have
\[
T_l \big(\!\operatorname{qTr}_{V_2} \big(\! \stackrel{V_2}{M}{}^{\!\! (i)}\stackrel{V_2}{M}{}^{\!\! (i+1)}\cdots \stackrel{V_2}{M}{}^{\!\! (i + k)}\big) \big) = \operatorname{qTr}_{T_l(V_2)} \big(\!\stackrel{T_l(V_2)}{M}{}^{\!\! (i)}\stackrel{T_l(V_2)}{M}{}^{\!\! (i+1)}\cdots \stackrel{T_l(V_2)}{M}{}^{\!\! (i + k)}\big),
\]
where $T_l(V_2)$ is the virtual representation in the Grothendieck ring of $U_A^{\rm res}$-modules, obtained by plugging $V_2$ in the $l$-th Chebyshev polynomial $T_l$. Indeed
\[
T_l\big(\!\operatorname{qTr}\big(\!\stackrel{V_2}{M}\big)\big)
= \operatorname{qTr}\big(\!\stackrel{\! T_l(V_2)}{M}\!\big)
\]
because (using, e.g., the second picture in Figure~\ref{fig6.2})
\[
\operatorname{qTr}\big(\!\!\stackrel{\! V_2^{\otimes k}}{M}\!\!\big) = \big(\!\operatorname{qTr}\big(\!\stackrel{V_2}{M}\big)\big)^k.
\]
This and the computations in the first half of the proof of Proposition \ref{Thcentral} prove our claim.
\end{remk}

\begin{remk}
Recall that in Theorem \ref{teoOBS} we have set $\zeta:={\rm i}q^{1/2}$. Let $\epsilon'$ be such that $(\epsilon')^2 = -\epsilon$, and define $K_{\epsilon'}(\Sigma) := K_\zeta(\Sigma)\otimes_{\mz[\zeta,\zeta^{-1}]} \mc_{\epsilon'}$, where $\mc_{\epsilon'}$ is the $\mz\big[\zeta,\zeta^{-1}\big]$-module $\mc$, where $\zeta$ acts by multiplication by $\epsilon'$. By a result of~\cite{FKL}, the set $\{T(\gamma), \gamma\in\mathcal{S}_l(\Sigma)\}$ is a $\mc$-basis of the center of $K_{\e'}(\Sigma)$. Therefore Theorem \ref{teoOBS} and the proof of~(1) show that $\mathcal{Z}(\Ll_{0,n}^\e)^\Gg$ and $\eta$ generate $\mathcal{Z}\big(\big(\Ll_{0,n}^A\big)^{U_A}_\e\big)$.
\end{remk}

\begin{remk}\label{addX}
As in the proof of (1), denote by $\partial_1,\dots,\partial_{n+1}$ the boundary components of $\Sigma$, and fix the generators of $\pi_1(\Sigma)$ to be $\partial_1,\dots,\partial_n$ (for some choice of basepoints). Consider the fiber product (with the notations of Remark \ref{tildeG})
\[
\tilde{G}^n\times_{G/\!/G} G/\!/G= \big\lbrace ((\tilde{g}_1,\dots,\tilde{g}_n),[t])\in \tilde{G}^n\times G/\!/G\ |\ p(g_1\cdots g_n) = p_l([t])\big\rbrace .
\]
As before Corollary \ref{equivcor}, one can identify $\tilde{G}^n\times_{G/\!/G} G/\!/G$ with an $l$-fold branched covering space $\tilde R_G(\Sigma)$ of $\tilde R_G'(\Sigma)$, whose points are given by representations $\rho\colon \pi_1(\Sigma)\ra G$ endowed with a choice of solution $x\in \mc$ of the equation $T_l(x) = \operatorname{Tr}(\rho(\partial_i))$, for every $i=1,\dots,n$ and $i = n+1$. Taking algebraic quotients yields a branched covering map $\tilde X_G(\Sigma)\ra X_G(\Sigma)$ of degree $l^{n+1}$, and one can lift $\widetilde{\rm Fr}\circ \mathfrak{c}^{*-1}$ to an isomorphism from $\mathcal{O}(\tilde X_G(\Sigma))$ onto the algebra generated by $\eta$ and~$\mathcal{Z}(\Ll_{0,n}^\e)^\Gg$.
\end{remk}

\section{Applications to skein algebras}\label{appskein}
As usual denote by $K_\zeta(\Sigma)$ the skein algebra of the sphere with $n+1$ punctures, $n>1$. Recall~$K_\zeta(\Sigma)$ is an algebra over $\mz\big[\zeta,\zeta^{-1}\big]$, where $\zeta$ is an indeterminate. In this section we use Theorem~\ref{teoOBS} to reformulate some of our results on $\big(\Ll_{0,n}^A\big)^{U_A}$ in the case of $\mathfrak{g} = {\mathfrak{sl}}(2)$ in terms of~$K_\zeta(\Sigma)$.

In \cite{BW0,PS} it was proved that $K_\zeta(\Sigma)$ does not have non trivial zero divisors, and its center was computed by topological means. On the contrary the proof we give below of these two facts is purely algebraic, based on properties of $\big(\Ll_{0,n}^A\big)^{U_A}$ proved in Section \ref{Lgnalg}, that hold true for any complex finite dimensional simple Lie algebra $\mathfrak{g}$.

\begin{cor}
The skein algebra $K_\zeta(\Sigma)$ does not have non trivial zero divisors, and its center is the polynomial algebra over $\mz\big[\zeta,\zeta^{-1}\big]$ generated by the classes $\partial_1,\dots,\partial_{n+1}$ of the boundary components of $\Sigma$.
\end{cor}

\begin{proof}The claims are direct consequences of Theorem \ref{teoOBS}, the formulas \eqref{formcentq}, and Pro\-po\-si\-tion~\ref{nozeroq} and the last claim of Theorem~\ref{CenterLinv} in the case of $\mathfrak{g}={\mathfrak{sl}}(2)$.
\end{proof}

Next we deduce from the results of Section~\ref{specialization} some properties of the center of the spe\-ci\-a\-li\-za\-ti\-ons of $K_\zeta(\Sigma)$ at roots of unity of order $4l$, $l\geq 3$ odd. Let $\epsilon'$ be such that $(\epsilon')^2 = -\epsilon$. Define
\[
K_{\epsilon'}(\Sigma) := K_\zeta(\Sigma)\otimes_{\mz[\zeta,\zeta^{-1}]} \mc_{\epsilon'},
\]
where $\mc_{\epsilon'}$ is the $\mz\big[\zeta,\zeta^{-1}\big]$-module $\mc$, where $\zeta$ acts by multiplication by~$\epsilon'$.
The sets $\mathcal{S}(\Sigma)$ and $\{T(\gamma), \gamma\in \mathcal{S}(\Sigma)\}$ are $\mz\big[\zeta,\zeta^{-1}\big]$-basis of $K_\zeta(\Sigma)$ as well as $\mc$-basis of $K_{\epsilon'}(\Sigma)$. By Theorem~\ref{teoOBS} we get an isomorphism of algebras
\[
\mathcal{W}\colon\ K_{\epsilon'}(\Sigma) \ra \big(\Ll_{0,n}^A\big)^{U_A}_\e.
\]
Denote by $\mathcal{Z}(K_{\epsilon'}(\Sigma))$ the center of $K_{\epsilon'}(\Sigma)$. Theorem~\ref{centerprop1}(1) implies immediately:

\begin{cor}
The set $\{T(\gamma), \gamma\in\mathcal{S}_l(\Sigma)\}$ is a $\mc$-basis of a central subalgebra $\mathcal{Z}'(K_{\epsilon'}(\Sigma))$ of~$K_{\epsilon'}(\Sigma)$, which is generated by $\eta$ and the image of the embedding
\[
{\rm Ch}_{\mathcal{W}} :=\mathcal{W}^{-1}\circ \widetilde{\rm Fr} \circ \tilde{\mathfrak{c}}^{*-1}\colon\ \mathcal{O}(\tilde{X}_G'(\Sigma)) \ra \mathcal{Z}(K_{\epsilon'}(\Sigma)).
\]
Therefore $\mathcal{Z}'(K_{\epsilon'}(\Sigma))$ is endowed with a natural Poisson bracket, the image of $\{\, ,\, \}_{\rm Gold}$, which extends to an action by derivations of $\mathcal{Z}'(K_{\epsilon'}(\Sigma))$ on $K_{\epsilon'}(\Sigma)$.
\end{cor}

In fact $\mathcal{Z}'(K_{\epsilon'}(\Sigma)) = \mathcal{Z}(K_{\epsilon'}(\Sigma))$ (see \cite{FKL}). One can check that ${\rm Ch}_{\mathcal{W}}$ is a version of the threading map ${\rm Ch}\colon \mathcal{K}_{\epsilon^{'l^2}}(\Sigma) \ra \mathcal{Z}_{\epsilon'}(\Sigma)$ of Bonahon--Wong (see \cite{BW1}, and also \cite{FKL2}). The point is that it affords an explicit realization of $Ch$ in classical invariant theory terms (via the for\-mu\-las~\eqref{exprWT}--\eqref{exprWT2}), and pulls the geometric tools of Section \ref{specialization} onto $K_{\epsilon'}(\Sigma)$. A more symmetric statement is obtained by extending ${\rm Ch}_{\mathcal{W}}$ to the ring $\mathcal{O}(\tilde{X}_G(\Sigma))$ of Remark \ref{addX}, so that $\partial_{n+1}$ belongs to its image.

\subsection*{Acknowledgements} We thank our colleagues of the work group on moduli spaces at IMAG for discussions, especially Paul-Emile Paradan and Damien Calaque. We also thank Pavel~Etingof, Matthieu Faitg, Charles~Frohman, and Catherine~Meusburger for their interest and exchanges related to the present work.
Finally, we also thank the referees for their suggestions, which greatly improved the exposition of the paper.

\pdfbookmark[1]{References}{ref}
\LastPageEnding

\end{document}